\pdfoutput=1
\RequirePackage{ifpdf}
\ifpdf 
\documentclass[pdftex]{sigma}
\else
\documentclass{sigma}
\fi

\numberwithin{equation}{section}

\newtheorem{Theorem}{Theorem}[section]
\newtheorem{Corollary}[Theorem]{Corollary}
\newtheorem{Lemma}[Theorem]{Lemma}
\newtheorem{Proposition}[Theorem]{Proposition}
 { \theoremstyle{definition}
\newtheorem{Definition}[Theorem]{Definition}

\newtheorem{Example}[Theorem]{Example}
\newtheorem{Remark}[Theorem]{Remark} }

\usepackage{mathtools}
\usepackage{mathrsfs}

\usepackage[all,cmtip]{xy}
\usepackage{bbm}
\usepackage[vcentermath]{youngtab}

\DeclareMathOperator{\Hom}{Hom}
\DeclareMathOperator{\End}{End}

\newcommand{\vp}{\varphi}

\makeatletter
\newsavebox{\@brx}
\newcommand{\llangle}[1][]{\savebox{\@brx}{\(\m@th{#1\langle}\)}%
 \mathopen{\copy\@brx\kern-0.5\wd\@brx\usebox{\@brx}}}
\newcommand{\rrangle}[1][]{\savebox{\@brx}{\(\m@th{#1\rangle}\)}%
 \mathclose{\copy\@brx\kern-0.5\wd\@brx\usebox{\@brx}}}
\makeatother

\begin{document}

\allowdisplaybreaks

\newcommand{\arXivNumber}{1712.08575}

\renewcommand{\thefootnote}{}

\renewcommand{\PaperNumber}{040}

\FirstPageHeading

\ShortArticleName{Local Moduli of Semisimple Frobenius Coalescent Structures}

\ArticleName{Local Moduli of Semisimple Frobenius Coalescent\\ Structures}

\Author{Giordano COTTI~$^\dag$, Boris DUBROVIN~$^\ddag$ and Davide GUZZETTI~$^\ddag$}

\AuthorNameForHeading{G.~Cotti, B.~Dubrovin and D.~Guzzetti}

\Address{$^\dag$~Max-Planck-Institut f\"ur Mathematik, Vivatsgasse 7 - 53111 Bonn, Germany}
\EmailD{\href{mailto:gcotti@sissa.it}{gcotti@sissa.it}, \href{mailto:G.Cotti.1@bham.ac.uk}{G.Cotti.1@bham.ac.uk}}
\URLaddressD{\url{https://www.birmingham.ac.uk/staff/profiles/maths/cotti-giordano.aspx}}

\Address{$^\ddag$~SISSA, Via Bonomea 265 - 34136 Trieste, Italy}
\EmailD{\href{mailto:guzzetti@sissa.it}{guzzetti@sissa.it}}

\ArticleDates{Received June 18, 2019, in final form April 13, 2020; Published online May 07, 2020}

\Abstract{We extend the analytic theory of Frobenius manifolds to semisimple points with coalescing eigenvalues of the operator of multiplication by the Euler vector field. We clarify which freedoms, ambiguities and mutual constraints are allowed in the definition of mono\-dromy data, in view of their importance for conjectural relationships between Frobenius manifolds and derived categories. Detailed examples and applications are taken from singularity and quantum cohomology theories. We explicitly compute the monodromy data at points of the Maxwell Stratum of the $A_3$-Frobenius manifold, as well as at the small quantum cohomology of the Grassmannian $\mathbb G_2\big(\mathbb C^4\big)$. In {the latter} case, we analyse in details the action of the braid group on the monodromy data. This proves that these data can be expressed in terms of characteristic classes of mutations of Kapranov's exceptional 5-block collection, as conjectured by one of the authors.}


\Classification{34M56; 53D45; 18E30}

\tableofcontents

\renewcommand{\thefootnote}{\arabic{footnote}}
\setcounter{footnote}{0}

\section{Introduction and results}

There is a conjectural relation, formulated by the second author (\cite{dubro0,dubro4}, see also \cite{bm4,hmt} and references therein), between the enumerative geometry of a wide class of smooth projective va\-rieties and their derived category of coherent sheaves. In particular, there is an increasing interest for an explicit description of certain local invariants, called monodromy data, of {\it semisimple} quantum cohomologies in terms of characteristic classes of exceptional collections in the derived categories \cite{dubro4,gamma1}. Being intentioned to address this problem, which, to our opinion, is still not well understood, we have realized that some issues in the theory of Frobenius manifolds need to be preliminarily clarified, and that an extension of the theory itself is necessary, in view of the fact that quantum cohomologies of certain classes of homogeneous spaces may show a~{\it coalescence phenomenon}.

In this paper, after reviewing the definition of the monodromy data, such as the {\it Stokes matrix} and the {\it central connection matrix}, we clarify their mutual constraints, the freedom and the natural operations allowed when we associate the data to a chart of the manifold. See Theorem~\ref{6-07-17-1} and \ref{6-07-17-2} in the Introduction and Sections~\ref{sec1} and~\ref{freedom}.
This issue does not
 seem to be sufficiently clear in the existing literature (some minor imprecisions are found also in \cite{dubro1,dubro0,dubro2}, especially concerning the central connection matrix), and it is fundamental in order to study the above mentioned conjectures.

Then, we extend the analytic theory of Frobenius manifolds in order to take into account a~\emph{coalescence phenomenon}, which occurs already for simple classes of varieties (e.g., for almost all complex Grassmannians~\cite{cotti0}). By this
 we mean that the operator of multiplication by the {\it Euler vector field} does not have a simple spectrum at some points where nevertheless the Frobenius algebra is semisimple. We call these points \emph{semisimple coalescence points} (see Definition~\ref{3agosto2016-9}). Such a phenomenon forbids an immediate application of the analytic theory of Frobenius mani\-folds to the computation of monodromy data. On the other hand, typically, the Frobenius structure is explicitly known only at the locus of semisimple coalescence points. Thus, we need to prove that the monodromy data associated with each region of the manifold can be computed starting only from the knowledge of the manifold at coalescence points.
 From the analytic point of view, coalescence implies that we have to deal with isomonodromic linear differential systems which violate one of the main assumptions\footnote{See \cite[p.~312]{JMU1}, assumption that the eigenvalues of $A_{\nu-r_\nu}$ are distinct. See also condition~(2) at p.~133 of~\cite{painkapaev}.} of the monodromy preserving deformations theory of M.~Jimbo, T.~Miwa and K.~Ueno~\cite{JMU1}.  Applying the results of~\cite{CDG0}, where the isomonodromy deformation theory has been extended to coalescence loci, we will show that the monodromy data computed at a semisimple coalescence point are the data associated with a whole neighbourhood of the point. The result is in Theorem \ref{6-07-17-3} of the Introduction and in Theorem \ref{mainisoth}. Moreover, by an action of the braid group, these data suffice to compute the data of the whole manifold (see Section \ref{30luglio2016-3}).

We give two explicit examples of the above procedure. One is from singularity theory in Section \ref{A_3case}, where we compute the monodromy data at points of the Maxwell Stratum of the $A_3$-Frobenius manifold. The monodromy data of the $A_3$-Frobenius manifold are well known, but here we present for the first time in the literature their computation at a coalescence point, where the isomonodromic system highly simplifies and can be solved in terms of Hankel functions. The validity of the result of the computation for the whole $A_3$-Frobenius manifold is justified by Theorem \ref{mainisoth}.

 The second example is new: in Section \ref{29novembre2016-1}, we explicitly compute the Stokes matrix and the central connection matrix of the quantum cohomology of the Grassmannian $\mathbb G_2\big(\mathbb C^4\big)$. The computation can be done only at the locus of small cohomology, which is the only locus where the structure of the manifold is known. This is a locus of coalescence points.
 Thus, we will have to do the computation {\it at a coalescence point}, explicitly obtaining Stokes and central connection matrices there. Theorem \ref{mainisoth} becomes crucial in order to give geometrical meaning to our computation: it guarantees that the result obtained {\it at a coalescence point} provides the monodromy data in a whole neighbourhood of the point. Consequently the result can be extended {\it to the whole manifold} by an action of the braid group (explained in the paper). In this way, we will obtain the monodromy data of $QH^\bullet\big(\mathbb G_2\big(\mathbb C^4\big)\big)$ and explicitly relate them to characteristic classes of objects of an exceptional collection in $\mathcal D^b\big(\mathbb G_2\big(\mathbb C^4\big)\big)$, establishing a correspondence between each region of the quantum cohomology and a full exceptional collection. See Theorem~\ref{6-07-17-4} in the Introduction and Theorem \ref{resultg24}. To our best knowledge, it seems that such an {\it explicit} description has not been done in the literature.

 The results of the example in Section \ref{29novembre2016-1} are important, because they make it evident that both formulations in \cite{dubro4} and \cite{gamma1}, refining the original conjecture of \cite{dubro0}, require more investigations.
 A refinement of the conjecture of \cite{dubro0,dubro4}, its proof for the case of all complex Grassmannians, and the relation to the version given in \cite{gamma1}, will be the content of our paper \cite{CDG1}.

 In conclusion, the example in Section \ref{29novembre2016-1} yields an explicit result of great theoretical significance in the theory of Frobenius manifolds and quantum cohomology, because it clarifies the conjecture explained above and because it shows the crucial role of Theorem \ref{mainisoth}. From another point of view, it may also be considered as a non-trivial example of analysis of a linear differential system with coalescing eigenvalues, showing a high level of complexity. As such, it may serve as a useful example in the field of linear systems and isomonodromy deformations.

Before explaining the results of the paper in more detail, we briefly recall preliminary basic facts.
A Frobenius manifold $M$ is a complex manifold, with finite dimension
\begin{gather*} n:= {\rm dim}_\mathbb{C}(M),\end{gather*}
 endowed with a structure of associative, commutative algebra with product $\circ_p$ and unit on each tangent space $T_pM$, analytically depending on the point~$p\in M$; in order to be \emph{Frobenius} the algebra must also satisfy an invariance property with respect to a symmetric non degenerate bilinear form $\eta$ on $TM$, called \emph{metric}, invariant w.r.t.\ the product~$\circ_p$, i.e.,
\[\eta(a\circ_pb,c)=\eta(a,b\circ_pc)\qquad\text{for all }a,b,c\in T_pM,\ p\in M,
\]
whose corresponding Levi-Civita connection $\nabla$ is \emph{flat} and, moreover, the unit vector field is flat. The above structure is required to be compatible with a $\mathbb C^*$-action on~$M$ (the so-called quasihomogeneity assumption). A precise definition and description will be given in Section~\ref{sec1}, including the potentiality (namely, the fact that the structure constants of the algebra product~$\circ_p$ are third derivatives of a function at~$p$).

The geometry of a Frobenius manifold is (almost) equivalent to the flatness condition for an extended connection $\widehat\nabla$ defined on the pull-back $\pi^*TM$ {of the tangent vector bundle} along the projection map $\pi\colon \mathbb C^*\times M \to M$. Consequently, we can look for $n$ holomorphic functions $\tilde t^1,\dots,\tilde t^n\colon\mathbb C^*\times M\to\mathbb C$ such that $\big(z,\tilde t^1,\dots,\tilde t^n\big)$ are $\widehat\nabla$-flat coordinates. In $\nabla$-flat coordinates $t=\big(t^1,\dots ,t^n\big)$, the $\widehat\nabla$-flatness condition $\widehat\nabla {\rm d}\tilde t(z,t)=0$ for a single function $\tilde t$
reads
\begin{gather}\label{intro1}
\frac{\partial \zeta}{\partial z} =\left(\mathcal U(t)+\frac{1}{z}\mu\right)\zeta,\\
\label{intro1.1}
\frac{\partial \zeta}{\partial t^\alpha} =z\mathcal C_\alpha(t) \zeta,\qquad\alpha=1,\dots, n,
\end{gather}
where the entries of the column vector $\zeta(z,t)$ are the components of the $\eta$-gradient of $\tilde t$
\begin{gather}\label{30nov2016-1}
\operatorname{grad}\tilde t:=\zeta^\alpha(z,t)\frac{\partial}{\partial t^\alpha},\qquad
 \zeta^\alpha(z,t):=\eta^{\alpha\nu}\frac{\partial\tilde t}{\partial t^\nu},\qquad \eta_{\alpha\beta}:=\eta\left(\frac{\partial}{\partial t^\alpha},\frac{\partial}{\partial t^\beta}\right),
\end{gather}
and $\mathcal C_\alpha(t)$,
$\mathcal U(t)$ and \begin{gather*}\mu:=\operatorname{diag}(\mu_1,\dots,\mu_n)\end{gather*} are $n\times n$ matrices described in Section~\ref{sec1}, satisfying $\eta \mathcal U=\mathcal U^{\rm T}\eta$ and $\eta\mu+\mu^{\rm T}\eta=0$.

A fundamental matrix solution of \eqref{intro1}--\eqref{intro1.1} provides $n$ independent $\widehat\nabla$-flat coordinates $\big(\tilde t^1,\dots,\tilde t^n\big)$. For fixed~$t$, the equation~\eqref{intro1} is an ordinary linear differential system with rational coefficients, with \emph{Fuchsian} singularity at $z=0$ and an \emph{irregular} singularity of
 Poincar\'e rank~1 at $z=\infty$.

A point $p\in M$ is called {\it semisimple} if the Frobenius algebra $T_pM$ is semisimple, i.e.,
 without nilpotents. A Frobenius manifold is {\it semisimple} if it contains an open dense subset $M_{ss}$ of semisimple points. If the matrix $\mathcal U$ is diagonalizable at $p$ with
\emph{pairwise distinct} eigenvalues, then $p\in M_{ss}$, see \cite{dubro1,dubro2}; this is a simple consequence of the definition of $\mathcal U$, given in the paper in expression~(\ref{26gen2020-1}). The condition is not necessary: there exist semisimple
points $p\in M_{ss}$ where $\mathcal U$ does not have a simple spectrum. In this case, if we move in $M_{ss}$ along a~curve terminating at~$p$ then some eigenvalues of $\mathcal U(t)$ coalesce.

The eigenvalues $u:=(u_1,\dots, u_n)$ of the operator $\mathcal U$, with chosen labelling, define a local system of coordinates $p\mapsto u=u(p)$ in a neighborhood of any semisimple point~$p$, called \emph{canonical}. In canonical coordinates, we set
\begin{gather}\label{intro2.0}
\operatorname{grad}\tilde t^\alpha(u,z)\equiv \sum_i Y^i_\alpha (u,z)f_i(u),\qquad f_i(u):=\frac{1}{\eta\big(\frac{\partial}{\partial u_i}|_u,\frac{\partial}{\partial u_i}|_u\big)^\frac{1}{2}}\left.\frac{\partial}{\partial u_i} \right|_u
\end{gather}
for some choice of the square roots.
Then the equations \eqref{intro1}, \eqref{intro1.1}, are equivalent to the following system
\begin{gather}\label{intro2}
\frac{\partial Y}{\partial z} =\left(U+\frac{V(u)}{z}\right)Y,\\
\label{intro2.1}\frac{\partial Y}{\partial u_k} =(zE_k+V_k(u))Y,\qquad 1\leq k\leq n,
\end{gather}
where $(E_k)_{ij}:=\delta_{ik}\delta_{jk}$, $U=\operatorname{diag}(u_1,\dots, u_n)$, $V$ is skew-symmetric and
\[U:=\Psi\mathcal U\Psi^{-1},\qquad V:=\Psi\mu\Psi^{-1},\qquad V_k(u):=\frac{\partial\Psi(u)}{\partial u_k}\Psi(u)^{-1}.
\]
Here, $\Psi(u)$ is a matrix defined by the change of basis between $\big(\frac{\partial}{\partial t^1},\dots, \frac{\partial}{\partial t^n}\big)$ and the normalized canonical vielbein $(f_1,\dots, f_n)$
\[\quad \frac{\partial}{\partial t^\alpha}=\sum_{i=1}^n \Psi_{i\alpha}f_i.
\] The compatibility conditions of the equations \eqref{intro2}--\eqref{intro2.1} are
\begin{gather}\label{compatib1}
[U,V_k]=[E_k,V],
\qquad
\frac{\partial V}{\partial u_k}=[V_k,V].
\end{gather}
When $u_i\neq u_j$ for $i\neq j$, equations~\eqref{compatib1} coincide with the Jimbo--Miwa--Ueno {\it isomonodromy deformation equations} for system~\eqref{intro2}, with deformation parameters $(u_1,\dots, u_n)$ \cite{JM2,JM3, JMU1}. This isomonodromic property allows to classify semisimple Frobenius manifolds by {\it locally constant monodromy data} of~\eqref{intro2}. Conversely, such \emph{local invariants} allow to reconstruct the Frobenius structure by means of an inverse Riemann--Hilbert problem \cite{dubro1,dubro2,guzzetti2}. Below, we briefly recall how they are defined in \cite{dubro1,dubro2}.

In \cite{dubro1,dubro2} it was shown that system \eqref{intro2} has a fundamental solution near $z=0$ in {\it Levelt normal form}
\begin{gather}
\label{3luglio2016-1}
Y_0(z,u)=\Psi(u)\Phi(z,u) z^\mu z^R,\qquad \Phi(z,u):=\mathbbm 1+\sum_{k=1}^\infty\Phi_k(u)z^k,
\end{gather}
satisfying the orthogonality condition
\begin{gather}\label{3luglio2016-1-nov}
 \Phi(-z,u)^{\rm T}\eta  \Phi(z,u)=\eta \qquad \text{for all $z\in \mathcal R,  u\in M$}.
\end{gather}
Here, \begin{gather*}\mathcal R:= \widetilde{\mathbb{C}\backslash \{0\}}\end{gather*} denotes the universal cover of $\mathbb C\setminus\left\{0\right\}$ and $R$ is a certain nilpotent matrix, which is non-zero only if $\mu$ has some eigenvalues differing by non-zero integers.
Since $z=0$ is a regular singularity, $\Phi(z,u)$ is convergent.

If $u=(u_1,\dots, u_n)$ are pairwise distinct, so that $U$ has distinct eigenvalues, then the system~\eqref{intro2} admits a formal solution of the form
\begin{gather}
Y_{\rm formal}(z,u)=G(z,u){\rm e}^{z U},\nonumber\\ G(z,u)=\mathbbm 1+\sum_{k=1}^\infty G_k(u) \frac{1}{z^k},\qquad G(-z,u)^{\rm T}G(z,u)=\mathbbm 1.\label{intro3}
\end{gather}
Although $Y_{\rm formal}$ in general does not converge, it always defines the asymptotic expansion of a~unique genuine solution on any sectors in the universal covering~$\mathcal R$, having central opening angle $\pi+\varepsilon$ for~$\varepsilon>0$ sufficiently small.

The choice of a ray $\ell_+(\phi):= \{z\in\mathcal R\colon\arg z=\phi \}$ with directional angle $\phi\in\mathbb R$ induces a~decomposition of the Frobenius manifold into disjoint chambers.\footnote{This definition does not appear in \cite{dubro1,dubro2}. See also Remark~\ref{3dicembre2017-1}.} An {\it $\ell$-chamber} is defined (see Definition~\ref{ellecella}) to be any connected component of the open dense subset of points $p\in M$ such that the eigenvalues of $\mathcal U$ at $p$ are all distinct (so, in particular, they are points of~$M_{ss}$), and the ray $\ell_+(\phi)$ does not coincide with any Stokes rays at~$p$, namely $\Re (z(u_i(p)-u_j(p)))\neq 0$ for $i\neq j$ and $z\in \ell_+(\phi)$.

Let $p$ belong to an $\ell$-chamber, and let $u=(u_1,\dots ,u_n)$ be the canonical coordinates in a~neighbourhood of $p$ contained in the chamber. Then, there exist unique solutions $Y_{\rm left/right}(z,u)$ such that
\[
Y_{\rm left/right}(z,u)\sim Y_{\rm formal}(z,u) \qquad \text{for $z\to \infty$},
\]
respectively in the sectors
\begin{gather*}
\Pi_{\rm right}^\varepsilon(\phi) := \{z\in\mathcal R\colon \phi-\pi-\varepsilon < \arg z< \phi+\varepsilon\},\\
\Pi_\text{left}^\varepsilon(\phi):=\{z\in\mathcal R\colon \phi-\varepsilon<\arg z<\phi+\pi+\varepsilon\}.
\end{gather*}The two solutions $Y_{\rm left/right}(z,u)$ are connected by the multiplication by an invertible matrix $S$, called \emph{Stokes matrix}:
\[
Y_\text{left}(z,u)=Y_{\rm right}(z,u)S, \qquad \text{for all $z\in\mathcal R$}.
\]
$S$ has the ``triangular structure'' described in Theorem \ref{proprstok}. Namely, $S_{ij}\neq 0$ implies $S_{ji}=0$. In particular, diag$(S)=(1,\dots,1)$ and $S_{ij}=S_{ji}=0$ whenever $u_i=u_j$.
Moreover, there exists a \emph{central connection matrix} $C$, whose properties will be described later, such that
\[
Y_{\rm right}(z,u)=Y_0(z,u)C,\qquad \text{ for all $z\in \mathcal R$.}
\]
 In \cite{dubro1} and \cite{dubro2} it is shown that the coefficients $\Phi_k$'s and $G_k$'s are holomorphic at any point of every $\ell$-chamber and that the
{\it monodromy data } $\mu, R, S, C$ are {\it constant} over a $\ell$-chamber (the isomonodromy Theorem I and~II of \cite{dubro2}, cf.\ Theorems~\ref{isomonod1} and~\ref{isoth2} below). They define local invariants of the semisimple Frobenius manifold $M$. In this sense, there is a local identification of a semisimple Frobenius manifold with the space of isomonodromy deformation parameters $(u_1,\dots ,u_n)$ of the equation~\eqref{intro2}.

\subsection{Results}
We now describe the results of the paper at points 1, 2 and 3 below.
\begin{itemize}\itemsep=0pt
\item[{\bf 1.}] {\bf Ambiguity in associating monodromy data with a point of the manifold} (cf.\ Sections~\ref{sec1} and~\ref{freedom}).
\end{itemize}

 From the above discussion, we see that with a point $p\in M_{ss}$ such that $u_1(p),\dots ,u_n(p)$ are pairwise distinct, we associate the monodromy data $(\mu, R, S, C)$. These data are constant on the whole $\ell$-chamber containing $p$. Nevertheless, there is not a unique choice of $(\mu, R, S, C)$ at $p$. The understanding of this issue is crucial in order to undertake a meaningful and well-founded study of the
 conjectured relationships of the monodromy data coming from quantum cohomology of smooth projective varieties with derived categories of coherent sheaves on these varieties.

The starting point is the observation that a normal form \eqref{3luglio2016-1} is not unique because of some freedom in the choice of $\Phi$ and $R$ (in particular, even for a fixed $R$ there is a freedom in $\Phi$). The description of this freedom was given in \cite{dubro2}, with a minor imprecision, to be corrected below. Let us identify\footnote{See Section \ref{3agosto2016-7} for the precise definition.} all tangent spaces $T_pM$, for $p\in M$, using the Levi-Civita connection on~$M$, with a $n$-dimensional complex vector space $V$, so that $\mu\in\operatorname{End}(V)$ is a linear operator antisymmetric w.r.t.\ the bilinear form $\eta$. Let  $\mathcal G(\eta,\mu)$ be the complex $(\eta,\mu)$-parabolic orthogonal Lie group, consisting of all  endomorphisms $G\colon V\to V$ of the form $G=\mathbbm 1_V+\Delta$, with $\Delta$ a~$\mu$-nilpotent endomorphism, and such that  $\eta\big({\rm e}^{{\rm i}\pi \mu}Ga,Gb\big)=\eta\big({\rm e}^{{\rm i}\pi \mu}a,b\big)$ for any $a,b\in V$ (see Section~\ref{3agosto2016-7} and Definition~\ref{3agosto2016-1}). We denote by $\mathfrak g(\eta,\mu)$ its Lie algebra.

\begin{Theorem}[Section \ref{3agosto2016-7}\footnote{In the description of the monodromy phenomenon of solutions of the system~\eqref{intro2} near $z=0$ the assumption of semisimplicity is not used. This will be crucial only for the description of solutions near $z=\infty$. Theorem~\ref{6-07-17-1} can be formulated for system (\ref{intro1}), having the fundamental solution $\Xi_0=\Psi^{-1}Y_0$.}]\label{6-07-17-1}
Given a fundamental matrix solution of system \eqref{intro2} in Levelt form \eqref{3luglio2016-1} near $z=0$ holomorphically depending on $(u_1,\dots, u_n)$ and satisfying the orthogonality condition~\eqref{3luglio2016-1-nov}, with $\mu=\Psi (u)^{-1} V(u)\Psi(u)$ constant and diagonal, then the holomorphic function $R=R(u)$ takes values in the Lie algebra $\mathfrak g(\eta,\mu)$.
Moreover,
\begin{enumerate}\itemsep=0pt
\item[$1.$] All other solutions in Levelt form near $z=0$ are $Y_0(z,u) G(u)$, where $G$ is a holomorphic function with values in $\mathcal G(\eta,\mu)$; the Levelt normal form of $Y_0(z,u)G(u)$ has again the structure~\eqref{3luglio2016-1} with $R(u)$ replaced with $\widetilde{R}(u):=G(u)R(u)G(u)^{-1}$ $($cf.\ Theorem~{\rm \ref{3agosto2016-5})}.
\item[$2.$]
 Because of the compatibility of \eqref{intro2} and \eqref{intro2.1}, $G(u)$ can be chosen so that $\widetilde R$ is independent of $u$ $($isomonodromy Theorem I in {\rm \cite{dubro2}}, Theorem~{\rm \ref{isomonod1})}.
\item[$3.$]
 For a fixed $R\in\mathfrak g(\eta,\mu)$, the isotropy subgroup $\mathcal G(\eta,\mu)_R$ of transformations $G\in \mathcal G(\eta,\mu)$, such that $GRG^{-1}=R$, can be identified with the group
 \begin{gather}\label{3dicembre2017-2}
\mathcal C_0(\eta,\mu, R):= \left\{
G\in {\rm GL}(V)\colon \ \begin{aligned}
  &P_G(z):=z^\mu z^R Gz^{-R}z^{-\mu}\in\operatorname{End}(V)[z],\\
  & P_G(0)\equiv\mathbbm 1_V, \\
  &\eta(P_G(-z)v_1,P_G(z)v_2)=\eta(v_1,v_2),\\
  & \text{for all }v_1,v_2\in V
 \end{aligned}
\right\}.
\end{gather}
\end{enumerate}
\end{Theorem}

The definition (\ref{3dicembre2017-2}) can be re-written in coordinates as follows
\begin{gather*}
\mathcal C_0(\eta,\mu, R):= \left\{
 \begin{aligned}
 & G\in  {\rm GL}(n,\mathbb C)\colon~ P_G(z):=z^\mu z^R Gz^{-R}z^{-\mu}\text{ is a matrix-valued} \\
 &\text{polynomial such that } P_{G}(0)=\mathbbm 1,\ \text{and }P_G(-z)^{\rm T}\eta P_G(z)=\eta
 \end{aligned}
\right\}.
\end{gather*}
If $G\in\mathcal C_0(\eta,\mu, R)$ and $Y_0(z,u)=\Psi(u)\Phi(z,u) z^\mu z^R$, then $Y_0(z,u)G=\Psi(u)\Phi(z,u) P_G(z) z^\mu z^R$. The refinement introduced here consists in the restriction of the group ${\mathcal C}_0(\mu, R)$ introduced in~\cite{dubro2} by adding the condition
 \begin{gather*}P_G(-z)^{\rm T}\eta P_G(z)=\eta,
 \end{gather*}
which does not appear in \cite{dubro2}. In \cite{dubro2} neither the $\eta$-orthogonality conditions appeared in the definition of the group $\mathcal C_0(\mu, R)$, nor this group was identified with the isotropy subgroup of $R$ w.r.t.\ the adjoint action of $\mathcal G(\eta,\mu)$ on its Lie algebra $\mathfrak g(\eta,\mu)$. These $\eta$-orthogonality conditions are crucial for preserving~(\ref{3luglio2016-1-nov}) and the constraints~(\ref{29giugno2017-1}) of all monodromy data $(\mu, R,S,C)$ (see also Theorem~\ref{constraint}).

Let us now summarize the freedom in assigning the monodromy data $(\mu,R,S,C)$ to a given semi-simple point~$p$ of the Frobenius manifold. It has various origins: it can come from a re-ordering of the canonical coordinates $u_1(p), \dots, u_n(p)$, from changing signs of the normalized idempotents, from changing the Levelt fundamental solution at $z=0$ and, last but not least, from changing the slope of the oriented line $\ell_+(\phi)$. Taking into account all these possibilities,
we have the following

 \begin{Theorem}[Section \ref{freedom}] \label{6-07-17-2}  Let $p\in M_{ss}$ be such that $(u_1(p),\dots ,u_n(p))$ are pairwise distinct. If $(\mu, R,S,C)$ is a set of monodromy data computed at $p$, then with a different labelling of the eigenvalues, different signs, different choice of $Y_0(z,u)$ and different $\phi$, another set of monodromy data can be computed at {\rm the same} $p$, which lies in the orbit of $(\mu, R,S,C)$ under the following actions:
 \begin{itemize}\itemsep=0pt
 \item the action of the group of permutations $\mathfrak S_n$
\begin{gather*}
S\longmapsto PSP^{-1},\qquad C \longmapsto CP^{-1},
\end{gather*}
which corresponds to a relabelling $(u_1,\dots, u_n)\mapsto(u_{\tau(1)},\dots,u_{\tau(n)})$, where $\tau\in \mathfrak S_n$ and
the invertible matrix $P$ has entries $P_{ij}=\delta_{j\tau (i)}$. For a suitable choice of the permutation, $PSP^{-1}$ is in upper-triangular form;
\item the action of the group $(\mathbb Z/2\mathbb Z)^{\times n}$ \[
S\longmapsto \mathcal{I}S\mathcal{I},\qquad C\longmapsto C~\mathcal I,
\]
where $\mathcal I$ is a diagonal matrix with entries equal to $1$ or $-1$, which corresponds to a change of signs of the square roots in~\eqref{intro2.0};
\item the action of the group $\mathcal C_0(\eta,\mu, R)$
\[
S\longmapsto S,\qquad C\longmapsto G C,\qquad G\in \mathcal C_0(\eta,\mu, R),
\]
which corresponds to a change $Y_0(z,u)\mapsto Y_0(z,u) G$ as in Theorem~{\rm \ref{6-07-17-1}}.

\item the action of the braid group, as in formulae~\eqref{connbraid2-29luglio} and \eqref{connbraid2},
 \[
S\mapsto A^\beta(S)\cdot S\cdot \big(A^\beta(S)\big)^{\rm T}, \qquad C\mapsto C\cdot \big(A^\beta(S)\big)^{-1},
\]
 where $\beta$ is a specific braid associated with a translation of $\phi$, corresponding to a rotation of $\ell_+(\phi)$. More details are in Section~{\rm \ref{freedom}}.
\end{itemize}

Any representative of $\mu$, $R$, $S$, $C$ in the orbit of the above actions satisfies the monodromy identity
\begin{gather*}CS^{\rm T}S^{-1}C^{-1}={\rm e}^{2\pi{\rm i}\mu}{\rm e}^{2\pi{\rm i}R},\end{gather*}
and the constraints
\begin{gather}\label{29giugno2017-1}
S=C^{-1}{\rm e}^{-\pi{\rm i} R}{\rm e}^{-\pi{\rm i} \mu}\eta^{-1}\big(C^{\rm T}\big)^{-1},
\qquad S^{\rm T}=C^{-1}{\rm e}^{\pi{\rm i} R}{\rm e}^{\pi{\rm i} \mu}\eta^{-1}\big(C^{\rm T}\big)^{-1}.
\end{gather}
\end{Theorem}

We stress again that the freedoms in Theorem~\ref{6-07-17-2} must be taken into account when we want to investigate the relationship between monodromy data and similar objects in the theory of derived categories

\begin{itemize}\itemsep=0pt
\item[{\bf 2.}] {\bf Isomonodromy theorem at semisimple coalescence points} (cf.\ Section~\ref{isomonocoal}).
\end{itemize}

\begin{Definition}\label{3agosto2016-9}A point $p\in M_{ss}$ such that the eigenvalues of $\mathcal U$ at $p$ are not pairwise distinct is called a {\it semisimple coalescence point}.
\end{Definition}

The isomonodromy deformations results presented above apply if $U$ has distinct eigenvalues. If two or more eigenvalues coalesce, as it happens at semisimple coalescence points of Definition~\ref{3agosto2016-9}, then a priori solutions $Y_{\rm left/right}(z,u)$ are expected to have a singular behaviour (branching at $u_i-u_j= 0$ and/or divergent limits for $u_i-u_j\to 0$ along any direction), the coefficients of the formal expansion of $Y_\text{formal}(z,u)$ may have poles at $u_i-u_j= 0$, and monodromy data must be redefined.

In almost all studied cases of quantum cohomology the structure of the manifold is explicitly known only at the locus of small quantum cohomology defined in terms of the three-points genus zero Gromov--Witten invariants. Along this locus the coalescence phenomenon may occur (for example, coalescence occurs in case of the quantum cohomology of almost all Grassmannians~\cite{cotti0}). Therefore, if we want to compute the monodromy data, we can only rely on the information available at coalescence points. Thus, we need to extend the analytic theory of Frobenius manifolds, in order to include this case, showing that the monodromy data are well defined at a semisimple coalescence point, and locally constant. Moreover, from these data we must be able to reconstruct the data for the whole manifold. We stress that this extension of the theory is essential in order to study the conjectural links to derived categories.

The extension is based on the observation that the matrix $\Psi(u)$ is holomorphic at semisimple points including those of coalescence (see Lemma~\ref{proppsi}). Consequently, the matrices~$V_k$'s and~$V$ are holomorphic at any semisimple point, and~$V$ is holomorphically similar to~$\mu$. These are exactly the sufficient conditions allowing the application of the general results obtained in~\cite{CDG0}, which yield the following

\begin{Theorem}[cf.\ Theorem \ref{mainisoth} below]\label{6-07-17-3}
Let $p_0$ be a semisimple coalescence point with canonical coordinates $u(p_0)=\big(u^{(0)}_1,\dots, u^{(0)}_n\big)$.\footnote{\label{dldhlhod}Up to permutation, these coordinates can be arranged as \begin{gather*}
u^{(0)}_1=\dots=u^{(0)}_{r_1},\\
u^{(0)}_{r_1+1}=\dots=u^{(0)}_{r_1+r_2},\\
 \cdots\cdots\cdots\cdots\cdots\cdots\cdots\cdots\\
u^{(0)}_{r_1+\dots+r_{s-1}+1}=\dots =u^{(0)}_{r_1+\dots +r_{s-1}+r_s},
\end{gather*}
where $r_1,\dots , r_s$ ($r_1+\dots +r_{s-1}+r_s=n$) are the multiplicities of the eigenvalues of $U\big(u^{0}\big)=\operatorname{diag} \big(u^{(0)}_1,\dots, u^{(0)}_n\big)$.} Moreover,
\begin{itemize}\itemsep=0pt
\item Let $\phi\in \mathbb{R}$ be fixed so that $\ell_+(\phi)$ does not coincide with any Stokes ray at~$p_0$, namely $\Re \big(z\big(u_i^{(0)}-u_j^{(0)}\big)\big)\neq 0$ for $u_i^{(0)}\neq u_j^{(0)}$ and $z\in \ell_+(\phi)$.

\item Consider the closed polydisc centered at $u^{(0)}=\big(u^{(0)}_1,\dots, u^{(0)}_n\big)$ and of size $\epsilon_0>0$
\begin{gather*}
\mathcal{U}_{\epsilon_0}\big(u^{(0)}\big):=
\big\{u\in\mathbb{C}^n\colon \max_{1\leq i\leq n}\big|u_i-u_i^{(0)}\big|\leq \epsilon_0 \big\}.
\end{gather*}
 Let $\epsilon_0>0$ be sufficiently small, so that $\mathcal{U}_{\epsilon_0}\big(u^{(0)}\big)$ is homeomorphic by the coordinate map to a neighbourhood of $p_0$ contained in $M_{ss}$. An additional upper bound for $\epsilon_0$ will be specified in Section~{\rm \ref{isomonocoal}}, see equation~\eqref{2agosto2016-1}.

\item Let $\Delta\subset \mathcal{U}_{\epsilon_0}\big(u^{(0)}\big)$ be the locus in the polydisc $\mathcal{U}_{\epsilon_0}\big(u^{(0)}\big)$ where some eigenvalues of $U(u)=\operatorname{diag} (u_1,\dots, u_n)$ coalesce.\footnote{Namely, $u_i=u_j$ for some $1\leq i \neq j \leq n$ whenever $u\in\Delta$. The bound on $\epsilon_0$, to be clarified later, implies that, with the arrangement of footnote~\ref{dldhlhod} the sets
\begin{gather*}
\{u_1, \dots  ,u_{r_1}\},\quad \{u_{r_1+1},~\dots ~,u_{r_1+r_2}\},\quad \dots  ,\quad
\{u_{r_1+\cdots+r_{s-1}+1}, \dots ,u_{r_1+\cdots+r_{s-1}+r_s}\}
\end{gather*}
do not intersect for any $u\in \mathcal{U}_{\epsilon_0}\big(u^{(0)}\big)$. In particular, $u^{(0)}\in\Delta$ is a point of ``maximal coalescence''.}
\end{itemize}

Then, the following results hold:
\begin{enumerate}\itemsep=0pt
\item[$1.$] System \eqref{intro2} at the {\rm fixed} value $u=u^{(0)}$ admits a unique formal solution, which we denote with $\mathring{Y}_{\rm formal}(z)$, having the structure~\eqref{intro3}, namely
 $\mathring{Y}_{\rm formal}(z)=\Bigl(\mathbbm 1+\sum\limits_{k=1}^\infty \mathring{G}_k z^{-k}\Bigr){\rm e}^{zU}$; moreover, it admits unique fundamental solutions, which we denote with $\mathring{Y}_{\rm left/right}(z)$, ha\-ving asymptotic representation $\mathring{Y}_{\rm formal}(z)$ in sectors $\Pi_{\rm left/right}^\varepsilon(\phi)$, for suitable $\varepsilon>0$.\footnote{A more precise characterisation of the angular amplitude of the sectors will be given later.} Let~$\mathring{S}$ be the Stokes matrix such that
 \begin{gather*}\mathring{Y}_{\rm left}(z)=\mathring{Y}_{\rm right}(z)\mathring{S}.\end{gather*}

\item[$2.$] The coefficients $G_k(u)$, $k\geq 1$, in \eqref{intro3} are holomorphic over $\mathcal{U}_{\epsilon_0}\big(u^{(0)}\big)$, and $G_k\big(u^{(0)}\big)=\mathring{G}_k$; moreover $Y_{\rm formal}\big(z,u^{(0)}\big)=\mathring{Y}_{\rm formal}(z)$.

\item[$3.$] For fixed $z$, $Y_{\rm left}(z,u)$, $Y_{\rm right}(z,u)$, computed in a neighbourhood of a point $u\in \mathcal{U}_{\epsilon_0}\big(u^{(0)}\big)\backslash \Delta$, can be $u$-analytically continued as single-valued holomorphic functions on the whole \linebreak $\mathcal{U}_{\epsilon_0}\big(u^{(0)}\big)$. Moreover
\begin{gather*}
Y_{\rm left/right}\big(z,u^{(0)}\big)=\mathring{Y}_{\rm left/right}(z).
\end{gather*}

\item[$4.$] For any $\epsilon_1<\epsilon_0$ the asymptotic relations
\begin{gather*}
Y_{\rm left/right}(z,u){\rm e}^{-zU}\sim I+\sum_{k=1}^\infty G_k(u)z^{-k},\qquad \hbox{for $z\to\infty$ in $\Pi_{\rm left/right}^\varepsilon(\phi)$},
\end{gather*}
hold uniformly in $u\in \mathcal{U}_{\epsilon_1}\big(u^{(0)}\big)$. In particular they also hold at $u\in \Delta$.

\item[$5.$] Denote by $\mathring{Y}_0(z)$ a solution of system \eqref{intro2} with the {\rm fixed} value $u=u^{(0)}$, in Levelt form $\mathring{Y}_0(z)=\Psi\big(u^{(0)}\big) (\mathbbm 1+O(z) )z^\mu z^{\mathring{R}}$, having monodromy data $\mu$ and $\mathring{R}$. For any such $\mathring{Y}_0(z)$ there exists a fundamental solution $Y_0(z,u)$ in Levelt form \eqref{3luglio2016-1}, \eqref{3luglio2016-1-nov} holomorphic in~$\mathcal{U}_{\epsilon_0}\big(u^{(0)}\big)$, such that its monodromy data $\mu$ and $R$ are independent of $u$ and
\begin{gather}\label{2luglio2017-2}
Y_0\big(z,u^{(0)}\big)=\mathring{Y}_0(z),\qquad R=\mathring{R}.
\end{gather}
Let $\mathring{C}$ be the central connection matrix for $\mathring{Y}_0$ and $\mathring{Y}_{\rm right}$; namely \begin{gather*}\mathring{Y}_{\rm right}(z)=\mathring{Y}_0(z)\mathring{C}.\end{gather*}

\item[$6.$] For any $\epsilon_1<\epsilon_0$ the monodromy data $\mu$, $R$, $S$, $C$ of system \eqref{intro2} are well defined and constant in the whole $ \mathcal{U}_{\epsilon_1}\big(u^{(0)}\big)$, so that the system is isomonodromic in $\mathcal{U}_{\epsilon_1}\big(u^{(0)}\big)$. They coincide with the data
associated with the fundamental solutions $\mathring{Y}_{\rm left/right}(z)$ and $\mathring{Y}_0(z)$ above, namely
\begin{gather*} R=\mathring{R},\qquad S=\mathring{S},\qquad C=\mathring{C}.
\end{gather*}
 The entries of $S=(S_{ij})_{i,j=1}^n$ with indices corresponding to coalescing canonical coordinates vanish:
\begin{gather}\label{9dicembre2016-1}
S_{ij}=S_{ji}=0 \qquad \text{for all $i\neq j$ such that} \quad u_i^{(0)}=u_j^{(0)}.
\end{gather}
\end{enumerate}
\end{Theorem}

Theorem \ref{6-07-17-3} implies that, in order to compute the monodromy data $\mu$, $R$, $S$, $C$ in the whole $\mathcal{U}_{\epsilon_0}\big(u^{(0)}\big)$, it suffices to compute $\mu$, $\mathring{R}$, $\mathring{S}$, $\mathring{C}$ {\it at} $u^{(0)}$. These can be used to obtain the monodromy data at any other point of $M_{ss}$ (including semisimple coalescence points, by Theorem~\ref{6-07-17-3}), by the action of the braid group $\mathcal B_n$ introduced in~\cite{dubro1} and~\cite{dubro2}
 \begin{gather} \label{2luglio2017-5}
S\longmapsto A^\beta  S  \big(A^\beta\big)^{\rm T}, \qquad C\longmapsto C\big(A^\beta\big)^{-1},
\end{gather}
 as in formulae (\ref{connbraid2-29luglio}) and~(\ref{connbraid2}). This action is well defined whenever $u_1,\dots ,u_n$ are pairwise distinct. It allows to obtain the monodromy data associated with all $\ell$-chambers. Therefore, the action can be applied to $S$, $C$ as defined by the above theorem starting from a point of~$\mathcal{U}_{\epsilon_0}\big(u^{(0)}\big)$ where $u_1,\dots ,u_n$ are pairwise distinct.

We will give two detailed applications of the above theorem. The first example, in Section~\ref{A_3case}, is the analysis of the monodromy data at the points of one of the two irreducible components of the bifurcation diagram of the Frobenius manifold associated with the Coxeter group~$A_3$. This is the simplest polynomial Frobenius structure in which semisimple coalescence points appear. The whole structure is globally and explicitly known, and the system~\eqref{intro2} at generic points is solvable in terms of oscillatory integrals. At semisimple points of coalescence, however, the system considerably simplifies, and it reduces to a Bessel equation. Thus, the asymptotic analysis of its solutions can be easily completed using Hankel functions, and $S$ and $C$ can be immediately computed. By Theorem~\ref{6-07-17-3} above, these are monodromy data of points in a whole neighbourhood of the coalescence point. We explicitly verify that the fundamental solutions expressed by means of oscillatory integrals converge to those expressed in terms of Hankel functions at a coalescence point, and that the computation done away from the coalescence point provides the same~$S$ and~$C$, as Theorem~\ref{6-07-17-3} predicts. In particular, the Stokes matrix~$S$ computed invoking Theorem~\ref{6-07-17-3} is in agreement with both the well-known results of~\cite{dubro1,dubro2}, stating that $S+S^{\rm T}$ coincides with the Coxeter matrix of the group $W(A_3)$ (group of symmetries of the regular tetrahedron), and with the analysis of~\cite{dubrmazz} for monodromy data of the algebraic solutions of PVI$_\mu$ corresponding to $A_3$ (see also \cite{CDG0} for this last point).

The second example is the quantum cohomology $QH^{\bullet}\big(\mathbb G_2\big(\mathbb C^4\big)\big)$ of the Grassmannian $\mathbb G_2\big(\mathbb C^4\big)$, which sheds new light on the conjecture mentioned in the beginning.
\begin{itemize}\itemsep=0pt
\item[{\bf 3.}] {\bf Quantum cohomology of the Grassmannian $\mathbb G_2\big(\mathbb C^4\big)$} (cf.\ Section~\ref{29novembre2016-1}).
\end{itemize}

We consider the Frobenius structure on $QH^{\bullet}\big(\mathbb G_2\big(\mathbb C^4\big)\big)$. The \emph{small quantum ring} -- or {\it small quantum cohomology} -- of Grassmannians has been one of the first cases of quantum cohomology rings to be studied both in physics \cite{vafa, wit1} and mathematical literature \cite{bert, siebti}, so that a~quantum extension of the classical Schubert calculus has been obtained \cite{buc1}. However, the ring structure of the big quantum cohomology is not explicitly known, so that the computation of the monodromy data can only be done at the small quantum cohomology locus. It happens that the small quantum locus of almost all Grassmannians $\mathbb G_k\big(\mathbb C^n\big)$ is made of semisimple coalescence points, see~\cite{cotti0}; the case of $\mathbb G_2\big(\mathbb C^4\big)$ is the simplest case where this phenomenon occurs. Therefore, in order to compute the monodromy data, we invoke Theorem~\ref{6-07-17-3} above.

In Section \ref{29novembre2016-1}, we carry out the asymptotic analysis of the system (\ref{intro2}) {\it at the coalescence locus}, corresponding to $t=0\in QH^\bullet\big(\mathbb G_2\big(\mathbb C^4\big)\big)$. We explicitly compute the monodromy data $\mu$ and $R$ (see~(\ref{2luglio2017-1}) and~(\ref{28luglio2016-10})) and the Stokes matrix, receiving the matrix $S$ in expression~\eqref{29luglio2016-1} (with $v=6$). For the computation of $S$, we take an admissible\footnote{Namely, $\ell_+(\phi)$ defined above is an admissible ray.} line $\ell:=\big\{z\in\mathbb{C}\colon z={\rm e}^{{\rm i}\phi}\big\}$ with the slope $0<\phi<\frac{\pi}{4}$. The signs in the square roots in~\eqref{intro2.0} and the labelling of $(u_1,\dots,u_6)$ are chosen in Section~\ref{smcohg24}.  In order to compute the central connection matrix, we choose a~specific fundamental solution (\ref{2luglio2017-2}) of \eqref{intro2} with fixed $t=0$, namely the {\it enumerative-topological} fundamental solution\footnote{This is the solution $\Psi(0)Y(z,0)=\Psi(0)H(z,0)z^\mu z^R$ in Proposition~\ref{2luglio2017-3}, where $\Phi$ is called~$H$.} $Y_0(z):= \Psi\big|_{t=0} \Phi(z) z^\mu z^R$,  whose coefficients are the genus 0 Gromov--Witten invariants with descendants
\[\Phi(z)^\alpha_\beta=\delta^\alpha_\beta+\sum_{n=0}^\infty\sum_\lambda\sum_{\nu\in\textnormal{Eff}(\mathbb G)\setminus\left\{0\right\}}\langle\tau_nT_\beta, T_\lambda\rangle_{0,2,\nu}^{\mathbb G}\eta^{\lambda\alpha}z^{n+1},\]
with
\[ \langle\tau_nT_\beta, T_\lambda\rangle_{0,2,\nu}^{\mathbb G}:=\int_{[\mathbb G_{0,2,\nu}]^\text{vir}}\psi_1^n\cup\text{ev}_1(T_\beta)\cup\text{ev}_2(T_\nu),
\]
and $(\eta^{\mu\nu})$ the inverse of Poincar\'e metric. This solution will be precisely described in Section~\ref{topsolfano} (cf.\ Proposition~\ref{2luglio2017-3}). The computation yields the connection matrix~$C$ reported in Appendix~\ref{appc} (with $v=6$).

Given $S$ and $C$ computed as explained above, let us denote by $S^\prime$ and $C^\prime$
 the data obtained by the action
 \begin{gather}
 S \longmapsto  \mathcal{I}P S (\mathcal{I}P)^{-1}=:S^\prime,\nonumber\\
\label{c'}
C \longmapsto  G C  (\mathcal{I}P)^{-1}=:C^\prime,
\end{gather}
of the groups of Theorem~\ref{6-07-17-2}, where $P$ and $\mathcal I$ will be explicitly written in the proof of Theo\-rem~\ref{resultg24}, while $G=A$ or $G=AB\in \mathcal C_0(\eta,\mu, R)$, for certain matrices $A$ and $B$ given in~(\ref{6luglio2017-10}), (\ref{6luglio2017-11}) below. Geometrically, $P$, $\mathcal I$ and $G$ correspond respectively to
 \begin{itemize}\itemsep=0pt
\item an appropriate re-ordering of the canonical coordinates $u_1,\dots, u_6$ near $0\in QH^\bullet(\mathbb G)$, yielding the Stokes matrix in upper-triangular form.
\item another determination of signs in the square roots of~\eqref{intro2.0} of the normalized idempotents vector fields $(f_i)_i$
\item another choice of the fundamental solution of the equation \eqref{intro2} in Levelt-normal form \eqref{3luglio2016-1}, obtained from the enumerative-topological solution by the action $Y_0\mapsto Y_0G$ of $\mathcal C_0(\eta,\mu, R)$. \end{itemize}
Given these explicit data, we prove Theorem~\ref{6-07-17-4} below, which clarifies for $\mathbb G_2\big(\mathbb C^4\big)$ the conjecture, formulated by the second author in~\cite{dubro0} (see also \cite{bm4,hmt} and references therein) and then refined in~\cite{dubro4}, relating the enumerative geometry of a Fano manifold with its derived category (see also Remark~\ref{17-dec-2016}). More details and new more general results about this conjecture are the contents of our paper~\cite{CDG1} (see also Remark~\ref{17-dec-2016}).

For brevity, let $\mathbb G:=\mathbb G_2\big(\mathbb C^4\big)$. Fix the Schubert basis \begin{gather*}(T_0,T_1,T_2,T_3,T_4,T_5)=(1,\sigma_1,\sigma_{2},\sigma_{1,1},\sigma_{2,1},\sigma_{2,2})\end{gather*} of $H^\bullet(\mathbb G;\mathbb C)$. Let $c\in\mathbb C^*$ be defined by
\[\int_{\mathbb G}\sigma_{2,2}=c.
\]
Denote by $\mathcal S$ the tautological bundle on $\mathbb G$ and by $\mathbb S^\lambda$ the Schur functor associated with the Young diagram $\lambda$.
Let $(E_1,\dots, E_6)$ be the 5-block\footnote{This means that $\chi(E_3,E_4)=\chi(E_4, E_3)=0$ and thus that both $(E_1,E_2,E_3,E_4,E_5,E_6)$ and $(E_1,E_2,E_4,E_3,\allowbreak E_5,E_6)$ are exceptional collections: we will write \begin{gather*}\left(E_1,\quad E_2,\quad\begin{matrix}
E_3\\
E_4
\end{matrix},\quad E_5,\quad E_6\right)\end{gather*}if we consider the exceptional collection with an unspecified order.} exceptional collection, obtained from the Kapranov exceptional $5$-block collection
\[\left(\mathbb S^0\mathcal S^*,\quad \mathbb S^{1}\mathcal S^*,\quad
\begin{matrix}\mathbb S^{2}\mathcal S^*\\ \mathbb S^{1,1}\mathcal S^*\end{matrix},\quad\mathbb S^{2,1}\mathcal S^*,\quad\mathbb S^{2,2}\mathcal S^*\right),
\]
 by mutation\footnote{The definition of the action of the braid group on the set of exceptional {collections} will be given in Section~\ref{30luglio2016-7}, slightly modifying (by a shift) the classical definitions that the reader can find, e.g., in~\cite{helix}. Our convention for the composition of {action} of braids is the following: braids act on an exceptional collection/monodromy datum \emph{on the right}. The braid $\beta_{34}$ acts on the 5-block collection $(E_1,\dots, E_6)$ above just as a permutation of the third and fourth elements of the block. }
 under the {\rm inverse} of any of the following braids\footnote{Curiously, these braids show a mere mirror symmetry: notice that they are indeed equal to their specular reflection. Any contingent geometrical meaning of this fact deserves further investigations.} in $\mathcal B_6$:
 \begin{figure}[h!]\centering
\def\svgscale{1}
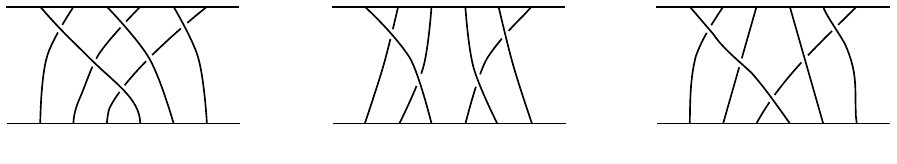
\end{figure}

Consider the Taylor expansion of the $\Gamma$-function
\[
\Gamma(1-x)=\exp\left\{ \gamma x +\sum_{n=2}^\infty \frac{\zeta(n)}{n}x^n \right\},\]
where
\[ \gamma:=\lim\limits_{n\to\infty}\left( - \log n+\sum\limits_{k=1}^n\frac{1}{k}\right)\] denotes the Euler--Mascheroni constant and $\zeta(n)$ is the Riemann zeta function. Let us introduce the following characteristic classes:
\begin{gather*} \widehat\Gamma^\pm(\mathbb G):=\prod_j \Gamma(1\pm\delta_j)\qquad\text{where }\delta_j\text{'s are the Chern roots of }T\mathbb G,\\
\operatorname{Ch}(V):=\sum_k {\rm e}^{2\pi{\rm i} x_k},\qquad x_k\text{'s are the Chern roots of a vector bundle $V$ on $\mathbb G$},\\ \operatorname{Ch}(V^\bullet):=\sum_j(-1)^j\operatorname{Ch}\big(V^j\big)\qquad\text{for a bounded complex $V^\bullet$ of vector bundles on $\mathbb G$.}
\end{gather*}

\begin{Theorem}[monodromy data of $QH^\bullet(\mathbb G)$ cf.~Theorem \ref{resultg24}]\label{6-07-17-4}
The Stokes matrix and the central connection matrix at $t=0\in QH^\bullet(\mathbb G)$ are related to the full exceptional collection $(E_1,\dots, E_6)$ in the following way.

The central connection matrix $C^\prime$ in \eqref{c'} equals the matrix $($one for both choices of sign $\pm)$ associated with the $\mathbb C$-linear morphism
\begin{align*}\nonumber\mathfrak X_{\mathbb G}^\pm\colon \ K_0(\mathbb G)\otimes_{\mathbb Z}\mathbb C& \to H^\bullet(\mathbb G;\mathbb C),\\
[E]&\mapsto \frac{1}{(2\pi)^2c^{\frac{1}{2}}}\widehat \Gamma^\pm(\mathbb G)\cup\operatorname{Ch}(E)
\end{align*}
 computed w.r.t.\ the basis $([E_1],\dots, [E_6])$ of $K_0(\mathbb G)$, and the Schubert basis $(T_0,T_1,T_2,T_3,T_4,T_5)$. In both cases $(\pm)$, the Stokes matrix $S^\prime$ coincides with the inverse of the Gram matrix of the Grothendieck--Euler--Poincar\'e pairing $ (\chi(E_i,E_j))_{i,j=1}^n$.
\end{Theorem}

The matrix $C^\prime$ in the above theorem is obtained by formula~\eqref{c'} for $G=A$ and $G=AB$ as follows:
\begin{itemize}\itemsep=0pt
\item the matrix representing $\mathfrak X_{\mathbb G}^-$ is equal to $C^\prime$ computed w.r.t.\ the solution $Y_0(z)\cdot A$, where $A\in\mathcal C_0(\eta,\mu, R)$ is
\begin{gather}\label{6luglio2017-10}
 A=\left(
\begin{matrix}
 1 & 0 & 0 & 0 & 0 & 0 \\
 2 {\rm i}\pi & 1 & 0 & 0 & 0 & 0 \\
 -2 \pi ^2 & 2 {\rm i}\pi & 1 & 0 & 0 & 0 \\
 -2 \pi ^2 & 2 {\rm i}\pi & 0 & 1 & 0 & 0 \\
 -\frac{1}{3} \big(8 {\rm i}\pi ^3\big) & -4 \pi ^2 & 2 {\rm i}\pi & 2 {\rm i}\pi & 1 & 0 \\
 \frac{4 \pi ^4}{3} & -\frac{1}{3} \big(8 {\rm i}\pi ^3\big) & -2 \pi ^2 & -2 \pi ^2 & 2 {\rm i}\pi & 1
\end{matrix}
\right);
\end{gather}

\item the matrix representing $\mathfrak X_{\mathbb G}^+$ is equal $C^\prime$ computed w.r.t.\ the solution $Y_0(z)\cdot A\cdot B$, where $B \in\mathcal C_0(\eta,\mu, R)$ is
\begin{gather}\label{6luglio2017-11}
B=\left(
\begin{matrix}
 1 & 0 & 0 & 0 & 0 & 0 \\
 -8 \gamma & 1 & 0 & 0 & 0 & 0 \\
 32 \gamma ^2 & -8 \gamma & 1 & 0 & 0 & 0 \\
 32 \gamma ^2 & -8 \gamma & 0 & 1 & 0 & 0 \\
 \frac{8}{3} \big(\zeta (3)-64 \gamma ^3\big) & 64 \gamma ^2 & -8 \gamma & -8 \gamma & 1 & 0 \vspace{1mm}\\
 \frac{64}{3} \big(16 \gamma ^4-\gamma \zeta (3)\big) & \frac{8}{3} \big(\zeta (3)-64 \gamma ^3\big) & 32 \gamma ^2 & 32 \gamma ^2 & -8 \gamma & 1
\end{matrix}
\right).
\end{gather}
\end{itemize}

A few more words about Theorem~\ref{6-07-17-4} are in order. The theorem provides the monodromy data at $t=0\in QH^\bullet(\mathbb G)$. It is worth stressing that the result acquires its geometrical significance because of Theorem~\ref{6-07-17-3} (i.e., Theorem \ref{mainisoth}), which guarantees that the data computed at the coalescence point $t=0$ are the monodromy data of the manifold in any chamber whose closure contains the point $t=0$. Without Theorem~\ref{6-07-17-3}, the results of the computations of Section~\ref{29novembre2016-1} would have little geometrical meaning for the theory of Frobenius manifolds.

From the data computed at $t=0$, the Stokes and the central connection matrices {\it at all other points of the small quantum cohomology} and/or w.r.t.\ other possible admissible lines $\ell$, will be computed in Section~\ref{recSC}, by a careful and pedagogical application of the action of the braid group.
They satisfy the same properties as in Theorem~\ref{6-07-17-4} w.r.t.\ other full exceptional $5$-block collections, obtained from $(E_1,\dots, E_6)$ by alternate mutation under the braids
\[\omega_1:=\beta_{12}\beta_{56},\qquad\omega_2:=\beta_{23}\beta_{45}\beta_{34}\beta_{23}\beta_{45},\qquad \hat{\omega}_1:=\beta_{12}\beta_{34}\beta_{56}.
\]
See Section \ref{recSC} for details.

It is important to remark that the Kapranov $5$-block exceptional collection itself appears {\it neither} at $t=0$ {\it nor} anywhere else along the locus of the small quantum cohomology. See Corollary~\ref{7febbraio2020-7}.

The monodromy data in any other chamber of $QH^\bullet(\mathbb G)$ can be obtained from the data~$S^\prime$,~$C^\prime$ computed at $0\in QH^\bullet(\mathbb G)$ (or from $PSP^{-1}$ and $CP^{-1}$), by the action~(\ref{2luglio2017-5}) of the braid group.

It is also worth noticing that the coalescence phenomenon should be taken into account in order to understand the appearance of the two exceptional collections $(E_1,E_2,E_3,E_4,E_5,E_6)$ and $(E_1,E_2,E_4, E_3,E_5,E_6)$ in Theorem~\ref{6-07-17-4}. In Section~\ref{isomonocoal} we will explain that, as a consequence of coalescence, a small neighbourhood of $0\in QH^\bullet(\mathbb G)$ is divided into two $\ell$-cells (see Definition~\ref{ellcell}). The passage from one exceptional collection to the other reflects the passage from one $\ell$-cell to the other.

\begin{Remark}In Theorem \ref{6-07-17-4}, we have two morphisms $\mathfrak X_{\mathbb G}^\pm$; the sign $(-)$ is the one taken in~\cite{dubro4}, whereas $(+)$ is the one taken in \cite{gamma1}.
\end{Remark}

\begin{Remark}\label{17-dec-2016}
Our explicit results suggest that the conjecture formulated in \cite{dubro4} and \cite{gamma1} requires some refinements, at least as far as the central connection matrix $C$ is concerned. Indeed, the connection matrix $C^\prime$ in Theorem \ref{6-07-17-4} belongs to the $\mathcal C_0(\eta,\mu, R)$-orbit of the connection matrix obtained from the topological-enumerative fundamental solution, but {\it is not } the connection matrix w.r.t.\ the topological-enumerative solution. These refinements are fully discussed in our paper \cite{CDG1}.
\end{Remark}

\begin{Remark}
In \cite{bay4} it is shown that the class of smooth projective varieties admitting generically semisimple quantum cohomology is closed w.r.t.\ the operation of blowing up at a finite number of points. Since this holds true also for the class of varieties which admit full exceptional collections in their derived categories, it is tempting to conjecture that the mentioned relationship between monodromy data and exceptional collections can be extended also for \emph{non-Fano} varieties. This is already suggested in \cite{bay4}. To the best of our knowledge, no explicit computations of the monodromy data have been done in the non-Fano case. The computations of the monodromy data for the $\frac{1}{2}K3$-surface, the rational elliptic surface obtained by blowing up 9 points in $\mathbb P^2$, could represent a significant step in this direction. This will represent a future research project of the authors.
\end{Remark}

Remarkably, our results suggest the validity of a constraint on the kind of exceptional collections associated {with} the monodromy data in a neighborhood of a semisimple coalescing point of the quantum cohomology $QH^\bullet(X)$ of a smooth projective variety $X$. If the eigenvalues $u_i$'s coalesce, at some semisimple point $t_0$, to $s< n$ values $\lambda_1,\dots,\lambda_s$ with multiplicities $p_1,\dots, p_s$ (with $p_1+\dots +p_s=n$, here $n$ is the sum of the Betti numbers of $X$), then the correspon\-ding monodromy data can be expressed in terms of Gram matrices and characteristic classes of objects of a \emph{full $s$-block exceptional collection}, i.e., a collection of the type
\[\mathcal E:=(\underbrace{E_1,\dots,E_{p_1}}_{\mathcal B_1},\underbrace{E_{p_1+1},\dots, E_{p_1+p_2}}_{\mathcal B_2},\dots, \underbrace{E_{p_1+\dots+p_{s-1}+1},\dots, E_{p_1+\dots+p_s}}_{\mathcal B_s}),\]
$ E_j\in\operatorname{Obj}\big(\mathcal D^b(X)\big)$,
where for each pair $(E_i, E_j)$ in a same block $\mathcal B_k$ the orthogonality conditions hold
\[\operatorname{Ext}^\ell(E_i, E_j)=0,\qquad \text{for any }\ell.
\]In particular, any reordering of the objects inside a single block $\mathcal B_j$ preserves the exceptiona\-li\-ty of~$\mathcal E$. More results about the nature of exceptional collections arising in this context and about their dispositions in the locus of small quantum cohomology for the class of complex Grassmannians will be explained in our paper~\cite{CDG1}.

\subsection{Plan of the paper} In Section \ref{sec1}, we review the analytic theory of Frobenius manifolds, their monodromy data and the isomonodromy theorems, according to~\cite{dubro1,dubro0,dubro2}. In particular, we characterise the freedom in the choice of the central connection matrix $C$, introducing the group $\mathcal C_0(\eta,\mu, R)$. We define a chamber-decomposition of the manifold, which depends on the choice of an oriented line $\ell$ in the complex plane: this is a natural structure related to the local invariance of the monodromy data (isomonodromy Theorems~\ref{isomonod1} and \ref{isoth2}), as well as of their discontinuous jumps from one chamber to another one, encoded in the action of the braid group, as a wall-crossing phenomenon.

In Section \ref{freedom} we review all freedoms and all other natural transformations on the monodromy data.

In Section \ref{isomonocoal} we extend the isomonodromy theorems and give a complete description of monodromy data in a neighborhood of semisimple coalescence points, specialising the result of \cite{CDG0} in Theorem \ref{mainisoth}.

 In Section \ref{A_3case}, we study the $A_3$ Frobenius manifold near the Maxwell stratum. We compute monodromy data both invoking Theorem~\ref{mainisoth} above and using oscillatory integrals. We compare the two approaches, so providing an explicit example of how Theorem~\ref{mainisoth} works. We also show how monodromy data mutate along a loop inside the Maxwell stratum.

In Section \ref{qcohg24}, we explicitly compute all monodromy data of the quantum cohomology of the Grassmannian $\mathbb G_2\big(\mathbb C^4\big)$. The result allows us to explicitly {verify} the conjecture of \cite{dubro0,dubro4} relating {the} monodromy data to characteristic classes of objects of an exceptional collection in~$\mathcal D^b\big(\mathbb G_2\big(\mathbb C^4\big)\big)$.

In Section \ref{topsolfano} we give an analytic characterisation of the {\it enumerative-topological solution}, in a~different way with respect to~\cite{gamma1}.

\section{Moduli of semisimple Frobenius manifolds}\label{sec1}
We denote with $\bigodot$ the symmetric tensor product of vector bundles, and with $(-)^\flat$ the standard operation of lowering the index of a $(1,k)$-tensor using a fixed inner product.

\begin{Definition}A \emph{Frobenius manifold} structure on a complex manifold $M$ of dimension $n$ is defined by giving
\begin{enumerate}\itemsep=0pt
\item[(FM1)] a symmetric {nondegenerate} $\mathcal O(M)$-bilinear tensor $\eta\in\Gamma\big(\bigodot^2T^*M\big)$, {called metric}, whose corresponding Levi-Civita connection $\nabla$ is flat;
\item[(FM2)] a $(1,2)$-tensor $c\in\Gamma\big(TM\otimes\bigodot^2T^*M\big)$ such that
\begin{itemize}\itemsep=0pt
\item the induced multiplication of vector fields $X\circ Y:=c(-,X,Y)$, for $X,Y\in\Gamma(TM)$, is \emph{associative},
\item $c^\flat\in\Gamma\big(\bigodot^3T^*M\big)$,
\item $\nabla c^\flat\in\Gamma\big(\bigodot^4T^*M\big)$;
\end{itemize}
\item[(FM3)] a vector field $e\in\Gamma(TM)$, called the \emph{unity vector field}, such that
\begin{itemize}\itemsep=0pt
\item the bundle morphism $c(-,e,-)\colon TM\to TM$ is the identity morphism,
\item $\nabla e=0$;
\end{itemize}
\item[(FM4)] a vector field $E\in\Gamma(TM)$, called the \emph{Euler vector field}, such that
\begin{itemize}\itemsep=0pt
\item $\mathfrak L_Ec=c$,
\item $\mathfrak L_E\eta=(2-d)\cdot \eta$, where $d\in\mathbb C$ is called the \emph{charge} of the Frobenius manifold.
\end{itemize}
\end{enumerate}
\end{Definition}

For simplicity it will be assumed that the tensor $\nabla E\in TM\otimes T^*M$ is diagonalizable.

Since the connection $\nabla$ is flat, there exist local flat coordinates that we denote $ \big(t^1,\dots ,t^n\big)$, w.r.t.\ which the metric $\eta$ is constant and the connection $\nabla$ coincides with the partial derivatives $\partial_\alpha=\partial/\partial t^\alpha$, $\alpha=1,\dots ,n$. Because of flatness and the conformal Killing condition, the Euler vector field is affine, i.e.,
\begin{gather*}\nabla\nabla E=0,\qquad\text{so that}\quad E=\sum_{\alpha=1}^n\big((1-q_\alpha)t^\alpha+r_\alpha\big)\frac{\partial}{\partial t^\alpha},\qquad q_\alpha,r_\alpha\in\mathbb C.
\end{gather*}

Following \cite{dubro1, dubro0, dubro2}, we choose flat coordinates so that $\frac{\partial}{\partial t^1}\equiv e$ and $r_\alpha\neq 0$ only if \mbox{$q_\alpha=1$} (this can always be done, up to an affine change of coordinates). In flat coordinates, let $\eta_{\alpha\beta}=\eta(\partial_\alpha,\partial_\beta)$, and $c_{\alpha\beta}^\gamma=c\big(dt^\gamma,\partial_\alpha,\partial_\beta\big)$, so that $\partial_\alpha\circ\partial_\beta=c_{\alpha\beta}^\gamma \partial_\gamma$. Condition (FM2) means that $c_{\alpha\beta\gamma}:=\eta_{\alpha\rho}c_{\beta\gamma}^\rho$ and $\partial_{\alpha}c_{\beta\gamma\delta}$ are symmetric in all indices. This implies the local existence of a~function~$F$ such that
\[c_{\alpha\beta\gamma}=\partial_\alpha\partial_\beta\partial_\gamma F.
\]
The associativity of the algebra is equivalent to the following conditions for~$F$, called WDVV-equations
\[{\partial_\alpha\partial_\beta\partial_\gamma}F \eta^{\gamma\delta}{\partial_\delta\partial_\epsilon\partial_\nu}F ={\partial_\nu\partial_\beta\partial_\gamma}F \eta^{\gamma\delta}{\partial_\delta\partial_\epsilon\partial_\alpha}F,
\]while axiom (FM4) is equivalent to
\[\eta_{\alpha\beta}=\partial_1\partial_\alpha\partial_\beta F,\qquad \mathfrak L_EF=(3-d)F+Q(t),
\]with $Q(t)$ a quadratic expression in $t_\alpha$'s. Conversely, given a solution of the WDVV equations, satisfying the quasi-homogeneity conditions above, a structure of Frobenius manifold is naturally defined on {an} open subset of the space of parameters $t^\alpha$'s.

Let us consider the canonical projection $\pi\colon\mathbb P^1_{\mathbb C}\times M\to M$, and the pull-back of the tangent bundle~$TM$:
\[\xymatrix{
\pi^*TM\ar[r]\ar[d]&TM\ar[d]\\
\mathbb P^1_{\mathbb C}\times M\ar[r]^{\quad\pi}&M.}
\]
We will denote by
\begin{enumerate}\itemsep=0pt
\item[1)] $\mathscr T_M$ the sheaf of sections of $TM$,
\item[2)] $\pi^*\mathscr T_M$ the pull-back sheaf, i.e., the sheaf of sections of~$\pi^*TM$,
\item[3)] $\pi^{-1}\mathscr T_M$ the sheaf of sections of~$\pi^*TM$ constant on the fibers of~$\pi$.
\end{enumerate} Introduce two $(1,1)$-tensors $\mathcal U$, $\mu$ on $M$ defined by
\begin{gather}
\label{26gen2020-1}\mathcal U(X):=E\circ X,\qquad\mu(X):=\frac{2-d}{2}X-\nabla_XE
\end{gather}
for all $X\in\Gamma(TM)$. In flat coordinates $(t^\alpha)_{\alpha=1}^n$ chosen as above, the operator~$\mu$ is \emph{constant} and in \emph{diagonal} form
\[
\mu=\text{ diag}(\mu_1,\dots ,\mu_n),\qquad \mu_\alpha=q_\alpha-\frac{d}{2}\in\mathbb C.
\]
All the tensors $\eta$, $e$, $c$, $E$, $\mathcal U$, $\mu$ can be lifted to $\pi^*TM$, and their lift will be denoted with the same symbol. So, also the Levi-Civita connection $\nabla$ is lifted on $\pi^*TM$, and it acts so that
\[\nabla_{\partial_z}Y=0\qquad \text{for }Y\in\big(\pi^{-1}\mathscr T_M\big)(M).
\]
Let us now twist this connection by using the multiplication of vectors and the operators~$\mathcal U$,~$\mu$.

\begin{Definition}
Let $\widehat{M}:=\mathbb C^*\times M$. The {\it deformed connection} $\widehat{\nabla}$ on the vector bundle $\pi^*TM|_{\widehat M}\allowbreak \to\widehat M$ is defined by
\begin{gather*} \widehat\nabla_XY=\nabla_XY+z\cdot X\circ Y,
\\
\widehat\nabla_{\partial_z}Y=\nabla_{\partial_z}Y+\mathcal U(Y)-\frac{1}{z}\mu(Y)
\end{gather*} for $X,Y\in(\pi^*\mathscr T_M)\big(\widehat M\big)$.
 \end{Definition}
 The crucial fact is that the deformed extended connection~$\widehat \nabla$ is \emph{flat}.

 \begin{Theorem}[\cite{dubro1,dubro2}]The flatness of $\widehat\nabla$ is equivalent to the following conditions on $M$
 \begin{itemize}\itemsep=0pt
 \item $\nabla c^\flat$ is completely symmetric,
 \item the product on each tangent space of $M$ is associative,
 \item $\nabla\nabla E=0$,
 \item $\mathfrak L_Ec=c$.
 \end{itemize}
 \end{Theorem}

 Because of this {integrability} condition, we can look for \emph{deformed flat coordinates} $\big(\tilde t^1,\dots ,\tilde t^n\big)$, with $\tilde t^\alpha=\tilde t^\alpha(t,z)$. These coordinates are defined by $n$ independent solutions of the equation
 \[
 \widehat\nabla {\rm d}\tilde t=0.
 \]Let $\xi$ denote a column vector of components of the differential ${\rm d}\tilde t$. The above equation becomes the linear system
 \begin{gather} \label{28luglio2016-4}
 \begin{cases}
 \partial_\alpha\xi=z\mathcal C_\alpha^{\rm T}(t)\xi,\\
 \partial_z\xi=\left(\mathcal U^{\rm T}(t)-\dfrac{1}{z}\mu^{\rm T}\right)\xi,
 \end{cases}
 \end{gather}
 where $\mathcal C_\alpha$ is the matrix $(\mathcal C_\alpha)^\beta_\gamma=c^\beta_{\alpha\gamma}$.
 We can rewrite the system in the form
 \begin{gather}\label{sistemaor}
 \begin{cases}
 \partial_\alpha\zeta=z\mathcal C_\alpha\zeta,\\
 \partial_z\zeta=\left(\mathcal U+\dfrac{1}{z}\mu\right)\zeta,
 \end{cases}
 \end{gather}
 where $\zeta:=\eta^{-1}\xi$.  In order to obtain~\eqref{sistemaor}, we have also used the invariance of the product, encoded in the relations
 \begin{gather}
 \eta^{-1}\mathcal C_\alpha^{\rm T}\eta=\mathcal C_\alpha,\nonumber\\
\label{usym} \mathcal U^{\rm T}\eta=\eta \mathcal U,
\end{gather}
and the $\eta$-skew-symmetry of $\mu$
 \begin{gather}\label{mantisym}
 \mu^{\rm T} \eta+\eta\mu=0.
\end{gather}
 Geometrically, $\zeta$ is the $\eta$-gradient of a deformed flat coordinate as in~\eqref{30nov2016-1}.
Monodromy data of system \eqref{sistemaor} define local invariants of the Frobenius manifold, as explained below.

\subsection[Spectrum of a Frobenius manifold and its monodromy data at $z=0$]{Spectrum of a Frobenius manifold and its monodromy data at $\boldsymbol{z=0}$} \label{3agosto2016-7}

 Let us fix a point $t$ of the Frobenius manifold $M$ and let us focus on the associated equation
\begin{gather}
\label{monoor1}\partial_z\zeta=\left(\mathcal U(t)+\frac{1}{z}\mu(t)\right)\zeta.
\end{gather}

\begin{Remark}\label{24novembre2016-1}If $\zeta_1$, $\zeta_2$ are {solutions} of the equation \eqref{monoor1}, then the two products
\[\langle\zeta_1,\zeta_2\rangle_\pm:=\zeta_1^{\rm T}\big({\rm e}^{\pm\pi{\rm i}}z\big)\eta\zeta_2(z)
\]
are independent of $z$. Indeed we have
\begin{align*}\partial_z\big(\zeta_1^{\rm T}\big({\rm e}^{\pm\pi{\rm i}}z\big)\eta\zeta_2(z)\big)&=\partial_z\big(\zeta_1^{\rm T}\big({\rm e}^{\pm\pi{\rm i}}z\big)\big)\eta\zeta_2(z)+\zeta_1^{\rm T}\big({\rm e}^{\pm\pi{\rm i}}z\big)\eta\partial_z\zeta_2(z)\\
&=\zeta_1^{\rm T}\big({\rm e}^{\pm\pi{\rm i}}z\big)\left[\eta \mathcal U-\mathcal U^{\rm T}\eta+\frac{1}{z}\left(\mu \eta+\eta\mu\right)\right]\zeta_2(z)\\
&=0\quad\text{ by \eqref{usym} and \eqref{mantisym}.}
\end{align*}
\end{Remark}

In order to give an intrinsic description of the structure of the \emph{normal forms} of solutions of equation \eqref{monoor1}, as well as a geometric characterization of the ambiguity and freedom up to which they are defined, we introduce the concept of the \emph{spectrum} of a Frobenius manifold (see also \cite{dubro2, dubro3}). Let $(V,\eta,\mu)$ be the datum of
\begin{itemize}\itemsep=0pt
\item an $n$-dimensional complex vector space $V$,
\item a bilinear symmetric non-degenerate form $\eta$ on $V$,
\item a diagonalizable endomorphism $\mu\colon V\to V$ which is $\eta$-antisymmetric
\[\eta(\mu a,b)+\eta(a,\mu b)=0\qquad\text{for any }a,b\in V.
\]
\end{itemize}
Let $\operatorname{spec}(\mu)=(\mu_1,\dots ,\mu_n)$ and let $V_{\mu_\alpha}$ be the eigenspace of a $\mu_\alpha$.

\begin{Definition}\label{munilp} Let $(V,\eta,\mu)$ as above. We say that an endomorphism $A\in\End(V)$ is $\mu$-\emph{nilpotent} if
\[AV_{\mu_\alpha}\subseteq\bigoplus_{m\in\mathbb Z_{\geq 1}}V_{\mu_\alpha+m}\quad\text{for any }\mu_\alpha\in\operatorname{spec}(\mu).
\]In particular such an operator is nilpotent in the usual sense. We can decompose a $\mu$-nilpotent operator $A$ in components $A_k$, $k\geq 1$, such that
\[A_kV_{\mu_\alpha}\subseteq V_{\mu_\alpha+k}\qquad\text{for any }\mu_\alpha\in\operatorname{spec}(\mu),
\]so that the following identities hold:
\[z^\mu Az^{-\mu}=A_1z+A_2z^2+A_3z^3+\cdots,\qquad [\mu,A_k]=kA_k\qquad\text{for }k=1,2,3,\dots.
\]
\end{Definition}

\begin{Definition}\label{3agosto2016-1} Let $(V,\eta,\mu)$ as above. Let us define on $V$ a second non-degenerate bilinear form $\left\{\cdot,\cdot\right\}$ by the equation
\[\left\{a,b\right\}:=\eta\big({\rm e}^{{\rm i}\pi\mu}a,b\big),\qquad\text{for all }a,b\in V.
\]
The set of all $\left\{\cdot,\cdot\right\}$-isometries $G\in\End(V)$ of the form
\[G=\mathbbm 1_{V}+\Delta,
\]with $\Delta$ a $\mu$-nilpotent operator, is a Lie group $\mathcal G(\eta,\mu)$, called $(\eta,\mu)$-\emph{parabolic orthogonal group}. Its Lie algebra $\mathfrak g(\eta,\mu)$ coincides with the set of all $\mu$-nilpotent operators $R$ which are also $\left\{\cdot,\cdot\right\}$-\emph{skew-symmetric} in the sense that
\[ \{Rx,y \}+ \{x,Ry \}=0.
\]In particular, any such matrix $R$ commutes with the operator ${\rm e}^{2\pi{\rm i}\mu}$.
\end{Definition}

The following result gives a description, in coordinates, of both $\mu$-nilpotents operators and elements of $\mathfrak g(\eta,\mu)$ and also describes some of their properties.

\begin{Lemma}\label{19.10.17-1}
Let $(V,\eta,\mu)$ as above, and let us fix a basis $(v_i)_{i=1}^n$ of eigenvectors of $\mu$.
\begin{enumerate}\itemsep=0pt
\item[$1.$] The operator $A\in\End(V)$ is $\mu$-nilpotent if and only if its associate matrix w.r.t.\ the basis $(v_i)_{i=1}^n$ satisfies the condition
\[(A)^\alpha_\beta=0\qquad \text{unless }\mu_\alpha-\mu_\beta\in\mathbb N^*.
\]
\item[$2.$] If $A\in\End(V)$ is a $\mu$-nilpotent operator, then the matrices associated with its components $(A_k)_{k\geq 1}$ w.r.t.\ the basis $(v_i)_{i=1}^n$ satisfy the condition
\begin{gather}\label{30nov2016-3}(A_k)^\alpha_\beta=0\qquad \text{unless }\mu_\alpha-\mu_\beta=k,\quad k\in\mathbb N^*.
\end{gather}
\item[$3.$] A $\mu$-nilpotent operator $A\in\End(V)$ is an element of $\mathfrak g(\eta,\mu)$ if and only if the matrices of its components $(A_k)_{k\geq 1}$ w.r.t.\ $(v_i)_{i=1}^n$ satisfy the further conditions
\begin{gather}\label{30nov2016-2}A_k^{\rm T}=(-1)^{k+1}\eta A_k\eta^{-1},\qquad k\geq 1.
\end{gather}
\item[$4.$] If $A\in\mathfrak g(\eta,\mu)$, then the following identity holds
\begin{gather}\label{relazioneRemu}z^{A^{\rm T}}\eta {\rm e}^{\pm {\rm i}\pi\mu}z^{A}=\eta {\rm e}^{\pm {\rm i}\pi\mu},
\end{gather}for any $z\in\mathbb C^*$.
\end{enumerate}
\end{Lemma}

\begin{proof}
The proof for points (1), (2), (3) can be found in \cite{dubro2}. For the identity~\eqref{relazioneRemu}, notice that~\eqref{30nov2016-2} implies
\[z^{A^{\rm T}}=\eta\big(z^{A_1-A_2+A_3-A_4+\cdots}\big)\eta^{-1}.
\]Moreover, from \eqref{30nov2016-3} we deduce that
\[{\rm e}^{\mp {\rm i}\pi\mu}A_k{\rm e}^{\pm {\rm i}\pi\mu}=(-1)^kA_k.
\]So, we conclude that
\begin{align*}
z^{A^{\rm T}}\eta {\rm e}^{\pm {\rm i}\pi\mu}z^A&=\eta\big(z^{A_1-A_2+A_3-A_4+\cdots}\big)\big({\rm e}^{\pm {\rm i}\pi\mu}z^A{\rm e}^{\mp {\rm i}\pi\mu}\big){\rm e}^{\pm {\rm i}\pi\mu}\\
&=\eta\big(z^{A_1-A_2+A_3-A_4+\cdots}\big)\big(z^{-A_1+A_2-A_3+A_4-\cdots}\big){\rm e}^{\pm {\rm i}\pi\mu}=\eta  {\rm e}^{\pm {\rm i}\pi\mu}.\tag*{\qed}
\end{align*}\renewcommand{\qed}{}
\end{proof}

The parabolic orthogonal group $\mathcal G(\eta,\mu)$ acts canonically on its Lie algebra $\mathfrak g(\eta,\mu)$ by the adjoint representation $\operatorname{Ad}\colon \mathcal G(\eta,\mu)\to \operatorname{Aut}(\mathfrak g(\eta,\mu))$:
\[\operatorname{Ad}_G(R):=G\cdot R\cdot G^{-1},\qquad\text{for all }G\in\mathcal G(\eta,\mu),\quad R\in\mathfrak g(\eta,\mu).
\]Such an action, in general is not free.

\begin{Definition}\label{29.05.17-1}
Let $R\in\mathfrak g(\eta,\mu)$. We define the group $\mathcal C_0(\eta,\mu, R)$ as the isotropy group of $R$ for the adjoint representation $\operatorname{Ad}\colon \mathcal G(\eta,\mu)\to \operatorname{Aut}(\mathfrak g(\eta,\mu))$.
\end{Definition}

The following Lemma can be easily directly proved from Definitions \ref{munilp}, \ref{3agosto2016-1}, \ref{29.05.17-1} and from results of Lemma \ref{19.10.17-1}.

\begin{Lemma}
Let $(V,\eta,\mu)$ a triple as above. If $G\in\mathcal G(\eta,\mu)$, and $R\in\mathfrak g(\eta,\mu)$ then
\[z^\mu z^{{\rm Ad}_GR}Gz^{-R} z^{-\mu}
\]
is an element of $\operatorname{End}(V)[z]$, i.e., it is polynomial in the indeterminate~$z$. Furthermore, the following is an equivalent characterization of the isotropy subgroup $\mathcal C_0(\eta,\mu, R)$:
\[\mathcal C_0(\eta,\mu, R)=\left\{
 G\in \operatorname{GL}(V)\colon \ \begin{aligned}
  & P_G(z):=z^\mu z^R Gz^{-R}z^{-\mu}\in\operatorname{End}(V)[z],\\
& P_G(0)\equiv\mathbbm 1_V, \\
 & \eta(P_G(-z)v_1,P_G(z)v_2)=\eta(v_1,v_2),\\
  &\text{for all }v_1,v_2\in V
 \end{aligned}
\right\}.
\]
\end{Lemma}

\begin{Definition}If $(V_1,\eta_1,\mu_1)$, $(V_2,\eta_2,\mu_2)$ are two triples as above, a morphism of triples $f\colon (V_1,\eta_1,\mu_1)\to (V_2,\eta_2,\mu_2)$, is the datum of a linear morphism $f\colon V_1\to V_2$, compatible with the metrics and the operators $\mu_1$, $\mu_2$, i.e.,
\begin{gather*}
\eta_1(v,w) =\eta_2(f(v),f(w)),\qquad v,w\in V_1,\qquad
\mu_2\circ f =f\circ\mu_1.
\end{gather*}
\end{Definition}

Given a Frobenius manifold $M$ (not necessarily semisimple), we can canonically associate to it an isomorphism class $[(V,\eta,\mu)]$ of triples as above, which will be called the \emph{spectrum} of $M$. Attached with any point $p\in M$, indeed, we have a triple $(T_pM,\eta_p, \mu_p)$. Given $p_1,p_2\in M$, the two triples are (non-canonically) isomorphic: using the Levi-Civita connection, for any path $\gamma\colon [0,1]\to M$ with $\gamma(0)=p_1$ and $\gamma(1)=p_2$, the parallel transport along $\gamma$ provides an isomoprhism of the triples at~$p_1$ and~$p_2$.

\begin{Definition}[\cite{dubro2,dubro3}]
A Frobenius manifold $M$ is called \emph{resonant} if $\mu_\alpha-\mu_\beta\in\mathbb Z\backslash\{0\}$ for some $\alpha\neq\beta$; otherwise, $M$ is called \emph{non-resonant}.
\end{Definition}

We can now give a complete (componentwise) description of normal forms of solutions of the system \eqref{monoor1}.

\begin{Theorem}[\cite{dubro1,dubro2}]\label{normalth}
Let $M$ be a Frobenius manifold (not necessarily semisimple). The system \eqref{monoor1} admits fundamental matrix solutions of the form
\begin{gather}\label{28luglio2016-9}Z(z,t)=\Phi(z,t)\cdot z^\mu z^{R(t)},\qquad \Phi(z,t) =\sum_{k\in\mathbb N}\Phi_k(t)z^k,\qquad\Phi_0(t)\equiv\mathbbm 1,\\ \label{30nov2016-4}\Phi(-z,t)^{\rm T}\cdot \eta\cdot \Phi(z,t) =\eta,
\end{gather}
where $\Phi_k\in\mathcal O(M)\otimes \mathfrak{gl}_n(\mathbb C)$, and $R\in\mathcal O(M)\otimes\mathfrak g(\eta,\mu)$. A solution of such a form will be said to be in \emph{Levelt normal form} at $z=0$.
\end{Theorem}

\begin{Remark}\label{ossmunondiag} In the general case, although not related to Frobenius manifolds, when $\mu$ is not diagonalizable and has a non-trivial nilpotent part, analogous results can be proved. However, the normal form becomes a little more complicated: e.g., it is no more defined by requiring that some entries of matrices~$R_k$ are nonzero, but that some \emph{blocks} are. For a detailed analysis of such case, we recommend the book by F.R.~Gantmacher~\cite{gant}.
\end{Remark}

Because of the Fuchsian character of the singularity $z=0$, the power series~$\Phi$ of Theorem~\ref{normalth} is convergent, and defines a genuine analytic solution. In general, solutions in Levelt normal form are not unique. As the following result shows, the \emph{freedom} in the choice of solutions in normal form are suitably quantified by the Lie groups $\mathcal G(\eta,\mu)$ and its isotropic subgroups $\mathcal C_0(\eta,\mu, R)$.

\begin{Theorem}[\cite{dubro1,dubro2}]\label{3agosto2016-5}
Let $M$ be a Frobenius manifold $($not necessarily semisimple$)$. Solutions of \eqref{monoor1} in normal form are not unique. Given two of them
\begin{gather*}
Z(z,t)=\Phi(z,t)\cdot z^\mu z^{R(t)},\qquad
\widetilde Z(z,t)=\widetilde\Phi(z,t)\cdot z^\mu z^{\widetilde R(t)},
\end{gather*}
there exists a unique holomorphic $ \mathcal G(\eta,\mu)$-valued function \begin{gather*}G(t)=\mathbbm 1+\Delta(t)\end{gather*} on $M$ such that
\begin{gather*} \widetilde Z(z,t)=Z(z,t)\cdot G(t),
\qquad
\widetilde R(t)=G(t)^{-1}\cdot R(t)\cdot G(t),\qquad \widetilde\Phi(z,t)=\Phi(z,t)\cdot P_G(z,t),
\end{gather*} where
\begin{gather*}P_G(z,t): =z^\mu\cdot G(t)\cdot z^{-\mu}
 =\mathbbm 1+z\Delta_1(t)+z^2\Delta_2(t)+\cdots,
\end{gather*}
$(\Delta_k)_{k\geq 1}$ being the components of~$\Delta$. In particular, if $\widetilde R= R$, then $G$ is $\mathcal C_0(\eta,\mu, R)$-valued.
\end{Theorem}

\begin{Remark}A first description of the freedom and ambiguities in the definition of the monodromy data was given in \cite{dubro1,dubro2}. In particular, a complex Lie group $\mathcal C_0(\mu, R)$ was introduced in order to describe the freedom of normal forms of solutions of~\eqref{monoor1}. Such a group \emph{is too big}, and in particular does not preserve the orthogonality condition~\eqref{30nov2016-4}. It must be replaced by $\mathcal C_0(\eta,\mu, R)$ of Definition~\ref{29.05.17-1}, which is the correct one.
\end{Remark}

For non-resonant Frobenius manifolds the corresponding $(\eta,\mu)$-parabolic orthogonal group $\mathcal G(\eta,\mu)$ together with all its subgroups $\mathcal C_0(\eta,\mu, R)$ are trivial. Since these groups are the responsible of a certain freedom in the choice of a normal form for solutions of~\eqref{monoor1} (according to Theorem~\ref{3agosto2016-5}), it follows that for non-resonant Frobenius manifolds such a choice is unique.

So far, we have focused on the system~\eqref{monoor1} at a fixed point of the manifold. Now let us vary the point $t$ in system \eqref{monoor1}, so that a fundamental solution $\Phi(z,t)z^\mu z^{R(t)}$, as in \eqref{28luglio2016-9}, depends on $t$. If instead of considering only the equation~\eqref{monoor1}, we focus on the whole system \eqref{sistemaor}, then the previous results can be further refined: namely, a $t$-independent choice for the exponent $R$ is allowed. Again, even for a fixed exponent $R$, solutions on normal forms are not unique, and they are parametrized by the isotropy group $\mathcal C_0(\eta,\mu, R)$.

\begin{Theorem}[isomonodromy Theorem I, \cite{dubro1,dubro2}]\label{isomonod1} Let $M$ be a Frobenius manifold $($not necessarily semisimple$)$.
\begin{enumerate}\itemsep=0pt
\item[$1.$] The system \eqref{sistemaor} admits fundamental matrix solutions of the form
\begin{gather*}Z(z,t) =\Phi(z,t)\cdot z^\mu z^{R},\\ \Phi(z,t)=\sum_{k\in\mathbb N}\Phi_k(t)z^k,\qquad\Phi_0(t) \equiv\mathbbm 1,\qquad \Phi(-z,t)^{\rm T}\cdot \eta\cdot \Phi(z,t)=\eta,
\end{gather*}
where $\Phi_k\in\mathcal O(M)\otimes \mathfrak{gl}_n(\mathbb C)$, and $R\in \mathfrak g(\eta,\mu)$ is \emph{independent} of~$t$. In particular the monodromy $M_0=\exp(2\pi{\rm i} \mu)\exp(2\pi{\rm i} R)$ at $z=0$ does not depend on~$t$.
\item[$2.$]
Solutions of the whole system \eqref{sistemaor} in normal form are not unique. Given two of them
\begin{gather*}
Z(z,t) =\Phi(z,t)\cdot z^\mu z^{R},\qquad
\widetilde Z(z,t)=\widetilde\Phi(z,t)\cdot z^\mu z^{\widetilde R},
\end{gather*}
there exists a unique matrix $G\in \mathcal G(\eta,\mu)$, say $G=\mathbbm 1+\Delta$, such that
\begin{gather*}\widetilde Z(z,t)=Z(z,t)\cdot G,
\qquad \widetilde R=G^{-1}\cdot R\cdot G,\qquad \widetilde\Phi(z,t)=\Phi(z,t)\cdot P_G(z,t),
\end{gather*} where
\begin{gather*}P_G(z,t): =z^\mu\cdot G\cdot z^{-\mu}
 =\mathbbm 1+z\Delta_1+z^2\Delta_2+\cdots,
\end{gather*}
$(\Delta_k)_{k\geq 1}$ being the components of $\Delta$. In particular, if $\widetilde R= R$, then $G\in\mathcal C_0(\eta,\mu, R)$.
\end{enumerate}
\end{Theorem}
\begin{proof} Let $Z(z,t)$ be a solution of \eqref{sistemaor}, and let $M_0(t)$ be the monodromy of $Z(\cdot, t)$ at $z=0$:
\[Z\big({\rm e}^{2\pi{\rm i}}z,t\big)=Z(z,t)\cdot M_0(t).
\]
The coefficients of the equations \begin{gather*}\partial_\alpha Z(z,t)=z\mathcal C_\alpha(t)\cdot Z(z,t),\qquad\alpha=1,\dots,n\end{gather*} being holomorphic in $z$, we have that
\begin{align*}
\partial_\alpha Z(z,t)\cdot Z(z,t)^{-1}&=\partial_\alpha Z\big({\rm e}^{2\pi{\rm i}}z,t\big)\cdot Z({\rm e}^{2\pi{\rm i}}z,t)^{-1}\\
&=\partial_\alpha  (Z(z,t)\cdot M_0(t) )\cdot ( Z(z,t)\cdot M_0(t) )^{-1}\\
&=\partial_\alpha Z(z,t)\cdot Z(z,t)^{-1} + Z(z,t)\cdot \partial_\alpha M_0(t)\cdot M_0(t)^{-1}\cdot Z(z,t)^{-1},
\end{align*}for any $\alpha$. Hence
\[\partial_\alpha M_0(t)=0,\qquad \alpha=1,\dots, n.
\]
By Theorem \ref{3agosto2016-5}, we necessarily conclude that $R$ is $t$-independent.
\end{proof}

\begin{Definition}[\cite{dubro1,dubro2}]Given a Frobenius manifold $M$, we will call \emph{monodromy data of $M$ at $z=0$} the data $(\mu, [R])$, where $[R]$ denotes the $\mathcal G(\eta,\mu)$-class of exponents of formal solutions in Levelt normal form of the system \eqref{sistemaor} as in Theorem \ref{normalth}. According to Theorem~\ref{isomonod1}, a~representative $R$ can be chosen independent of the point $t\in M$.
\end{Definition}

We conclude this section with a result giving sufficient conditions on solutions of the system~\eqref{sistemaor} for \emph{resonant} Frobenius manifolds in order that they satisfy the $\eta$-orthogonality condition~\eqref{30nov2016-4}. In its essence, this result is stated and proved in~\cite{gamma1}, in the specific case of quantum cohomologies of Fano manifolds.

\begin{Proposition}\label{6marzo2017}
Let $M$ be a \emph{resonant} Frobenius manifold, and $t_0\in M$ a fixed point.
\begin{enumerate}\itemsep=0pt
\item[$1.$] Suppose that there exists a fundamental solution of~\eqref{sistemaor} of the form
\[Z(z,t)=\Phi(z,t)z^\mu z^R,\qquad \Phi(t)=\mathbbm 1+\sum_{j=1}^\infty \Phi_j(t)z^j,
\]with $R$ satisfying all the properties of the Theorem~{\rm \ref{normalth}}, such that
\[H(z):=z^{-\mu} \Phi(z,t_0)z^\mu
\]is a holomorphic function at $z=0$ and $H(0)\equiv \mathbbm 1$. Then $\Phi (z,t)$ satisfies the constraint \begin{gather*}\Phi(-z,t)^{\rm T}\eta \Phi(z,t)=\eta\end{gather*} for all points $t\in M$.

\item[$2.$] If a solution with the properties above exists, then it is unique.
\end{enumerate}
\end{Proposition}

\begin{proof} From Remark \ref{24novembre2016-1}, we already know that the following bracket must be independent of~$z$:
\begin{align*}\langle Z(z,t_0), Z(z,t_0)\rangle_+&= \big(\Phi(-z,t_0)\big({\rm e}^{{\rm i}\pi}z\big)^\mu \big({\rm e}^{{\rm i}\pi}z\big)^R\big)^{\rm T}\eta\big(\Phi(z,t_0)z^\mu z^R\big)\\
&= \big(\big({\rm e}^{{\rm i}\pi}z\big)^\mu H(-z) \big({\rm e}^{{\rm i}\pi}z\big)^R\big)^{\rm T}\eta\big(z^\mu H(z)z^R\big)\\
&={\rm e}^{{\rm i}\pi R^{\rm T}}z^{R^{\rm T}}H(-z)^{\rm T}{\rm e}^{{\rm i}\pi\mu}z^\mu\eta z^\mu H(z)z^R\\
&={\rm e}^{{\rm i}\pi R^{\rm T}}z^{R^{\rm T}}H(-z)^{\rm T}{\rm e}^{{\rm i}\pi\mu}\eta H(z)z^R.
\end{align*}
By taking the first term of the Taylor expansion in $z$ of the r.h.s., and using~\eqref{relazioneRemu}, we get
\[\langle Z(z,t_0), Z(z,t_0)\rangle_+={\rm e}^{{\rm i}\pi R^{\rm T}}{\rm e}^{{\rm i}\pi\mu}\eta.
\]So, using again the equation $z^{\mu^{\rm T}}\eta z^\mu=\eta$ and \eqref{relazioneRemu}, we can conclude that
\[\Phi(-z,t_0)^{\rm T}\eta \Phi(z,t_0)=\big(\big({\rm e}^{{\rm i}\pi}z\big)^\mu\big({\rm e}^{{\rm i}\pi}z\big)^R \big)^{-{\rm T}}\langle Z(z,t_0), Z(z,t_0)\rangle_+\big(z^\mu z^R\big)^{-1}=\eta.
\]Because of \eqref{sistemaor} and the property of $\eta$-compatibility of the Frobenius product, we have that
\[\frac{\partial}{\partial t^\alpha}\big(\Phi(-z,t)^{\rm T}\eta \Phi(z,t)\big)=z\cdot \Phi(-z,t)^{\rm T}\cdot\big(\eta\mathcal C_\alpha-\mathcal C_\alpha^{\rm T}\eta\big)\cdot\Phi(z,t)=0.
\]This concludes the proof of (1). Let us now suppose that there are two solutions
\[\Phi_1(z,t)z^\mu z^R,\qquad \Phi_2(z,t)z^\mu z^R
\]such that
\begin{align}\label{10-06-16}z^{-\mu}\Phi_1(z,t_0)z^\mu&=\mathbbm 1+z K_1+z^2 K_2+\cdots,\\
\label{10-06-16.bis}z^{-\mu}\Phi_2(z,t_0)z^\mu&=\mathbbm 1+z K'_1+z^2 K'_2+\cdots.
\end{align}The two solutions must be related by
\[ \Phi_2(z,t)z^\mu z^R=\Phi_1(z,t)z^\mu z^R\cdot C
\]for some matrix $C\in\mathcal C_0(\eta,\mu, R)$. This implies that $\Phi_2(z,t)=\Phi_1(z,t)\cdot P(z)$, where $P(z)$ is a~matrix valued polynomial of the form
\begin{gather*}
P(z)=\mathbbm 1+z\Delta_1+z^2\Delta_2+\cdots,\qquad\text{with} \  (\Delta_k)^\alpha_\beta=0 \ \text{unless} \ \mu_\alpha-\mu_\beta=k, \ \text{and} \ P(1)\equiv C.
\end{gather*}
We thus have $z^{-\mu}\Phi_1^{-1}\Phi_2 z^\mu=z^{-\mu}P(z)z^\mu$, and
\[\left(z^{-\mu}P(z)z^\mu\right)^\alpha_\beta=\delta^\alpha_\beta+\sum_{k}\left(\Delta_k\right)^\alpha_\beta z^{k-\mu_\alpha+\mu_\beta}=\delta^\alpha_\beta+\sum_{k}\left(\Delta_k\right)^\alpha_\beta\equiv C.
\]Then, from 
\eqref{10-06-16}, \eqref{10-06-16.bis} it immediately follows that $C=\mathbbm 1$, which proves that \mbox{$\Phi_1=\Phi_2$}.
\end{proof}

\subsection{Semisimple Frobenius manifolds}\label{semisimplesection}
\begin{Definition}A finite dimensional commutative and associative $\mathbb K$-algebra $A$ with unit is called \emph{semisimple} if there is no nilpotent element, i.e., an element $a\in A\setminus\left\{0\right\}$ such that $a^k=0$ for some $k\in\mathbb N$.
\end{Definition}

In what follows we will always assume that the ground field is $\mathbb C$.

\begin{Lemma}\label{ssfrobalg}Let $A$ be a $\mathbb C$-Frobenius algebra of dimension $n$. The following are equivalent:
\begin{enumerate}\itemsep=0pt
\item[$1)$] $A$ is semisimple;
\item[$2)$] $A$ is isomorphic to $\mathbb C^{\oplus n}$;
\item[$3)$] $A$ has a basis of idempotents, i.e., elements $\pi_1,\dots,\pi_n$ such that
\[\pi_i\circ \pi_j=\delta_{ij}\pi_i,
\qquad \eta(\pi_i,\pi_j)=\eta_{ii}\delta_{ij},
\]for a suitable non-degenerate multiplication invariant pairing~$\eta$ on~$A$;
\item[$4)$] there is a vector $\mathcal E\in A$ such that the multiplication operator $\mathcal E\circ \colon A\to A$ has $n$ pairwise distinct eigenvalues.
\end{enumerate}
\end{Lemma}

\begin{proof}
All these equivalences are well known. For the equivalence of (1) and (2), see for example \cite[Chapter~V, Section~6, Proposition 5, p.~A.V.34]{boualg4-7}. Another elementary proof can be found in the lectures notes~\cite{dubro2}. The fact that (2) and (3) are equivalent is trivial. Let us prove that~(3) and (4) are equivalent. If~(3) holds it is sufficient just to take
\[\mathcal E=\sum_kk\pi_k.
\]
So $\mathcal E\circ$ has spectrum $\left\{1,\dots,n\right\}$. Let us now suppose that~(4) holds. The commutativity of the algebra implies that all operators $a\circ\colon A\to A$ commute. Therefore they preserve the one dimensional eigenspaces of $\mathcal E\circ$ and thus they are all diagonalizable. It follows that all operators~$a\circ$, being commuting and diagonalizable linear operators on a finite dimensional vector space, are simultaneously diagonalizable. Thus idempotents are constructed by suitably re\-sca\-ling the eigenvectors of $\mathcal E\circ$.
\end{proof}

\begin{Definition}[semisimple Frobenius manifolds] A point $p$ of a Frobenius manifold $M$ is \emph{semisimple} if the corresponding Frobenius algebra $T_pM$ is semisimple. If there is an open dense subset $M_{ss}\subset M$ of semisimple points, then $M$ is called a \emph{semisimple Frobenius manifold}.
\end{Definition}

 It is evident from point (4) of Lemma~\ref{ssfrobalg} that semisimplicity is an \emph{open property}: if $p$ is semisimple, then all points in a neighborhood of $p$ are semisimple.

\begin{Definition}[caustic and bifurcation set]\label{caustic}Let $M$ be a semisimple Frobenius manifold. We call \emph{caustic} the set
\[\mathcal K_M:=M\setminus M_{ss}= \{p\in M\colon T_pM\text{ is not a semisimple Frobenius algebra}\}.
\]
We call \emph{bifurcation set} of the Frobenius manifold the set
\[
\mathcal B_M:= \{p\in M\colon \operatorname{spec}\left(E\circ_p\colon T_pM\to T_pM\right)\text{ is not simple} \}.
\]
By Lemma \ref{ssfrobalg}, we have $\mathcal{K}_M\subseteq\mathcal B_M$. Semisimple points in $\mathcal B_M\setminus\mathcal K_M$ are called \emph{ semisimple coalescence points}.
\end{Definition}

The bifurcation set $\mathcal B_M$ and the caustic $\mathcal K_M$ are either empty or {a} hypersurface {(in general a singular one)}, invariant w.r.t.\ the unit vector field~$e$ (see~\cite{hertling}). For Frobenius manifolds defined on the base space of semiuniversal unfoldings of a singularity, these sets coincide with the \emph{bifurcation diagram} and the \emph{caustic} as defined in the classical setting of singularity theory \cite{arnold1, singularity1}. In this context, the set $\mathcal B_M\setminus\mathcal K_M$ is called \emph{Maxwell stratum}. Remarkably, all these subsets typically admit a naturally induced Frobenius submanifold structure~\cite{strachan, strachan1}. In what follows we will assume that the semisimple Frobenius manifold $M$ admits nonempty bifurcation set $\mathcal B_M$, and set of semisimple coalescence points $\mathcal B_M\setminus\mathcal K_M$.

At each point $p$ in the open dense semisimple subset $M_{ss}\subseteq M$, there are $n$ idempotent vectors
 \begin{gather*}
 \pi_1(p),\dots,\pi_n(p)\ \in T_pM,
 \end{gather*}
 unique up to a permutation. By Lemma~\ref{ssfrobalg} there exists a suitable local vector field $\mathcal E$ such that $\pi_1(p),\dots,\pi_n(p)$ are eigenvectors of the multiplication $\mathcal E\circ$, with simple spectrum at $p$ and consequently in a whole neighborhood of~$p$. Using the results exposed in \cite{kato} about analytic deformation of operators with simple spectrum w.r.t.\ one complex parameter, in particular the results stating analyticity of eigenvectors and eigenprojections, and extending them to the case of more parameters using Hartogs' theorem, we deduce the following

\begin{Lemma}\label{4luglio2017-1} The idempotent vector fields are holomorphic at a semisimple point $p$, in the sense that, chosen an ordering $\pi_1(p),\dots,\pi_n(p)$, there exist a neighborhood of $p$ where the resulting local vector fields are holomorphic.
\end{Lemma}

Notice that, although the idempotents are defined (and unique up to a permutation) at each point of $M_{ss}$, it is not true that there exist $n$ globally well-defined holomorphic idempotent vector fields. Indeed, the caustic $\mathcal K_M$ is in general a locus of algebraic branch points: if we consider a semisimple point $p$ and a close loop $\gamma\colon [0,1]\to M$, with base point $p$, encircling $\mathcal K_M$, along which a coherent ordering is chosen, then
 \[\big(\pi_1(\gamma(0)),\dots,\pi_n(\gamma(0))\big)\qquad\text{and}\qquad \big(\pi_1(\gamma(1)),\dots,\pi_n(\gamma(1))\big)
\]may differ by a permutation. Thus, the idempotent vector fields are holomorphic and single-valued on simply connected open subsets not containing points of the caustic.

 \begin{Remark}\label{26.05.17-3} More generally, under the assumption $\mathcal K_M\neq \varnothing$, the {idempotent} vector fields define single-valued and holomorphic local sections of the tangent bundle $TM$ on any connected open set $\Omega\subseteq M\setminus\mathcal K_M=M_{ss}$ satisfying the following property: for any $z\in \Omega$ the inclusions
\begin{gather*}
\Omega\xhookrightarrow{\quad\alpha\quad} M_{ss} \xhookrightarrow{\quad\beta\quad} M
\end{gather*}
 induce morphisms in homotopy
\begin{gather*}
\xymatrix{\pi_1(\Omega,z)\ar[r]^{\alpha_*\ }& \pi_1(M_{ss},z)\ar[r]^{\ \beta_* }&\pi_1(M,z)}
\end{gather*}
such that $\operatorname{im}(\alpha_*)\cap\ker(\beta_*)=\left\{0\right\}$. Moreover, this means that the structure group of the tangent bundle of $M_{ss}$ is reduced to the symmetric group $\mathfrak S_n$, and that the \emph{local} isomorphism of $\mathcal O_{M_{ss}}$-algebras
\[\mathscr T_{M_{ss}}\underset{\rm loc}{\cong}\mathcal O_{M_{ss}}^{\oplus n},
\]existing everywhere, can be replaced by a global one by considering a Frobenius structure prolonged to an unramified covering of degree at most $n!$ (see \cite{manin}).
\end{Remark}

\begin{Theorem}[\cite{dubronapoli,dubro1,dubro2}]\label{coorddubr}
Let $p\in M_{ss}$ be a semisimple point, and $(\pi_i(p))_{i=1}^n$ a basis of idempotents in $T_pM$. Then \[[\pi_i,\pi_j]=0;\]as a consequence there exist local coordinates $u_1,\dots, u_n$ such that
\[\pi_i=\frac{\partial}{\partial u_i}.
\]
\end{Theorem}

\begin{Definition}[canonical coordinates \cite{dubro1,dubro2}] Let $M$ a Frobenius manifold and $p\in M$ a~semisimple point. The coordinates defined in a neighborhood of~$p$ of Theorem~\ref{coorddubr} are called \emph{canonical coordinates}.
\end{Definition}

Canonical coordinates are defined only up to permutations and shifts. They are holomorphic local coordinates in a simply connected neighbourhood of a semisimple point not containing points of the caustic~$\mathcal{K}_M$, or more generally on domains with the property of Remark~\ref{26.05.17-3}. Holomorphy holds also at semisimple coalescence points.

\begin{Theorem}[\cite{dubro2}] If $u_1,\dots,u_n$ are canonical coordinates near a semisimple point of a~Frobenius manifold $M$, then $($up to shifts$)$ the following relations hold
\[\frac{\partial}{\partial u_i}\circ\frac{\partial}{\partial u_i}=\delta_{ij}\frac{\partial}{\partial u_i}, \qquad e=\sum_{i=1}^n\frac{\partial}{\partial u_i}, \qquad
E=\sum_{i=1}^nu_i\frac{\partial}{\partial u_i}.
\]
\end{Theorem}

In this paper we will fix the shifts of canonical coordinates so that they coincide with the eigenvalues of the $(1,1)$-tensor $E\circ$.

\begin{Definition}[matrix $\Psi$]Let $M$ be a semisimple Frobenius manifold, $t^1,\dots,t^n$ be local flat coordinates such that $\frac{\partial}{\partial t^1}=e$ and $u_1,\dots, u_n$ be canonical coordinates. Introducing the orthonormal basis
\begin{gather}\label{5luglio2017-1}
f_i:=\frac{1}{\eta\big(\frac{\partial}{\partial u_i},\frac{\partial}{\partial u_i}\big)^{\frac{1}{2}}}\frac{\partial}{\partial u_i}
\end{gather}
 for {an} arbitrary {choice} of signs in the square roots, we define a matrix $\Psi$ (depending on the point of the Frobenius manifold) whose elements $\Psi_{i\alpha}$ ($i$-th row, $\alpha$-th column) are defined by the relation
\[\frac{\partial}{\partial t^\alpha}=\sum_{i=1}^n\Psi_{i\alpha}f_i,\qquad\alpha=1,\dots,n.
\]
\end{Definition}

\begin{Lemma}\label{proppsi}The matrix $\Psi$ is a single-valued holomorphic function on any simply connected open subset not containing points of the caustic $\mathcal{K}_M$, or more generally on any open domain $\Omega$ as in Remark~{\rm \ref{26.05.17-3}}. Moreover, it satisfies the following relations:
\begin{gather*} \Psi^{\rm T}\Psi=\eta,\qquad
\Psi_{i1}=\eta\left(\frac{\partial}{\partial u_i},\frac{\partial}{\partial u_i}\right)^{\frac{1}{2}},\\
\frac{\partial}{\partial u_i}=\Psi_{i1}\sum_{\alpha,\beta=1}^n\Psi_{i\alpha}\eta^{\alpha\beta}\frac{\partial}{\partial t^\beta},
\qquad
c_{\alpha\beta\gamma}=\sum_{i=1}^n\frac{\Psi_{i\alpha}\Psi_{i\beta}\Psi_{i\gamma}}{\Psi_{i1}}.
\end{gather*}
If $\mathcal U$ is the operator of multiplication by the Euler vector field, then $\Psi$ diagonalizes it:
\[\Psi\mathcal U\Psi^{-1}=U:=\operatorname{diag}(u_1,\dots, u_n).
\]
\end{Lemma}
\begin{proof} The first assertion is a direct consequence of the analogous property of the idempotents vector fields, as in Lemma~\ref{4luglio2017-1}. All the other relations follow by computations (see~\cite{dubro2}, notice a misprint in formula~(3.16) there, where~$f_i$ must be replaced by~$\pi_i$). \end{proof}

We stress that $\Psi$ and the coordinates $u_i$'s are holomorphic also at semisimple coalescence points, due to the same property of the idempotents.

\subsection{Monodromy data for a semisimple Frobenius manifold}\label{monossfm}
Monodromy data at $z=\infty$ are defined in \cite{dubro1,dubro0,dubro2} at {a} point of a semisimple Frobenius manifold \emph{not belonging to the bifurcation set}. In the present section we review these issues, and we enlarge the definition to \emph{all semisimple points}, including the bifurcation ones, namely the semisimple coalescence points of Definition~\ref{3agosto2016-9}.

In this section, we fix an open subset $\Omega\subseteq M_{ss}$ satisfying the property of Remark~\ref{26.05.17-3}, so that we can choose and fix on~$\Omega$
\begin{itemize}\itemsep=0pt
\item an ordering {of} idempotent vector fields and canonical local coordinates $p\mapsto u(p)$, $p\in\Omega$,
\item a {choice of} the square roots in the definition of normalized idempotent vector fields $f_i$'s, and hence a determination of the matrix $\Psi$.
\end{itemize}
In this way, system \eqref{sistemaor} and system \eqref{idemsist} below, are determined.
In the idempotent frame
\begin{gather}\label{30nov2016-5}
y=\Psi\zeta,
\end{gather}
system \eqref{sistemaor} becomes
\begin{gather}\label{idemsist}
\begin{cases}
\partial_iy=(zE_i+V_i)y,\\
\partial_zy=\left(U+\dfrac{1}{z}V\right)y,
\end{cases}
\end{gather}
where $(E_i)^\alpha_\beta=\delta^\alpha_i\delta^\beta_i$ and
\begin{gather}\label{4luglio2017-2}
V:=\Psi\mu\Psi^{-1},\qquad V_i:=\partial_i\Psi\cdot\Psi^{-1},
\\
U:=\Psi\mathcal U\Psi^{-1}=\operatorname{diag}(u_1,\dots,u_n),\nonumber
\end{gather}
with not necessarily $u_i\neq u_j $ when $ i\neq j$. By Lemma~\ref{proppsi}, $\Psi(u)$, $V(u)$ and $V_i(u)$'s are holomorphic on~$\Omega$.

\begin{Lemma}\label{30nov2016-9}The matrix $V=\Psi\mu\Psi^{-1}$ is antisymmetric, i.e., $V^{\rm T}+V=0$. Moreover,
\[
\text{if $u_i=u_j$, then } V_{ij}=V_{ji}=0.
\]
\end{Lemma}

\begin{proof}Antisymmetry is an easy consequence of \eqref{mantisym} and the $\eta$-orthogonality of $\Psi$ (see~\cite{dubro2}). Moreover, compatibility conditions of the system \eqref{idemsist} imply that
\[[E_i,V]=[V_i,U].
\]Reading this equation for entries at place $(i,j)$, we find that
\[V_{ij}=(u_j-u_i)(V_i)_{ij}.
\]
Now, $(V_i)_{ij}$ is holomorphic, by Lemma~\ref{proppsi} and~(\ref{4luglio2017-2}), so that if $i\neq j$, but $u_i=u_j$, then $V_{ij}=0$.
\end{proof}

We focus on the second linear system
\begin{gather}\label{semiseq}\partial_zy=\left(U+\frac{1}{z}V\right)y,
\end{gather}
and study it at {\it a fixed} point $p\in\Omega$.

\begin{Theorem}\label{formalredu} Let $\Omega\subseteq M_{ss}$ as in Remark~{\rm \ref{26.05.17-3}}. At a $($fixed$)$ point $p\in \Omega$, there exists a unique formal $($in general divergent$)$ series
\[F(z):=\mathbbm 1+\sum_{k=1}^\infty\frac{A_k}{z^k}\]
with
\[F^{\rm T}(-z)F(z)=\mathbbm 1,
\]such that the transformation $\tilde{y}=F(z)y$
 reduces the corresponding system~\eqref{semiseq} at $p$ to the one with constant coefficients
\[\partial_z\tilde y=U\tilde y.
\]
Hence, system \eqref{semiseq} has a unique {\rm formal solution}
\begin{gather}\label{28luglio2016-7}
Y_{\rm formal}(z)=G(z){\rm e}^{zU},\qquad G(z):=F(z)^{-1}=\mathbbm 1+\sum_{k=1}^\infty\frac{G_k}{z^k}.
\end{gather}
\end{Theorem}

\begin{proof}By a direct substitution, one finds the following recursive equations for the coefficients~$A_k$:
\[
[U,A_1]=V,\qquad [U,A_{k+1}]=A_kV-kA_k,\qquad k=1,2,\dots.
\]
If $(i,j)$ is such that $u_i\neq u_j$ then we can determine $(A_{k+1})^i_j$ by the second equation in terms of entries of~$A_k$; if $u_i=u_j$ then we can determine $(A_{k+1})^i_j$ from the successive equation:
\[[U,A_{k+2}]=A_{k+1}V-(k+1)A_{k+1}.
\]
{Indeed}, the $(i,j)$-entry of the l.h.s.\ is~0 and, by Lemma~\ref{30nov2016-9} $(A_{k+1}V)^i_j$ is a~linear combination of {already determined entries} $(A_{k+1})^i_h$, with $u_i\neq u_h$. In such a way we can construct~$F(z)$. Let us now prove that $F^{\rm T}(-z)F(z)=\mathbbm 1$.
Let us take any (formal or analytic) solution $Y$ of the original system, and pose
\[A:=Y\big({\rm e}^{-{\rm i}\pi}z\big)^{\rm T}Y(z).
\]$A$ is a constant matrix, since it does not depend on~$z$. Thus, for an appropriate constant matrix~$C$ we have
\[F(z)Y(z)={\rm e}^{zU}C,
\]from which we deduce that
\[F(z)^{-1}=Y(z)C^{-1}{\rm e}^{-zU},\qquad
F(-z)^{-{\rm T}}={\rm e}^{zU}C^{-{\rm T}}Y({\rm e}^{-{\rm i}\pi}z)^{\rm T}.
\]So
\[F(-z)^{-{\rm T}}F(z)^{-1}={\rm e}^{zU}C^{-{\rm T}}AC^{-1}{\rm e}^{-zU}.
\]
Comparing the constant terms of the expansion of the r.h.s and the l.h.s.\ we conclude that $C^{-{\rm T}}AC^{-1}=\mathbbm1$.
\end{proof}

Notice in the above proof that {the equation} $[U,A_{k+1}]=A_kV-kA_k$, that is $(u_i-u_j)(A_{k+1})_j^i=(A_kV-kA_k)_j^i$, implies that, if we let $p$ vary in $\Omega$ then the $G_k$'s define holomorphic matrix valued functions $G_k(u)$ at points $u$, lying in $u(\Omega)$, such that $u_i\neq u_j$ for $i\neq j$. Accordingly, the formal matrix solution
\begin{gather}
\label{7luglio2017-1}
Y_{\rm formal}(z,u)=G(z,u){\rm e}^{zU},\qquad G(z,u)=\mathbbm 1+\sum_{k=1}^\infty\frac{G_k(u)}{z^k},
\end{gather}
 is well defined and holomorphic w.r.t.\ $u=u(p)$ away from semisimple coalescence points in $\Omega$. In Theorem~\ref{mainisoth} below, we will show that $Y_{\rm formal}(z,u)$ extends holomorphically also at semisimple coalescence points.

\begin{Remark}The proof of Theorem \ref{formalredu} is based on a simple computation, which holds both at a coalescence and a non-coalescence semisimple point. The statement can also be deduced from the more general results of~\cite{BJL2} (see also~\cite{CDG0}). A similar computation can be found also in \cite{tele} and \cite{gamma1}. Notice however that this computation does not provide any information about the analiticity of $G(u)$ in case of coalescence $u_i\to u_j$, $i\neq j$. The analiticity of $Y_{\rm formal}(z,u)$ -- and of actual fundamental solutions -- at a semisimple coalescence point follows from the results proved in \cite{CDG0}, and will be the content of Theorem \ref{mainisoth} below.
 \end{Remark}

In order to study actual solutions at $p\in\Omega$, we introduce Stokes rays. In what follows, we denote by ${\rm pr}\colon \mathcal R \to \mathbb{C} \backslash \{0\}$ the covering map. For pairs $(u_i, u_j)$ such that $u_i\neq u_j$, we {locally choose arguments}~$\alpha_{ij}$ of $\arg(u_i-u_j)$ {within} the interval $[0;2\pi{[}$, and we let
\[\tau_{ij}:=\frac{3\pi}{2}-\alpha_{ij}.
\]

\begin{Definition}[Stokes rays]We call \emph{Stokes rays} of the system \eqref{semiseq} the rays in the universal covering $\mathcal R$ defined by
\[R_{ij,k}:=\left\{z\in\mathcal R\colon \arg z=\tau_{ij}+2k\pi\right\},\qquad k\in\mathbb Z.
\]
\end{Definition}

The characterisation of Stokes rays is as follows:
 $z\in R_{ij,k}$ if and only if
\[\operatorname{Re}((u_i-u_j)z)=0,\qquad
\operatorname{Im}((u_i-u_j)z)<0,\qquad z\in \mathcal{R}.
\]
For given $1\leq i\neq j \leq n$, the projection of the rays $R_{ij,k}$, $k\in \mathbb{Z}$, on the $\mathbb C$-plane
\[
R_{ij}:=\operatorname{pr} (R_{ij,k} )
\]
 {does not depend on $k$} and is also called a \emph{Stokes ray}. It coincides with the ray defined in~\cite{dubro2}, namely
\begin{gather}\label{28luglio2016-13}
R_{ij}= \{z\in\mathbb C\colon z=-{\rm i}\rho(\overline{u_i}-\overline{u_j}),\ \rho>0 \}.
\end{gather}
 Stokes rays have a natural orientation from 0 to $\infty$. For $z\in \mathbb C$ we have
\begin{gather*}
\big|{\rm e}^{zu_i}\big| =\big|{\rm e}^{zu_j}\big|\qquad\text{if }z\in R_{ij},\\
\big|{\rm e}^{zu_i}\big| >\big|{\rm e}^{zu_j}\big|\qquad\text{if }z\ \text{is on the left of }R_{ij},\\
\big|{\rm e}^{zu_i}\big| <\big|{\rm e}^{zu_j}\big|\qquad\text{if }z\ \text{is on the right of }R_{ij}.
\end{gather*}

\begin{Definition}[admissible rays and line]\label{5luglio2017-7}
Let $\phi\in\mathbb R$ and let us define the rays in $\mathcal{R}$
\begin{gather*}
\ell_+(\phi):= \{z\in\mathcal R\colon \arg z=\phi \},\\
\ell_-(\phi):= \{z\in\mathcal R\colon \arg z=\phi-\pi \}.
\end{gather*}
We will say that these rays are \emph{admissible at $u$}, for the system~\eqref{semiseq}, if they do not coincide with any Stokes rays $R_{ij,k}$ for any~$i$,~$j$ s.t.\ $u_i\neq u_j$ and any $k\in\mathbb Z$. Moreover, a line $\ell(\phi):= \big\{z=\rho {\rm e}^{{\rm i}\phi}, \allowbreak \rho\in\mathbb{R}\big\}$ of the complex plane, with the orientation induced by~$\mathbb{R}$, is called \emph{admissible at $u$} for the system~\eqref{semiseq} if
\[\operatorname{Re}z(u_i-u_j)|_{z\in\ell\setminus 0}\neq 0
\]for any $i$, $j$ s.t.\ $u_i\neq u_j$. In other words, a line is admissible if it does not contain (projected) Stokes {ray} $R_{ij}$.
\end{Definition}
Notice that {the rays} $\operatorname{pr}{(\ell_\pm(\phi))}$ are contained in the line $\ell(\phi)=\big\{z=\rho {\rm e}^{{\rm i}\phi},~\rho\in\mathbb{R}\big\}$, and that the orientation induced by $\mathbb{R}$ is such that the positive part of $\ell(\phi)$ is $\operatorname{pr}{ (\ell_+(\phi))}$.

\begin{Definition}[$\ell$-chambers]\label{ellecella} Given a semisimple Frobenius manifold $M$, and fixed an oriented line $\ell(\phi)=\{z=\rho {\rm e}^{{\rm i}\phi},~\rho\in\mathbb{R}\}$ in the complex plane, consider the open dense subset of points $p\in M$ such that{\samepage
\begin{itemize}\itemsep=0pt
\item the eigenvalues of $U$ at $p$ are pairwise distinct,
\item the line $\ell$ is admissible at $u(p)=(u_1(p),\dots ,u_n(p))$.
\end{itemize}
We call {\it $\ell$-chamber} any connected component $\Omega_\ell$ of this set.}
\end{Definition}

The definition is well posed, since it does not depend on the ordering of the idempotents (i.e., the labelling of the canonical coordinates) and on the signs in the square roots defining~$\Psi$.
Any $\ell$-chamber satisfies the property of Remark \ref{26.05.17-3}: hence, idempotent vector fields and canonical coordinates are single-valued and holomorphic on any $\ell$-chamber.
The topology of an $\ell$-chamber in $M$ can be highly non-trivial (it should not be confused with the simple topology in $\mathbb{C}^n$ of an {\it $\ell$-cell} of Definition \ref{ellcell} below). For example, in \cite{guzzettikon} the analytic continuation of the Frobenius structure of the Quantum Cohomology of $\mathbb P^2$ is studied: it is shown that there exist points $(u_1,u_2, u_3)\in\mathbb C^3$ with $u_i\neq u_j$, which do not correspond to any true geometric point of the Frobenius manifold. This is due to singularities of the change of coordinates $u\mapsto t$.

\begin{Remark}\label{3dicembre2017-1}In \cite[Section 3.4]{dubro0}, the second author introduced a strictly related notion of \emph{charts} of semisimple Frobenius manifolds. Although both definitions of chambers and charts are subordinate to the choice of an oriented line $\ell$, notice some differences between the two concepts. Basically, $\ell$-chambers are a non-coordinatized version of charts. Given a semisimple Frobenius manifold, its decomposition is intrinsically defined and it depends on the spectrum of~$\mathcal U$ as a~\emph{set}, without particular reference to any ordering of canonical coordinates.

Conversely, adopting an \emph{inverse-problem} point of view, as in Section 3.4 of \cite{dubro0}, charts are identified with open sets of $n$-tuples $(u_1,\dots, u_n)\in\mathbb C^n$ with pairwise distinct values of $u_i$'s in $\ell$-\emph{lexicographical order} (see Definition~\ref{30luglio2016-1}), and in correspondence to which a suitable Riemann-Hilbert problem is solvable, so that the local Frobenius structure can be reconstructed. Furthermore, it is also required a condition guaranteeing that the changes of coordinates $t\mapsto u$, $u\mapsto t$ are not singular. Note that in both cases (charts or chambers), semisimple coalescence points are not considered: hence, despite of their name, charts do not really constitute an atlas of the Frobenius manifold.
\end{Remark}

For a fixed $\phi\in\mathbb R$, we define the sectors
\begin{gather*} \Pi_{\rm right}(\phi):=\left\{z\in\mathcal R\colon \phi-\pi< \arg z< \phi\right\},\\
\Pi_{\rm left}(\phi):=\left\{z\in\mathcal R\colon \phi< \arg z< \phi+\pi\right\}.
\end{gather*}

\begin{Theorem}\label{fundamental} Let $\Omega\subset M_{ss}$ be as in Remark~{\rm \ref{26.05.17-3}} and let system \eqref{idemsist} be determined as in the beginning of this section. Let $\phi\in\mathbb R$ be fixed. Then the following statements hold.
\begin{enumerate}\itemsep=0pt
\item[$1.$] At any $p\in \Omega$ such that $\ell(\phi)$ is admissible at $u(p)$ and, for any $k\in\mathbb Z$ there exist two fundamental matrix solutions $Y^{(k)}_{{\rm left/right}}(z)$ uniquely determined by the asymptotic condition
\[
Y^{(k)}_{{\rm left/right}}(z)\sim Y_{\rm formal}(z),\qquad |z|\to\infty,\qquad z\in {\rm e}^{2\pi{\rm i} k}\Pi_{\rm left/right}(\phi).
\]

\item[$2.$] The above solutions $Y^{(k)}_{\rm left/right}$ satisfy
\begin{gather}
\label{othersolutions}Y^{(k)}_{\rm left/right}\big({\rm e}^{2\pi{\rm i} k}z\big)=Y^{(0)}_{\rm left/right}(z),\qquad z\in\mathcal R.
\end{gather}

\item[$3.$] In case $\Omega\equiv \Omega_\ell$ is an $\ell(\phi)$-chamber if $p$ varies in $\Omega_\ell$, then the solutions $Y^{(k)}_{{\rm left/right}}(z)$ define holomorphic functions
\begin{gather*}Y^{(k)}_{{\rm left/right}}(z, u)
\end{gather*}
 w.r.t.\ $u=u(p)$.
 Moreover, the asymptotic expansion
\begin{gather}
\label{5luglio2017-2}
Y^{(k)}_{{\rm left/right}}(z,u)\sim Y_{\rm formal}(z,u),\qquad |z|\to\infty,\qquad z\in {\rm e}^{2\pi{\rm i} k}\Pi_{\rm left/right}(\phi),
\end{gather}
 holds uniformly in~$u$ corresponding to~$p$ varying in $\Omega_\ell$. Here $Y_{\rm formal}(z,u)$ is the $u$-holo\-mor\-phic formal solution~\eqref{7luglio2017-1}.
\end{enumerate}
\end{Theorem}

\begin{proof}
The proof of (1) and (2) away from coalescence points is standard (see \cite{BJL1,dubro2,dubro3, Wasow}), while at coalescence points it follows from the results of~\cite{CDG0} and~\cite{BJL2}. Point (3) is stated in \cite{dubro2,dubro3}, though the name ``$\ell$-chamber'' does not appear there.
\end{proof}

\begin{Remark} The holomorphic properties at point (3) of Theorem \ref{fundamental} hold in a $\ell$-chamber, where there are no coalescence points.
In our Theorem \ref{mainisoth} below, we will see that point (3) actually holds in a set $\Omega\subset M_{ss}$ as in Remark \ref{26.05.17-3}, no matter whether it contains semisimple coalescence points or not. The only requirement is that $\ell(\phi)$ is admissible at $u=u(p)$ for any $p\in\Omega$.
\end{Remark}

\begin{figure}\centering
\def\svgscale{0.5}
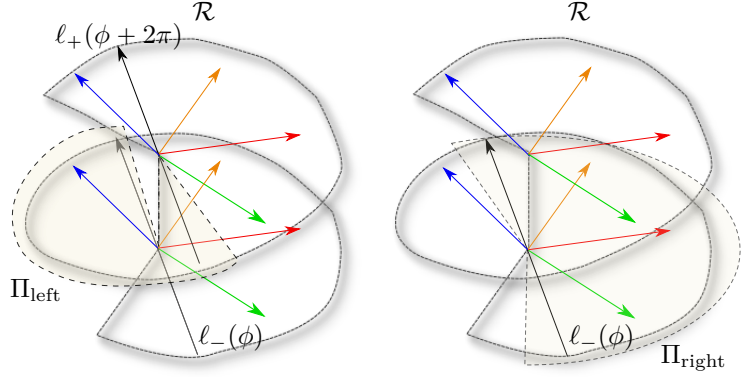
\caption{The figure shows $\Pi^\varepsilon_{\rm right}(\phi),\Pi^\varepsilon_{\rm left}(\phi)$ as dashed sectors, $\ell_\pm(\phi)$ in (black) and Stokes rays (in color).}
\end{figure}

\begin{Remark}\label{26.05.17-1}
The asymptotic relation~(\ref{5luglio2017-2}) means that for any compact $K\Subset\Omega_\ell$, for any $h\in\mathbb N$ and for any proper closed subsector $\overline{\mathcal S}\subsetneq {\rm e}^{2\pi{\rm i}k}\Pi_{\operatorname{right/left}} (\phi)$ there exists $C_{K,h,\overline{\mathcal S}}>0$ such that, if $z\in\overline{\mathcal S}\setminus \{0 \}$, then
\[
\sup_{u\in K}\left\| Y^{(k)}_{\operatorname{right/left}}(z,u)\cdot\exp(-z U)-\sum_{m=0}^{h-1}\frac{G_m(u)}{z^m}\right\|<\frac{C_{K,h,\overline{\mathcal S}}}{|z|^h}.
\]
 Actually, the solutions $Y^{(k)}_{\text{right/left}}(z,u)$ maintain their asymptotic {expansions} (\ref{5luglio2017-2}) in sectors
 wider than  ${\rm e}^{2\pi{\rm i} k}\Pi_{\text{right/left}}(\phi)$ {after} extending at least up to the nearest Stokes rays outside ${\rm e}^{2\pi{\rm i} k}\Pi_{\text{right/left}}(\phi)$. In particular, for any $p\in K\Subset\Omega_\ell$ and suitably small $\varepsilon=\varepsilon(K)>0$, then the asymptotics holds in ${\rm e}^{2\pi{\rm i} k}\Pi_{\text{right/left}}^\varepsilon(\phi)$, where
 \begin{gather*}\Pi^\varepsilon_{\rm right}(\phi): =\left\{z\in\mathcal R\colon \phi-\pi-\varepsilon< \arg  z< \phi+\varepsilon\right\},\\
\Pi^\varepsilon_{\text{left}}(\phi): =\left\{z\in\mathcal R\colon \phi-\varepsilon< \arg z< \phi+\pi+\varepsilon\right\}.
\end{gather*}
The positive number $\varepsilon$ is chosen small enough in such a way that, as~$p$ varies in the compact set~$K$, no Stokes ray is contained in the following sectors:
\begin{gather*}
\Pi^\varepsilon_{+}(\phi):= \{z\in\mathcal R\colon \phi-\varepsilon< \arg  z< \phi+\varepsilon \},
\\
\Pi^\varepsilon_{-}(\phi):= \{z\in\mathcal R\colon \phi-\pi-\varepsilon< \arg z< \phi-\pi+\varepsilon \}.
\end{gather*}
\end{Remark}

\begin{Lemma}\label{orthlemma}In the assumptions of Theorem~{\rm \ref{fundamental}}, for any $k\in\mathbb Z$ and any $z\in\mathcal R$ the following orthogonality relation holds:
\[Y^{(k)}_{\rm left}\big({\rm e}^{{\rm i}\pi}z\big)^{\rm T}Y^{(k)}_{\rm right}(z)=\mathbbm 1.
\]
\end{Lemma}

\begin{proof} From Remark~\ref{24novembre2016-1} we already know that the product above is independent of $z\in\mathcal R$. According to Remark~\ref{26.05.17-1}, if $\varepsilon>0$ is a sufficiently small positive number, then
\[Y^{(k)}_{{\rm left/right}}(z)\sim Y_{\rm formal}(z),\qquad |z|\to\infty,\qquad z\in {\rm e}^{2\pi{\rm i} k}\Pi_{\rm left/right}^\varepsilon(\phi).
\]
Consequently,
\[Y^{(k)}_{\rm left}\big({\rm e}^{{\rm i}\pi}z\big)\sim G(-z){\rm e}^{-zU},\qquad Y^{(k)}_{\rm right}(z)\sim G(z){\rm e}^{zU}, \qquad |z|\to\infty,\qquad z\in {\rm e}^{2\pi{\rm i} k}\Pi^\varepsilon_+(\phi).\]
Thus, $Y^{(k)}_{\rm left}\big({\rm e}^{{\rm i}\pi}z\big)^{\rm T}Y^{(k)}_{\rm right}(z)=\mathbbm 1$ for all $z\in {\rm e}^{2\pi{\rm i} k}\Pi^\varepsilon_+(\phi)$, and by analytic continuation for all $z\in\mathcal R$.
\end{proof}

Let $Y_0(z,u)$ be a fundamental solution of~\eqref{semiseq} near $z=0$ of the form~(\ref{3luglio2016-1}), i.e.,
\begin{gather}
Y_0(z,u)=\Psi(u)\Phi(z,u)z^\mu z^R,\qquad \Phi(z,u)=\mathbbm 1+\sum_{k=1}^\infty \Phi_k(u)z^k,
\nonumber\\
 \Phi(-z,u)^{\rm T}\eta  \Phi(z,u)=\eta,\label{30nov2016-6}
\end{gather}
 with $\Psi^{\rm T}\Psi=\eta$, obtained from \eqref{28luglio2016-9} and \eqref{30nov2016-4} through the constant gauge \eqref{30nov2016-5}. This solution is not affected by coalescence phenomenon and since $\mu$ and $R$ are independent of $p\in \Omega$, it is holomorphic w.r.t.\ $u$ (see \cite{dubro2,dubro3}). Recall that $Y_0(z,u)$ is not uniquely determined by the choice of~$R$.

\begin{Definition}[Stokes and central connection matrices]\label{5luglio2017-5}
 Let $\Omega\subset M_{ss}$ be as in Remark \ref{26.05.17-3} and let the system \eqref{idemsist} be determined as in the beginning of this section. Let $\phi\in\mathbb R$ be fixed. Let $p\in\Omega$ be such that $\ell(\phi)$ is admissible at $u(p)$. Finally, let $Y^{(0)}_{\text{right/left}}(z)$ be the fundamental solutions of Theorem \ref{fundamental} at $p$. The matrices $S$ and $S_{-}$ defined at $u(p)$ by the relations
\begin{gather}\label{matrstok1}Y^{(0)}_{\rm left}(z)=Y^{(0)}_{\rm right}(z) S,\qquad z\in\mathcal R,
\\
Y^{(0)}_{\rm left}\big({\rm e}^{2\pi{\rm i}}z\big)=Y^{(0)}_{\rm right}(z) S_-,\qquad z\in\mathcal R\nonumber
\end{gather}
are called \emph{Stokes matrices} of the system \eqref{semiseq} at the point $p$ w.r.t.\ the line $\ell(\phi)$. The matrix $C$ such that
\begin{gather*}
Y^{(0)}_{\rm right}(z)=Y_0(z,u(p)) C,\qquad z\in\mathcal R
\end{gather*} is called \emph{central connection matrix} of the system \eqref{idemsist} at~$p$, w.r.t.\ the line $\ell$ and the fundamental solution~$Y_0$.
\end{Definition}

\begin{Theorem}\label{proprstok} The Stokes matrices $ S$, $S_-$ and the central connection matrix $C$ of Definition~{\rm \ref{5luglio2017-5}} at a point $p\in \Omega$ satisfy the following properties, for all $k\in\mathbb Z$ and all $z\in\mathcal R$:
\begin{enumerate}\itemsep=0pt
\item[$1.$]
\begin{gather*}
Y^{(k)}_{\operatorname{left}}(z) =Y^{(k)}_{\operatorname{right}}(z)  S,\\
Y^{(k)}_{\operatorname{left}}(z) =Y^{(k+1)}_{\operatorname{right}}(z)  S_-,\\
Y^{(k)}_{\operatorname{right}}(z) =Y_0(z,u(p))~ M_0^{-k}  C,
\end{gather*}where $M_0=\exp(2\pi{\rm i} \mu)\exp(2\pi{\rm i} R)$;

\item[$2.$]
\begin{gather*}
Y^{(k)}_{\operatorname{right}}\big({\rm e}^{2\pi{\rm i}}z\big) =Y^{(k)}_{\operatorname{right}}(z)  S_-  S^{-1},\\
Y^{(k)}_{\operatorname{left}}\big({\rm e}^{2\pi{\rm i}}z\big) =Y^{(k)}_{\operatorname{right}}(z)  S^{-1}  S_-;
\end{gather*}
\item[$3.$]
\begin{gather*}
 S_- = S^{\rm T},\\
 S_{ii}=1, \qquad i=1,\dots, n,\\
 S_{ij}\neq 0\qquad \text{with }i\neq j\text{ only if } u_i\neq u_j \text{ and }R_{ij}\subset\operatorname{pr} (\Pi_{\operatorname{left}}(\phi) ).
\end{gather*}
\end{enumerate}
\end{Theorem}

\begin{proof} The first and second identities of (1) follow from equation \eqref{othersolutions}. For the third note that
\begin{gather*}
Y^{(k)}_{\rm right}(z)=Y^{(0)}_{\rm right}\big({\rm e}^{-2{\rm i}k\pi}z\big)=Y_{0}\big({\rm e}^{-2{\rm i}k\pi}z\big)C=Y_0(z)M_0^{-k}C.
\end{gather*}
Point (2) follows easily from the vanishing of the exponent of formal monodromy ($\operatorname{diag} V=0$). By definition of Stokes matrices we have that
\[Y^{(0)}_{\rm left}\big({\rm e}^{{\rm i}\pi}z\big)= Y^{(0)}_{\rm right}\big({\rm e}^{-{\rm i}\pi}z\big)S_-,\qquad Y^{(0)}_{\rm right}(z)= Y^{(0)}_{\rm left}(z)S^{-1},
\]
and by Lemma \ref{orthlemma}
\[S^{\rm T}_-\underbrace{Y^{(0)}_{\rm right}(z)^{\rm T}Y^{(0)}_{\rm left}\big({\rm e}^{{\rm i}\pi}z\big)}_{\mathbbm 1}S^{-1}\equiv\mathbbm1.
\]
We conclude $S^{\rm T}_-=S$. If we consider the sector $\Pi^\varepsilon_{+}(\phi)$ for sufficiently small $\varepsilon>0$ as in proof of Lemma \ref{orthlemma}, them from the relation
$Y^{(0)}_{\rm left}(z)= Y^{(0)}_{\rm right}(z)S$,
we deduce that
\[{\rm e}^{z(u_i-u_j)}S_{ij}\sim\delta_{ij},\qquad |z|\to\infty,\qquad z\in\Pi^\varepsilon_{+}(\phi).
\] So, if $u_i=u_j$ we deduce $S_{ij}=\delta_{ij}$. If $i\neq j$ are such that $u_i\neq u_j$, then if $R_{ij}\subset \operatorname{pr} (\Pi_{\rm right}(\phi) )$
we have
 \[\big|{\rm e}^{z(u_i-u_j)}\big|\to\infty\qquad\text{for }|z|\to\infty,\qquad z\in\Pi^\varepsilon_{+}(\phi),
\] and hence necessarily $S_{ij}=0$. For the opposite ray $R_{ji}\subset\operatorname{pr}\left(\Pi_{\rm left}\right)$ we have
 \[\big|{\rm e}^{z(u_i-u_j)}\big|\to0\qquad\text{for }|z|\to\infty,\qquad z\in\Pi^\varepsilon_{+}(\phi),
\]so $S_{ij}$ need not to be~0. This proves~(3).
\end{proof}

The monodromy data must satisfy some important constraints, summarised in the following theorem, whose proof is {omitted} in~\cite{dubro0,dubro2}.

\begin{Theorem}\label{constraint} The monodromy data $\mu$, $R$, $S$, $C$ at a point $p\in\Omega$ as in Definition~{\rm \ref{5luglio2017-5}} satisfy the identities:
\begin{enumerate}\itemsep=0pt
 \item[$1)$] $CS^{\rm T}S^{-1}C^{-1}=M_0={\rm e}^{2\pi{\rm i}\mu}{\rm e}^{2\pi{\rm i}R}$,
\item[$2)$] $S=C^{-1}{\rm e}^{-\pi{\rm i} R}{\rm e}^{-\pi{\rm i} \mu}\eta^{-1}\big(C^{\rm T}\big)^{-1}$,
\item[$3)$] $S^{\rm T}=C^{-1}{\rm e}^{\pi{\rm i} R}{\rm e}^{\pi{\rm i} \mu}\eta^{-1}\big(C^{\rm T}\big)^{-1}.$
\end{enumerate}
\end{Theorem}
\begin{proof}
The first identity has a simple topological motivation: loops around the origin in $\mathbb C^*$ are homotopic to loops around infinity. So, one easily obtains the relation using Theorem \ref{proprstok}, and the definition of central connection matrix. Using the orthogonality relations for solutions, equation \eqref{relazioneRemu} and the fact that \[z^{\mu^{\rm T}}\eta  z^{\mu}=\eta\] ($\mu$ being diagonal and $\eta$-antisymmetric), we can now prove the identities (2) and (3). By Lemma~\ref{orthlemma} we have that
\begin{gather*}\mathbbm1 =Y^{(0)}_{\rm right}(z)^{\rm T}Y^{(0)}_{\rm left}\big({\rm e}^{{\rm i}\pi}z\big)
 =Y^{(0)}_{\rm right}(z)^{\rm T}Y^{(0)}_{\rm right}\big({\rm e}^{{\rm i}\pi}z\big)S=C^{\rm T}Y_0(z)^{\rm T}Y_0\big({\rm e}^{{\rm i}\pi}z\big)CS.
\end{gather*}
Now we have
\begin{align*}
Y_0(z)^{\rm T}Y_0\big({\rm e}^{{\rm i}\pi}z\big)&
=
z^{R^{\rm T}}z^{\mu^{\rm T}}\big(\Phi(z)^{\rm T}\Psi^{\rm T}\Psi \Phi(-z)\big)\big({\rm e}^{{\rm i}\pi}z\big)^\mu
\big({\rm e}^{{\rm i}\pi}z\big)^R\\
&=z^{R^{\rm T}}z^{\mu^{\rm T}}\eta z^\mu {\rm e}^{{\rm i}\pi\mu}z^R{\rm e}^{{\rm i}\pi R}
 =\eta {\rm e}^{{\rm i}\pi\mu}{\rm e}^{{\rm i}\pi R}.
\end{align*}
This shows the first identity. For the second one, we have that
\begin{gather*}
\mathbbm1 =Y^{(0)}_{\rm right}(z)^{\rm T}Y^{(0)}_{\rm left}\big({\rm e}^{{\rm i}\pi}z\big)
 =Y^{(0)}_{\rm right}(z)^{\rm T}Y^{(0)}_{\rm right}\big({\rm e}^{-{\rm i}\pi}z\big)S^{\rm T}
 =C^{\rm T}Y_0(z)^{\rm T}Y_0\big({\rm e}^{-{\rm i}\pi}z\big)CS^{\rm T}.
\end{gather*}
Again, we have
\begin{align*}Y_0(z)^{\rm T}Y_0\big({\rm e}^{-{\rm i}\pi}z\big)&=z^{R^{\rm T}}z^{\mu^{\rm T}}\big(\Phi(z)^{\rm T}\Psi^{\rm T}\Psi \Phi(-z)\big)\big({\rm e}^{-{\rm i}\pi}z\big)^\mu
\big({\rm e}^{-{\rm i}\pi}z\big)^R\\
&=z^{R^{\rm T}}z^{\mu^{\rm T}}\eta z^\mu {\rm e}^{-{\rm i}\pi\mu}z^R{\rm e}^{-{\rm i}\pi R}
=\eta {\rm e}^{-{\rm i}\pi\mu}{\rm e}^{-{\rm i}\pi R}.\tag*{\qed}
\end{align*}\renewcommand{\qed}{}
\end{proof}

It follows from point (3) of Theorem \ref{fundamental} that $S$ and $C$ depend holomorphically on $p$ varying in an $\ell$-chamber $\Omega_\ell$, namely they define analytic matrix valued functions $S(u)$ and $C(u)$, $u=u(p)$. Moreover, due to the compatibility conditions $[E_i,V]=[V_i,U]$ and $\partial_i\Psi =V_i\Psi$, the system \eqref{idemsist} is isomonodromic. Therefore $\partial_iS=\partial_i C=0$. Indeed, the following holds:

\begin{Theorem}[isomonodromy Theorem, II, \cite{dubro1, dubro0, dubro2}]\label{isoth2}
 The Stokes matrix $S$ and the central connection matrix $C$, computed w.r.t.\ a line $\ell$, are independent of $p$ varying in an $\ell$-chamber. The values of $S$, $C$ in two different $\ell$-chambers are related by an action of the braid group of Section~{\rm \ref{freedom}}.
\end{Theorem}

\section[Ambiguity in definition of monodromy data and braid group action]{Ambiguity in definition of monodromy data\\ and braid group action}\label{freedom}

In associating the data $(\mu, R, S, C)$ to $p\in M$ several choices have been done, all preserving the constraints of Theorem \ref{constraint}
\begin{gather*}
S=C^{-1}{\rm e}^{-{\rm i}\pi R}{\rm e}^{-{\rm i}\pi \mu}\eta^{-1}\big(C^{-1}\big)^{\rm T},\\
S^{\rm T}=C^{-1}{\rm e}^{{\rm i}\pi R}{\rm e}^{{\rm i}\pi\mu}\eta^{-1}\big(C^{-1}\big)^{\rm T}.
\end{gather*}

 While the operator $\mu$ is completely fixed by the choice of flat coordinates as in Section \ref{sec1}, $R$~is determined only up to conjugacy class of the $(\eta,\mu)$-parabolic orthogonal group $\mathcal G(\eta,\mu)$ as in Theorem~\ref{3agosto2016-5}. Suppose now that $R$ has been chosen in this class.
The remaining local invariants~$ S$,~$C$ are subordinate to the following choices:
\begin{enumerate}\itemsep=0pt
\item[1)] an oriented line $\ell(\phi)=\big\{z=\rho {\rm e}^{{\rm i}\phi},~\rho\in\mathbb{R}\big\}$ in the complex plane;
\item[2)] for given $\phi\in\mathbb{R}$, the change $\phi\mapsto \phi-2k\pi$, $k\in \mathbb{Z}$, or dually, for fixed $\phi$, the change $Y^{(0)}_{\rm left/right}(z)\mapsto Y^{(k)}_{\rm left/right}(z)$;
\item[3)] the choice of {an} ordering of canonical coordinates on each $\ell$-chamber $\Omega_\ell$;
\item[4)] the choice of the {branches} of the square roots~(\ref{5luglio2017-1}) defining the matrix $\Psi$ on each $\ell$-chamber $\Omega_\ell$;
\item[5)] the choice of {solution} $Y_0$ in the Levelt normal form corresponding to the same exponent~$R$.
\end{enumerate}
The transformations of the data depending on the choice of $\ell$ in~(1) will be studied in the next Section. Here we describe how the freedoms~(2), (3), (4) and (5) affect the data $(S, C)$:
\begin{itemize}\itemsep=0pt
\item \emph{Action of the additive group $\mathbb Z$}: according to formula (\ref{othersolutions}), $S$ remains invariant and
\[C\mapsto M_0^{-k}\cdot C,\qquad k\in\mathbb Z,\qquad M_0={\rm e}^{2\pi{\rm i}\mu}{\rm e}^{2\pi{\rm i}R},\qquad t\in\Omega_\ell.
\]

\item \emph{Action of the group of permutations} $\mathfrak S_n$: if $\tau$ is a permutation,
we can reorder the canonical coordinates:\[(u_1,\dots, u_n)\mapsto(u_{\tau(1)},\dots,u_{\tau(n)}).\]
The system (\ref{semiseq}) is changed to $U\mapsto P UP^{-1}=\operatorname{diag}(u_{\tau(1)},\dots,u_{\tau(n)})$, $V\mapsto P VP^{-1}$.
The fundamental matrices change as follows: $Y^{(0)}_{\rm left/right}\mapsto P Y^{(0)}_{\rm left/right}P^{-1}$ and $Y_0\longmapsto P Y_0$.
 Therefore
\begin{gather}\label{29luglio2016-9}
S\mapsto P SP^{-1},\qquad C\mapsto CP^{-1}.
\end{gather}

\item \emph{Action of the group} $(\mathbb Z/2\mathbb Z)^{\times n}$: by changing signs of the normalized idempotents (matrix~$\Psi$) we change the signs of the entries of the matrices $S$ and $C$. If $\mathcal I$ is a diagonal matrix with~$1$'s or~$(-1)$'s on the diagonal, the system (\ref{semiseq}) is changed to $U\mapsto \mathcal{I} U\mathcal{I}\equiv U$, $V\mapsto \mathcal{I}V\mathcal{I}$. Correspondingly, $Y_{\rm left/right}\mapsto \mathcal{I} Y_{\rm left/right}\mathcal I$, $Y_0\mapsto \mathcal{I} Y_0$. Therefore
\[S\mapsto \mathcal IS\mathcal I,\qquad C\mapsto C\mathcal I.
\]

\item \emph{Action of the group }$\mathcal C_0(\eta,\mu, R)$: for chosen $R$, the choice of a fundamental system at the origin having the form~\eqref{30nov2016-6} is {defined} up to $Y_0\mapsto Y_0G$, where $G\in \mathcal C_0(\eta,\mu, R)$ of Definition~\ref{29.05.17-1}. The corresponding left action on $C$ is
 \begin{gather*}
 C\longmapsto GC,\qquad G\in \mathcal C_0(\eta,\mu, R).
 \end{gather*}
\end{itemize}

Among all possible {orderings} of the canonical coordinates, a particularly useful one is the \emph{lexicographical order} w.r.t.\ an admissible line $\ell(\phi)$, defined as follows. Let us introduce the following rays in the complex plane:
\begin{gather}\label{3giugno2019-1}
L_j(\phi):=\big\{u_j+\rho {\rm e}^{{\rm i}\left(\frac{\pi}{2}-\phi\right)}\colon \rho\in\mathbb R_+\big\},\qquad j=1,\dots, n.
\end{gather}
Each of them starts from the point $u_j$ and is considered to be {\it oriented} from $u_j$ to $\infty$.

\begin{Definition}[lexicographical order]\label{30luglio2016-1} The canonical coordinates $u_j$'s are in \emph{$\ell$-lexicogra\-phi\-cal order} if $L_j(\phi)$ is to the left of $L_k(\phi)$ (w.r.t.\ the above orientation), for any $1\leq j<k\leq n$.
\end{Definition}

If $u_1,\dots, u_n$ are in lexicographical order w.r.t.\ the admissible line $\ell(\phi)$, then:
\begin{enumerate}\itemsep=0pt
\item[1)] the Stokes matrix is in upper triangular form,
\item[2)] $R_{i,j}\subseteq \operatorname{pr}\left(\Pi_{\rm left}(\phi)\right)$ if and only if $i<j$,
\item[3)] the nearest Stokes rays to the positive half-line $\operatorname{pr}(\ell_+(\phi))$ are of the form
\[ R_{i,i+1}\subseteq \operatorname{pr} (\Pi_{\rm left}(\phi) ),
\qquad R_{j,j-1}\subseteq \operatorname{pr} (\Pi_{\rm right}(\phi) ),
\]
where $1\leq i\leq n-1$ and $2\leq j\leq n$.
\end{enumerate}

In general, condition $(1)$ alone does not imply that the canonical coordinates are in lexicographical order: it does if and only if the number of nonzero entries of the Stokes matrix $S$ is maximal (and equal to~$\frac{n(n+1)}{2}$). In this case, by Theorem~\ref{proprstok}, necessarily $u_i\neq u_j$ for $i\neq j$. On the other hand, if there are some  vanishing entries $S_{ij}=S_{ji}=0$ for $i\neq j$, and $S$ is upper triangular, then also $PSP^{-1}$ in~\eqref{29luglio2016-9} is upper triangular for any permutation exchanging $u_i$ and $u_j$ corresponding to $S_{ij}=S_{ji}=0$.
For example, this happens at a coalescence point: by Theorem~\ref{proprstok}, the entries $S_{ij}$ with $i\neq j$ are $0$ corresponding to coalescing values $u_i=u_j$, $i\neq j$.

\begin{Definition}[triangular order]\label{30luglio2016-2}
We say that $u_1,\dots ,u_n$ are in triangular order w.r.t.\ the line $\ell$ whenever $S$ is upper triangular.
\end{Definition}

It follows from the preceding discussion that at a semisimple coalescence point there are more than one triangular orders. Moreover, {\it any} of them is also lexicographical. For further comments, see Remark~\ref{9dicembre2016-3}.

\subsection[Action of the braid group $\mathcal B_n$]{Action of the braid group $\boldsymbol{\mathcal B_n}$}

In this section, canonical coordinates are pairwise distinct, corresponding to a non-coalescence semisimple points lying in $\ell$-chambers. The braid group is
\[
\mathcal B_n=\pi_1\big(\big(\mathbb C^n\setminus\Delta\big)/\mathfrak S_n\big),
\]
where $\Delta$ stands for the union of all diagonals in $\mathbb C^n$. It is generated by $n-1$ elementary braids $\beta_{12},\beta_{23},\dots,\beta_{n-1,n}$, with the relations
\begin{gather*} \beta_{i,i+1}\beta_{j,j+1}=\beta_{j,j+1}\beta_{i,i+1}\qquad\text{for }i+1\neq j,j+1\neq i,
\\
\beta_{i,i+1}\beta_{i+1,i+2}\beta_{i,i+1}=\beta_{i+1,i+2}\beta_{i,i+1}\beta_{i+1,i+2}.
\end{gather*}
Any braid in $\mathcal B_n$ is a product of the generators $\beta_{12},\beta_{23},\dots,\beta_{n-1,n}$ and their inverses.

The action of the braid group $\mathcal B_n$ on the monodromy data manifests whenever some Stokes ray and the chosen line~$\ell$ cross under rotation. This can happen in two ways:
\begin{itemize}\itemsep=0pt
\item First: we let vary the point of the Frobenius manifold at which we compute the data, keeping fixed the line~$\ell$; this is the case if, starting from the data computed in {an} $\ell$-chamber we want to compute the data in a~neighboring $\ell$-chamber, or even more in general if we want to analyze properties of the analytic continuation of the whole Frobenius structure by letting varying the coordinates $(u_1,\dots, u_n)$ on the universal cover $\widetilde{\mathbb C^n\setminus\Delta}$.
\item Second: we fix the point at which we compute the data and change the admissible line~$\ell$ by a rotation.
\end{itemize}
In the first case the $\ell$-chambers are fixed, in the second case they change: indeed, the {given} point of the Frobenius manifold is in two different chambers before and after the rotation of~$\ell$. In both cases, we will always label the canonical coordinates $(u_1,\dots,u_n)$ in lexicographical order w.r.t.~$\ell$ both before and after the transformation (so that, in particular, any Stokes matrix is always in upper triangular form).

Any continuous deformation of the $n$-tuple $(u_1,\dots,u_n)$, represented as a deformation of $n$ points in $\mathbb C$ never colliding, can be decomposed into \emph{elementary} ones. If we restrict to the case of a continuous deformation which ends exactly with the same {initially} ordered pattern of points, then we can identify an elementary deformation with a generator of the \emph{pure} braid group, i.e., $\pi_1\big(\mathbb C^n\setminus\Delta\big)$. Otherwise, by allowing permutations, we can identify an elementary deformation with a generator of the braid group $\mathcal B_n$. In particular, an elementary deformation which will be denoted by $\beta_{i,i+1}$ consists in a counter-clockwise rotation of $u_i$ w.r.t.~$u_{i+1}$, so that the two exchange. All other points $u_j$'s are subjected to a sufficiently small perturbation, so that the corresponding Stokes' rays almost do not move. $\beta_{i,i+1}$ corresponds to
\begin{itemize}\itemsep=0pt
\item clockwise rotation of the Stokes' ray $R_{i,i+1}$ crossing the line $\ell$,
\item or, dually, counter-clockwise rotation of the line $\ell$ crossing the Stokes' ray $R_{i,i+1}$.
\end{itemize}
This determines the following mutations of the monodromy data, as shown in \cite{dubro1} and \cite{dubro2}.
For the deformation of $u_i$,
$u_{i+1}$ relatively moving anticlockwise, associated with $\beta_{i,i+1}$, we have
\begin{gather}
\label{stokesbraid1}
S\mapsto S^{\beta_{i,i+1}}:= A^{\beta_{i,i+1}}(S) S A^{\beta_{i,i+1}}(S)^{\rm T},
\end{gather}
where
\begin{gather*}
\big(A^{\beta_{i,i+1}}(S)\big)_{hh}=1,\qquad h=1,\dots, n\quad h\neq i,i+1,\\
\big(A^{\beta_{i,i+1}}(S)\big)_{i+1,i+1}=-s_{i,i+1}, \qquad
\big(A^{\beta_{i,i+1}}(S)\big)_{i,i+1}=\big(A^{\beta_{i,i+1}}(S)\big)_{i+1,i}=1,
\end{gather*}
and all the other entries are zero. For the inverse braid $\beta_{i,i+1}^{-1}$ ($u_i$ and $u_{i+1}$ move
 clockwise) the mutation is $S^{\beta_{i,i+1}^{-1}}:= A^{\beta_{i,i+1}^{-1}}(S) S  A^{\beta_{i,i+1}^{-1}}(S)^{\rm T}$, where
\begin{gather*}
\big(A^{\beta_{i,i+1}^{-1}}(S)\big)_{hh}=  1, \qquad h=1,\dots ,n, \qquad h\neq i, i+1,\\
\big(A^{\beta_{i,i+1}^{-1}}(S) \big)_{i,i}=-s_{i,i+1},\qquad
 \big(A^{\beta_{i,i+1}^{-1}}(S) \big)_{i,i+1}=
\big(A^{\beta_{i,i+1}^{-1}}(S) \big)_{i+1,i}=1,
\end{gather*}
and all the other entries are zero. For a generic braid $\beta=\beta_{i_1,i_1+1}^{\pm 1}\cdot\beta_{i_2,i_2+1}^{\pm 1} \cdots \beta_{i_N,i_N+1}^{\pm 1}$, which is a product of $N\geq 1$ elementary braids or their inverses, the action is
 \begin{gather}\label{connbraid2-29luglio}
 S\mapsto S^{\beta}=  A^{\beta}(S)S\big[A^{\beta}(S)\big]^{\rm T}.
\end{gather}
We remark that $S^{\beta}$ is still upper triangular. The action on the central connection matrix (in lexicographical order) is
\begin{gather}\label{connbraid2}
C\mapsto C^\beta:= C A^\beta(S)^{-1}.
\end{gather}
The matrix $A^{\beta}(S)$ is constructed successively applying \eqref{stokesbraid1}, as follows:
\begin{gather*}
 A^{\beta}(S)=
 A^{\beta_{i_N,i_N+1}^{\pm 1}}\big(S^{
 \beta_{i_1,i_1+1}^{\pm 1}\cdot \beta_{i_2,i_2+1}^{\pm 1}\cdots \beta_{i_{N-1},i_{N-1}+1}^{\pm 1} }\big)
\cdots A^{\beta_{i_2,i_2+1}^{\pm 1}}\big(S^{\beta_{i_1,i_1+1}^{\pm 1}}\big) \cdot A^{\beta_{i_1,i_1+1}^{\pm 1}}(S),
\end{gather*}
where
\begin{gather*}S^{\beta_{i_1,i_1+1}^{\pm 1}}=A^{\beta_{i_1,i_1+1}^{\pm 1}}(S) S A^{\beta_{i_1,i_1+1}^{\pm 1}}(S)^{\rm T},\\
S^{\beta_{i_1,i_1+1}^{\pm 1}\cdot \beta_{i_2,i_2+1}^{\pm 1}} =\big[A^{\beta_{i_2,i_2+1}^{\pm 1}}\big(S^{\beta_{i_1,i_1+1}^{\pm 1}}\big)\big] S^{\beta_{i_1,i_1+1}^{\pm 1}}\big[A^{\beta_{i_2,i_2+1}^{\pm 1}}\big(S^{\beta_{i_1,i_1+1}^{\pm 1}}\big)\big]^{\rm T},\\
S^{\beta_{i_1,i_1+1}^{\pm 1}\cdot \beta_{i_2,i_2+1}^{\pm 1}\cdot \beta_{i_3,i_3+1}^{\pm 1}} =\big[A^{\beta_{i_3,i_3+1}^{\pm 1}}\big(S^{\beta_{i_1,i_1+1}^{\pm 1}\cdot \beta_{i_2,i_2+1}^{\pm 1}}\big)\big]
 S^{\beta_{i_1,i_1+1}^{\pm 1}\cdot \beta_{i_2,i_2+1}^{\pm 1}}\\
 \hphantom{S^{\beta_{i_1,i_1+1}^{\pm 1}\cdot \beta_{i_2,i_2+1}^{\pm 1}\cdot \beta_{i_3,i_3+1}^{\pm 1}}=}{} \times
\big[A^{\beta_{i_3,i_3+1}^{\pm 1}}\big(S^{\beta_{i_1,i_1+1}^{\pm 1}\cdot \beta_{i_2,i_2+1}^{\pm 1}}\big)\big]^{\rm T},
\end{gather*}
and so on.

Now, let us consider a complete counter-clockwise $2\pi$-rotation of the admissible line $\ell$, and observe the following:
\begin{enumerate}\itemsep=0pt
\item In the generic case (i.e., when the canonical coordinates $u_j$'s are in general position) there are $n(n-1)$ distinct projected Stokes' rays~$R_{jk}$. An elementary braid acts any time the line $\ell$ crosses a Stokes ray. So, in total, we expect that a complete rotation of $\ell$ correspond to the product of $n(n-1)$ elementary braids $\beta_{i,i+1}$'s.
\item Since the formal monodromy is vanishing, the effect of the rotation of~$\ell$ on the Stokes matrix is trivial, while the central connection matrix~$C$ is transformed to~$M_0^{-1}C$, $M_0$~being the monodromy at the origin (point~(1) of Theorem~\ref{proprstok}). As a consequence, the complete rotation of the line~$\ell$ can be viewed as a deformation of points $u_j$'s \emph{commuting} with any other braid.
\end{enumerate}
From point $(2)$ we deduce that the braid corresponding to the complete rotation of $\ell$ is an element of the center
\[Z(\mathcal B_n)=\big\{(\beta_{12}\beta_{23}\cdots\beta_{n-1,n})^{kn}\colon k\in\mathbb Z\big\}.
\]From point $(1)$ and from the fact that $\ell$ rotates counter-clockwise we deduce the following
\begin{Lemma}\label{centerbraidlemma}The braid corresponding to a complete counter-clockwise $2\pi$-rotation of $\ell$ is
\[(\beta_{12}\beta_{23}\cdots\beta_{n-1,n})^{n},
\]and its acts on the monodromy data as follows:
\begin{itemize}\itemsep=0pt
\item trivially on Stokes matrices,
\item the central connection matrix is transformed as $C\mapsto M_0^{-1}C$.
\end{itemize}
\end{Lemma}

\section{Isomonodromy theorem at coalescence points}\label{isomonocoal}

So far the monodromy data, $S$ and $C$ have been defined pointwise and then the deformation theo\-ry has been described at point (3) of Theorem~\ref{fundamental} and in Theorem~\ref{isoth2}, away from coalescence points. In particular, $S$ and $C$ are constant in any $\ell$-chamber, and the matrices $Y^{(k)}_{\rm left/right}(z,u)$ are $u$-holomorphic in all $\ell$-chambers. In this section we generalize the deformation theory to semisimple coalescence points. We show that {the} monodromy data, which are well defined at a coalescence point, actually provide the monodromy data in a neighborhood of the point, and can be extended to the whole {Frobenius} manifold through the action of the braid group. In this section we will use the following notation for objects computed at a coalescence point: a~matrix~$Y$,~$S$ or~$C$ will be denoted $\mathring{Y}$, $\mathring{S}$ or $\mathring{C}$.

Let $p_0\in\mathcal B_M\setminus\mathcal K_M$ be a semisimple coalescence point. Consider a neighbourhood $\Omega\subseteq M\setminus\mathcal K_M$ of~$p_0$, satisfying the property of Remark~\ref{26.05.17-3}. An ordering for canonical coordinates $(u_1,\dots, u_n)$ and a holomorphic branch of the function $\Psi\colon\Omega\to {\rm GL}_n(\mathbb C)$ can be chosen in $\Omega$. We denote by $u(p):=(u_1(p),\dots, u_n(p))$ the value of the canonical coordinate map $u\colon \Omega\to\mathbb C^n$, and we define
 \[\Delta_\Omega:=\big\{u(p)= (u_1(p),\dots , u_n(p) )\in\mathbb{C}^n\, \big|\, p\in \Omega\cap\mathcal B_M\big\}.
 \]
Therefore, if $u\in \Delta_\Omega$, then $u_i=u_j$ for some $i\neq j$. The coordinates $u(p_0)$ of~$p_0$ will be denoted $u^{(0)}=\big(u_1^{(0)},\dots ,u_n^{(0)}\big)$. $\Delta_\Omega$ is not empty and contains $u^{(0)}$. Let $r_1$, \dots , $r_s$ be the multiplicities of the eigenvalues of $U\big(u^{(0)}\big)=\operatorname{diag}\big(u_1^{(0)},\dots ,u_n^{(0)}\big)$, with $s< n$, $r_1+\cdots +r_s=n$. By a permutation of $(u_1,\dots ,u_n)$, there is no loss in generality (cf.\ Section~\ref{freedom}) if we assume that the entries of $u^{(0)}$ are
\begin{gather}
u_1^{(0)}=\cdots=u_{r_1}^{(0)}=:\lambda_1,\nonumber\\
u_{r_1+1}^{(0)}=\cdots=u_{r_1+r_2}^{(0)}=:\lambda_2,\nonumber\\
 \cdots\cdots\cdots\cdots\cdots\cdots\cdots\cdots\cdots\nonumber\\
u_{r_1+\cdots+r_{s-1}+1}^{(0)}=\cdots=u_{r_1+\cdots+r_{s-1}+r_s-1}^{(0)}=u_n^{(0)}=:\lambda_s,\label{30gen2016-1}
\end{gather}
$\lambda_k\neq \lambda_l$ for $k\neq l$.
Let
\begin{gather*}
\delta_i:=\frac{1}{2}\min \big\{ \big|\lambda_i-\lambda_j +\rho {\rm e}^{\sqrt{-1}\left(\frac{\pi}{2}-\phi\right)}\big|, ~ j\neq i, ~ \rho\in\mathbb{R}\big\},
\end{gather*}
and let $\epsilon_0$ be a small positive number such that
\begin{gather}\label{2agosto2016-1}
\epsilon_0<\min_{1\leq i\leq s} \delta_i.
\end{gather}
Notice that $\delta_i>0$, because $\phi$ is admissible.
We will assume that $\epsilon_0$ is sufficiently small so that the polydisc at $u^{(0)}$, defined by\footnote{ Here $\overline{B}(\lambda_i;\epsilon_0)$ is the closed ball in $\mathbb{C}$ with center $\lambda_i$ and radius $\epsilon_0$. Note that if the uniform norm $|u|=\max_i |u_i|$ is used, as in \cite{CDG0}, then
$\mathcal{U}_{\epsilon_0}\big(u^{(0)}\big)=\big\{u\in \mathbb{C}^n\,\bigr|\,\big|u-u^{(0)}\big|\leq \epsilon_0\big\}.$}
\begin{gather*}
\mathcal{U}_{\epsilon_0}\big(u^{(0)}\big):=\bigtimes_{i=1}^s\overline{B}(\lambda_i;\epsilon_0)^{\times r_i},
\end{gather*}
 is completely contained in the image $u(\Omega)$ of the chart $\Omega$.

\begin{Lemma} For $\epsilon_0$ satisfying \eqref{2agosto2016-1}, if $u$ varies in $ \mathcal{U}_{\epsilon_0}\big(u^{(0)}\big)$, the sets
 \begin{gather}
I_1:=\{u_1, \dots  ,u_{r_1}\},\quad I_2:= \{u_{r_1+1}, \dots  ,u_{r_1+r_2}\},\quad \dots ,\nonumber\\
 I_s:=\{u_{r_1+\cdots+r_{s-1}+1}, \dots ,u_{r_1+\cdots+r_{s-1}+r_s}\}\label{19-05-17.1}
\end{gather}
do never intersect. Thus, $u^{(0)}$ is a point of {\it maximal coalescence} in $\mathcal{U}_{\epsilon_0}\big(u^{(0)}\big)$. We will say that a coordinate $u_a$ {\it is close to $\lambda_j$} if it belongs to $I_j$, which is to say that $u_a\in\overline{B}(\lambda_j;\epsilon_0)$.
\end{Lemma}

Let us fix $\phi\in\mathbb{R}$ so that the line $\ell_+(\phi)$, $\ell_-(\phi)$, $\ell(\phi)$ are admissible at~$p_0$ (Definition~\ref{5luglio2017-7}).
 For $u\in \mathcal{U}_{\epsilon_0}\big(u^{(0)}\big)$, with $\epsilon_0$ as in~(\ref{2agosto2016-1}), consider the subset $\mathfrak{R}(u)$ of Stokes rays $R_{ab,k}$ in the universal covering $\mathcal{R}$ which are associated with all couples of eigenvalues~$u_a$ and~$u_b$ such that~$u_a$ is close to a~$\lambda_i$ and~$u_b$ is closed to $\lambda_j$ {for some} $i\neq j$. Then, the following holds:

 \begin{Lemma} Let $\epsilon_0$ be as in~\eqref{2agosto2016-1}. If $u_a$ varies in $\overline{B}(\lambda_i;\epsilon_0)$ and~$u_b$ in $\overline{B}(\lambda_j;\epsilon_0)$, then the rays $R_{ab,k}\in\mathfrak{R}(u)$ continuously rotate, but they never cross $\ell_+(\phi)$ and $\ell_{-}(\phi)$. In other words, the projections $R_{ab}=\operatorname{pr}(R_{ab,k})$ never cross $\ell(\phi)$ in~$\mathbb{C}$.
\end{Lemma}

The choice of the line $\ell$, admissible at~$p_0$, induces a cell decomposition of $\mathcal{U}_{\epsilon_0}\big(u^{(0)}\big)$, according to the following

\begin{Definition}\label{ellcell}Let $\ell$ be admissible at $u^{(0)}$. An $\ell$-cell of $\mathcal{U}_{\epsilon_0}\big(u^{(0)}\big)$ is any connected component of the open dense subset of points $u\in \mathcal{U}_{\epsilon_0}\big(u^{(0)}\big)$ such that $u_1,\dots ,u_n$ are pairwise distinct and~$\ell$ is admissible at~$u$.
\end{Definition}

\begin{Proposition}[\cite{CDG0}]An $\ell$-cell is homeomorphic to a ball.
\end{Proposition}

\begin{figure}[h]\centering
\def\svgscale{0.7}
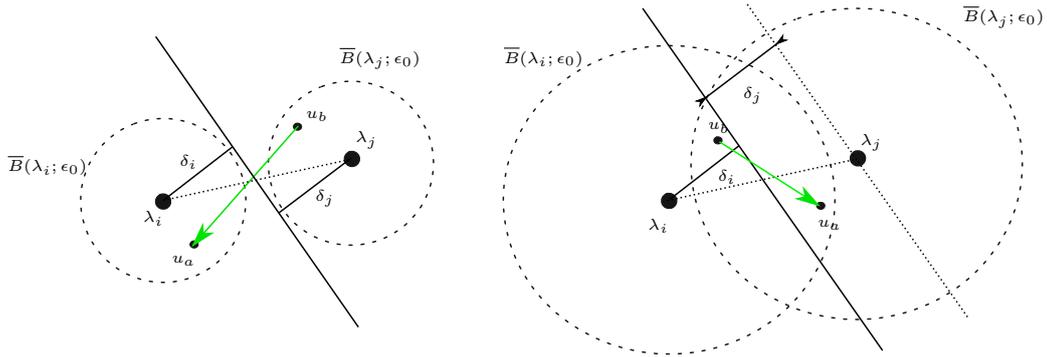
\caption{Points $\lambda_i$, $\lambda_j$ and $u_a$, $u_b$ are represented on the same complex plane. The thick line has slope $\pi/2-\phi$. As $u$ varies, for values of $\epsilon_0$ sufficiently small (left figure) the Stokes rays $R_{ab}$ and $R_{ba}$ associated {with} $u_a$ in the disk $\overline{B}(\lambda_i;\epsilon_0)$ and $u_b$ in the disk of $ \overline{B}(\lambda_j;\epsilon_0)$ do not cross the line $\ell$. If the disks have radius exceeding $\min_{1\leq i\leq s} \delta_i$ as in~(\ref{2agosto2016-1}) (see right figure) then the Stokes rays $R_{ab}$, $R_{ba}$ cross the line~$\ell$.}
\end{figure}

We notice that, if $u(p)$ is in a $\ell$-cell then $p$ lies in {an} $\ell$-chamber. Thus, if $\mathcal{D}$ is an open subset whose closure is contained in a cell of $\mathcal{U}_{\epsilon_0}\big(u^{(0)}\big)$, according to Theorem~\ref{fundamental}(3), the system
\begin{gather}\label{23gen2016-1}
\frac{dY}{dz}=\left(U+\frac{V(u)}{z}\right)Y,
\end{gather}
for $u\in \mathcal{D}$ admits two fundamental solutions $Y^{(0)}_{\rm right/left}(z,u)$ uniquely determined by the canonical asymptotic representation $Y^{(0)}_{\rm right/left}(z,u)\sim Y_{\rm formal}(z,u)$ as in~\eqref{7luglio2017-1} valid in the sectors $\Pi_{\rm left/right}(\phi)$ respectively. It follows from the proof of Theorem~\ref{formalredu} that $Y_{\rm formal}(z,u)$ is
 $u$-holomorphic in $\mathcal{U}_{\epsilon_0}\big(u^{(0)}\big)\backslash\Delta_\Omega$.
By Remark~\ref{26.05.17-1} actually the asymptotic representation is valid in wider sectors $\mathcal S_{\rm left/right}(u)$, defined as the sectors which contain $\Pi_{\rm left/right}(\phi)$ and extends up to the nearest Stokes rays. By Theorem~\ref{isoth2} the {above} system with $u\in\mathcal D$ is isomonodromic, so that the Stokes matrix~$S$ defined in formula~\eqref{matrstok1} is constant.

Let us now turn our attention to the coalescence point $u^{(0)}$. From the results of \cite{CDG0} -- and more generally in~\cite{BJL2}~-- it follows that there are a {\it unique} formal solution at $u^{(0)}$,
\begin{gather}\label{7febbraio2020-3}
\mathring{Y}_{\rm formal}(z)=\left(\mathbbm 1+\sum_{k=1}^\infty\frac{ \mathring{G}_k}{ z^{k}}\right){\rm e}^{zU},
\end{gather}
and unique actual solutions $\mathring{Y}^{(0)}_{\rm left}(z)$ and $\mathring{Y}^{(0)}_{ \rm right}(z)$, with asymptotic representation given by $\mathring{Y}_{\rm formal}(z)$ in $\Pi_{\rm left/right}$, and
 in wider sectors $\mathcal{S}_{\rm left}\big(u^{(0)}\big)$ and $\mathcal{S}_{ \rm right}\big(u^{(0)}\big)$ respectively. The Stokes matrices of $\mathring{Y}^{(0)}_{\rm right}(z)$ and $\mathring{Y}^{(0)}_{\rm left}(z)$ are defined by
\begin{gather*}
\mathring{Y}^{(0)}_{\rm left}(z)=\mathring{Y}^{(0)}_{\rm right}(z)\mathring{S},
\qquad \mathring{Y}^{(0)}_{\rm left}\big({\rm e}^{2\pi{\rm i}}z\big)=\mathring{Y}^{(0)}_{\rm right}(z)\mathring{S}_{-},
\qquad\mathring{S}_{-}=\mathring{S}^{\rm T}.
\end{gather*}
A priori, the following problems could emerge.
\begin{enumerate}\itemsep=0pt
\item The asymptotic representations \begin{gather*}Y^{(0)}_{\rm left/right}(z,u)\sim Y_{\rm formal}(z,u),\qquad \text{for } |z|\to \infty \qquad \text{and}\\
    z\in\bigcap_{u\in\mathcal D}\mathcal S_{\rm left/right}(u)\supsetneq \Pi_{\rm left/right}(\phi)\end{gather*}
{do} no longer hold for $u$ outside the cell containing $\mathcal D$.
\item The coefficients $G_k(u)$'s of \eqref{7luglio2017-1} may be divergent at~$\Delta_\Omega$.
\item The locus $\Delta_\Omega$ is expected to be a locus of singularities for the solutions $Y_{\rm formal}(z,u)$ in~(\ref{7luglio2017-1}) and $Y^{(0)}_{\rm left/right}(z,u)$.
$Y_{\rm formal}(z,u)$.
\item The Stokes {matrix $S$} may differ from $\mathring{S}$.
\end{enumerate}

We notice that the system (\ref{23gen2016-1}) at $u^{(0)}$ also has a fundamental solution in Levelt form at $z=0$,
\begin{gather}\label{4agosto2016-6}
\mathring{Y}_0(z)= \Psi\big(u^{(0)}\big)(I+\mathcal{O}(z))z^\mu z^{\mathring{R}},
\end{gather}
with a certain exponent $\mathring{R}$. Hence, a central connection matrix $\mathring{C}$ is defined by
\begin{gather*}
 \mathring{Y}^{(0)}_{\rm right}(z)=\mathring{Y}_0(z)\mathring{C}.
\end{gather*}

We recall that $\Psi(u)$ is holomorphic in the whole $\mathcal{U}_{\epsilon_0}\big(u^{(0)}\big)$, so that $V_{ij}(u)$ vanishes along~$\Delta_\Omega$ whenever $u_i=u_j$ (see Lemma~\ref{30nov2016-9}). These are sufficient conditions to apply the main theorem of~\cite{CDG0}, adapted and particularised to the case of Frobenius manifolds, which becomes the following:

\begin{Theorem}\label{mainisoth}Let $M$ be a semisimple Frobenius manifold, $p_0\in\mathcal B_M\setminus\mathcal K_M$ and $\Omega\subseteq M_{ss}=M\setminus\mathcal K_M$ an open connected neighborhood of~$p_0$ with the property of Remark~{\rm \ref{26.05.17-3}} on which a~holomorphic branch for canonical coordinates $u\colon\Omega\to\mathbb C^n$ and $\Psi\colon\Omega\to {\rm GL}_n(\mathbb C)$ has been fixed. Let~$\epsilon_0$ be a~real positive number as above, and consider the corresponding neighborhood $\mathcal{U}_{\epsilon_0}\big(u^{(0)}\big)$ of $u^{(0)}=u(p_0)$.
 Then
\begin{enumerate}\itemsep=0pt
\item[$1.$] The coefficients $G_k(u)$, $k\geq 1$, in~\eqref{7luglio2017-1} are holomorphic over $\mathcal{U}_{\epsilon_0}\big(u^{(0)}\big)$,
\begin{gather*}G_k\big(u^{(0)}\big)=\mathring{G}_k
\qquad \text{and} \qquad
Y_{\rm formal}\big(z,u^{(0)}\big)=\mathring{Y}_{\rm formal}(z).
\end{gather*}
\item[$2.$] $Y^{(0)}_{\rm left}(z,u)$, $Y^{(0)}_{\rm right}(z,u)$, can be $u$-analytically continued as single-valued holomorphic functions\footnote{Hence, they are holomoprhic on $\mathcal{R}\times\mathcal{U}_{\epsilon_0}\big(u^{(0)}\big)$.} on $\mathcal{U}_{\epsilon_0}\big(u^{(0)}\big)$. Moreover
\begin{gather*}
Y^{(0)}_{\rm left/right}\big(z,u^{(0)}\big)=\mathring{Y}^{(0)}_{\rm left/right}(z).
\end{gather*}

\item[$3.$] For any solution $\mathring{Y}_0(z)$ as in~\eqref{4agosto2016-6} there exists a fundamental solution $Y_0(z,u)$ in Levelt form~\eqref{30nov2016-6} such that
\begin{gather*}
Y_0\big(z,u^{(0)}\big)=\mathring{Y}_0(z), \qquad R=\mathring{R}.
\end{gather*}

\item[$4.$] For any {positive} $\epsilon_1<\epsilon_0$, the asymptotic relations
\begin{gather}\label{27.05.17-2}
Y^{(0)}_{\rm left/right}(z,u)\sim \left(\mathbbm 1+\sum_{k=1}^\infty \frac{G_k(u)}{z^k}\right){\rm e}^{zU},\qquad z\to\infty \hbox{ in }\Pi_{\rm left/right}(\phi),
\end{gather}
hold uniformly in $u\in \mathcal{U}_{\epsilon_1}\big(u^{(0)}\big)$. In particular they hold also at points of $\Delta_\Omega\cap \mathcal{U}_{\epsilon_1}\big(u^{(0)}\big)$ and at~$u^{(0)}$.

\item[$5.$] For any $u\in\mathcal{U}_{\epsilon_0}\big(u^{(0)}\big)$ consider the sectors $\widehat{\mathcal{S}}_{\rm right}(u)$ and $\widehat{\mathcal{S}}_{\rm left}(u)$ which contain the sectors $\Pi_{\rm right}(\phi)$ and $\Pi_{\rm left}(\phi)$ respectively, and extend up to the nearest Stokes rays in the set $\mathfrak{R}(u)$ defined above. Let
\begin{gather*}
\widehat{\mathcal{S}}_{\rm left/right}=\bigcap_{u \in \mathcal{U}_{\epsilon_0}\big(u^{(0)}\big)} \widehat{\mathcal{S}}_{\rm left/right}(u).
\end{gather*}
Observe that for sufficiently small $\varepsilon>0$ the sectors \begin{gather*}
\Pi^\varepsilon_{\rm{right}}(\phi):= \{z\in\mathcal R\colon \phi-\pi-\varepsilon< \arg z< \phi+\varepsilon \},\\
\Pi^\varepsilon_{\rm{left}}(\phi):= \{z\in\mathcal R\colon \phi-\varepsilon< \arg z< \phi+\pi+\varepsilon \}
\end{gather*} are strictly contained in $ \widehat{\mathcal{S}}_{\rm right}$ and $ \widehat{\mathcal{S}}_{\rm left}$ respectively. Then, the asymptotic relations~\eqref{27.05.17-2} actually hold in the sectors $\widehat{\mathcal{S}}_{\rm left/right}$.

\item[$6.$] The monodromy data $\mu$, $R$, $C$, $S$ of system~\eqref{23gen2016-1}, defined and constant in an open subset~$\mathcal{D}$ of a cell of $\mathcal{U}_{\epsilon_0}\big(u^{(0)}\big)$, are actually defined and constant at any $u\in \mathcal{U}_{\epsilon_1}\big(u^{(0)}\big)$, namely the system is isomonodromic in $\mathcal{U}_{\epsilon_1}\big(u^{(0)}\big)$. They coincide with the data $\mu$, $\mathring{R}$, $\mathring{C}$, $\mathring{S}$ associated with the fundamental solutions $\mathring{Y}_{\rm left/right}(z)$ and $\mathring{Y}_0(z)$ of system~\eqref{23gen2016-1} at $u^{(0)}$. The entries of $S=(S_{ij})_{i,j=1}^n$ satisfy the vanishing condition~\eqref{9dicembre2016-1}, namely
\begin{gather}\label{9dicembre2016-2}
S_{ij}=S_{ji}=0 \qquad \text{for all $i\neq j$ such that} \quad u_i^{(0)}=u_j^{(0)}.
\end{gather}
\end{enumerate}
\end{Theorem}

This Theorem allows us to obtain the monodromy data $\mu$, $R$, $C$, $S$ in a~neighbourhood of a~coalescence point just by computing them \emph{at} the coalescence point, namely just by computing $\mu$, $\mathring{R}$, $\mathring{C}$, $\mathring{S}$. Its importance has been explained in the Introduction and will be illustrated in subsequent sections.

\begin{Remark}\label{9dicembre2016-3}Suppose that $S$ is upper triangular. By formula~\eqref{9dicembre2016-2} it follows that in any $\ell$-cell of $\mathcal{U}_{\epsilon_0}\big(u^{(0)}\big)$ the order of the canonical coordinates is triangular, according to Definition \ref{30luglio2016-2}, and at most in one cell the order is lexicographical (Definition~\ref{30luglio2016-1}). Moreover, if $\lambda_1$, \dots , $\lambda_s$ are in lexicographical order, then the order of canonical coordinates is lexicographical in exactly one cell.
\end{Remark}

\subsection{Reconstruction of monodromy data of the whole manifold}\label{30luglio2016-3}
The monodromy data of the Frobenius manifold can be obtained from those computed in Theorem~\ref{mainisoth} around $u^{(0)}$. Without loss of generality, let us suppose that the ordering \eqref{30gen2016-1} is such that $\lambda_1,\dots,\lambda_s$ are in $\ell$- lexicographical order. Then, the matrix~$ S$ computed at the coalescence point~$u^{(0)}$ is upper triangular. Therefore, by Theorem~\ref{mainisoth}, the matrix is constant and  upper triangular in the whole polydisc $\mathcal{U}_{\epsilon_0}\big(u^{(0)}\big)$. In particular, it is upper triangular in every  cell of $\mathcal{U}_{\epsilon_1}\big(u^{(0)}\big)$. This means that $u_1,\dots ,u_n$ are in triangular order  (Definition~\ref{30luglio2016-2}) in each such cell, and in particular they are in lexicographical order in only one  of these cells (Definition~\ref{30luglio2016-1}, Remark~\ref{9dicembre2016-3}). Note that any permutation of canonical coordinates preserving the  sets $I_1,\dots, I_s$ of~\eqref{19-05-17.1} maintains the upper triangular structure of $S$, namely the triangular  order of $u_1,\dots ,u_n$ in each cell of $\mathcal{U}_{\epsilon_o}\big(u^{(0)}\big)$. The permutation changes the cell where  the order is lexicographical. Now, each cell of the polydisc $\mathcal{U}_{\epsilon_0}\big(u^{(0)}\big)$ is contained in a chamber of the manifold (identifying coordinates with points of the manifold, which is possible because of the holomorphy of canonical coordinates near \emph{semisimple} coalescent points). Let us start from the cell of $\mathcal{U}_{\epsilon_0}\big(u^{(0)}\big)$ where $u_1,\dots ,u_n$ are in lexicographical order. The monodromy data of Theorem~\ref{mainisoth} in this cell are the constant data of the chamber containing the cell (Theorems~\ref{isomonod1} and~\ref{isoth2}). Since in this chamber $u_1,\dots ,u_n$ are in lexicographical order (and distinct!), we can apply the action of the braid group to $S$ and $C$, as dictated by formulae~\eqref{stokesbraid1},~\eqref{connbraid2}. In this way, the monodromy data for any other chamber of the manifold are obtained, as explained in Section~\ref{freedom}.

\section[First detailed example of application of Theorem \ref{mainisoth}: the $A_3$ Frobenius\\ manifold. Stokes phenomenon for Pearcey-type oscillating integrals\\ from Hankel functions]{First detailed example of application of Theorem \ref{mainisoth}:\\ the $\boldsymbol{A_3}$ Frobenius manifold. Stokes phenomenon\\ for Pearcey-type oscillating integrals from Hankel functions}\label{A_3case}

 With the example of $A_3$ Frobenius manifold below, we show how Theorem~\ref{mainisoth} allows the computation of monodromy data in an elementary way, by means of Hankel {\it special functions}. Moreover, we apply the results of Section~\ref{freedom}, especially showing how the braid group can be used to reconstruct the data for the whole manifold, starting from a coalescence point. The reader not interested in a general introduction to Frobenius manifolds associated with singularity theory may skip Sections~\ref{28luglio2016-1} and~\ref{28luglio2016-2} and go directly to Section~\ref{28luglio2016-3}.

 \subsection{Singularity theory and Frobenius manifolds} \label{28luglio2016-1}
 Let $f$ be a quasi-homogeneous polynomial on~$\mathbb C^m$ with an isolated simple singularity at $0\in\mathbb C^m$. According to V.I.~Arnol'd~\cite{arn} simple singularities are classified by simply-laced Dynkin diagrams~$A_n$ (with $n\geq 1$), $D_n$ (with $n\geq 4$), $E_6$, $E_7$, $E_8$. Denoting by $(x_1\dots, x_m)$ the coordinates in~$\mathbb C^m$ (for singularities of type~$A_n$ we consider $m=1$), the classification of simple singularities is summarized in Table~\ref{singADE}. Let~$\mu$ be the Milnor number of~$f$ (note that $\mu=n$ for~$A_n$, $D_n$ and~$E_n$), and
 \[
 f(x,a):=f(x)+\sum_{i=1}^\mu a_i\phi_i(x)
 \]
 be a miniversal unfolding of $f$, where $a$ varies in a ball $B\subseteq \mathbb C^\mu$, and $(\phi_1(x),\dots,\phi_\mu(x))$ is a~basis of the Milnor ring.  Using K.~Saito's theory of \emph{primitive forms}~\cite{Saito2}, a~flat metric and a~Frobenius manifold structure can be defined on the base space~$B$ \cite{BV}. See also~\cite{Saito1}. For any fixed $a\in B$, let the critical points be $x_i(a)=\big(x^{(1)}_i,\dots, x^{(m)}_i\big)$, $i=1,\dots, \mu$, defined by the condition $\partial_{x_\alpha}f(x_i, a)=0$ for any $\alpha=1,\dots, m$, with critical values $u_i(a):=f(x_i(a),a)$. The open ball~$B$ can be stratified as follows:
 \begin{enumerate}\itemsep=0pt
 \item[1)] the stratum of \emph{generic points}, i.e., points where both critical points $x^{(i)}$'s and critical values $u_i$'s are distinct;
 \item[2)] the \emph{Maxwell stratum}, which is the closure of the set of points with distinct critical points $x^{(i)}$'s but some coalescing critical values $u_i$'s;
 \item[3)] the \emph{caustic}, where some critical points coalesce.
 \end{enumerate}
The union of the Maxwell stratum and the caustic is called \emph{function bifurcation diagram} $\Xi$ of the singularity (see~\cite{singularity1} and~\cite{AGZV-I,AGZV-II}).

 The complement of the caustic consists exclusively of semisimple points of the Frobenius manifold, and the critical values $u_i(a):=f(x_i(a),a)$ are the canonical coordinates.

 In this section we want to show how one can reconstruct local information near semisimple points in the Maxwell stratum, by invoking Theorem~\ref{mainisoth}. We will focus on the simplest example of $A_3$.

\begin{Remark}For simplicity of exposition, here we focus only on the case of Frobenius structures associated with simple singularities. The existence of primitive forms for arbitrary singularities was proved by M.~Saito in~\cite{msaito}.
\end{Remark}

\begin{table}\centering
\begin{tabular}{|c||l|l|}\hline
& \multicolumn{1}{|c|}{Singularity} & \multicolumn{1}{|c|}{Versal Deformation}\\
\hline
\hline
$A_n$& $f(x)=x^{n+1}$ & $f(x,a)=x^{n+1}+a_{n-1}x^{n-1}+\dots+a_1 x+a_0$\tsep{2pt}\bsep{2pt}\\
\hline
$D_n$& $f(x)= x_1^{n-1}+x_1 x_2^2$ & $f(x,a)=x_1^{n-1}+x_1 x_2^2+a_{n-1}x_1^{n-2}+\dots+a_1+a_0 x_2$\tsep{2pt}\bsep{2pt} \\
\hline
$E_6$ & $f(x)=x_1^4+x_2^3$ & $f(x,a)=x_1^4+x_2^3+a_6x_1^2x_2 +a_5x_1x_2$\tsep{2pt} \\
 & & $\hphantom{f(x,a)=}{}+a_4x^2_1 +a_3x_2 +a_2x_1 +a_1$\bsep{2pt} \\
\hline
$E_7$ & $f(x)=x_1^3x_2+x_2^3$ & $f(x,a)=x^3_1x_2 +x^3_2 +a_7x_2 +a_6x^2_1 +a_5x_1 $\tsep{2pt}\\
 & & $\hphantom{f(x,a)=}{} +a_4x_1x_2+a_3x_1x_2+a_2x_2+a_1$\bsep{2pt}\\
\hline
$E_8$ & $f(x)=x^5_1 +x^3_2$ & $f(x,a)=x^5_1 +x^3_2 +a_8x^3_1x_2 +a_7x^2_1x_2 +a_6x^3_1$\tsep{2pt}\\
 & & $\hphantom{f(x,a)=}{}+a_5x_1x_2 +a_4x^2_1 +a_3x_2 +a_2x_1 +a_1$\bsep{2pt} \\
\hline
\end{tabular}
\caption{Arnol'd's classification of simple singularities, and their corresponding miniversal deformations.}\label{singADE}
\end{table}

\subsection[Frobenius structure of type $A_n$]{Frobenius structure of type $\boldsymbol{A_n}$}\label{28luglio2016-2}

General references for this section are \cite{DVV,dubro1,dubro5,dubro2}. Let us consider the affine space $M\cong\mathbb C^n$ of all polynomials
\[f(x,a)=x^{n+1}+a_{n-1}x^{n-1}+\dots+a_1x+a_0,
\]where $(a_0,\dots,a_{n-1})\in M$ are used as coordinates. We call \emph{bifurcation diagram} $\Xi$ of the singularity $A_n$ the set of polynomials in $M$ with some coalescing critical values. The bifurcation diagram $\Xi$ is an algebraic subvariety in $M$ which consists of two irreducible components (the derivative w.r.t.\ the variable $x$ will be denoted by $(\cdot)'$):

\begin{itemize}\itemsep=0pt
\item the caustic $\mathcal K$, which is the set of polynomials with degenerate critical points (i.e., solutions of the system of equations $f'(x,a)=f''(x,a)=0$);\footnote{The equation of the caustic is $\Delta(f')=0$, where $\Delta(f'):=\operatorname{Res}(f',f'')$ is the discriminant of the polyno\-mial~$f'(x,a)$. The reader can consult the monograph \cite[Chapter~12]{GKZ}.}
\item The Maxwell stratum $\mathcal M$, defined as the closure of the set of polynomials with some coalescing critical values but different critical points.
\end{itemize} For more information about the topology and geometry of (the complement of) these strata, the reader can consult the paper \cite{Nekrasov}, and the monograph \cite{vassiliev}.
There is a naturally defined covering map $\rho\colon\widetilde{M}\to M$ of degree $n!$, whose fiber over a point $f(x,a)$ consists of total orderings of its critical points. On $\widetilde M$, $x_1,\dots, x_n$ are well defined functions such that
\[f'(x,\rho(w))=(n+1)\prod_{i=1}^n(x-x_i(w)),\qquad w\in\widetilde M.
\]
The caustic $\mathcal K$ is the ramification locus of the covering $\rho$. For any simply connected open subset $U\subseteq M\setminus \mathcal K$, we can choose a connected component $W$ of $\rho^{-1}(U)$. The restriction of the functions $x_1, \dots , x_n$ on~$W$ defines single-valued functions of $a\in U$, which are local branches of $x_1, \dots , x_n$. For further details see~\cite{manin}.

We define on $M$ the following structures:
\begin{enumerate}\itemsep=0pt
\item A free sheaf of rank $n$ of $\mathcal O_M$-algebras: this is the sheaf of Jacobi--Milnor algebras
\[\frac{\mathcal O_M[x]}{f'(x,a)\cdot\mathcal O_M[x]}.
\]For fixed $a\in M$, the fiber of this sheaf is the algebra $\mathbb C[x]/\langle f^\prime(x,a)\rangle$. We also define an $\mathcal O_M$-linear \emph{Kodaira--Spencer} isomorphism $\kappa\colon \mathscr T_M\to\mathcal O_M[x]/\langle f'(x,a)\rangle$ which associates to a vector field $\xi$ the class $\mathfrak L_\xi (f)=\xi (f) \mod f'$. In particular, for any $\alpha=0,\dots,n-1$ the class $\partial_{a_i}f$ is associated {with} the vector field $\partial_{a_i}$. In this way we introduce a product $\circ$ of vector fields defined by
\[\xi\circ\zeta:=\kappa^{-1}\left(\xi (f) \cdot \zeta (f)\mod f'\right).
\]
The product $\circ$ is associative, commutative and with unit $\partial_{a_0}$. We call {\it Euler vector field} the distinguished vector field~$E$ corresponding to the class $f\mod f'$ under the Kodaira--Spencer map $\kappa$. An elementary computation shows that
\[
E=\sum_{i=0}^{n-1}\frac{n+1-i}{n+1}a_i\frac{\partial}{\partial a_i},
\qquad
\mathfrak L_E(\circ)=\circ.
\]
\item A \emph{symmetric bilinear form} $\eta$, defined at a fixed point $a\in M$ as the Grothendieck residue
\begin{gather}\label{metricaAn}\eta_a(\xi,\zeta):=\frac{1}{2\pi{\rm i}}\int_{\Gamma_a}\frac{\xi(f)(u,a)\cdot\zeta(f)(u,a)}{f'(u,a)}\,{\rm d}u,
\end{gather}where $\Gamma_a$ is a circle, positively oriented, bounding a disc containing all the roots of $f'(u,a)$. It is a nontrivial fact that the bilinear form $\eta$ is non-degenerate (for a proof, see \cite{AGZV-I,AGZV-II}) and \emph{flat} (explicit flat coordinates can be found in \cite{SYS}: notice that the natural coordinates $a_i$'s are not flat). Notice that
\begin{gather*}\mathfrak L_E\eta=\frac{n+3}{n+1}\eta.\end{gather*}
\end{enumerate}

\begin{Theorem}\label{5agosto2016-1}The manifold $M$, endowed with the tensors $(\eta,\circ, \partial_{a_0},E)$, is a Frobenius manifold of charge $\frac{n-1}{n+1}$. The caustic $\mathcal{K}_M$, defined as in Definition~{\rm \ref{caustic}}, coincides with the caustic $\mathcal K$ of the singularity $A_n$ defined above. By analytic continuation, the semisimple Frobenius structures extends on the unramified covering space $\rho^{-1}(M\backslash\mathcal K)\subseteq\widetilde M$. Critical values define a system of canonical coordinates.
\end{Theorem}

The reader can find detailed proofs in \cite{dubro1,dubro2,manin,sabbah}. If $a$ is a given point of $M\setminus\mathcal K$, i.e., such that $f(x,a)$ has $n$ distinct Morse critical points ${x}_1,\dots,{x}_n$, then the elements
\[
\pi_i(a):=\kappa^{-1}\left(\frac{f'(x,a)}{ f''({x}_i,a)(x-x_i) }\right)\qquad\text{for }i=1,\dots, n
\]
are idempotents of $(T_aM,\circ_{a})$. This follows from the equality $f'(x,a)=(n+1)\prod\limits_{i=1}^n(x-x_i)$. Consider now {the critical points} $x_1(a),\dots,x_n(a)$ {locally well defined as functions of } $a$ varying in a simply connected open set away from the caustic. {The critical values} $u_i(a):=f(x_i(a),a)$ for $i=1,\dots, n$ can also be considered as functions on the same set. Since $\det\big(\frac{\partial u_i}{\partial a_j}\big)$ is the Vandermonde determinant of $x_i(a)$'s, the functions $u_i$'s define a system of local coordinates on~$M$. In order to see that $\pi_i\equiv\frac{\partial}{\partial u_i}$, it is sufficient to prove that $\kappa(\partial{u_i})(x_j)=\delta_{ij}$, i.e., $\frac{\partial f}{\partial u_i}(x_i)=\delta_{ij}$. This follows from the equalities
\[
\frac{\partial f (x_i(a),a)}{\partial a_j}=(x_i(a))^j,\qquad \frac{\partial }{\partial a_j}=\sum_i(x_i(a))^j\frac{\partial }{\partial u_i}.
\]

\subsection[The case of $A_3$: reduction of the system for deformed flat coordinates]{The case of $\boldsymbol{A_3}$: reduction of the system for deformed flat coordinates}\label{28luglio2016-3}

We consider the space $M$ of polynomials
\[f(x;a)=x^4+a_2x^2+a_1x+a_0,
\]where $a_0,a_1,a_2\in\mathbb C$ are ``natural'' coordinates on $M$. The residue theorem implies that the metric $\eta$, defined on $M$ as in~\eqref{metricaAn}, can be expressed as
\[
\eta_a(\xi,\zeta)=-\mathop{\operatorname{res}}\limits_{u=\infty}\frac{\xi(f)(u,a)\cdot\zeta(f)(u,a)}{f'(u,a)}{\rm d}u,
\]
and consequently
\[
\eta_{a}(\partial_i,\partial_j)=\mathop{\operatorname{res}}\limits_{v=0}\frac{v^{1-i-j}}{4+2a_2v^2+a_1v^3}{\rm d}v,
\]
where $\partial_i=\frac{\partial}{\partial a_i}$, $\partial_j=\frac{\partial}{\partial a_j}$. So we find that
\[\eta_a=\begin{pmatrix}
0&0&\dfrac{1}{4}\\
0&\dfrac{1}{4}&0\\
\dfrac{1}{4}&0&-\dfrac{a_2}{8}
\end{pmatrix}.
\]
Note that $a_0$, $a_1$, $a_2$ are not flat coordinates for $\eta$. The commutative and associative product defined on each tangent space $T_aM$, using the Kodaira--Spencer map, is given by the {structure} constants at a generic point $a\in M$:
\begin{gather*} \partial_0\circ\partial_i=\partial_i\qquad\text{for all }i,
\\
\partial_1\circ\partial_1=\partial_2,\qquad
\partial_1\circ\partial_2=-\frac{1}{2}a_2\partial_1-\frac{1}{4}a_1\partial_0,
\qquad
\partial_2\circ\partial_2=-\frac{1}{2}a_2\partial_2-\frac{1}{4}a_1\partial_1.
\end{gather*}
The Euler vector field is
 \begin{gather*}
 E:=\sum_{i=0}^2\frac{4-i}{4}a_i\partial_i=a_0\partial_0+\frac{3}{4}a_1\partial_1+\frac{1}{2}a_2\partial_2.
 \end{gather*}
 With such a structure $M$ is a Frobenius manifold. The (1,1)-tensor $\mathcal U$ of multiplication by $E$ is
\[\mathcal U(a)=\left(
\begin{matrix}
 a_0 & -\dfrac{1}{8} a_1 a_2 & -\dfrac{3}{16} a_1^2 \\
 \dfrac{3 a_1}{4} & a_0-\dfrac{a_2^2}{4} & -\dfrac{1}{2} a_1 a_2 \\
 \dfrac{a_2}{2} & \dfrac{3}{4}a_1 & a_0-\dfrac{a_2^2}{4}
\end{matrix}
\right).
\]
 Up to a multiplicative constant, the discriminant of the characteristic polynomial of $\mathcal U$ is equal to
\[a_1^2 \big(8 a_2^3+27 a_1^2\big)^3
\]and so the bifurcation set of the Frobenius manifold is the locus
\[
\mathcal{B}=\{a_1=0\}\cup\big\{8 a_2^3+27 a_1^2=0\big\}.
\]
The irreducible component $\{a_1=0\}$ is the Maxwell stratum, whereas the irreducible component $\big\{8 a_2^3+27 a_1^2=0\big\}$ is the caustic.

Let us focus on the set $\{a_1=0\}$, and let us look for semisimple points on it. It is enough to consider the multiplication by the vector field $\lambda\partial_1+\mu\partial_2$ ($\lambda,\mu\in\mathbb C$), and show that it has distinct eigenvalues. This is a (1,1)-tensor with components at points $(a_0,a_1,a_2)$ equal to
\[\left(
\begin{matrix}
 0 & -\dfrac{\mu}{4}a_1 & -\dfrac{\lambda}{2}a_1 \\
 \lambda & -\dfrac{\mu}{2} a_2 & -\dfrac{\lambda}{2} a_2-\dfrac{\mu}{4}a_1 \\
 \mu & \lambda & -\dfrac{\mu}{2} a_2
\end{matrix}
\right),
\]whose characteristic polynomial, at points $(a_0,0,a_2)$, has discriminant
\[-\frac{1}{8} \lambda^2 a_2^3 \big(2 \lambda^2+\mu^2 a_2\big)^2.
\]
So, the points $(a_0,0,a_2)$ with $a_2\neq 0$ are semisimple points of the bifurcation set, namely they belong to the Maxwell stratum. In view of Theorem~\ref{5agosto2016-1}, they are semisimple coalescence points of Definition~\ref{3agosto2016-9}. We would like to study deeper the behavior of the Frobenius structure near points $(a_0,a_1,a_2)=(0,0,h)$ of the Maxwell stratum, with fixed $a_0=0$ and with $h\in\mathbb C^*$.

\begin{Remark}The points $(a_0,0,0)$, instead, are not semisimple because we have evidently $\partial_2^2=0$ on them.
\end{Remark}

Let us introduce flat coordinates $t_1$, $t_2$, $t_3$ defined by
\[
\begin{cases}
a_0=t_1+\dfrac{1}{8}t^2_3,\\
a_1=t_2,\\
a_2=t_3,
\end{cases}\qquad
J=\left(\frac{\partial a_i}{\partial t_j}\right)_{i,j}=
\begin{pmatrix}
1&0&\dfrac{1}{4}t_3\\
0&1&0\\
0&0&1
\end{pmatrix}.
\]
In flat coordinates we have
\begin{gather*}
\eta=\left(
\begin{matrix}
 0 & 0 & \dfrac{1}{4} \\
 0 & \dfrac{1}{4} & 0 \\
 \dfrac{1}{4} & 0 & 0 \\
\end{matrix}
\right),\qquad \mathcal U(t_1,t_2,t_3)=\left(
\begin{matrix}
 t_1 & \dfrac{-5}{16} t_2 t_3 & -\dfrac{3}{16} t_2^2+\dfrac{1}{32} t_3^3 \\
 \dfrac{3 t_2}{4} & t_1-\dfrac{t_3^2}{8} & \dfrac{-5}{16} t_2 t_3 \\
 \dfrac{t_3}{2} & \dfrac{3 t_2}{4} & t_1
\end{matrix}
\right),\\
 \mu=\begin{pmatrix}
-\dfrac{1}{4}&&\\
&0&\\
&&\dfrac{1}{4}
\end{pmatrix}.
\end{gather*}
Thus, the second system in \eqref{28luglio2016-4}
\[\partial_z\xi=\left(\mathcal U^{\rm T}-\frac{1}{z}\mu\right)\xi
\]
reads
\begin{gather}\label{systemA3}
\begin{cases}
\partial_z\xi_1=\dfrac{3 }{4}\xi _2 t_2+\dfrac{1}{2}\xi _3 t_3+\xi _1 \left(t_1+\dfrac{1}{4 z}\right),\\
\partial_z\xi_2=-\dfrac{5}{16} \xi _1 t_2 t_3+\xi _2 \left(t_1-\dfrac{t_3^2}{8}\right)+\dfrac{3 }{4}\xi _3 t_2,\\
\partial_z \xi_3=\xi _1 \left(-\dfrac{3}{16} t_2^2+\dfrac{1}{32} t_3^3\right)-\dfrac{5}{16} \xi _2 t_2 t_3+\xi _3 \left(t_1-\dfrac{1}{4 z}\right).
\end{cases}
\end{gather}We know that, if $(t_1,t_2,t_3)$ is a semisimple point of the Frobenius manifold then the monodromy data are well defined, and that these are invariant under (small) deformations of $t_1$, $t_2$, $t_3$ by Theorems~\ref{isoth2} and~\ref{mainisoth}. The bifurcation set is now
\[
 \{t_2=0 \}\cup\big\{8 t_3^3+27 t_2^2=0\big\}.
\]
 Now, if we fix $a_0=0$, the tensor $\mathcal U$ at $(0,a_1,h)$, i.e., $(t_1,t_2,t_3)=\big({-}\frac{1}{8}h^2,t_2,h\big)$, is
\begin{gather*}
\mathcal U\left(-\frac{1}{8}h^2,t_2,h\right)=
\left(
\begin{matrix}
 -\dfrac{h ^2}{8} & -\dfrac{5h}{16} t_2 & \dfrac{1}{32} \big(h ^3-6 t_2^2\big) \\
 \dfrac{3 t_2}{4} & -\dfrac{h ^2}{4} & -\dfrac{5h}{16} t_2 \\
 \dfrac{h }{2} & \dfrac{3 t_2}{4} & -\dfrac{h ^2}{8}
\end{matrix}
\right).
\end{gather*}
The bifurcation locus is reached for $a_1=t_2=0$. At these points \begin{gather*} (t_1,t_2,t_3)=\left(-\frac{1}{8}h^2,0,h\right),\end{gather*}
we have
\[\mathcal U\left(-\frac{1}{8}h^2,0,h\right)=
\left(
\begin{matrix}
 -\dfrac{h ^2}{8} & 0 & \dfrac{h ^3}{32} \\
 0 & -\dfrac{h ^2}{4} & 0 \\
 \dfrac{h }{2} & 0 & -\dfrac{h ^2}{8}
\end{matrix}\right).
\]

Define the function
\[X(a):=\big[{-}9a_1+\sqrt{3}\big(27a_1^2+8a_2^3\big)^{\frac{1}{2}}\big]^{\frac{1}{3}},
\]
which has branch points along the caustic $\mathcal K= \{a_2=0 \}\cup\big\{27a_1^2+8a_2^3=0\big\}$. Fix a {branch} of~$X$ on a simply connected domain in $M \setminus \mathcal K$, that we also denote by~$X(a)$.
The critical points $x_1$, $x_2$, $x_3$ of $f(x,a)$ are equal to
\[x_i(a):=\frac{\vartheta_i \cdot a_2}{2\sqrt{3}\cdot X(a)}-\frac{\vartheta_i \cdot X(a)}{2\cdot 3^{2/3}},
\]where
\[\vartheta_1:=-1,\qquad\vartheta_2:=\frac{1-{\rm i}\sqrt{3}}{2},\qquad\vartheta_3:=\frac{1+{\rm i}\sqrt{3}}{2}
\]are the cubic roots of $(-1)$. Of course, different choices of {branches} of $X$ correspond to permutations of the $x_i$'s.  After some computations, we find the following expression for $\Psi$:
\begin{gather*}\Psi(t)=\left.
\left(
\begin{matrix}
 \frac{\sqrt{6 x_1^2+a_2}}{2 \sqrt{2} (x_1-x_2) (x_1-x_3)} & -\frac{(x_2+x_3) \sqrt{6 x_1^2+a_2}}{2 \sqrt{2} (x_1-x_2) (x_1-x_3)} & -\frac{\sqrt{6 x_1^2+a_2} (a_2-4 x_2 x_3)}{8 \sqrt{2} (x_1-x_2) (x_1-x_3)} \vspace{1mm}\\
 \frac{\sqrt{6 x_2^2+a_2}}{2 \sqrt{2} (x_1-x_2) (x_3-x_2)} & \frac{(x_1+x_3) \sqrt{6 x_2^2+a_2}}{2 \sqrt{2} (x_1-x_2) (x_2-x_3)} & \frac{\sqrt{6 x_2^2+a_2} (a_2-4 x_1 x_3)}{8 \sqrt{2} (x_1-x_2) (x_2-x_3)} \vspace{1mm}\\
 \frac{\sqrt{6 x_3^2+a_2}}{2 \sqrt{2} (x_1-x_3) (x_2-x_3)} & \frac{(x_1+x_2) \sqrt{6 x_3^2+a_2}}{2 \sqrt{2} (x_1-x_3) (x_3-x_2)} & \frac{(a_2-4 x_1 x_2) \sqrt{6 x_3^2+a_2}}{8 \sqrt{2} (x_1-x_3) (x_3-x_2)} \end{matrix}
\right)\right|_{a=a(t)},\end{gather*}
where
\[a_0=t_1+\frac{1}{8}t_3^2,\qquad a_1=t_2,\qquad a_2=t_3.
\]
The canonical coordinates are $u_i(t)=f(x_i(a(t)),a(t))$. In a neighbourhood of the point $(t_1,t_2,t_3)= \big({-}\frac{1}{8}h^2,0,h\big)$, i.e., for small $t_2$ and $h\neq 0$, we find
\begin{gather*}
u_1(t_2;h)=-\frac{t_2^2}{4 h}+\frac{t_2^4}{16 h^4}-\frac{t_2^6}{16 h^7}+\frac{3 t_2^8}{32 h^{10}}+O\big(|t_2|^{10}\big),\\
u_2(t_2;h)=-\frac{h^2}{4}+\frac{{\rm i} \sqrt{h} t_2}{\sqrt{2}}+\frac{t_2^2}{8 h}+\frac{{\rm i} t_2^3}{16 \sqrt{2} h^{5/2}}-\frac{t_2^4}{32 h^4}-\frac{21 {\rm i} t_2^5}{512 \sqrt{2} h^{11/2}}\\
\hphantom{u_2(t_2;h)=}{}+\frac{t_2^6}{32 h^7}+\frac{429 {\rm i} t_2^7}{8192 \sqrt{2} h^{17/2}}-\frac{3 t_2^8}{64 h^{10}}-\frac{46189 {\rm i} t_2^9}{524288 \sqrt{2} h^{23/2}}+O\big(|t_2|^{10}\big),\\
u_3(t_2;h)=-\frac{h^2}{4}-\frac{{\rm i} \sqrt{h} t_2}{\sqrt{2}}+\frac{t_2^2}{8 h}-\frac{{\rm i} t_2^3}{16 \sqrt{2} h^{5/2}}-\frac{t_2^4}{32 h^4}+\frac{21 {\rm i} t_2^5}{512 \sqrt{2} h^{11/2}}\\
\hphantom{u_3(t_2;h)=}{} +\frac{t_2^6}{32 h^7}-\frac{429 {\rm i} t_2^7}{8192 \sqrt{2} h^{17/2}}-\frac{3 t_2^8}{64 h^{10}}+\frac{46189 {\rm i} t_2^9}{524288 \sqrt{2} h^{23/2}}+O\big(|t_2|^{10}\big),\\
\Psi(t_2)=
\left(
\begin{matrix}
 \frac{1}{\sqrt{2} \sqrt{h}} & 0 & \frac{\sqrt{h}}{4 \sqrt{2}} \vspace{1mm}\\
 \frac{{\rm i}}{2 \sqrt{h}} & -\frac{1}{2 \sqrt{2}} & -\frac{1}{8} \big({\rm i}\sqrt{h}\big) \vspace{1mm}\\
 \frac{{\rm i}}{2 \sqrt{h}} & \frac{1}{2 \sqrt{2}} & -\frac{1}{8} \big({\rm i}\sqrt{h}\big)
\end{matrix}
\right)+t_2\left(
\begin{matrix}
 0 & -\frac{1}{2 \sqrt{2} h^{3/2}} & 0 \vspace{1mm}\\
 -\frac{3}{8 \sqrt{2} h^2} & -\frac{{\rm i}}{16 h^{3/2}} & -\frac{5}{32 \sqrt{2} h} \vspace{1mm}\\
 \frac{3}{8 \sqrt{2} h^2} & -\frac{{\rm i}}{16 h^{3/2}} & \frac{5}{32 \sqrt{2} h}
\end{matrix}
\right)\\
\hphantom{\Psi(t_2)=}{} +t_2^2\left(
\begin{matrix}
 -\frac{3}{4 \sqrt{2} h^{7/2}} & 0 & \frac{1}{16 \sqrt{2} h^{5/2}} \vspace{1mm}\\
 -\frac{39{\rm i}}{128 h^{7/2}} & \frac{15}{128 \sqrt{2} h^3} & -\frac{41{\rm i}}{512 h^{5/2}} \vspace{1mm}\\
 -\frac{39{\rm i}}{128 h^{7/2}} & -\frac{15}{128 \sqrt{2} h^3} & -\frac{41{\rm i}}{512 h^{5/2}}
\end{matrix}
\right)+O\big(t_2^3\big).
\end{gather*}
Hence, at points $(t_1,t_2,t_3)= \big({-}\frac{1}{8}h^2,0,h\big)$, canonical coordinates $u_i(0;h)$ are
\[
(u_1,u_2,u_3) =\left(0,-\frac{h^2}{2}, -\frac{h^2}{2}\right)
\]
and the system \eqref{systemA3} reduces to
\begin{gather}\label{sistemaA3}
\begin{cases}
\partial_z\xi_1=\left(-\dfrac{h^2}{8}+\dfrac{1}{4z}\right)\xi_1+\dfrac{h}{2}\xi_3,\vspace{1mm}\\
\partial_z\xi_2=-\dfrac{h^2}{4}\xi_2,\vspace{1mm}\\
\partial_z\xi_3=\dfrac{h^3}{32}\xi_1-\left(\dfrac{h^2}{8}+\dfrac{1}{4z}\right)\xi_3.
\end{cases}
\end{gather}
The second equation yields
\[\xi_2(z)=c\cdot {\rm e}^{-\frac{h^2}{4}z},\qquad c\in\mathbb C.
\]
From the first equation we find that
\begin{gather}\label{reconstructionxi3}\xi_3=\frac{2}{h}\left(\partial_z\xi_1+\frac{h^2}{8}\xi_1-\frac{1}{4z}\xi_1\right),
\end{gather}and so from the third equation we obtain
\[\frac{2}{h}\xi_1''(z)+\frac{h}{2}\xi_1'(z)+\frac{3}{8z^2h}\xi_1=0.
\]
Making the ansatz
\[\xi_1=z^{\frac{1}{2}}{\rm e}^{-\frac{h^2z}{8}}\Lambda(z),
\]
the equation for $\Lambda$ becomes the following Bessel equation:
\begin{gather}\label{besselA3}64z^2\Lambda''(z)+64z\Lambda'(z)-\big(4+z^2h^4\big)\Lambda(z)=0.
\end{gather}
Therefore, $\xi_1$ is of the form
\[\xi_1=z^{\frac{1}{2}}{\rm e}^{-\frac{h^2z}{8}}\left(c_1H^{(1)}_{\frac{1}{4}}\left(\frac{{\rm i}h^2}{8}z\right)+c_2H^{(2)}_{\frac{1}{4}}\left(\frac{{\rm i}h^2}{8}z\right)\right),\qquad c_1,c_2\in\mathbb C,
\]
 where $H^{(1)}_{\nu}(z)$, $H^{(2)}_{\nu}(z)$ stand for the Hankel functions of the first and second kind of parameter $\nu=1/4$. Notice that if $\Lambda(z)$ is a solution of equation~\eqref{besselA3}, then also $\Lambda\big({\rm e}^{\pm {\rm i}\pi}z\big)$ is a solution.

\subsection{Computation of Stokes and central connection matrices}
In order to compute the Stokes matrix, let us fix the line $\ell$ {coinciding} with the real axis. Such a line is admissible for all points $(t_1,t_2,t_3)= \big({-}\frac{1}{8}h^2,0,h\big)$ with
\[|\operatorname{Re} h|\neq |\operatorname{Im} h|,\qquad h\in \mathbb C^*.
\]
Indeed, the Stokes rays for $(u_1,u_2,u_3)=\big(0,-\frac{1}{4}h^2,-\frac{1}{4}h^2\big)$ are
\begin{gather*}
z=\pm {\rm i}\rho \overline{h}^2 \quad \Longrightarrow \quad \arg z= \frac{\pi}{2}-2\arg h \text{ (mod $\pi$)}.
\end{gather*}
Thus, admissibility corresponds to $\frac{1}{2}\pi-2\arg h\neq k\pi$, $k\in\mathbb Z$.
Let us compute the Stokes matrix in the case \[ -\frac{\pi}{4}<\arg h<\frac{\pi}{4}.\] The asymptotic expansion for fundamental solutions $\Xi_{\rm left}$, $\Xi_{\rm right}$ of the system~\eqref{sistemaA3}, is
\begin{align*}\eta\Psi^{-1}\left(\mathbbm 1+O\left(\frac{1}{z}\right)\right){\rm e}^{zU}&=\Psi^{\rm T}\left(\mathbbm 1+O\left(\frac{1}{z}\right)\right){\rm e}^{zU}\\
& =\left(\mathbbm 1+O\left(\frac{1}{z}\right)\right)
\left(
\begin{matrix}
 \frac{1}{\sqrt{2} \sqrt{h}} & \frac{{\rm i}}{2 \sqrt{h}}{\rm e}^{-\frac{1}{4} h^2 z} & \frac{{\rm i} }{2 \sqrt{h}} {\rm e}^{-\frac{1}{4} h^2 z}
 \\
 0 & -\frac{1}{2 \sqrt{2}} {\rm e}^{-\frac{1}{4} h^2 z}& \frac{1}{2 \sqrt{2}}{\rm e}^{-\frac{1}{4} h^2 z}
 \\
 \frac{\sqrt{h}}{4 \sqrt{2}} & -\frac{{\rm i}}{8} {\rm e}^{-\frac{1}{4} h^2 z} \sqrt{h} & -\frac{{\rm i}}{8} {\rm e}^{-\frac{1}{4} h^2 z} \sqrt{h}
\end{matrix}
\right),
\end{align*}
being
\[ U:=\Psi\mathcal U\Psi^{-1}=\operatorname{diag}(u_1,u_2,u_3)=\operatorname{diag}\left(0,-\frac{h^2}{4},-\frac{h^2}{4}\right).
\]
For the admissible line $\ell$ and for the above labelling of canonical coordinates the Stokes matrix must be of the form prescribed by Theorem \ref{proprstok}:
\begin{gather}
\label{stokesa3}S=\begin{pmatrix}
1&0&0\\
\alpha&1&0\\
\beta&0&1
\end{pmatrix}
\end{gather}for some constants $\alpha,\beta\in\mathbb C$ to be determined. This means that the last two columns of $\Xi_{\rm left}$ must be the analytic continuation of $\Xi_{\rm right}$.

\begin{Lemma}\label{lemmawatson}The following asymptotic expansions hold:
\begin{itemize}\itemsep=0pt
\item if $m\in\mathbb Z$, then\[H_{\frac{1}{4}}^{(1)}\left({\rm e}^{{\rm i} m\pi}\frac{{\rm i} h^2}{8}z\right)\sim \sqrt{\frac{2}{\pi}}\left({\rm e}^{{\rm i} m\pi}\frac{{\rm i} h^2}{8}z\right)^{-\frac{1}{2}}{\rm e}^{-\frac{3{\rm i}\pi}{8}}\exp\left(-{\rm e}^{{\rm i} m\pi}\frac{h^2}{8}z\right)
\]in the sector
\[-\frac{3}{2}\pi-m\pi-\arg\big(h^2\big)<\arg z<\frac{3}{2}\pi-m\pi-\arg\big(h^2\big);
\]

\item if $m\in\mathbb Z$, then
\[H_{\frac{1}{4}}^{(2)}\left({\rm e}^{{\rm i} m\pi}\frac{{\rm i} h^2}{8}z\right)\sim\sqrt{\frac{2}{\pi}}\left({\rm e}^{{\rm i} m\pi}\frac{{\rm i} h^2}{8}z\right)^{-\frac{1}{2}}{\rm e}^{\frac{3{\rm i}\pi}{8}}\exp\left({\rm e}^{{\rm i} m\pi}\frac{h^2}{8}z\right)
\]in the sector
\[-\frac{5}{2}\pi-m \pi-\arg\big(h^2\big)<\arg z<\frac{\pi}{2}-m\pi-\arg\big(h^2\big).
\]
\end{itemize}
\end{Lemma}

\begin{proof}These formulae easily follow from the following well-known asymptotic expansion of Hankel functions (see~\cite{watson}):
\[H^{(1)}_{\nu}(z)\sim\sqrt{\frac{2}{\pi z}}\exp\left({\rm i}\left(z-\frac{\nu}{2}\pi-\frac{\pi}{4}\right)\right)
,\qquad
-\pi+\delta\leq\arg z\leq2\pi-\delta,
\]$\delta$ being any positive acute angle. Analogously,
\begin{gather*} H^{(2)}_{\nu}(z)\sim\sqrt{\frac{2}{\pi z}}\exp\left(-{\rm i}\left(z-\frac{\nu}{2}\pi-\frac{\pi}{4}\right)\right),
\qquad
-2\pi+\delta\leq\arg z\leq\pi-\delta.\tag*{\qed}
\end{gather*}\renewcommand{\qed}{}
\end{proof}

Using Lemma \ref{lemmawatson}, we obtain
\begin{gather}\label{leftandright}
\Xi_{\rm left}(z)=\begin{pmatrix}
\xi^{L}_{(1),1}&\xi^{L}_{(2),1}&\xi^{L}_{(3),1}\vspace{1mm}\\
0&-\frac{{\rm e}^{-\frac{1}{4} h^2 z}}{2 \sqrt{2}} & \frac{{\rm e}^{-\frac{1}{4} h^2 z}}{2 \sqrt{2}}\vspace{1mm}\\
*&*&*
\end{pmatrix},\qquad\!
\Xi_{\rm right}(z)=\begin{pmatrix}
\xi^{R}_{(1),1}&\xi^{R}_{(2),1}&\xi^{R}_{(3),1}\vspace{1mm}\\
0&-\frac{{\rm e}^{-\frac{1}{4} h^2 z}}{2 \sqrt{2}} & \frac{{\rm e}^{-\frac{1}{4} h^2 z}}{2 \sqrt{2}}\vspace{1mm}\\
*&*&*
\end{pmatrix},
\end{gather}
where
\[\xi^{L}_{(2),1}(z)=\xi^{L}_{(3),1}(z)=\xi^{R}_{(2),1}(z)=\xi^{R}_{(3),1}(z) =\frac{{\rm i}\sqrt{\pi}}{8}h^{\frac{1}{2}}{\rm e}^{{\rm i}\frac{5}{8}\pi}z^{\frac{1}{2}}{\rm e}^{-\frac{zh^2}{8}}H^{(1)}_{\frac{1}{4}}\left(\frac{{\rm i}h^2}{8}z\right),
\]
with the required asymptotic expansion in the following sector containing both $\Pi_{\rm left}$ and $\Pi_{\rm right}$
\[\left\{z\in\mathcal R\colon -\frac{3}{2}\pi-\arg\big(h^2\big)<\arg z<\frac{3}{2}\pi-\arg\big(h^2\big)\right\},
\]
and
\begin{gather*}
\xi^{L}_{(1),1}(z) =\frac{\sqrt{\pi}}{4\sqrt{2}}h^{\frac{1}{2}}{\rm e}^{{\rm i}\frac{\pi}{8}}z^{\frac{1}{2}}{\rm e}^{-\frac{zh^2}{8}}H^{(1)}_{\frac{1}{4}}\left({\rm e}^{-{\rm i}\pi}\frac{{\rm i}h^2}{8}z\right),\\
\xi^{R}_{(1),1}(z) =\frac{\sqrt{\pi}}{4\sqrt{2}}h^{\frac{1}{2}}{\rm e}^{-{\rm i}\frac{\pi}{8}}z^{\frac{1}{2}}{\rm e}^{-\frac{zh^2}{8}}H^{(2)}_{\frac{1}{4}}\left(\frac{{\rm i}h^2}{8}z\right),
\end{gather*}
with the required expansion respectively in the sectors
\begin{gather*}
\left\{z\in\mathcal R\colon -\frac{\pi}{2}-\arg\big(h^2\big)<\arg z<\frac{5}{2}\pi-\arg \big(h^2\big)\right\}\supseteq \Pi_{\rm left},
\\
\left\{z\in\mathcal R\colon -\frac{5}{2}\pi-\arg\big(h^2\big)<\arg z<\frac{\pi}{2}-\arg\big(h^2\big)\right\}\supseteq\Pi_{\rm right}.
\end{gather*}
The entries of $\Xi_{\rm left}$, $\Xi_{\rm right}$ denoted by~$*$ are reconstructed from the first rows, by applying equation~\eqref{reconstructionxi3}.

From the second rows of $\Xi_{\rm left}$, $\Xi_{\rm right}$ we can immediately say that the entries $\alpha$, $\beta$ of~\eqref{stokesa3} must be equal. Specializing the following well-known connection formula for Hankel special functions
\begin{gather}\label{connhankel1}\sin(\nu\pi)H^{(1)}_\nu\big(z{\rm e}^{m\pi{\rm i}}\big)=-\sin((m-1)\nu\pi)H^{(1)}_\nu(z)-{\rm e}^{-\nu\pi{\rm i}}\sin(m\nu\pi)H^{(2)}_\nu(z),\qquad m\in\mathbb Z,\!\!\!
\end{gather} to the case $m=-1$, $\nu=\frac{1}{4}$, we easily obtain
\[
\xi^{L}_{(1),1}(z)=\xi^{R}_{(1),1}(z)-\xi^{R}_{(2),1}(z)-\xi^{R}_{(3),1}(z),
\]
which means that $\alpha=\beta=-1$. So, we have obtained that, at points $(t_1,t_2,t_3)= \big({-}\frac{1}{8}h^2,0,h\big)$ with \begin{gather*}|\operatorname{Re} h|>|\operatorname{Im} h|,\qquad-\frac{\pi}{4}<\arg h<\frac{\pi}{4}\end{gather*} (and consequently in their neighborhood, by Theorem \ref{mainisoth}) the Stokes matrix is
\begin{gather*}
S=\begin{pmatrix}
1&0&0\\
-1&1&0\\
-1&0&1
\end{pmatrix}.
\end{gather*}

In order to compute the central connection matrix, we observe that the $A_3$ Frobenius manifold structure is \emph{non-resonant}, i.e., the components of the tensor $\mu$ are such that $\mu_\alpha-\mu_\beta\notin\mathbb Z$ for $\alpha \neq \beta$. This implies that the $(\eta,\mu)$-parabolic orthogonal group is trivial, and that the fundamental system of \eqref{sistemaA3} near the origin $z=0$ can be uniquely chosen in such a way that
\begin{gather}\label{originA3}\Xi_0(z)=(\eta+O(z))z^\mu.
\end{gather}
Now, let us recall the following Mellin--Barnes integral representations of Hankel functions (see~\cite{watson})
\begin{gather*}
H^{(1)}_\nu(z)=-\frac{\cos(\nu\pi)}{\pi^{\frac{5}{2}}}{\rm e}^{{\rm i}(z-\pi\nu)}(2z)^\nu\int_{-\infty {\rm i}}^{\infty {\rm i}}\Gamma(s)\Gamma(s-2\nu)\Gamma\left(\nu+\frac{1}{2}-s\right)(-2{\rm i}z)^{-s}\,{\rm d}s,\\
H^{(2)}_\nu(z)=\frac{\cos(\nu\pi)}{\pi^{\frac{5}{2}}}{\rm e}^{-{\rm i}(z-\pi\nu)}(2z)^\nu\int_{-\infty {\rm i}}^{\infty {\rm i}}\Gamma(s)\Gamma(s-2\nu)\Gamma\left(\nu+\frac{1}{2}-s\right)(2{\rm i}z)^{-s}\,{\rm d}s,
\end{gather*}
which are valid for
\begin{itemize}\itemsep=0pt
\item $2\nu\notin2\mathbb Z+1$,
\item respectively in the sectors $|\arg(\mp {\rm i}z)|<\frac{3}{2}$,
\item and where the integration path separates the poles of $\Gamma(s)\Gamma(s-2\nu)$ from the poles of $\Gamma\big(\nu+\frac{1}{2}-s\big)$.
\end{itemize}
Specializing these integral forms to $\nu=\frac{1}{4}$, and deforming the integration path so that it reduces to positively oriented circles around the poles
\[s\in\frac{1}{2}-\frac{1}{2}\mathbb N,
\]we immediately obtain the following expansion of the solution $\xi^{(1)}_{1,R}$, $\xi^{(2)}_{1,R}$, $\xi^{(3)}_{1,R}$ for the points $(t_1,t_2,t_3)= \big({-}\frac{1}{8}h^2,0,h\big)$,  with $-\frac{\pi}{4}<\arg h<\frac{\pi}{4}$, valid for small values of~$|z|$:
\begin{Lemma}At the points $(t_1,t_2,t_3)= \big({-}\frac{1}{8}h^2,0,h\big)$,  with $-\frac{\pi}{4}<\arg h<\frac{\pi}{4}$ the following expansions hold:
\begin{gather*}
\xi^{R}_{(1),1}(z)= \frac{(1+{\rm i}) \Gamma \left(\frac{5}{4}\right)}{\sqrt{\pi }}z^{\frac{1}{4}}+\frac{\left(\frac{1}{4}-\frac{{\rm i}}{4}\right) h \Gamma \left(\frac{3}{4}\right)}{\sqrt{\pi }}z^{3/4}\\
\hphantom{\xi^{R}_{(1),1}(z)=}{} -\frac{\left(\frac{1}{32}+\frac{{\rm i}}{32}\right) h^2 \Gamma \left(\frac{1}{4}\right)}{\sqrt{\pi }}z^{5/4} -\frac{\left(\frac{1}{32}-\frac{{\rm i}}{32}\right) h^3 \Gamma \left(\frac{3}{4}\right)}{\sqrt{\pi }}z^{7/4} +O\big(|z|^{9/4}\big),\\
\xi^{R}_{(2),1}(z)=\xi^{R}_{(3),1}(z)=\frac{{\rm i} \Gamma \left(\frac{5}{4}\right)}{\sqrt{\pi }}z^\frac{1}{4}-\frac{4 {\rm i} h \Gamma \left(\frac{11}{4}\right)}{21 \sqrt{\pi }}z^{3/4}-\frac{{\rm i} h^2 \Gamma \left(\frac{5}{4}\right)}{8 \sqrt{\pi }}z^{5/4}+\frac{{\rm i} h^3 \Gamma \left(\frac{11}{4}\right)}{42 \sqrt{\pi }}z^{7/4}\\
\hphantom{\xi^{R}_{(2),1}(z)=}{} +O\big(|z|^{9/4}\big).
\end{gather*}
Moreover, using equation \eqref{reconstructionxi3} we find that
\begin{gather*}
\xi^{R}_{(1),3}(z)=\frac{\left(\frac{1}{4}-\frac{{\rm i}}{4}\right) \Gamma \left(\frac{3}{4}\right)}{\sqrt{\pi } }z^{-\frac{1}{4}}-\frac{\left(\frac{1}{32}-\frac{{\rm i}}{32}\right) h^2 \Gamma \left(\frac{3}{4}\right)}{\sqrt{\pi }}z^{3/4}+O\big(|z|^{5/4}\big),
\\\
\xi^{R}_{(2),3}(z)=\xi^{R}_{(3),3}(z)=-\frac{4 {\rm i} \Gamma \left(\frac{11}{4}\right)}{21 \sqrt{\pi }}z^{-\frac{1}{4}}+\frac{{\rm i} h^2 \Gamma \left(\frac{11}{4}\right)}{42 \sqrt{\pi }}z^{3/4}+O\big(|z|^{5/4}\big).
\end{gather*}
\end{Lemma}

\begin{proof} These expansions are the first terms of the expressions
\begin{gather*}
\xi^{R}_{(1),1}(z)= \frac{\sqrt{\pi}}{4\sqrt{2}}h^{\frac{1}{2}}{\rm e}^{-{\rm i}\frac{\pi}{8}}z^{\frac{1}{2}}{\rm e}^{-\frac{zh^2}{8}} \frac{{\rm e}^{-{\rm i} \left(-\frac{\pi }{4}+\frac{{\rm i} h^2 z}{8}\right)} \big({\rm i} h^2 z\big)^{\frac{1}{4}}}{2 \pi ^{5/2}}
\\
\hphantom{\xi^{R}_{(1),1}(z)=}{} \times 2\pi{\rm i}\sum_{n=0}^\infty\underset{s=\frac{1}{2}-\frac{n}{2}}{\operatorname{res}} \left(\Gamma(s)\Gamma(s-2\nu)\Gamma\left(\nu+\frac{1}{2}-s\right)\left({\rm e}^{{\rm i}\pi}\frac{h^2z}{4}\right)^{-s}\right),
\end{gather*}
and
\begin{gather*}\xi^{R}_{(2),1}(z)= \xi^{R}_{(3),1}(z)=\frac{{\rm i}\sqrt{\pi}}{8}h^{\frac{1}{2}}{\rm e}^{{\rm i}\frac{5}{8}\pi}z^{\frac{1}{2}}{\rm e}^{-\frac{zh^2}{8}} \left(-\frac{{\rm e}^{{\rm i} \left(-\frac{\pi }{4}+\frac{1}{8} {\rm i} h^2 z\right)} \big({\rm i} h^2 z\big)^\frac{1}{4}}{2 \pi ^{5/2}}\right)
\\
\hphantom{\xi^{R}_{(2),1}(z)=}{} \times 2\pi{\rm i}\sum_{n=0}^\infty\underset{s=\frac{1}{2}- \frac{n}{2}}{\operatorname{res}}\left(\Gamma(s)\Gamma(s-2\nu)\Gamma\left(\nu+\frac{1}{2}-s\right)\left(-{\rm e}^{{\rm i}\pi}\frac{h^2z}{4}\right)^{-s}\right).\tag*{\qed}
\end{gather*}\renewcommand{\qed}{}
\end{proof}

By a direct comparison between these expansion of solution $\Xi_{\rm right}(z)$ of~\eqref{leftandright} and the dominant term of~\eqref{originA3}, namely
\[\left(
\begin{matrix}
 0 & 0 &\dfrac{z^\frac{1}{4}}{4}\\
 0 & \dfrac{1}{4} & 0 \\
 \dfrac{z^{-\frac{1}{4}}}{4}& 0 & 0
\end{matrix}
\right),
\]
we obtain the central connection matrix
\[C=\frac{1}{\pi^\frac{1}{2}}\left(
\begin{matrix}
 (1-{\rm i}) \Gamma \left(\dfrac{3}{4}\right) & -{\rm i} \Gamma \left(\dfrac{3}{4}\right) & -{\rm i} \Gamma \left(\dfrac{3}{4}\right)\vspace{1mm}\\
 0 & -\sqrt{2 \pi } & \sqrt{2 \pi } \vspace{1mm}\\
 (1+{\rm i}) \Gamma \left(\dfrac{1}{4}\right) & {\rm i} \Gamma \left(\dfrac{1}{4}\right) & {\rm i} \Gamma \left(\dfrac{1}{4}\right)
\end{matrix}
\right).
\]
Notice that such a matrix satisfies all the constraints of Theorem \ref{constraint}.

We can put the Stokes matrix in triangular form using two different permutations of the canonical coordinates $\big(0,-h^2/4,-h^2/4\big)$, namely
\begin{itemize}\itemsep=0pt
\item re-labelling $(u_1,u_2,u_3)\mapsto(u_2,u_3,u_1)$, corresponding to the permutation matrix
\[P=\begin{pmatrix}
0&1&0\\
0&0&1\\
1&0&0
\end{pmatrix},
\]
\item or re-labelling $(u_1,u_2,u_3)\mapsto(u_3,u_2,u_1)$, corresponding to the permutation matrix
\[P=\begin{pmatrix}
0&0&1\\
0&1&0\\
1&0&0
\end{pmatrix}.
\]
\end{itemize}
In \looseness=-1 both cases these are the lexicographical orders of two different $\ell$-cells which divide any sufficiently small neighborhood of the point $(t_1,t_2,t_3)= \big({-}\frac{1}{8}h^2,0,h\big)$, with $|\operatorname{Re} h|>|\operatorname{Im} h|$ and $-\frac{\pi}{4}<\arg h<\frac{\pi}{4}$, in which Theorem~\ref{mainisoth} applies.
Using both permutations, the Stokes matrix becomes
\begin{gather}\label{27giugno2017-1}
S_{\text{lex}}=PSP^{-1}=
\begin{pmatrix}
1&0&-1\\
0&1&-1\\
0&0&1
\end{pmatrix},
\end{gather}
which can be thought as in {the} lexicographical form in one of the $\ell$-cells. The central connection matrix, instead, has the following lexicographical forms in the two $\ell$-cells:
\begin{gather}
\label{29-Feb-16-1}
C_{\text{lex}}=\frac{1}{\pi^\frac{1}{2}}\left(
\begin{matrix}
 -{\rm i} \Gamma \left(\dfrac{3}{4}\right) & -{\rm i} \Gamma \left(\dfrac{3}{4}\right) &(1-{\rm i}) \Gamma \left(\dfrac{3}{4}\right) \vspace{1mm}\\
 \mp\sqrt{2 \pi } & \pm\sqrt{2 \pi }&0 \vspace{1mm}\\
 {\rm i} \Gamma \left(\dfrac{1}{4}\right) & {\rm i} \Gamma \left(\dfrac{1}{4}\right)& (1+{\rm i}) \Gamma \left(\dfrac{1}{4}\right) \\
\end{matrix}
\right),
\end{gather}where we take the first sign if the lexicographical order is the relabeling $(u_1,u_2,u_3) \mapsto(u_2,u_3,u_1)$, the second if it is the re-labeling $(u_1,u_2,u_3)\mapsto(u_3,u_2,u_1)$.

\subsection[A ``tour'' in the Maxwell stratum: reconstruction of neighboring monodromy data]{A ``tour'' in the Maxwell stratum:\\ reconstruction of neighboring monodromy data}

\begin{figure}\centering
 \includegraphics[width=0.45 \textwidth]{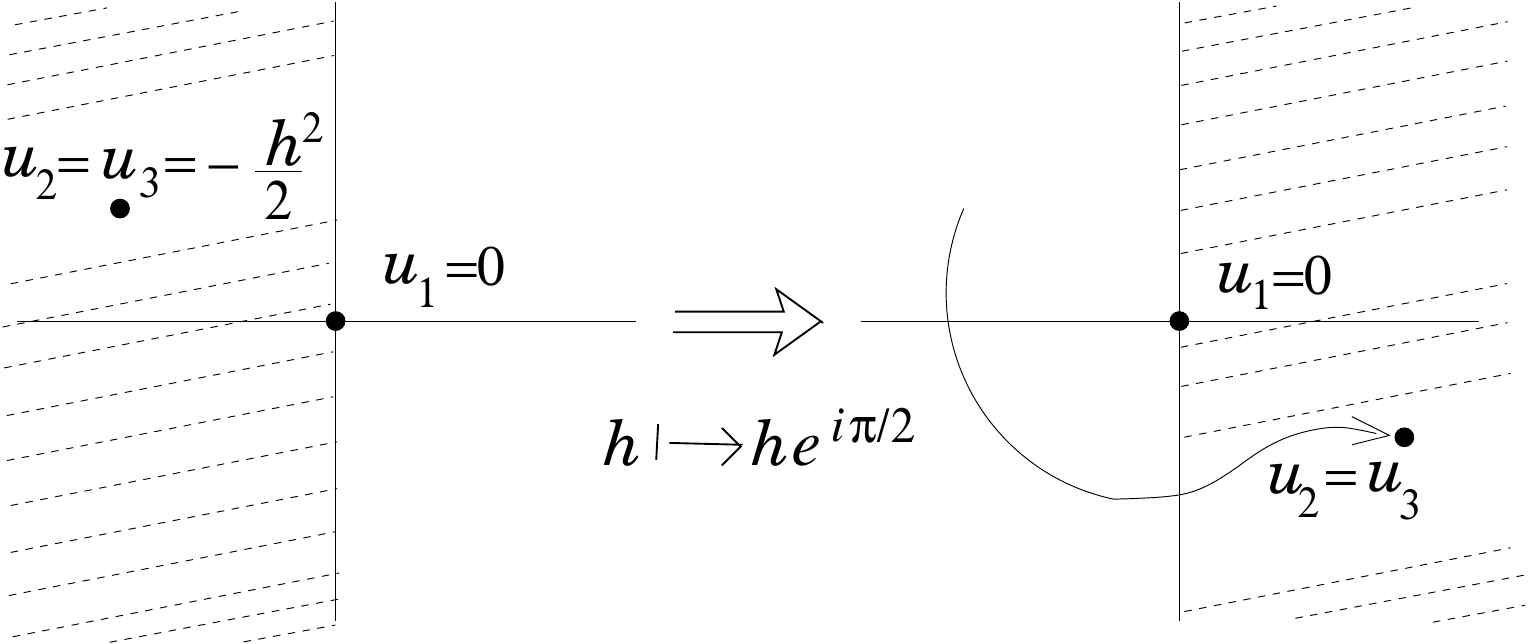}
\caption{The triple $(u_1,u_2,u_2)$ is represented by three points $u_1$, $u_2$, $u_3$ in $\mathbb{C}$. We move along $h\mapsto h{\rm e}^{{\rm i}\frac{\pi}{2}}$, starting from $-\frac{\pi}{4}<\arg h<\frac{\pi}{4}$. The two dashed regions in the left and right figures correspond respectively to $-\frac{\pi}{4}<\arg h<\frac{\pi}{4}$ and $\frac{\pi}{4}<\arg h<\frac{3\pi}{4}$.}\label{27giugno2017-5}
\end{figure}

From the data (\ref{27giugno2017-1}) and (\ref{29-Feb-16-1}), by an action of the braid group, we can compute $S$ and $C$ in the neighborhood of all other points $(t_1,t_2,t_3)= \big({-}\frac{1}{8}h^2,0,h\big)$ with $|\operatorname{Re} h|\neq|\operatorname{Im} h|$.
As an example, let us determine the Stokes matrix for points
 \[(t_1,t_2,t_3)= \left(-\frac{1}{8}h^2,0,h\right),\qquad \text{with } \frac{\pi}{4}<\arg h<\frac{3}{4}\pi.\]
Starting from a point in the region $-\frac{\pi}{4}<\arg h<\frac{\pi}{4}$ and moving counter-clockwise towards the region $\frac{\pi}{4}<\arg h<\frac{3}{4}\pi$, the two coalescing canonical coordinates $u_2=u_3=-\frac{1}{2}h^2$ move in the $u_i$'s-plane counter-clockwise w.r.t.\ $u_1=0$. For example, in Fig.~\ref{27giugno2017-5} we move along a curve $h\mapsto h{\rm e}^{{\rm i}\frac{\pi}{2}}$, starting in $-\frac{\pi}{4}<\arg h<\frac{\pi}{4}$. At $\arg h=\frac{\pi}{4}$, the Stokes rays $R_{12}=\big\{z=-{\rm i}\rho\overline{h}^2,~\rho>0\big\}$ and $R_{21}=\big\{z={\rm i}\rho\overline{h}^2,~\rho>0\big\}$ cross the real line~$\ell$, and a braid must act on the monodromy data. In order to determine the braid and the transformed monodromy data, we proceed according to the prescription of Section~\ref{30luglio2016-3}, as follows.

\begin{figure}[th]\centering
\def\svgscale{.45}
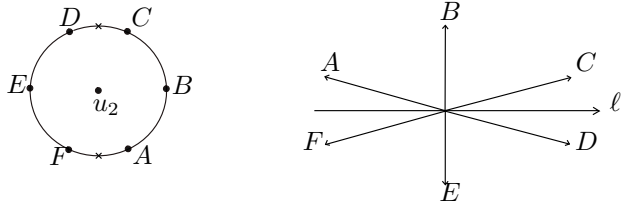
\caption{In the left picture we represent relative positions of $u_3$ w.r.t.~$u_2$ such that the real line~$\ell$ is admissible. On the right, we represent the corresponding positions of the Stokes ray~$R_{23}$. Notice that if we let vary~$u_3$, by a deformation of the parameter~$\varphi$, starting from $A$, going through~$B$ up to~$C$, the corresponding Stokes ray does not cross the line~$\ell$, and no braids act. If we continue the deformation of~$\varphi$ from $C$ to $D$, an elementary braid acts on the monodromy data.}\label{raggiA3}
\label{5agosto2016-4}
\end{figure}

(1) We split the coalescing canonical coordinates, for example by considering the point
\begin{gather}\label{5agosto2016-3}
(t_1,t_2,t_3)=\left(-\frac{1}{8}h^2,\varepsilon {\rm e}^{{\rm i}\varphi},h\right), \qquad \text{with } -\frac{\pi}{4}<\arg h<\frac{\pi}{4}
\end{gather}
for chosen $\varphi$ and $\epsilon$, being $\varepsilon$ small (so that $\varepsilon^2 \ll \varepsilon$). The corresponding canonical coordinates
\begin{gather}\label{27giugno2017-8}
u_1=O\big(\varepsilon^2\big),\\
\label{27giugno2017-9}
u_2=-\frac{h^2}{4}+\varepsilon |h|^{\frac{1}{2}}\exp\left[{\rm i}\left(\frac{\arg h}{2}+\varphi+\frac{\pi}{2}\right)\right]+O\big(\varepsilon^2\big),\\
\label{27giugno2017-10}
u_3 =-\frac{h^2}{4}+\varepsilon |h|^{\frac{1}{2}}\exp\left[{\rm i}\left(\frac{\arg h}{2}+\varphi-\frac{\pi}{2}\right)\right]+O\big(\varepsilon^2\big),
\end{gather}
give a point $(u_1,u_2,u_3)$ which lies in one of the two cells (Definition~\ref{ellcell}) which divide a polydisc centred at $(u_1,u_2,u_3)=\big(0,-\frac{1}{2}h^2, -\frac{1}{2}h^2\big)$. The Stokes rays are
\begin{gather*}
R_{12}=\big\{z=-{\rm i}\rho\overline{h}^2+O(\varepsilon),~\rho>0\big\}, \qquad
R_{13}=\big\{z=-{\rm i}\rho\overline{h}^2+O(\varepsilon),~\rho>0\big\},
\\ 
R_{23}=\left\{z=\rho\exp\left[-{\rm i}\left(\frac{\arg h}{2}+\varphi+\pi\right)\right]+O\big(\varepsilon^2\big),~\rho>0\right\},
\end{gather*}
and opposite ones $R_{21}$, $R_{31}$, $R_{32}$. Notice that in order for the real line $\ell$ to remain admissible, we choose $\varphi\neq k\pi -\frac{1}{2}\arg h$, $k\in\mathbb{Z}$, $-\frac{\pi}{4}<\arg h<\frac{\pi}{4}$. The position of $R_{23}$ w.r.t.\ the real line $\ell$ is determined by the sign of $\cos\big(\frac{\arg h}{2}+\varphi+\frac{\pi}{2}\big)$. As long as $\varphi
$ varies in such a way that $\operatorname{sgn}\cos\big(\frac{\arg h}{2}+\varphi+\frac{\pi}{2}\big)$
does not change, then~$R_{23}$ does not cross~$\ell$. See Fig.~\ref{5agosto2016-4}. This means that $(u_1,u_2,u_3)$ remains inside the same cell, i.e., the point corresponding to coordinates~\eqref{5agosto2016-3} remains inside an $\ell$-chamber, where the isomonodromy Theorem~\ref{isoth2} applies.

(2) The Stokes matrix must be put in triangular form $S_\text{lex}$ (\ref{27giugno2017-1}).
 In particular,
\begin{itemize}\itemsep=0pt
\item if $\cos\big(\frac{\arg h}{2}+\varphi+\frac{\pi}{2}\big)<0$, then $R_{23}$ is on the \emph{left} of~$\ell$, and the lexicographical order is given by the permutation $(u_1,u_2,u_3)\mapsto (u_1^\prime ,u_2^\prime,u_3^\prime)=(u_2,u_3,u_1)$;
\item if $\cos\big(\frac{\arg h}{2}+\varphi+\frac{\pi}{2}\big)>0$, then $R_{23}$ is on the \emph{right} of $\ell$, and the lexicographical order is given by the permutation $(u_1,u_2,u_3)\mapsto (u_1^\prime ,u_2^\prime,u_3^\prime)=(u_3,u_2,u_1)$.
\end{itemize}
 We choose the cell where the triangular order coincides with the lexicographical order. The passage to the other $\ell$-cell is obtained by a counter-clockwise rotation of $u_1^\prime$ w.r.t.\ $u_2^\prime$, which corresponds to the action of the elementary braid~$\beta_{12}$. Its action~(\ref{stokesbraid1}) is a permutation matrix, since $(S_\text{lex})_{12}= 0$; it is a trivial action on~$S_\text{lex}$, but not on~$C_\text{lex}$, as~(\ref{29-Feb-16-1}) shows.

\begin{figure}[t]
\centerline{\includegraphics[width=.37 \textwidth]{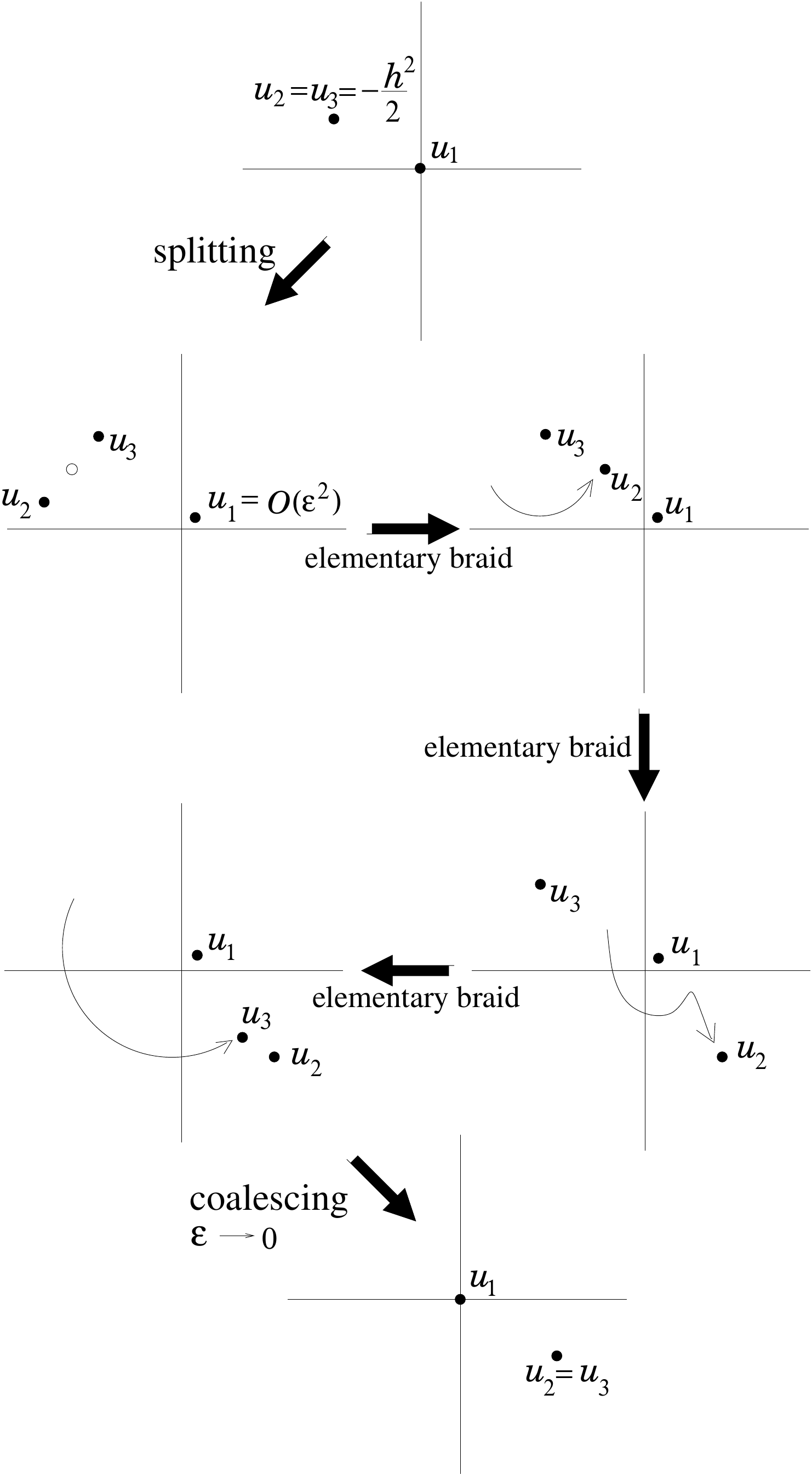}}
\caption{The transition in Fig.~\ref{27giugno2017-5} by splitting and elementary steps. After the splitting, we obtain a point $(u_1,u_2,u_3)$, as in (\ref{27giugno2017-8})--(\ref{27giugno2017-10}), lying in an $\ell$-cell of the polydisc centred at $(u_1,u_2,u_3)=\big(0,-\frac{1}{2}h^2,-\frac{1}{2}h^2\big)$ of the left part of Fig.~\ref{27giugno2017-5}. The transformation of Fig.~\ref{27giugno2017-5} is obtained by successive steps following the arrows. The final step is the right part of Fig.~\ref{27giugno2017-5}. The first elementary braid is $\beta_{12}$ (because $u_1^\prime=u_2$, $u_2^\prime=u_3$ in the upper left figure). The second is $\beta_{23}$ (after relabelling in lexicographical order, $u_2^\prime=u_2$ and $u_3^\prime=u_1$ in the upper right figure). The third is $\beta_{12}$.}\label{27giugno2017-6}
\end{figure}

\begin{figure}[th]\centering
\def\svgscale{.8}
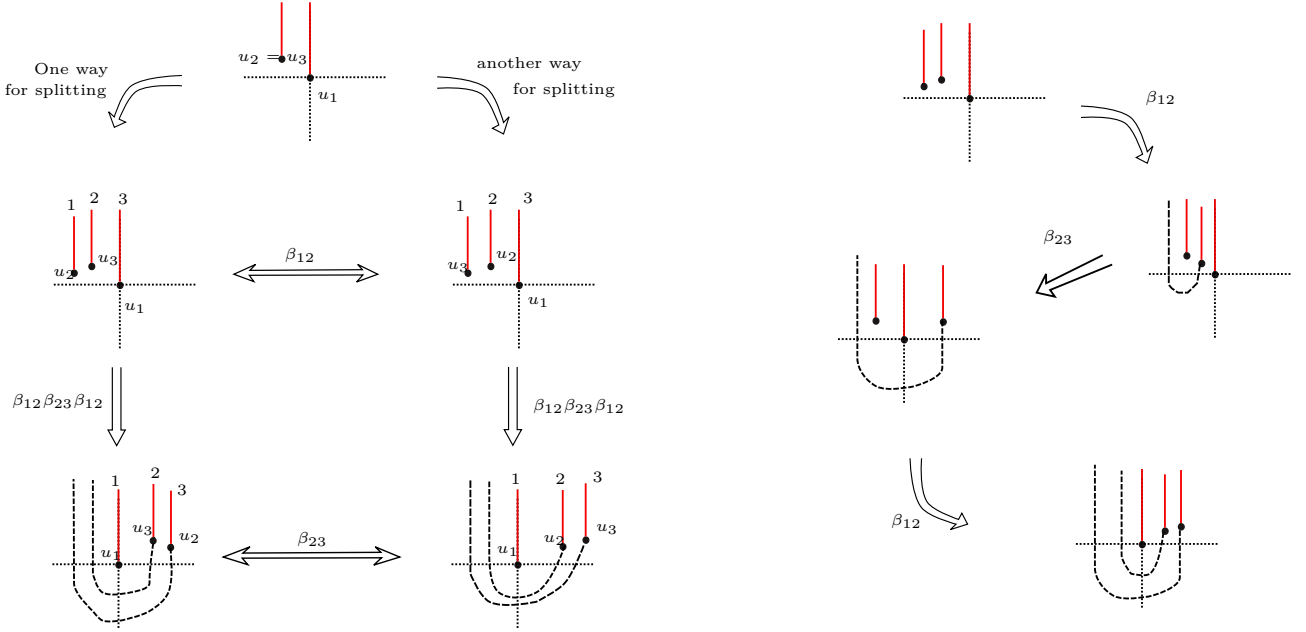
\caption{In the picture we represent $u_1$, $u_2$, $u_3$ as points in $\mathbb C$. On the left we describe all the braids necessary to pass from a neighborhood of $(t_1,t_2,t_3)= \big({-}\frac{1}{8}h^2,0,h\big)$ with $-\frac{\pi}{4}<\arg h<\frac{\pi}{4}$ to one with $\frac{\pi}{4}<\arg h<\frac{3}{4}\pi$. Different columns of this diagram correspond to different $\ell$-cells of the same neighborhood. The passage from such one cell to the other is through an action of an elementary braid ($\beta_{12}$ or $\beta_{23}$) acting as a permutation matrix. In the picture on the right, we show the decomposition of the global transformation in elementary ones.}\label{braidA3}
\end{figure}

(3) We move along a curve $h\mapsto h {\rm e}^{{\rm i}\frac{\pi}{2}}$ in the $h$-plane from a point (\ref{5agosto2016-3}) up to a point
\[
(t_1,t_2,t_3)= \left(-\frac{1}{8}h^2,\varepsilon {\rm e}^{{\rm i}\varphi^\prime},h\right),\qquad \text{with} \quad \frac{\pi}{4}<\arg h<\frac{3}{4}\pi,
\]
for some $\varphi^\prime \neq k\pi -\frac{1}{2}\arg h$, $k\in\mathbb{Z}$, $\frac{\pi}{4}<\arg h<\frac{3\pi}{4}$. The transformation in Fig.~\ref{27giugno2017-5}, due to the splitting, can substituted by the sequence of transformations in Fig.~\ref{27giugno2017-6}, each step corresponding to an elementary braid.  Each elementary braid corresponds to a Stokes ray crossing clock-wise the real line $\ell$ as $h$ varies along the curve $h\mapsto h {\rm e}^{{\rm i}\frac{\pi}{2}}$.\footnote{Notice that the ray $R_{23}$ rotates slower than $R_{12}$, $R_{13}$: namely, the angular velocity of $R_{23}$ is approximately (i.e., modulo negligible corrections in powers of~$\varepsilon$) equal to $\frac{1}{4}$ the one of~$R_{12}$,~$R_{13}$.}
The total braid is then factored into the product of the elementary braids as in Fig.~\ref{braidA3}, namely
\[\beta_{12}\beta_{23}\beta_{12},\qquad\text{or}\qquad\beta_{12}\beta_{23}\beta_{12}\beta_{23}.
\]
 Applying formulae \eqref{stokesbraid1}, \eqref{connbraid2}, we obtain
\begin{gather}\label{stokesA3-2}
S_\text{lex}^{\beta_{12}\beta_{23}\beta_{12}}=S_\text{lex}^{\beta_{12}\beta_{23}\beta_{12}\beta_{23}}=
\begin{pmatrix}
1&1&1\\
0&1&0\\
0&0&1
\end{pmatrix}.
\end{gather}
These are the monodromy data in the two $\ell$-cells of a polydisc centred at the point
\[(t_1,t_2,t_3)= \left(-\frac{1}{8}h^2,0,h\right),\qquad \text{with} \quad \frac{\pi}{4}<\arg h<\frac{3}{4}\pi.
\]
The braid $\beta_{23}$ is responsible for the passage from one cell to the other. Its action $A^{\beta_{23}}\big(S_\text{lex}^{\beta_{12}\beta_{23}\beta_{12}}\big)$ is a permutation matrix, since $\big(S_\text{lex}^{\beta_{12}\beta_{23}\beta_{12}}\big)_{23}=0$, which explains the equality in~(\ref{stokesA3-2}). By the action~(\ref{connbraid2}), the central connection matrix~\eqref{29-Feb-16-1}, instead, assumes the following two forms (differing for a permutation of the second and third column)
\begin{gather*}
C_{\text{lex}}^{\beta_{12}\beta_{23}\beta_{12}}=
\frac{1}{\pi^\frac{1}{2}}\left(
\begin{matrix}
 (1+ {\rm i}) \Gamma \left(\dfrac{3}{4}\right) & - {\rm i} \Gamma \left(\dfrac{3}{4}\right) & -{\rm i} \Gamma \left(\dfrac{3}{4}\right) \vspace{1mm}\\
 0 & \pm\sqrt{2 \pi } & \mp\sqrt{2 \pi } \vspace{1mm}\\
 (1-{\rm i}) \Gamma \left(\dfrac{1}{4}\right) & {\rm i} \Gamma \left(\dfrac{1}{4}\right) & {\rm i} \Gamma \left(\dfrac{1}{4}\right)
\end{matrix}
\right),
\\
C_{\text{lex}}^{\beta_{12}\beta_{23}\beta_{12}\beta_{23}}=
\frac{1}{\pi^\frac{1}{2}}\left(
\begin{matrix}
 (1+{\rm i}) \Gamma \left(\dfrac{3}{4}\right) & -{\rm i} \Gamma \left(\dfrac{3}{4}\right) & -{\rm i} \Gamma \left(\dfrac{3}{4}\right) \vspace{1mm}\\
 0 & \mp\sqrt{2 \pi } & \pm\sqrt{2 \pi } \vspace{1mm}\\
 (1-{\rm i}) \Gamma \left(\dfrac{1}{4}\right) & {\rm i} \Gamma \left(\dfrac{1}{4}\right) & {\rm i} \Gamma \left(\dfrac{1}{4}\right)
\end{matrix}
\right).
\end{gather*}

In Table~\ref{tabellaA3} we show the monodromy data for other values of $\arg h$, with the corresponding braid. In Fig.~\ref{braidA3} we represent the braid corresponding to the passage from $-\frac{\pi}{4}<\arg h <\frac{\pi}{4}$ to $\frac{5}{4}\pi <\arg h< \frac{7}{4}\pi$.

\begin{figure}[t]\centering
\def\svgscale{0.8}
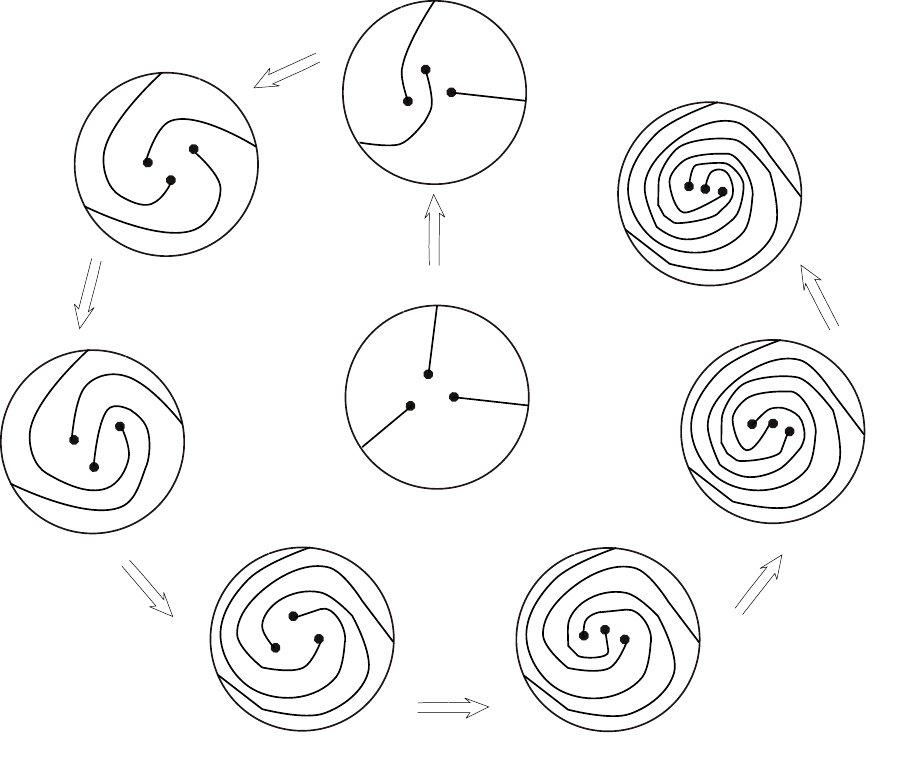
\caption{Using the diagram representation of the braid group as mapping class group of the punctured disk, we draw the braids acting along a curve $h\mapsto {\rm e}^{\frac{3\pi{\rm i}}{2}}h$, starting from the chambers close to $(t_1,t_2,t_3)= \big({-}\frac{1}{8}h^2,0,h\big)$ with $-\frac{\pi}{4}<\arg h<\frac{\pi}{4}$, and reaching the ones with $\frac{5}{4}\pi<\arg h<\frac{7}{4}\pi$. The braids in red describe mutations of the split pair $u_2$, $u_3$: their action on the monodromy data is a permutation matrix. In the central disk, the blue numbers refer to the lexicographical order w.r.t.\ the real axis $\ell$ (i.e., from the left to the right). The braids are the same for both cases $(a,b)=(2,3)$ and $(3,2)$.}\label{27giugno2017-2}
\end{figure}

\begin{Remark}The reader can re-obtain this result by direct computation
 observing that, for points
 \begin{gather*}(t_1,t_2,t_3)= \left(-\frac{1}{8}h^2,0,h\right), \qquad \text{with}\quad \frac{\pi}{4}<\arg h<\frac{3}{4}\pi,\end{gather*}
 the left and right solutions of~\eqref{sistemaA3} defining the Stokes matrix\footnote{Notice that for the points with $\frac{\pi}{4}<\arg h<\frac{3}{4}\pi$ the original labelling of canonical coordinates $(u_1,u_2,u_3)=\big(0,-\frac{h^2}{4},-\frac{h^2}{4}\big)$ already put the Stokes matrix in upper triangular form.} are of the form \eqref{leftandright} with:
\begin{gather*}
\xi^{L}_{(1),1}
=\xi^{R}_{(1),1}
=\frac{\sqrt{\pi}}{4\sqrt{2}}h^{\frac{1}{2}}{\rm e}^{{\rm i}\frac{\pi}{8}}z^{\frac{1}{2}}{\rm e}^{-\frac{zh^2}{8}}H^{(1)}_{\frac{1}{4}}\left({\rm e}^{-{\rm i}\pi}\frac{{\rm i}h^2}{8}z\right),
\\
\xi^{L}_{(2),1}(z)
=\xi^{L}_{(3),1}(z)=\frac{{\rm i}\sqrt{\pi}}{8}h^{\frac{1}{2}}{\rm e}^{{\rm i}\frac{3}{8}\pi}z^{\frac{1}{2}}{\rm e}^{-\frac{zh^2}{8}}H^{(2)}_{\frac{1}{4}}\left({\rm e}^{-3{\rm i}\pi}\frac{{\rm i}h^2}{8}z\right),
\\
\xi^{R}_{(2),1}(z)
=\xi^{R}_{(3),1}(z)=\frac{{\rm i}\sqrt{\pi}}{8}h^{\frac{1}{2}}{\rm e}^{{\rm i}\frac{5}{8}\pi}z^{\frac{1}{2}}{\rm e}^{-\frac{zh^2}{8}}H^{(1)}_{\frac{1}{4}}\left(\frac{{\rm i}h^2}{8}z\right),
\end{gather*} having the expected asymptotic expansions in suitable sectors containing $\Pi_{{\rm left}}$ and/or $\Pi_{\rm right}$ by Lemma~\ref{lemmawatson}. Thus, by some manipulation of formulae~\eqref{connhankel1} and
\[\sin(\nu\pi)H^{(2)}_{\nu}\big(z{\rm e}^{m\pi{\rm i}}\big)={\rm e}^{\nu\pi{\rm i}}\sin(m\nu\pi)H^{(1)}_\nu(z)+\sin((m+1)\nu\pi)H^{(2)}_\nu(z),\qquad m\in\mathbb Z,
\]one sees that
\[\xi^{L}_{(2),1}(z)=\xi^{R}_{(1),1}(z)+\xi^{R}_{(2),1}(z), \qquad\xi^{L}_{(3),1}(z)=\xi^{R}_{(1),1}(z)+\xi^{R}_{(3),1}(z),
\]which are equivalent to~\eqref{stokesA3-2}. For the computation of the central connection matrix, one can use analogous Puiseux series expansions of the solution~$\Xi_{\rm right}(z)$, obtained from the integral representation of Hankel functions given above.
\end{Remark}

\begin{table}[t]\centering
\begin{tabular}{|c||c|c|c|}
\hline
 $\arg h$ & $S_\text{lex}$ & $C_\text{lex}$&Braid\\
\hline
\hline
& & &\\[-2ex]
${\big]}{-}\frac{\pi}{4},\frac{\pi}{4}{\big[} $&$\begin{pmatrix}
1&0&-1\\
0&1&-1\\
0&0&1
\end{pmatrix}$ &$\frac{1}{\pi^\frac{1}{2}}\left(
\begin{matrix}
 -{\rm i} \Gamma \left(\frac{3}{4}\right) & -{\rm i} \Gamma \left(\frac{3}{4}\right) &(1-{\rm i}) \Gamma \left(\frac{3}{4}\right) \vspace{1mm}\\
 \mp\sqrt{2 \pi } & \pm\sqrt{2 \pi }&0 \vspace{1mm}\\
 {\rm i} \Gamma \left(\frac{1}{4}\right) & {\rm i} \Gamma \left(\frac{1}{4}\right)& (1+{\rm i}) \Gamma \left(\frac{1}{4}\right)
\end{matrix}
\right)$&$\color{red}{\beta_{12}}$\\
& &&\\[-2ex]
\hline
&&&\\[-2ex]
$\big]\frac{\pi}{4},\frac{3\pi}{4}\big[$&$\begin{pmatrix}
1&1&1\\
0&1&0\\
0&0&1
\end{pmatrix}$ &$\frac{1}{\pi^\frac{1}{2}}\left(
\begin{matrix}
 (1+{\rm i}) \Gamma \left(\frac{3}{4}\right) & -{\rm i} \Gamma \left(\frac{3}{4}\right) & -{\rm i} \Gamma \left(\frac{3}{4}\right) \vspace{1mm}\\
 0 & \pm\sqrt{2 \pi } & \mp\sqrt{2 \pi } \vspace{1mm}\\
 (1-{\rm i}) \Gamma \left(\frac{1}{4}\right) & {\rm i} \Gamma \left(\frac{1}{4}\right) & {\rm i} \Gamma \left(\frac{1}{4}\right)
\end{matrix}
\right)$&$\beta_{12}\beta_{23}\beta_{12}\color{red}{\beta_{23}}$\\
& &&\\[-2ex]
\hline
&&&\\[-2ex]
$\big]\frac{3\pi}{4},\frac{5\pi}{4}\big[$&$\left(
\begin{matrix}
 1 & 0 & -1 \\
 0 & 1 & -1 \\
 0 & 0 & 1
\end{matrix}
\right)$&$\frac{1}{\pi^\frac{1}{2}}\left(
\begin{matrix}
 \Gamma \left(\frac{3}{4}\right) & \Gamma \left(\frac{3}{4}\right) & (1+{\rm i}) \Gamma \left(\frac{3}{4}\right) \vspace{1mm}\\
 \mp \sqrt{2 \pi } & \pm \sqrt{2 \pi } & 0 \vspace{1mm}\\
 \Gamma \left(\frac{1}{4}\right) & \Gamma \left(\frac{1}{4}\right) & (1-{\rm i}) \Gamma \left(\frac{1}{4}\right) \end{matrix}
\right)$&$(\beta_{12}\beta_{23})^3\color{red}{\beta_{12}}$\\
&&&\\[-2ex]
\hline
&&&\\[-2ex]
$\big]\frac{5\pi}{4},\frac{7\pi}{4}\big[$&$\left(
\begin{matrix}
 1 & 1 & 1 \\
 0 & 1 & 0 \\
 0 & 0 & 1
\end{matrix}
\right)$&$\frac{1}{\pi^\frac{1}{2}}\left(
\begin{matrix}
 (-1+{\rm i}) \Gamma \left(\frac{3}{4}\right) & \Gamma \left(\frac{3}{4}\right) & \Gamma \left(\frac{3}{4}\right) \vspace{1mm}\\
 0 & \pm \sqrt{2 \pi } & \mp \sqrt{2 \pi } \vspace{1mm}\\
 (-1-{\rm i}) \Gamma \left(\frac{1}{4}\right) & \Gamma \left(\frac{1}{4}\right) & \Gamma \left(\frac{1}{4}\right)
\end{matrix}
\right)$& $\begin{matrix}(\beta_{12}\beta_{23})^3
\\
\times \beta_{12} \beta_{23}\beta_{12}\color{red}{\beta_{23}}\end{matrix}$\\
&&&\\[-2ex]
\hline
&&&\\[-2ex]
$\big]\frac{7\pi}{4},\frac{9\pi}{4}\big[$&
$\left(
\begin{matrix}
 1 & 0 & -1 \\
 0 & 1 & -1 \\
 0 & 0 & 1
\end{matrix}
\right)$&$\frac{1}{\pi^\frac{1}{2}}
\left(
\begin{matrix}
{\rm i} \Gamma \left(\frac{3}{4}\right) & {\rm i} \Gamma \left(\frac{3}{4}\right) & (-1+{\rm i}) \Gamma \left(\frac{3}{4}\right) \vspace{1mm}\\
  \mp\sqrt{2 \pi } &\pm \sqrt{2 \pi } & 0\\
 -{\rm i} \Gamma \left(\frac{1}{4}\right) & -{\rm i} \Gamma \left(\frac{1}{4}\right) & (-1-{\rm i}) \Gamma \left(\frac{1}{4}\right)
\end{matrix}
\right)$&
$(\beta_{12}\beta_{23})^6\color{red}{\beta_{12}}$\bsep{20pt}\\
\hline
\end{tabular}

\caption{For different values of $\arg h$, in the open angular intervals in the left column, we tabulate the monodromy data $(S_{\text{lex}},C_{\text{lex}})$, in lexicographical order, in the two $\ell$-cells which divide a sufficiently small neighborhood of the point $(t_1,t_2,t_3)= \big({-}\frac{1}{8}h^2,0,h\big)$. The difference of the data in the two $\ell$-cells (just a permutation of two columns in the central connection matrix) is obtained by applying the braid written in red: if it is not applied the sign to be read is the first one, the second one otherwise. Notice that the central element $(\beta_{12}\beta_{23})^3$ acts trivially on the Stokes matrices, and by a left multiplication by $M_0^{-1}=\operatorname{diag}({\rm i},1,-{\rm i})$ on the central connection matrix.}\label{tabellaA3}
\end{table}

\subsection{Monodromy data as computed outside the Maxwell stratum}
In this section, we compute the Stokes matrix $S$ at non-coalescence points in a neighbourhood of a coalescence one, by means of oscillatory integrals. We show that $S$ coincides with that obtained at the coalescence point in the previous section. Moreover, we explicitly show that the fundamental matrices converge to those computed at the coalescence point, exactly as prescribed by our Theorem~\ref{mainisoth}.

The system \eqref{systemA3} admits solutions given in terms of oscillating integrals,
\begin{gather}\label{oscillint1}
\xi_1(z,t) =z^\frac{1}{2}\int_\gamma\exp \{z\cdot f(x,t) \}\,{\rm d}x,\\
\label{oscillint2}
\xi_2(z,t) =z^\frac{1}{2}\int_\gamma x\exp \{z\cdot f(x,t) \}\,{\rm d}x,\\
\label{oscillint3}
\xi_3(z,t) =z^\frac{1}{2}\int_\gamma \left(x^2+\frac{1}{4}t_3\right)\exp \{z\cdot f(x,t) \}\,{\rm d}x,
\end{gather}
where $f(x,t)=x^4+t_3x^2+t_2 x+t_1+\frac{1}{8}t_3^2$. Here $\gamma$ is any cycle along which $\operatorname{Re}(z\cdot f(x,t))\to-\infty$ for $|x|\to+\infty$, i.e., a relative cycle in $H_1(\mathbb C,\mathbb C_{T,z,t})$, with
\[\mathbb C_{T,z,t}:=\left\{x\in\mathbb C\colon \operatorname{Re}(zf(x,t))< -T\right\},\qquad \text{with $T$ very large positive number.}
\]

\begin{figure}\centering
\begin{minipage}[b]{0.44\textwidth}
\centering
\def\svgscale{.5}
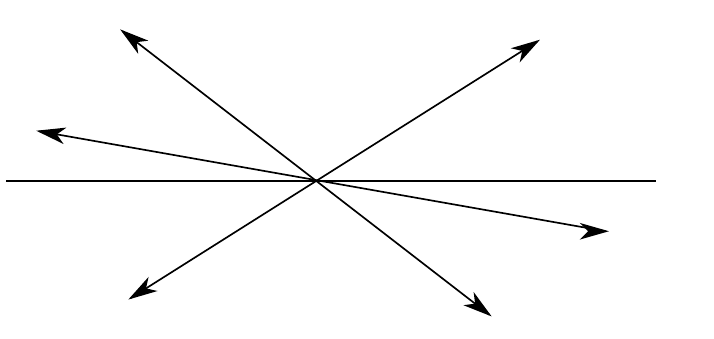
\caption{Disposition of the Stokes rays for a~point in the chosen $\ell$-chamber.}
\label{raggiA3bis}
\end{minipage}\qquad
\begin{minipage}[b]{0.45\textwidth}
\centering
\def\svgscale{.3}
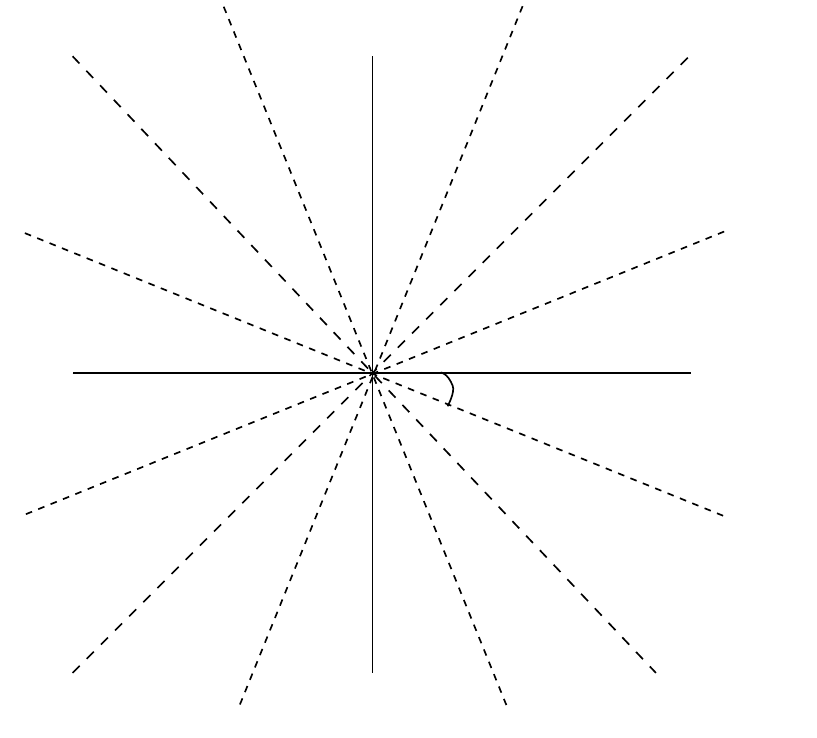
\caption{Integration contours $\mathcal I_i$ which define the functions $\mathfrak I_i$'s.}
\label{integrationcont}
\end{minipage}
\end{figure}

First, we show that the Stokes matrix at points in $\ell$-chambers near the coalescence point $(t_1,t_2,t_3)=\big({-}\frac{1}{8}h^2,0,h\big)$ coincide with the one previously computed, in accordance with Theorem \ref{mainisoth}. In what follows we will focus on the $\ell$-chamber made of points $(t_1,t_2,t_3)\!=\!\big(\!{-}\frac{1}{8}h^2\!,\varepsilon {\rm e}^{{\rm i}\phi}\!,h\big)$, where $-\frac{\pi}{4}<\arg h<\frac{\pi}{4}$, and $\varepsilon$, $\phi$ are small positive numbers. For points in this $\ell$-chamber, the Stokes rays are disposed as described in Fig.~\ref{raggiA3bis}.

Notice that in order to compute the Stokes matrix at a semisimple point with distinct canonical coordinates it {suffices} to know the first rows of $\Xi_{{\rm left/right}}$.
 Assuming that $z\in\mathbb R_+$, we define the following three functions obtained from the integrals \eqref{oscillint1} with integration cycles $\mathcal I_i$ as in Fig.~\ref{integrationcont}:
\begin{gather}\label{oscintLef}
\mathfrak I_i(z,t_2):=\int_{\mathcal I_i}\exp\big(z\big(x^4+hx^2+t_2 x\big)\big)\,{\rm d}x,\qquad i=1,2,3.
\end{gather}
For the specified integration cycles, the integrals $\mathfrak I_i(z,t_2)$ are convergent in the half-plane $|\arg z|<\frac{\pi}{2}$. A continuous deformation of a path~${\mathcal I_i}$, which maintains its asymptotic directions in the shaded sectors, yields a convergent integral and defines the analytic continuation of $\mathfrak I_i(z,t_2)$ on the whole sector $|\arg z| <\pi$. If we vary~$z$ (excluding $z=0$), the shaded regions continuously rotate clockwise or counterclockwise. In order to obtain the analytic continuation of the functions $\mathfrak I_i(z,t_2)$ to the whole universal cover $\mathcal R$, we can simply rotate the integration contours $\mathcal I_i$. This procedure also makes it clear that the functions~$\mathfrak I_i$ have monodromy of order~$4$: indeed as~$\arg z$ increases or decreases by $2\pi$, the shaded regions are cyclically permuted.

\begin{figure}[h]\centering
\def\svgscale{.7}
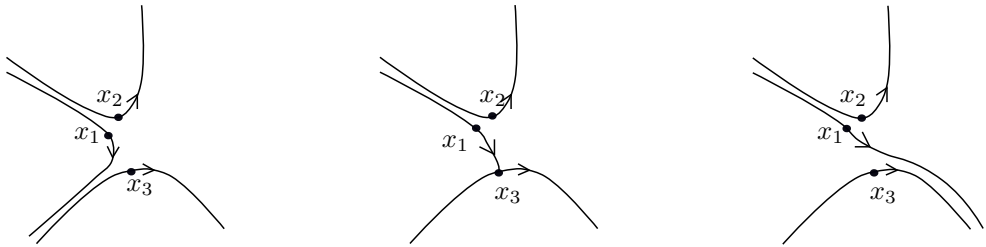
\caption{Discontinuous change of a Lefschetz thimbles. As $z$ varies in $\mathcal R$, we pass from the configuration on the left to the one on the right. The middle configuration is realized when $z$ is on a~Stokes ray: in this case there is a downward geodesic-flow line connecting two critical points~$x_1$ and~$x_3$.}\label{lefschetz}
\end{figure}

In order to obtain information about the asymptotic expansions of the functions $\mathfrak I_i$, we associate to any critical point $x_i$ a relative cycle $\mathcal L_i$, called \emph{Lefschetz thimble}, defined as the set of points of $\mathbb C$ which can be reached along the downward geodesic-flow
\begin{gather}\label{flow}\frac{{\rm d}x}{{\rm d}\tau}=-\bar{z}\frac{\partial\bar{f}}{\partial\bar{x}},\qquad \frac{{\rm d}\bar{x}}{{\rm d}\tau}=-z\frac{\partial f}{\partial{x}}
\end{gather}starting at the critical point $x_i$ for $\tau\to-\infty$. Morse and Picard--Lefschetz theory guarantees that the cycles $\mathcal L_i$ are smooth {one-dimensional} submanifolds of $\mathbb C$, \emph{piecewise smoothly} dependent on the parameters $z,t$, and they represent a basis for the relative homology groups $H_1(\mathbb C,\mathbb C_{T,z,t})$. Moreover, the Lefschetz thimbles are steepest descent paths: namely, $\operatorname{Im}(zf(x,t))$ is constant on each connected component of $\mathcal L_i\setminus\left\{x_i\right\}$ and $\operatorname{Re}(zf(x,t))$ is strictly decreasing along the flow~\eqref{flow}. Thus, after choosing an orientation, the paths of integration defining the functions~$\mathfrak I_i$ can be expressed as integer combinations of the thimbles $\mathcal L_i$ for any value of~$z$: \begin{gather}\label{lincomblefsch}\mathcal I_i=n_1\mathcal L_1+n_2\mathcal L_2+n_3\mathcal L_3,\qquad n_i\in\mathbb Z.\end{gather}
If we let $z$ vary, the Lefschetz thimbles change. When $z$ crosses a Stokes ray, Lefschetz thimbles jump discontinuously, as shown in Fig.~\ref{lefschetz}. In particular, for $z$ on a Stokes ray there exists a~flow line of \eqref{flow} connecting two critical points~$x_i$'s.

\begin{figure}[t]\centering
\def\svgscale{0.8}
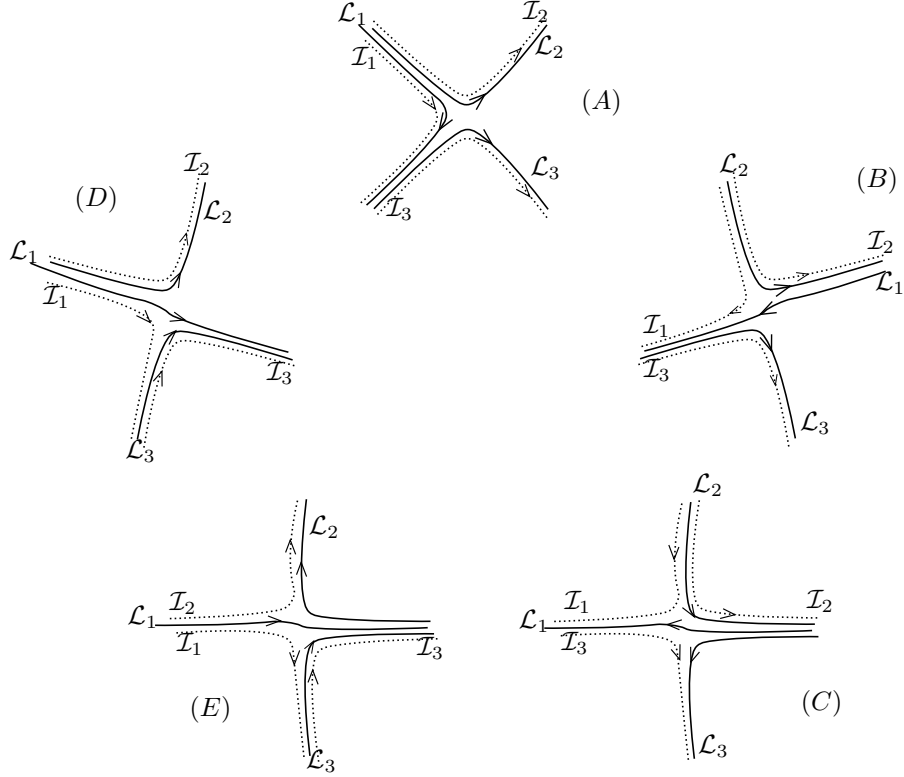
\caption{In this figure it is shown how the Lefschetz thimbles $\mathcal L_i$'s (continuous lines), and the integrations contours $\mathcal I_i$'s (dotted lines) change by analytic continuation with respect to the variable~$z$. The configuration $(A)$ corresponds to the case $\arg z=0$. Increasing $\arg z$ the configuration $(B)$ and $(C)$ are reached after crossing the Stokes rays $R_{31}$, and $R_{21}$ respectively. Decreasing $\arg z$, we obtain the configurations $(D)$ and $(E)$ after crossing the rays $R_{12}$ and $R_{13}$ respectively. Note that when $z$ crosses the Stokes rays $R_{32}$ and $R_{23}$ no Lefschetz thimble changes, coherently with the detailed analysis done in~\cite{CDG0}.}\label{24-11-16}
\end{figure}

This discontinuous change of the thimbles implies a discontinuous change of the integer coefficients $n_i$ in~\eqref{lincomblefsch}, and a discontinuous change of the leading term of the asymptotic expansions of the functions $\mathfrak I_i$'s. Using the notations introduced in Fig.~\ref{24-11-16}, in each configuration the following identities hold:
\begin{gather*} (A)\colon \ \begin{cases}
\mathcal I_1=\mathcal L_1,\\
\mathcal I_2=\mathcal L_2,\\
\mathcal I_3=\mathcal L_3,
\end{cases}\qquad
(B)\colon \ \begin{cases}
\mathcal I_1=\mathcal L_1+\mathcal L_2,\\
\mathcal I_2=\mathcal L_2,\\
\mathcal I_3=\mathcal L_3,
\end{cases}\qquad
(C)\colon \ \begin{cases}
\mathcal I_1=\mathcal L_1+\mathcal L_2,\\
\mathcal I_2=\mathcal L_2,\\
\mathcal I_3=-\mathcal L_1+\mathcal L_3,
\end{cases}
\\
(D)\colon \ \begin{cases}
\mathcal I_1=\mathcal L_1-\mathcal L_3,\\
\mathcal I_2=\mathcal L_2,\\
\mathcal I_3=\mathcal L_3,
\end{cases}\qquad
(E)\colon \ \begin{cases}
\mathcal I_1=\mathcal L_1-\mathcal L_3,\\
\mathcal I_2=\mathcal L_1+\mathcal L_2,\\
\mathcal I_3=\mathcal L_3.
\end{cases}
\end{gather*}

By a {straightforward} application of the Laplace method we find that, at least for sufficiently small positive values of $\arg z$, the following asymptotic expansions hold
\[\mathfrak I_i(z,t_2)= \pi^\frac{1}{2}{\rm i}z^{-\frac{1}{2}}\big(6x_i^2+h\big)^{-\frac{1}{2}}{\rm e}^{zu_i}\left(1+O\left(\frac{1}{z}\right)\right).
\]
Since the deformations of the thimbles $\mathcal I_2$, $\mathcal I_3$ happen for values of $z$ for which the exponent ${\rm e}^{zu_1}$ is subdominant, we immediately conclude that the functions
\begin{gather}\label{soluz2leftright}\xi^L_{(2),1}(z,t_2)=\xi^R_{(2),1}(z,t_2)=\pm {\rm i}\pi^{-\frac{1}{2}}z^{\frac{1}{2}}\frac{6x_2^2+h}{2\sqrt{2}(x_1-x_2)(x_3-x_2)}\mathfrak I_2(z,t_2),
\\
\label{soluz3leftright}\xi^L_{(3),1}(z,t_2)=\xi^R_{(3),1}(z,t_2)=\pm {\rm i}\pi^{-\frac{1}{2}}z^{\frac{1}{2}}\frac{6x_3^2+h}{2\sqrt{2}(x_1-x_3)(x_2-x_3)}\mathfrak I_3(z,t_2)
\end{gather}
have asymptotic expansions
\[\Psi_{21}{\rm e}^{zu_2}\left(1+O\left(\frac{1}{z}\right)\right),\qquad \Psi_{31}{\rm e}^{zu_3}\left(1+O\left(\frac{1}{z}\right)\right),
\]respectively, both in $\Pi_{\rm left}$ and $\Pi_{\rm right}$. Thus, we can immediately say that the Stokes matrix computed at a point $(t_1,t_2,t_3)=\big({-}\frac{1}{8}h^2,\varepsilon {\rm e}^{{\rm i}\phi},h\big)$ is of the form
\[S=\begin{pmatrix}
1&0&0\\
*&1&0\\
*&0&1
\end{pmatrix}.
\]
Note that the arbitrariness of the orientations of the Lefschetz thimbles can be incorporated in the {choice of the entries} of the~$\Psi$ matrix, and hence it will affect the monodromy data by the action of the group $(\mathbb Z/2\mathbb Z)^3$.

After a careful analysis of the deformations of the Lefschetz thimbles, one finds that the solutions $\xi_{(1),1}^L(z,t_2)$, $\xi_{(1),1}^R(z,t_2)$ are respectively given by
\begin{gather}\label{soluz1right}\xi_{(1),1}^R(z,t_2) =\pm {\rm i}\Psi_{11}\pi^{-\frac{1}{2}}z^\frac{1}{2}\big(6x_1^2+h\big)^\frac{1}{2}\left(\mathfrak I_1(z,t_2)+\mathfrak I_3(z,t_2)\right),\\
\label{soluz1left}\xi_{(1),1}^L(z,t_2) =\pm {\rm i}\Psi_{11}\pi^{-\frac{1}{2}}z^\frac{1}{2}\big(6x_1^2+h\big)^\frac{1}{2}\left(\mathfrak I_1(z,t_2)-\mathfrak I_2(z,t_2)\right),
\end{gather}
having the asymptotic expansion
\[\Psi_{11}{\rm e}^{zu_1}\left(1+O\left(\frac{1}{z}\right)\right)
\]
in $\Pi_{\rm right}$ and $\Pi_{\rm left}$ respectively. This immediately allows one to compute the remaining entries of the Stokes matrix
\begin{gather*}S_{21} =\frac{\Psi_{11}\big(6x_1^2+h\big)^\frac{1}{2}}{\Psi_{21}\big(6x_2^2+h\big)^\frac{1}{2}}
=\pm\frac{\big(6x_1^2+h\big)(x_3-x_2)}{(x_1-x_3)\big(6x_2^2+h\big)}\equiv \pm 1,\\
S_{31} =\frac{\Psi_{11}\big(6x_1^2+h\big)^\frac{1}{2}}{\Psi_{31}\big(6x_3^2+h\big)^\frac{1}{2}}
=\pm\frac{\big(6x_1^2+h\big)(x_2-x_3)}{(x_1-x_2)\big(6x_3^2+h\big)}\equiv \pm 1.
\end{gather*}
This result is independent on the point $(t_1,t_2,t_3)=\big({-}\frac{1}{8}h^2,\varepsilon {\rm e}^{{\rm i}\phi},h\big)$ of the chosen $\ell$-chamber. It coincides with the Stokes matrix obtained at the coalescence point $(t_1,t_2,t_3)=\big({-}\frac{1}{8}h^2,0,h\big)$, in complete accordance with our Theorem~\ref{mainisoth}.

\begin{Remark}
It is interesting to note that the isomonodromy condition in this context is equivalent to the condition
\[\frac{f''(x_1)}{f''(x_2)}=-\frac{x_1-x_3}{x_2-x_3},
\]
a relation that the reader can easily show to be valid for any polynomial $f(x)$ of fourth degree with three non-degenerate critical points $x_1$, $x_2$, $x_3$.
\end{Remark}

\begin{figure}[th]\centering
\begin{minipage}[b]{0.45\textwidth}
\centering
\def\svgscale{.6}
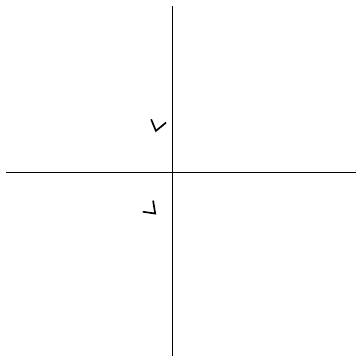
\caption{Integration contour $L$ used~in the integral representation of the Weber parabolic cylinder functions. \newline {} \newline {} }\label{pathL}
\end{minipage}\qquad
\begin{minipage}[b]{0.45\textwidth}
\centering
\def\svgscale{.6}
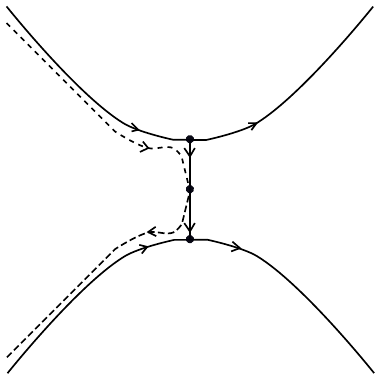
\caption{For $t_2=0$, we can decompose the integration cycle $\mathcal I_1$ into two pieces, $\mathcal I_1^1$, $\mathcal I_1^2$ used to define the functions $\mathfrak I_1^1$ and $\mathfrak I_1^2$. The continuous lines represent the Lefschetz thimbles through the critical points $x_i$'s.}\label{camminiintcoal}
\end{minipage}
\end{figure}

Our Theorem \ref{mainisoth} also states that, as $t_2\to 0$ the solutions \eqref{soluz2leftright}, \eqref{soluz3leftright}, \eqref{soluz1right}, \eqref{soluz1left} must converge to the ones computed in the previous section at the coalescence point. We show this explicitly below. In order to do this, it suffices to set $t_2=0$ in the integral \eqref{oscintLef}. With the change of variable $x=2^{-\frac{1}{4}}z^{-\frac{1}{4}}s^\frac{1}{2}$ , we obtain
\begin{align*}\mathfrak I_2(z,0)=\mathfrak I_3(z,0)&=2^{-\frac{5}{4}}z^{-\frac{1}{4}}\int_L\exp\left\{\frac{s^2}{2}+\left(\frac{hz^{\frac{1}{2}}}{\sqrt 2}\right)s\right\}\,{\rm d}s\\
&=2^{-\frac{5}{4}}z^{-\frac{1}{4}}(2\pi)^\frac{1}{2}{\rm e}^{-\frac{h^2z}{8}}D_{-\frac{1}{2}}\left(\frac{hz^{\frac{1}{2}}}{\sqrt 2}\right)\\
&=2^{-\frac{3}{2}}{\rm e}^{-\frac{h^2z}{8}}h^{\frac{1}{2}}K_{\frac{1}{4}}\left(\frac{h^2z}{8}\right)
 =\pi{\rm i}\cdot 2^{-\frac{5}{2}}h^\frac{1}{2}{\rm e}^{-\frac{h^2z}{8}}{\rm e}^{\frac{\pi{\rm i}}{8}}H^{(1)}_\frac{1}{4}\left(\frac{{\rm i}h^2z}{8}\right).
\end{align*}
Here $D_{\nu}(z)$ is the Weber parabolic cylinder function of order $\nu$, with integral representation \cite[p.~688]{abrasteg}
\[D_{-\frac{1}{2}}(z)=\pm\frac{{\rm e}^{\frac{1}{2}z^2}}{(2\pi)^{\frac{1}{2}}}\int_L s^{-\frac{1}{2}}\exp\left(\frac{s^2}{2}+zs\right)\,{\rm d}s,
\]
where
\[ \begin{cases}
(+) & \text{if }-\frac{3\pi}{2}+2k\pi<\arg s<-\frac{\pi}{2}+2k\pi,\\
(-) & \text{if }\frac{\pi}{2}+2k\pi<\arg s<\frac{3\pi}{2}+2k\pi,
\end{cases}
\]
the integration contour $L$ being the one represented in Fig.~\ref{pathL}, together with the identities
\[ D_{-\frac{1}{2}}(z)=\left(\frac{z}{2\pi}\right)^\frac{1}{2}K_{\frac{1}{4}}\left(\frac{1}{4}z^2\right),\qquad
K_{\nu} (z )=\begin{cases}
\frac{\pi{\rm i}}{2}{\rm e}^{\frac{\nu\pi{\rm i}}{2}}H^{(1)}_\nu\big(z{\rm e}^{\frac{\pi{\rm i}}{2}}\big),\\
-\frac{\pi{\rm i}}{2}{\rm e}^{-\frac{\nu\pi{\rm i}}{2}}H^{(2)}_\nu\big(z{\rm e}^{-\frac{\pi{\rm i}}{2}}\big).
\end{cases}
\]

It follows that
\begin{align*}\xi^L_{(2),1}(z,0)=\xi^R_{(2),1}(z,0)&=\pm {\rm i}\pi^{-\frac{1}{2}}z^{\frac{1}{2}}\frac{6x_2^2+h}{2\sqrt{2}(x_1-x_2)(x_3-x_2)}\mathfrak I_2(z,0)\\
&=\pm \frac{{\rm i}\sqrt{\pi}}{8}h^\frac{1}{2}{\rm e}^{\frac{5{\rm i}\pi}{8}}z^\frac{1}{2}{\rm e}^{-\frac{h^2z}{8}}H^{(1)}_\frac{1}{4}\left(\frac{{\rm i}h^2z}{8}\right),
\end{align*}
which coincides (up to an irrelevant sign) with the solution computed in the previous section at the coalescence point. The computations for $\xi_{(3),1}^L(z,0)=\xi_{(3),1}^R(z,0)$ are identical.

The computations for $\xi_{(1),1}^R$ and $\xi_{(1),1}^L$ are a bit more laborious. First of all let us observe that the integral
\[g(z):=\int_0^\infty\exp\left(-\frac{t^2}{2}-zt\right)t^{-\frac{1}{2}}\,{\rm d} t
\] is convergent for all $z\in\mathbb C$, defining an entire function.\footnote{This is in accordance with the expression of $g$ in terms of the modified Bessel function $K$, which gives
\[g\big({\rm e}^{\pm \pi{\rm i}}z\big)=2^{-\frac{1}{2}}{\rm e}^{\frac{z^2}{4}}{\rm e}^{\pm\frac{\pi{\rm i}}{2}}z^\frac{1}{2}K_\frac{1}{4}\left({\rm e}^{\pm 2\pi{\rm i}}\frac{z^2}{4}\right).
\]
From the symmetry $K_\frac{1}{4}\big({\rm e}^{4\pi{\rm i}}z\big)=-K_\frac{1}{4}(z)$ we deduce that $g\big({\rm e}^{-\pi{\rm i}}z\big)=g\big({\rm e}^{\pi{\rm i}}z\big)$.}
Moreover we have
\[
g(z)=\sqrt{\pi}{\rm e}^{\frac{z^2}{4}}D_{-\frac{1}{2}}(z)=2^{-\frac{1}{2}}{\rm e}^{\frac{z^2}{4}}z^\frac{1}{2}K_\frac{1}{4}\left(\frac{z^2}{4}\right).
\]
With a change of variable $t={\rm e}^{-{\rm i}\theta}\tau$ that rotates the half line $\mathbb R_+$ by $\theta$, we find the following identity
\begin{gather}\label{formulautile}
g(z)={\rm e}^{-{\rm i}\frac{\theta}{2}}\int_{{\rm e}^{{\rm i}\theta}\mathbb R_+}\exp\left(-{\rm e}^{-2{\rm i}\theta}\frac{\tau^2}{2}-{\rm e}^{-{\rm i}\theta}z\tau\right)\tau^{-\frac{1}{2}}\,{\rm d}\tau.
\end{gather}

For $t_2=0$ the integral $\mathfrak I_1(z,0)$ splits into two pieces:
\[\mathfrak I_1(z,0)=\mathfrak I_1^1(z)+\mathfrak I_1^2(z),\qquad\mathfrak I_1^i(z):=\int_{\mathcal I_1^i}\exp\big(z\big(x^4+hx^2\big)\big)\,{\rm d}x,\qquad i=1,2,
\]where the paths $\mathcal I_1^i$ are as in Fig.~\ref{camminiintcoal}. Setting $x=2^{-\frac{1}{4}}z^{-\frac{1}{4}}s^\frac{1}{2}$, the image of the paths $\mathcal I_1^i$ are in two different sheets of the Riemann surface with local coordinate~$s$. Keeping track of this, and of the orientations of the modified paths, using formula~\eqref{formulautile} for $\theta=\frac{3\pi{\rm i}}{2},\frac{5\pi{\rm i}}{2}$ and a small deformation of the paths of integration, we find that
\begin{align*}\mathfrak I_1^1(z)&=2^{-\frac{5}{4}}z^{-\frac{1}{4}}\left(-\int_{{\rm e}^{\frac{3\pi{\rm i}}{2}}\mathbb R_+}\exp\left\{\frac{s^2}{2}+\frac{hz^\frac{1}{2}}{\sqrt{2}}s\right\}s^{-\frac{1}{2}}\,{\rm d}s\right)
 =-2^{-\frac{5}{4}}z^{-\frac{1}{4}}{\rm e}^{\frac{3\pi{\rm i}}{4}}g\left({\rm e}^{\frac{\pi{\rm i}}{2}}\frac{hz^\frac{1}{2}}{\sqrt{2}}\right)\\
&=-2^{-\frac{5}{4}}z^{-\frac{1}{4}}{\rm e}^{\frac{3\pi{\rm i}}{4}}\cdot 2^{-\frac{1}{2}}{\rm e}^{-\frac{h^2z}{8}}\left({\rm e}^{\frac{\pi{\rm i}}{2}}\frac{hz^\frac{1}{2}}{\sqrt{2}}\right)^\frac{1}{2}K_{\frac{1}{4}}\left({\rm e}^{\pi{\rm i}}\frac{h^2z}{8}\right)
=\frac{1}{4}{\rm e}^{-\frac{h^2z}{8}}h^\frac{1}{2}K_{\frac{1}{4}}\left({\rm e}^{\pi{\rm i}}\frac{h^2z}{8}\right),
\end{align*}
and
\begin{align*}
\mathfrak I_1^2(z)&=2^{-\frac{5}{4}}z^{-\frac{1}{4}}\left(\int_{{\rm e}^{\frac{5\pi{\rm i}}{2}}\mathbb R_+}\exp\left\{\frac{s^2}{2}+\frac{hz^\frac{1}{2}}{\sqrt{2}}s\right\}s^{-\frac{1}{2}}\,{\rm d}s\right)
 =2^{-\frac{5}{4}}z^{-\frac{1}{4}}{\rm e}^{\frac{5\pi{\rm i}}{4}}g\left({\rm e}^{-\frac{\pi{\rm i}}{2}}\frac{hz^\frac{1}{2}}{\sqrt{2}}\right)\\
&=2^{-\frac{5}{4}}z^{-\frac{1}{4}}{\rm e}^{\frac{5\pi{\rm i}}{4}}\cdot 2^{-\frac{1}{2}}{\rm e}^{-\frac{h^2z}{8}}\left({\rm e}^{-\frac{\pi{\rm i}}{2}}\frac{hz^\frac{1}{2}}{\sqrt{2}}\right)^\frac{1}{2}K_{\frac{1}{4}}\left({\rm e}^{-\pi{\rm i}}\frac{h^2z}{8}\right) \\
& =-\frac{1}{4}{\rm e}^{-\frac{h^2z}{8}}h^\frac{1}{2}K_{\frac{1}{4}}\left({\rm e}^{-\pi{\rm i}}\frac{h^2z}{8}\right).
\end{align*}
Thus, in the limit $t_2=0$ we find that
\begin{align*}
\xi_{(1),1}^R(z,0)&=\pm {\rm i} \Psi_{11}\pi^{-\frac{1}{2}}z^\frac{1}{2}\big(6x_1^2+h\big)^\frac{1}{2}\big(\mathfrak I_1^1(z)+\mathfrak I_1^2(z)+\mathfrak I_3(z,0)\big)\\
&=\pm {\rm i} 2^{-\frac{5}{2}}\pi^{-\frac{1}{2}}z^\frac{1}{2}{\rm e}^{-\frac{h^2z}{8}}h^\frac{1}{2}\left\{K_\frac{1}{4}\left({\rm e}^{{\rm i}\pi}\frac{h^2z}{8}\right)-K_\frac{1}{4}\left({\rm e}^{-{\rm i}\pi}\frac{h^2z}{8}\right)+2^{\frac{1}{2}}K_\frac{1}{4}\left(\frac{h^2z}{8}\right)\right\}\\
&=\pm\pi^\frac{1}{2}z^\frac{1}{2}2^{-\frac{7}{2}}{\rm e}^{-\frac{h^2z}{8}}h^\frac{1}{2}{\rm e}^{-\frac{\pi{\rm i}}{8}}\\
& \quad{}
\times \left\{H^{(2)}_\frac{1}{4}\left({\rm e}^{\frac{{\rm i}\pi}{2}}\frac{h^2z}{8}\right)+{\rm e}^{\frac{\pi{\rm i}}{4}} H^{(1)}_\frac{1}{4}\left({\rm e}^{-\frac{{\rm i}\pi}{2}}\frac{h^2z}{8}\right)+2^{\frac{1}{2}}H^{(2)}_\frac{1}{4}\left({\rm e}^{-\frac{{\rm i}\pi}{2}}\frac{h^2z}{8}\right)\right\}\\
&=\pm \pi^\frac{1}{2}2^{-\frac{3}{2}}z^\frac{1}{2}{\rm e}^{-\frac{h^2z}{8}}h^\frac{1}{2}{\rm e}^{-\frac{\pi{\rm i}}{8}}H^{(2)}_\frac{1}{4}\left({\rm e}^{\frac{{\rm i}\pi}{2}}\frac{h^2z}{8}\right),
\end{align*}
which is exactly (modulo irrelevant signs) the solution at the coalescence point as computed in the previous section. We leave as an exercise for the reader to show that all other solutions $\xi_{(i),j}^{R/L}(z)$ converge to the ones computed at the coalescence point.

\section[Second example of application of Theorem \ref{mainisoth}: quantum cohomology\\ of the Grassmannian $\mathbb{G}_2\big(\mathbb C^4\big)$ and $\Gamma$-conjecture]{Second example of application of Theorem \ref{mainisoth}: quantum\\ cohomology of the Grassmannian $\boldsymbol{\mathbb{G}_2\big(\mathbb C^4\big)}$ and $\boldsymbol{\Gamma}$-conjecture}\label{29novembre2016-1}

In this section we prove Theorem~\ref{resultg24}, which is Theorem \ref{6-07-17-4} of the Introduction, we establish a~correspondence between each region of the quantum cohomology and a full exceptional collection, obtaining an {explicit} refinement of the original conjecture of~\cite{dubro0}. We also prove Proposition~\ref{29luglio2016-4}, which we believe to be an important characterisation of $\mathcal C_0(\eta,\mu, R)$.

The problem is to compute the monodromy data for the Frobenius manifold $QH^\bullet\big(\mathbb{G}_2\big(\mathbb C^4\big)\big)$, the quantum cohomology of the Grassmannian $\mathbb{G}_2\big(\mathbb C^4\big)$. This manifold has a locus of semisimple coalescent points, called small quantum cohomology. Moreover, the structure of the manifold is known {\it only} at the small quantum cohomology locus. Therefore, if we want an explicit computation of monodromy data, this can be done only at coalescence points. This is what we will do: the data will be calculated at the coalescence point $t=0$ in the small cohomology locus. The validity of the result for the whole Frobenius manifold $QH^\bullet\big(\mathbb{G}_2\big(\mathbb C^4\big)\big)$ is completely justified by our Theorem~\ref{mainisoth}, which is thus {\it crucial} to us. Without Theorem~\ref{mainisoth} our computations would be geometrically meaningless.\footnote{Without Theorem~\ref{mainisoth}, we would anyway have an interesting and non-trivial example of computation of monodromy data for a $6\times 6$ differential system with two coinciding eigenvalues at the irregular singularity. The purpose of this article goes beyond this; our goal is the study the monodromy data of a Frobenius manifold.}

 As a result of the computations, we prove Theorem \ref{resultg24} (Theorem \ref{6-07-17-4}), which clarifies and {verifies} a conjecture, formulated by the second author in \cite{dubro0} and then refined\footnote{The detailed comparison between the explicit computations of the monodromy data for complex Grassmannians and the $\Gamma$-classes proposed in \cite{gamma1}, is one of the contents of our paper \cite{CDG1}.} in \cite{dubro4} and \cite{gamma1}, in the case the quantum cohomology of the Grassmannian $\mathbb{G}_2\big(\mathbb C^4\big)$. Our explicit computations, using elementary analytic methods only, seems to be missing from the literature.

\subsection{Notations in Gromov--Witten theory}\label{qcohg24}
Let $X$ be a smooth complex projective variety with vanishing odd cohomology
\[H^{2k+1}(X;\mathbb C)=0,\qquad k\geq 0.
\]
Let us fix a homogeneous basis $(T_1,T_2,\dots, T_N)$ of $H^\bullet(X;\mathbb C)=\bigoplus_{k}H^{2k}(X;\mathbb C)$ such that
\begin{itemize}\itemsep=0pt
\item $T_1=1$ is the unity of the cohomology ring;
\item $\deg T_\alpha=:2q_\alpha$;
\item $T_2,\dots, T_r$ span $H^2(X;\mathbb C)$.
\end{itemize}
We will denote by $\eta\colon H^\bullet(X;\mathbb C)\times H^\bullet(X;\mathbb C)\to\mathbb C$ the Poincar\'e metric
\[\eta(\xi,\zeta):=\int_X\xi\cup\zeta,
\qquad
\eta_{\alpha\beta}:=\int_XT_\alpha\cup T_\beta.
\]
We also introduce the \emph{Novikov ring} $\Lambda:=\mathbb C[\![Q_2,\dots,Q_r]\!]$, and the symbol
\[Q^\beta:=Q_2^{\int_\beta T_2} \cdots  Q_r^{\int_\beta T_r}.
\]
Let $X_{g,n,\beta}$ be the moduli space of stable maps with target $X$, of genus $g$, with $n$ distinct marked points and of degree $\beta\in H_2(X;\mathbb Z)$. We will denote by
\[\langle\tau_{d_1}\gamma_1,\dots,\tau_{d_n}\gamma_{n}\rangle^X_{g,n,\beta}:=\int_{[X_{g,n,\beta}]^{\text{vir}}}\bigwedge_{i=1}^n\operatorname{ev}^*_i(\gamma_i)\cup\psi_i^{d_i}
\]the value of the \emph{Gromov--Witten invariant} (with \emph{gravitational descendants}, if some of the $d_i$'s is nonzero), where
\begin{itemize}\itemsep=0pt
\item $\gamma_1,\dots,\gamma_n\in H^\bullet(X;\mathbb C)$,
\item $(d_1,\dots, d_n)\in\mathbb N^n$,
\item $\psi_1,\dots,\psi_n\in H^2(X_{g,n,\beta};\mathbb Q)$ are the universal cotangent line classes,
\item $\operatorname{ev}_i\colon X_{g,n,\beta}\to X$ is the evaluation map at the $i$-th marked point,
\item $[X_{g,n,\beta}]^{\text{vir}}$ stands for the \emph{virtual fundamental class}. Recall that degree of the virtual cycle is equal to the virtual dimension (over $\mathbb R$)
\[\operatorname{vir\ dim}_{\mathbb R} X_{g,n,\beta}=2(1-g)\dim_\mathbb CX-2\int_{\beta}\omega_X+2(3g-3+n).
\]
\end{itemize}

It is convenient to collect Gromov--Witten invariants with descendants as coefficients of a~ge\-ne\-rating function, called \emph{genus~$g$ gravitational Gromov--Witten potential}, or simply \emph{genus~$g$ free energy}
\[
\mathcal F_g^X(\gamma):=\sum_{n=0}^\infty\sum_{\beta\in \operatorname{Eff}(X)}\frac{Q^\beta}{n!}\langle\underbrace{\gamma,\dots,\gamma}_{n\text{ times}}\rangle_{g,n,\beta}^X,
\]
the set $\operatorname{Eff}(X)\subseteq H_2(X;\mathbb Z)$ being the set of effective classes of $X$. Introducing (infinitely many) coordinates $\bold{t}:=(t^{\alpha,p})_{\alpha,p}$
\[\gamma=\sum_{\alpha,p}t^{\alpha,p}\tau_pT_{\alpha},
\] the free energy $\mathcal F^X_g\in\Lambda[\![{\bold t}]\!]$ can be seen as a function on the \emph{large phase-space}, and restricting the free energy to the \emph{small phase space} (naturally identified with $H^\bullet(X;\mathbb C)$),
\[F^X_g\big(t^{1,0},\dots, t^{N,0}\big):=\mathcal F^X_g(\bold t)\big|_{t^{\alpha,p}=0,\ p>0},
\]
one obtains the generating function of the Gromov--Witten invariants of genus $g$. It will also be convenient to introduce the \emph{genus g correlation functions} defined by the derivatives
\[\llangle\tau_{d_1}T_{\alpha_1},\dots,\tau_{d_n}T_{\alpha_n}\rrangle_g:=\frac{\partial}{\partial t^{\alpha_1,d_1}}\cdots\frac{\partial}{\partial t^{\alpha_n,d_n}}\mathcal F^X_g.
\]

 Let $t^\alpha:=t^{\alpha,0}$. By the Divisor axiom, the genus $0$ Gromov--Witten potential $F^X_0(t)$ can be seen as an element of the ring $\mathbb C\big[\!\big[t^1,Q_2{\rm e}^{t^2},\dots,Q_r{\rm e}^{t^r},t^{r+1},\dots, t^N\big]\!\big]$. In what follows we will be interested in the cases when $F^X_0$ is a convergent series
expansion
\begin{gather}\label{26gennaio2020-3}
F^X_0\in\mathbb C\big\{t^1,Q_2{\rm e}^{t^2},\dots,Q_r{\rm e}^{t^r},t^{r+1},\dots, t^N\big\}.
\end{gather}
 Without loss of generality we can put $Q_2=Q_3=\dots=Q_r=1$. Under the assumption~(\ref{26gennaio2020-3}), $F^X_0(t)$ defines an analytic function in an open neighbourhood $\Omega\subseteq H^\bullet(X;\mathbb C)$ of the point
\begin{gather}\label{classical1}
t^i=0,\qquad i=1,r+1,\dots, N;
\qquad
\operatorname{Re} t^i\to -\infty,\qquad i=2,3,\dots, r.
\end{gather}
The function $F^X_0$ is a solution of WDVV equations \cite{kon,manin,tian,voisin}, and thus it defines an analytic Frobenius manifold structure on $\Omega$. Using the canonical identifications of tangent spaces $T_p\Omega\cong H^\bullet(X;\mathbb C)\colon \partial_{t^\alpha}\mapsto T_\alpha$, the unit vector field is $e=\partial_{t^1}\equiv 1$, and the Euler vector field is
\[
E:=c_1(X)+\sum_{\alpha=1}^N\left(1-\frac{1}{2}\deg T_\alpha\right)t^\alpha T_\alpha.
\]
The resulting Frobenius structure is called \emph{quantum cohomology of }$X$, denoted $QH^\bullet (X)$. Notice that at the classical limit point \eqref{classical1} the algebra structure on the tangent spaces coincides with the classical cohomological algebra structure. Notice that, because of the divisor axiom, the Frobenius structure is $2\pi{\rm i}$-periodic in the 2-nd cohomology directions: the structure can be considered as defined on an open region of the quotient $H^\bullet(X;\mathbb C)/2\pi{\rm i} H^2(X;\mathbb Z)$.

There are no general results characterizing smooth projective varieties with semisimple quantum cohomology: however, for some classes of varieties such as
\begin{itemize}\itemsep=0pt
\item some Fano threefolds \cite{ciolli},
\item toric varieties \cite{iri7},
\item some homogeneous spaces \cite{cmp10},
\end{itemize}
it has been proved that the small quantum locus is made of semisimple points. Grassmannians are among these varieties.
Below, we focus on the small quantum cohomology of $\mathbb G_2\big(\mathbb C^4\big)$, namely the restriction to the locus $H^2\big(\mathbb G_2\big(\mathbb C^4\big);\mathbb C\big)$, with {coordinates} $\big(0,t^2,0,\dots ,0\big)$.

\subsection[Small quantum cohomology of $\mathbb{G}_2\big(\mathbb C^4\big)$]{Small quantum cohomology of $\boldsymbol{\mathbb{G}_2\big(\mathbb C^4\big)}$}\label{smcohg24}

\subsubsection{Generalities and proof of its semisimplicity}

For simplicity, let us use the notation $\mathbb G:= \mathbb{G}_2\big(\mathbb C^4\big)$.
From the general theory of Schubert calculus, it is known that $H^\bullet(\mathbb G;\mathbb C)$ is a complex vector space of dimension 6, and a basis is given by \emph{Schubert classes}:
\[\sigma_0:=1,\ \sigma_1,\ \sigma_2,\ \sigma_{1,1},\ \sigma_{2,1},\ \sigma_{2,2}.
\]
Each $\sigma_\lambda\in H^{2|\lambda|}(\mathbb G;\mathbb C)$. By posing
\[v_1:=\sigma_0,\qquad v_2:=\sigma_1,\qquad v_3:=\sigma_{2},\qquad v_4:=\sigma_{1,1},\qquad v_{5}:=\sigma_{2,1},\qquad v_6:=\sigma_{2,2},
\]
we will denote by $t^i$ the coordinate with respect to $v_i$. The coordinates in the small quantum cohomology are
\[
t=\big(0,t^2,0,\dots ,0\big).
\]
 By Pieri and Giambelli formulas one finds that the matrix of the Poincar\'e pairing
\[\eta(\alpha,\beta):=\int_{\mathbb G}\alpha\wedge\beta
\]
with respect to the {above} basis is given by
\[\eta=\begin{pmatrix}
0& 0 & 0 & 0 & 0 & c\\
0& 0 & 0 & 0 & c & 0\\
0& 0 & c & 0 & 0 & 0\\
0& 0 & 0 & c & 0 & 0\\
0& c & 0 & 0 & 0 & 0\\
c& 0 & 0 & 0 & 0 & 0
\end{pmatrix},\qquad c:=\int_{\mathbb{G}}\sigma_{2,2}.
\]
Using quantum Pieri--Bertram formula \cite{bert}, we deduce that the
matrix of the operator of multiplication by $\lambda\sigma_1+\mu\sigma_{1,1}$ is
\begin{gather}\label{matr1}\begin{pmatrix}
0&0&\mu q&0&\lambda q&0\\
\lambda&0&0&0&\mu q&\lambda q\\
0&\lambda&0&0&0&\mu q\\
\mu&\lambda&0&0&0&0\\
0&\mu&\lambda&\lambda&0&0\\
0&0&0&\mu&\lambda&0
\end{pmatrix},\qquad  q:={\rm e}^{t^2}.
\end{gather}
The discriminant of the characteristic polynomial of this matrix is
\[16777216 \lambda^4\mu^2q^8\big(\lambda^4+\mu^4q\big)^6
\]
and so, if $\lambda\neq 0$, $\mu\neq 0$ and $\lambda^4+q\mu^4\neq 0$, its eigenvalues are pairwise distinct. This is a sufficient condition to state that the quantum cohomology of $\mathbb G$ is semisimple.

 Notice that the value at the point $p$ of coordinates $\big(0,t^2,0,\dots,0\big)$ of the Euler field of quantum cohomology $QH^\bullet(\mathbb G)$ is\footnote{We identify $T_pH^\bullet(\mathbb G)$ with $H^\bullet(\mathbb G)$ in the canonical way.} given by the first Chern class $c_1(\mathbb G)=4\sigma_1$:
\[E|_p=4\frac{\partial}{\partial t^2}\equiv 4\sigma_1.
\]
The matrix $\mathcal U$ of multiplication by $E$ at the point $p$ is given by posing $\lambda=4$, $\mu=0$ in~\eqref{matr1}:
\[\mathcal U\big(0,t^2,0,\dots,0\big)\equiv 4 \mathcal C_2\big(0,t^2,0,\dots,0\big)=
\begin{pmatrix}
0& 0 & 0 & 0 & 4q & 0\\
4& 0 & 0 & 0 & 0 & 4q\\
0& 4 & 0 & 0 & 0 & 0\\
0& 4 & 0 & 0 & 0 & 0\\
0& 0 & 4 & 4 & 0 & 0\\
0& 0 & 0 & 0 & 4 & 0
\end{pmatrix}.
\]
 The characteristic polynomial is $p(z)=z^6-1024qz^2$, so that 0 is an eigenvalue with multipli\-ci\-ty~2. Therefore, the semisimple points with coordinates $\big(0,t^2,0,\dots,0\big)$ are \emph{semisimple coalescence points in the bifurcation set}.

\subsubsection[Idempotents at the points $\big(0,t^2,0,\dots,0\big)$]{Idempotents at the points $\boldsymbol{\big(0,t^2,0,\dots,0\big)}$}

 The multiplication by $\sigma_1+\sigma_{1,1}$ has pairwise distinct eigenvalues, at least at points for which $t^2\neq {\rm i}\pi(2k+1)$. Putting $\lambda=\mu=1$ in~\eqref{matr1}, we deduce that the characteristic polynomial of this operator is
\[p(z)=\big(q+z^2\big)\big({-}4q+q^2-8qz-2qz^2+z^4\big).
\]
So the six eigenvalues are
\begin{gather*}
{\rm i}q^{\frac{1}{2}},\qquad  -{\rm i}q^{\frac{1}{2}},\\
 \varepsilon_1:=-{\rm i}\sqrt2q^{\frac{1}{4}}-q^{\frac{1}{2}},
 \qquad
 \varepsilon_2:={\rm i}\sqrt2q^{\frac{1}{4}}-q^{\frac{1}{2}},
 \qquad
 \varepsilon_3:=-\sqrt2q^{\frac{1}{4}}+q^{\frac{1}{2}},
 \qquad
 \varepsilon_4:=\sqrt2q^{\frac{1}{4}}+q^{\frac{1}{2}},
\end{gather*}
and the corresponding eigenvectors are
\begin{gather*}
\pi_1:=-q-{\rm i}q^{\frac{1}{2}}\sigma_2+{\rm i}q^{\frac{1}{2}}\sigma_{1,1}+\sigma_{2,2},
\qquad
\pi_2:=-q+{\rm i}q^{\frac{1}{2}}\sigma_2-{\rm i}q^{\frac{1}{2}}\sigma_{1,1}+\sigma_{2,2},\\
\pi_{2+i}:=\big(q^2+q\varepsilon_i^2\big)+\big({-}q^2+2q\varepsilon_i+q\varepsilon_i^2\big)\sigma_1
+(2q+2q\varepsilon_i)\sigma_2+(2q+2q\varepsilon_i)\sigma_{1,1}\\
\hphantom{\pi_{2+i}:=}{} +\big({-}2q-q\varepsilon_i+\varepsilon_i^3\big)\sigma_{2,1}+\big(q+\varepsilon_i^2\big)\sigma_{2,2}.
\end{gather*}
Then,
\[\pi_i\cdot\pi_j=0\ \ \text{if }i\neq j,\qquad \pi_i^2=\lambda_i \pi_i\ \ \text{where }\lambda_i>0;
\]as a consequence, these vectors are orthogonal since, for $i\neq j$,
\[ \eta(\pi_i,\pi_j)=\eta(\pi_i\cdot\pi_j,1)=\eta(0,1)=0.\] Introducing the normalized eigenvectors
\[f_i:=\frac{\pi_i}{\eta(\pi_i,\pi_i)^\frac{1}{2}}
\]we obtain an orthonormal frame of normalized idempotent vectors, for any choice of the sign of the square roots.

Let us now introduce a matrix $\Psi=(\psi_{ij})$ such that
\[
\frac{\partial}{\partial t_\alpha}=\sum_i\psi_{i\alpha}f_i,\qquad \alpha=1,2,\dots ,n.
\]
Note that necessarily we have
\[\Psi^{\rm T}\Psi=\eta,\qquad\psi_{i1}=\frac{\eta(\pi_i,1)}{\eta(\pi_i,\pi_i)^{\frac{1}{2}}}.
\]
After some computations, we obtain
\[\Psi=\frac{c^{\frac{1}{2}}}{2}
\begin{pmatrix}
-{\rm i}q^{-\frac{1}{2}}& 0& -1& 1& 0& {\rm i}q^{\frac{1}{2}}\\
-{\rm i}q^{-\frac{1}{2}}& 0& 1& -1& 0& {\rm i}q^{\frac{1}{2}}\\
\frac{1}{\sqrt{2}q^{\frac{1}{2}}}& -\frac{{\rm i}}{q^\frac{1}{4}}& -\frac{1}{\sqrt 2}& -\frac{1}{\sqrt 2}& {\rm i}q^{\frac{1}{4}}& \frac{q^{\frac{1}{2}}}{\sqrt 2}\\
\frac{1}{\sqrt{2}q^{\frac{1}{2}}}& \frac{{\rm i}}{q^\frac{1}{4}}& -\frac{1}{\sqrt 2}& -\frac{1}{\sqrt 2}& -{\rm i}q^{\frac{1}{4}}& \frac{q^{\frac{1}{2}}}{\sqrt 2}\\
\frac{1}{\sqrt{2}q^{\frac{1}{2}}}& -\frac{1}{q^\frac{1}{4}}& \frac{1}{\sqrt 2}& \frac{1}{\sqrt 2}& -q^{\frac{1}{4}}& \frac{q^{\frac{1}{2}}}{\sqrt 2}\\
\frac{1}{\sqrt{2}q^{\frac{1}{2}}}& \frac{1}{q^\frac{1}{4}}& \frac{1}{\sqrt 2}& \frac{1}{\sqrt 2}& q^{\frac{1}{4}}& \frac{q^{\frac{1}{2}}}{\sqrt 2}\\
\end{pmatrix}.
\]
This matrix diagonalizes $\mathcal U$ as follows
\begin{gather}
U:=\Psi\mathcal U\Psi^{-1}=\big(\Psi^{\rm T}\big)^{-1}\widehat{\mathcal U}\Psi^{\rm T} \nonumber\\
\label{7febr2020-1}
\hphantom{U}{} =\begin{pmatrix}
u_1&&&&&\\
&u_2&&&&\\
&&u_3&&&\\
&&&u_4&&\\
&&&&u_5&\\
&&&&&u_6
\end{pmatrix}=4\sqrt{2}q^{\frac{1}{4}}\left(
\begin{matrix}
 0 & 0 & 0 & 0 & 0 & 0 \\
 0 & 0 & 0 & 0 & 0 & 0 \\
 0 & 0 & - {\rm i} & 0 & 0 & 0 \\
 0 & 0 & 0 & {\rm i} & 0 & 0 \\
 0 & 0 & 0 & 0 & -1& 0 \\
 0 & 0 & 0 & 0 & 0 & 1
\end{matrix}
\right).
\end{gather}
The eigenvalues $u_i$ stand for $u_i\big(0,t^2,\dots,0\big)$. Note that
\begin{gather}\label{28luglio2016-12}
u_i\big(0,t^2,0,\dots ,0\big)=q^{\frac{1}{4}} u_i(0,0,\dots ,0)={\rm e}^{\frac{t^2}{4}}u_i(0,0,\dots ,0).
\end{gather}

\subsection{Differential system for deformed flat coordinates}

Our goal is to obtain the monodromy data for the small quantum cohomology. Therefore, we need to consider system (\ref{28luglio2016-4}), which we rewrite as follows:
\begin{gather}
\label{28luglio2016-5-A}
\partial_z \xi=\left(\widehat{\mathcal U}-\frac{1}{z}\mu\right) \xi,\\
\label{28luglio2016-5-B}
\partial_2 \xi=z \widehat{\mathcal C_2} \xi,
\end{gather}
where $\xi$ is a column vector, whose components are $\xi_i=\partial_i \tilde{t}(t,z)$ (derivatives of a~deformed flat coordinate), and the matrix coefficients are
\begin{gather} \widehat{\mathcal U}:=\eta \mathcal U \eta^{-1}=
\begin{pmatrix}
0&4&0&0&0&0\\
0&0&4&4&0&0\\
0&0&0&0&4&0\\
0&0&0&0&4&0\\
4q&0&0&0&0&4\\
0&4q&0&0&0&0
\end{pmatrix},
\qquad
\widehat{\mathcal C_2}\equiv \frac{1}{4}~\widehat{\mathcal U},\nonumber\\
\label{2luglio2017-1}
\mu=\operatorname{diag} (-2,-1,0,0,1,2 ),
\end{gather}
with eigenvalues $\mu_\alpha=\frac{\deg (\partial/\partial_\alpha)-4}{2}$, $1\leq \alpha \leq 6$.
As it is customary in the analysis of differential systems, it is convenient to try a reduction to an equivalent scalar equation. To this purpose, we introduce the scalar function~$\phi$ defined by
\[\phi(t,z):=\frac{\xi_1(t,z)}{z^2}.
\]
In this way, the first equation of system~\eqref{28luglio2016-5-A} becomes a single {linear} differential equation
\begin{gather}\label{eq1}z^4\partial_z^5\phi +10z^3\partial_z^4\phi+25z^2\partial_z^3\phi+15z\partial_z^2\phi+\big(1-1024qz^4\big)\partial_z\phi-2048qz^3\phi=0.
\end{gather}
A column vector solution of \eqref{28luglio2016-5-A} can be reconstructed from the formulae
\begin{gather}
\xi_1 =z^2\phi,\nonumber\\
\xi_2 =\frac{1}{4}z^2\partial_z\phi,\nonumber\\
\xi_3 =\frac{1}{32}\big(z\partial_z\phi+z^2\partial_z^2\phi\big)+h,\nonumber\\
\xi_4 =\frac{1}{32}\big(z\partial_z\phi+z^2\partial_z^2\phi\big)-h,\nonumber\\
\xi_5 =\frac{1}{128}\big(\partial_z\phi+3z\partial_z^2\phi+z^2\partial_z^3\phi\big),\nonumber\\
\xi_6 =\frac{1}{512}\left(-512qz^2\phi+\frac{1}{z}\partial_z\phi+7\partial_z^2\phi +6z\partial_z^3\phi+z^2\partial_z^4\phi\right)\label{reconstruction1}
\end{gather} with a constant $h\in \mathbb{C}$ to be determined later (formulae \eqref{reconstruction2}--\eqref{reconstruction4}).

\begin{Remark}\label{7febbraio2020-2}The third and the fourth equalities \eqref{reconstruction1} follow from the fact that
\[\xi_3+\xi_4=\frac{1}{16}\big(z\partial_z\phi+z^2\partial_z^2\phi\big),\qquad\partial_z(\xi_3-\xi_4)=0, \qquad\partial_2(\xi_3-\xi_4)=0,
\]so that $\xi_3-\xi_4=2h$ is a constant.

This fact reflects the peculiarity of our systems~\eqref{28luglio2016-5-A} of being direct sum of a 1-dimensional and 5-dimensional systems. Indeed, by the change of variable $\widetilde{\xi}_j=\xi_j$, $j=1,2,5,6$, $\widetilde{\xi}_3=\xi_3+\xi_4$ and $\widetilde{\xi}_4=\xi_3-\xi_3$, form~\eqref{28luglio2016-5-A} we obtain the systems
\begin{gather*}
\partial_z \widetilde{\xi}_4=0 \qquad\!\hbox{and} \!\qquad \partial_z\left(\begin{matrix}
\widetilde{\xi}_1
\\
\widetilde{\xi}_2
\\
\widetilde{\xi}_3
\\
\widetilde{\xi}_5
\\
\widetilde{\xi}_6
\end{matrix}
\right)=
\left[
\left(\begin{matrix}
0 & 4 & 0 & 0 & 0
\\
0 & 0 & 4 & 0 & 0
\\
0 & 0 & 0 & 8 & 0
\\
4q & 0 & 0 & 0 & 4
\\
0 & 4q & 0 & 0 & 0
\end{matrix}\right)
+
\frac{1}{z}
\left(\begin{matrix}
2 & 0 & 0 & 0 & 0
\\
0 & 1 & 0 & 0 & 0
\\
0 & 0 & 0 & 0 & 0
\\
0 & 0 & 0 & -1 & 0
\\
0 & 0 & 0 & 0 & -2
\end{matrix}\right)
\right]
\left(\begin{matrix}
\widetilde{\xi}_1
\\
\widetilde{\xi}_2
\\
\widetilde{\xi}_3
\\
\widetilde{\xi}_5
\\
\widetilde{\xi}_6
\end{matrix}
\right).
\end{gather*}
The leading matrix has diagonal form $4\sqrt{2}q^{1/4}{\rm diag}(0,-{\rm i},{\rm i},-1,1)$, with {\it distinct} eigenvalues. Also system~\eqref{28luglio2016-5-B} decouples into a direct sum, and
\begin{gather*}
\partial_2 \widetilde{\xi}_4=0.
\end{gather*}
Notice that this decoupling into a direct sum also confirms, in our very particular case, the general structure \eqref{7febbraio2020-3} of the formal solution.
However, this decoupling is a peculiarity of the small cohomology locus of $\mathbb{G}$, but does not apply in general cases at coalescing eigenvalues of~$U$. The peculiarity is evident when we write system \eqref{28luglio2016-5-A}--\eqref{28luglio2016-5-B} as
\begin{gather*}
 y= \big(\Psi^{\rm T}\big)^{-1}\xi,\qquad \partial_z y=\left(U+\frac{V}{z}\right)y, \qquad
 \partial_2 y = \frac{1}{4} zU y,
 \end{gather*}
 where $U$ is the diagonal matrix \eqref{7febr2020-1} and
 \begin{gather*}
 V=\left(
 \begin{matrix}
 0 & 0 & {\rm i}/\sqrt{2} & {\rm i}/\sqrt{2} & {\rm i}/\sqrt{2} & {\rm i}/\sqrt{2}
 \\
 0 & 0 & {\rm i}/\sqrt{2} & {\rm i}/\sqrt{2} & {\rm i}/\sqrt{2} & {\rm i}/\sqrt{2}
 \\
 - {\rm i}/\sqrt{2}& - {\rm i}/\sqrt{2}& 0 & 0 & -{\rm i}/2 & {\rm i}/2
 \\
 - {\rm i}/\sqrt{2}&- {\rm i}/\sqrt{2} & 0 & 0 & {\rm i}/2 & -{\rm i}/2
 \\
 - {\rm i}/\sqrt{2} & -{\rm i}/\sqrt{2}& {\rm i}/2&-{\rm i}/2 & 0 & 0
 \\
 - {\rm i}/\sqrt{2} & -{\rm i}/\sqrt{2} & -{\rm i}/2&{\rm i}/2& 0 & 0
 \end{matrix}
 \right).
 \end{gather*}
Here $y=(y_1,\dots,y_6)^{\rm T}$ is a column vector. Now, the aforementioned peculiarity is that the first two columns of $V$ are equal. Hence, we can decouple into a direct sum using the variables $\widetilde{y}_1=y_1-y_2$, $ \widetilde{y}_2=y_1+y_2$ and $ \widetilde{y}_j=y_j$ for $j=3,4,5, 6$. Moreover, the first two rows are equal, which implies $\partial_z y_1=\partial_2 \widetilde{y}_1=0$.

We thank the anonymous referee for suggesting the above observations.
\end{Remark}

From (\ref{28luglio2016-5-B}) it follows that
\[\partial_2\phi=\frac{z}{4}\partial_z\phi,
\]
which implies the following functional form
\[\phi(t_2,z)=\vp\big(zq^{\frac{1}{4}}\big),
\]
for a scalar function $\vp$ of one variable. As a consequence, our problem~\eqref{eq1} reduces to the solution of a single scalar ordinary differential equation with independent variable $w=zq^{\frac{1}{4}}$ and dependent variable $\vp(w)$:
\[w^4\vp^{(5)}+10w^3\vp^{(4)}+25w^2\vp^{(3)}+15w\vp''+\big(1-1024w^4\big)\vp'-2048w^3\vp=0.
\]
Multiplying by $w\in\mathbb C^*$, we can rewrite this equation in a more compact form
\begin{gather}\label{eqip}\Theta^5\vp-1024w^4\Theta\vp-2048w^4\vp=0,
\end{gather}
where $\Theta$ is the Euler's differential operator $w\frac{d}{dw}$. The fact that we have reduced the problem to a fifth order scalar differential equation reflects the observation of Remark~\ref{7febbraio2020-2} that the system is a direct sum of a trivial 1-dimensional system and a non-trivial 5-dimensional one.

\subsubsection{Expected asymptotic expansions}

Let $\Xi$ be a fundamental matrix solution of system~(\ref{28luglio2016-5-A}), and let $Y$ be defined by
\begin{gather}\label{transformation}\Xi=\eta\Psi^{-1}Y.
\end{gather}
Then, $Y$ is a fundamental solution of system~(\ref{semiseq}). The asymptotic theory has been explained in Section~\ref{monossfm}, and Theorem~\ref{fundamental} applies. To the formal solution~(\ref{28luglio2016-7}) it corresponds a formal matrix solution
\[
\Xi_{\rm formal}=\eta\Psi^{-1}G(z)^{-1}{\rm e}^{zU}.
\]
To the fundamental solutions $Y_{\rm left/right}$, there correspond solutions $\Xi_{\rm left/right}$.
For fixed $t^2$, then $\Xi_{\rm left/right}\big(t^2,z\big)\equiv \Xi_{\rm left/right}\big({\rm e}^{\frac{t^2}{4}}z\big)$ has the following asymptotic expansion for $z\to \infty$
 \begin{gather} \Xi_{\rm left/right}\big(t^2,z\big)=\eta\Psi^{-1}\left(1+O\left(\frac{1}{z}\right)\right){\rm e}^{zU}
\nonumber\\
=\frac{c^{\frac{1}{2}}}{2} \left(I+O\left(\frac{1}{z}\right)\right) \left(
\begin{matrix}
 -\frac{{\rm i} {\rm e}^{z u_1}}{q^{\frac{1}{2}}}
 &
 -\frac{{\rm i} {\rm e}^{z u_2}}{q^{\frac{1}{2}}} & \frac{{\rm e}^{z u_3}}{\sqrt{2} q^{\frac{1}{2}}} & \frac{{\rm e}^{z u_4}}{\sqrt{2} q^{\frac{1}{2}}}
 & \frac{{\rm e}^{z u_5}}{\sqrt{2} q^{\frac{1}{2}}} & \frac{{\rm e}^{z u_6}}{\sqrt{2} q^{\frac{1}{2}}}
 \\
 0 & 0 & -\frac{{\rm i} {\rm e}^{z u_3}}{q^{\frac{1}{4}}} & \frac{{\rm i} {\rm e}^{z u_4}}{q^{\frac{1}{4}}} & -\frac{{\rm e}^{z u_5}}{q^{\frac{1}{4}}} & \frac{{\rm e}^{z u_6}}{q^{\frac{1}{4}}}
 \\
 -{\rm e}^{z u_1} & {\rm e}^{z u_2} & -\frac{{\rm e}^{z u_3}}{\sqrt{2}} & -\frac{{\rm e}^{z u_4}}{\sqrt{2}} & \frac{{\rm e}^{z u_5}}{\sqrt{2}} & \frac{{\rm e}^{z u_6}}{\sqrt{2}}
 \\
 {\rm e}^{z u_1} & -{\rm e}^{z u_2} & -\frac{{\rm e}^{z u_3}}{\sqrt{2}} & -\frac{{\rm e}^{z u_4}}{\sqrt{2}} & \frac{{\rm e}^{z u_5}}{\sqrt{2}} & \frac{{\rm e}^{z u_6}}{\sqrt{2}}
 \\
 0 & 0 & {\rm i} {\rm e}^{z u_3} q^{\frac{1}{4}} & -{\rm i} {\rm e}^{z u_4} q^{\frac{1}{4}} & -{\rm e}^{z u_5} q^{\frac{1}{4}} & {\rm e}^{z u_6} q^{\frac{1}{4}}
 \\
 {\rm i} {\rm e}^{z u_1} q^{\frac{1}{2}} & {\rm i} {\rm e}^{z u_2} q^{\frac{1}{2}} & \frac{{\rm e}^{z u_3} q^{\frac{1}{2}}}{\sqrt{2}} & \frac{{\rm e}^{z u_4} q^{\frac{1}{2}}}{\sqrt{2}} & \frac{{\rm e}^{z u_5} q^{\frac{1}{2}}}{\sqrt{2}} & \frac{{\rm e}^{z u_6} q^{\frac{1}{2}}}{\sqrt{2}}
 \\
\end{matrix}
\right).\label{28luglio2016-8}
\end{gather}

Now, one of our tasks is to explicitly compute the fundamental matrix solutions $\Xi_{\rm left/right}$ with behaviour~\eqref{28luglio2016-8} and the Stokes matrix connecting them, by means of the formulae \eqref{reconstruction1}. Therefore, we need fundamental systems of solutions of~\eqref{eq1} with the correct asymptotic behaviour such that formulae~\eqref{reconstruction1} will match with~\eqref{28luglio2016-8}. It will suffice to identify solutions $\phi\big(z,t^2\big)$ of \eqref{eq1} with the asymptotic behaviour given by the first row of \eqref{28luglio2016-8}. Such solutions will be identified in Section~\ref{26gennaio2020-2}, making use of the analysis of equation \eqref{eq1} developed in Section~\ref{30gennaio2020-1} below.

The behaviour in \eqref{28luglio2016-8} also allows to find the correct values of~$h$ in~(\ref{reconstruction1}). This values must be determined in order to match with the asymptotics of the third and fourth rows of~(\ref{28luglio2016-8}). We find
\begin{alignat}{3}\label{reconstruction2}
& h=-\frac{c^{\frac{1}{2}}}{2},\qquad&& \text{for the first column},& \\
\label{reconstruction3}
& h=\frac{c^{\frac{1}{2}}}{2},\qquad&& \text{for the second column},& \\
\label{reconstruction4}
& h=0, \qquad && \text{for the remaining columns}.&
\end{alignat}
Indeed, for $\phi$ corresponding to the first two columns we have respectively
\[\phi\big(z,t^2\big)=-\frac{c^{\frac{1}{2}}}{2}\frac{{\rm i} {\rm e}^{zu_1}}{z^2q^{\frac{1}{2}}}\left(1+O\left(\frac{1}{z}\right)\right)
\qquad
\text{or}
\qquad
\phi(z,t^2) =-\frac{c^{\frac{1}{2}}}{2}\frac{{\rm i} {\rm e}^{zu_2}}{z^2q^{\frac{1}{2}}}\left(1+O\left(\frac{1}{z}\right)\right).
\]
Since $u_1=u_2=0$, the above {expressions become}
\[\phi\big(z,t^2\big)
=-\frac{c^{\frac{1}{2}}}{2}\frac{{\rm i} }{z^2q^{\frac{1}{2}}}\left(1+O\left(\frac{1}{z}\right)\right).
\]
Then
\[ \frac{1}{32}\big(z \partial_z \phi + z^2 \partial^2_z\phi\big)=O\left(\frac{1}{z}\right).
\]
Comparing with the matrix elements $(3,1)$, $(4,1)$ and $(3,2)$, $(4,2)$ of~\eqref{28luglio2016-8} respectively, we obtain~(\ref{reconstruction2}) and~(\ref{reconstruction3}). For the remaining columns we proceed in the same way and find~\eqref{reconstruction4}.

\subsection{Solutions of the differential system}\label{30gennaio2020-1}

Defining $\varpi:=4w^4$ and writing $\vp(w)=\tilde{\vp}(\varpi)$, equation \eqref{eqip} becomes
\[
\Theta^5_{\varpi}\tilde{\vp}-\varpi \Theta_\varpi\tilde{\vp}-\frac{1}{2}\varpi\tilde\vp=0,
\qquad
\Theta_\varpi:=\varpi\frac{{\rm d}}{{\rm d}\varpi}=\frac{1}{4}\Theta_w.
\]
This belongs to the class of generalized hypergeometric differential equations (see \cite{abrasteg,luke1,nist,paris-kam} and references therein). By applying the Mellin transform $\mathfrak M$ to it, we obtain the finite difference equation
\begin{gather*}
s^5\tilde\tau(s)=\left(s+\frac{1}{2}\right)\tilde\tau(s+1),\qquad
\tilde\tau(s):=\mathfrak M(\tilde\vp)(s):=\int_0^\infty\tilde\vp(t)t^{s-1}\,{\rm d}t,
\end{gather*}
whose solutions are of the form
\[\tilde\tau(s)=\frac{\Gamma(s)^5}{\Gamma\big(s+\frac{1}{2}\big)}\psi(s),\qquad\psi(s)=\psi(s+1).
\]
Hence, we expect that solutions of \eqref{eqip} are of the form
\[\vp(w)=\frac{1}{2\pi{\rm i}}\int_\Lambda\frac{\Gamma(s)^5}{\Gamma\big(s+\frac{1}{2}\big)}\psi(s)4^{-s}w^{-4s}\,{\rm d}s,
\]for suitable chosen paths of integration $\Lambda$.  Indeed, the lemma below holds:

\begin{Lemma}\label{sol}The following functions are solutions of the generalized hypergeometric equation~\eqref{eqip}:
\begin{itemize}\itemsep=0pt
\item the function \[\vp_1(w):=\frac{1}{2\pi{\rm i}}\int_{\Lambda_1}\frac{\Gamma(s)^5}{\Gamma\big(s+\frac{1}{2}\big)}4^{-s}w^{-4s}\,{\rm d}s,\] defined for $-\frac{\pi}{2}<\arg w<\frac{\pi}{2}$, and where $\Lambda_1$ is any line in the complex plane from the point $\kappa-{\rm i}\infty$ to $\kappa+{\rm i}\infty$ for any $0<\kappa$;
\item the function
\[
\vp_2(w):=
\frac{1}{2\pi{\rm i}}\int_{\Lambda_2}\Gamma(s)^5\Gamma\left(\frac{1}{2}-s\right){\rm e}^{{\rm i}\pi s}4^{-s}w^{-4s}\,{\rm d}s,
\]
 defined for $-\frac{\pi}{2}< \arg w<\pi$, and where $\Lambda_2$ is any line in the complex plane from the point $\kappa-{\rm i}\infty$ to $\kappa+{\rm i}\infty$ for any $0<\kappa<\frac{1}{2}$.
\end{itemize}
\end{Lemma}

For the proof see Appendix~\ref{applemma}.

Note that solutions $\vp_1$ and $\vp_2$ are $\mathbb C$-linearly independent, since their Mellin transforms are. However we have the following identities
\begin{Lemma}By analytic continuation of the functions $\vp_1$ and $\vp_2$, we have
\begin{gather}\label{rel1}
\vp_2\big(w{\rm e}^{{\rm i}\frac{\pi}{2}}\big) =2\pi\vp_1(w)-\vp_2(w),\\
\label{rel2}\vp_2\big(w{\rm e}^{-{\rm i}\frac{\pi}{2}}\big)=2\pi\vp_1\big(w{\rm e}^{-{\rm i}\frac{\pi}{2}}\big)-\vp_2(w),\\
\label{rel3}\vp_2\big(w{\rm e}^{{\rm i}\frac{\pi}{2}}\big)=2\pi\vp_1(w)+\vp_2\big(w{\rm e}^{-{\rm i}\frac{\pi}{2}}\big)-2\pi\vp_1\big(w{\rm e}^{-{\rm i}\frac{\pi}{2}}\big).
\end{gather}
\end{Lemma}
\begin{proof}We have that
\[\Gamma\left(\frac{1}{2}+s\right)\Gamma\left(\frac{1}{2}-s\right)=\frac{\pi} {\sin\big(\pi\big(\frac{1}{2}+s\big)\big)}=\frac{2\pi {\rm e}^{\pm {\rm i}\pi s}}{ {\rm e}^{\pm 2{\rm i}\pi s}+1}
\]for a coherent choice of the sign. So
\[{\rm e}^{\pm 2{\rm i}\pi s}=\frac{2\pi {\rm e}^{\pm {\rm i}\pi s}}{\Gamma\big(\frac{1}{2}+s\big)\Gamma\big(\frac{1}{2}-s\big)}-1.
\]First let us choose the identity with $(-)$: we find that
\begin{align*}\vp_2\big(w{\rm e}^{{\rm i}\frac{\pi}{2}}\big)&=\frac{1}{2\pi{\rm i}}\int_{\Lambda_2}\Gamma(s)^5\Gamma\left(\frac{1}{2}-s\right){\rm e}^{{\rm i}\pi s}\left(\frac{2\pi {\rm e}^{- {\rm i}\pi s}}{\Gamma\big(\frac{1}{2}+s\big)\Gamma\big(\frac{1}{2}-s\big)}-1\right)4^{-s}w^{-4s}\,{\rm d}s\\
&=2\pi\vp_1(w)-\vp_2(w),
\end{align*}which is the first identity. The second one can be deduce analogously using the formula with~$(+)$ sign. Finally the third identity is the difference of~\eqref{rel1} and~\eqref{rel2}.
\end{proof}

The following lemma gives the asymptotic expansions of $\vp_1(w)$ and $\vp_1(w)$ when $w\to \infty$ in certain sectors. Later, in Lemma~\ref{asympt-revised} we will be able to prove the same asymptotics in larger sectors.
\begin{Lemma} \label{asympt} The following asymptotic expansions hold for $w\to \infty$ in the sectors specified below:
\begin{alignat*}{3}
& \vp_1(w)=(2\pi)^{\frac{3}{2}}\frac{{\rm e}^{-4\sqrt{2}w}}{4w^2}\left(1+O\left(\frac{1}{w}\right)\right),
\qquad && \text{for }-\frac{\pi}{2}<\arg w<\frac{\pi}{2},&\\
& \vp_2(w) = \frac{{\rm i}\pi^{\frac{5}{2}}}{2w^2}\left(1+O\left(\frac{1}{w}\right)\right),
\qquad && \text{for } -\frac{\pi}{2}<\arg w<\pi.&
\end{alignat*}
\end{Lemma}

For the proof, see Appendix \ref{applemma}. The sectors where the asymptotics is valid will be enlarged in Lemma \ref{asympt-revised}.

\subsection{Computation of monodromy data}\label{26gennaio2020-2}

\subsubsection[Solution at the origin and computation of $\mathcal C_0(\eta,\mu, R)$]{Solution at the origin and computation of $\boldsymbol{\mathcal C_0(\eta,\mu, R)}$}

Monodromy data at the origin $z=0$ are determined by the action of the first Chern class $c_1(\mathbb G)=4\sigma_1$ on the classical cohomology ring. So,
\begin{gather}\label{28luglio2016-10}
R=\begin{pmatrix}
0& 0 & 0 & 0 & 0 & 0\\
4& 0 & 0 & 0 & 0 & 0\\
0& 4 & 0 & 0 & 0 & 0\\
0& 4 & 0 & 0 & 0 & 0\\
0& 0 & 4 & 4 & 0 & 0\\
0& 0 & 0 & 0 & 4 & 0
\end{pmatrix}.
\end{gather}
By Theorems~\ref{normalth} and~\ref{isomonod1}, there exists a fundamental matrix solution~(\ref{28luglio2016-9})
\[
Y(z)=\Phi\big(t^2,z\big)z^{\mu}z^R,
\]
for some appropriate converging power series $\Phi\big(t^2,z\big)=\mathbbm 1+O(z)$ such that
\[
\Phi^{\rm T}\big(t^2,-z\big) \eta \Phi\big(t^2,z\big)=\eta.
\]
 Thus, a fundamental matrix for our problem is given by
\[\Xi_0(z)=\eta \Phi\big(t^2,z\big) z^{\mu}z^R =\Phi^{\rm T}\big(t^2,-z\big)^{-1}\eta  z^{\mu}z^R.
\]
By applying the iterative procedure in~\cite{dubro2} for the proof of Theorem~\ref{normalth}, at $t^2=0$ one finds the following fundamental solution
\begin{gather}\label{xi0}
\Xi_0(0,z)=S(0,z)\eta z^{\mu}z^R,
\\ S(0,z)=\left(
\begin{matrix}
 2 z^4+1 & 0 & 0 & 0 & 0 & 0 \\
 2 z^3 & 1-4 z^4 & 0 & 0 & 0 & 0 \\
 z^2 & - z^3 & 1 & 0 & 0 & 0 \\
 z^2 & - z^3 & 0 & 1 & 0 & 0 \\
 z & 0 & - z^3 & - z^3 & 4 z^4+1 & 0 \\
 z^4 & z & - z^2 & - z^2 & 2 z^3 & 1-2 z^4 \\
\end{matrix}
\right)+O\big(z^5\big).\nonumber
\end{gather}

Notice that the leading term of the solution $\Xi_0$ in \eqref{xi0} is exactly
\[ \eta z^{\mu}z^R=c\left(
\begin{matrix}
 \frac{64}{3} z^2 \log ^4(z) & \frac{64}{3} z^2 \log ^3(z) & 8 z^2 \log ^2(z) & 8 z^2 \log ^2(z) & 4 z^2 \log (z) & z^2 \vspace{1mm}\\
 \frac{64}{3} z \log ^3(z) & 16 z \log ^2(z) & 4 z \log (z) & 4 z \log (z) & z & 0 \vspace{1mm}\\
 8 \log ^2(z) & 4 \log (z) & 1 & 0 & 0 & 0 \vspace{1mm}\\
 8 \log ^2(z) & 4 \log (z) & 0 & 1 & 0 & 0 \vspace{1mm}\\
 \frac{4 \log (z)}{z} & \frac{1}{z} & 0 & 0 & 0 & 0 \vspace{1mm}\\
 \frac{1}{z^2} & 0 & 0 & 0 & 0 & 0
\end{matrix}
\right).
\]From the first row, we deduce that near $z=0$ any solution of the equation \eqref{eqip}, i.e.,
\[\Theta^5\vp-1024z^4\Theta\vp-2048z^4\vp=0
\] is of the form
\begin{gather}\label{near0}\vp(z)=\sum_{n\geq 0}z^n\big(a_n+b_n\log z+c_n\log^2 z+d_n\log^3 z+e_n\log^4 z\big),
\end{gather}where $a_0$, $b_0$, $c_0$, $d_0$, $e_0$ are arbitrary constants, and successive coefficients can be obtained recursively.

\begin{Proposition}\label{29luglio2016-4}Let $R$ be as in \eqref{28luglio2016-10}. Then, $\mathcal C_0(\eta,\mu, R)$ is the algebraic abelian group of complex dimension $3$ given by
\begin{gather*} \mathcal C_0(\eta,\mu, R)=\left\{
\begin{pmatrix}
1&0&0&0&0&0\\
\alpha_1&1&0&0&0&0\\
\alpha_2&\alpha_1&1&0&0&0\\
\alpha_3&\alpha_1&0&1&0&0\\
\alpha_4&\!\!\alpha_2+\alpha_3\!\!&\!\alpha_1\!&\!\alpha_1\!&1&0\\
\alpha_5&\alpha_4&\alpha_3&\!\alpha_2\!&\!\alpha_1\!&1
\end{pmatrix}\!\colon \alpha_i \in\mathbb C\text{ s.t.}\!\begin{cases}
\!\alpha_1^2-\alpha_2-\alpha_3=0,\\
\!\alpha_2^2+\alpha_3^2-2\alpha_1\alpha_4+2\alpha_5=0
\end{cases}\!\!\!\!\!\!\!
\right\}.
\end{gather*}
In particular, if $F(t)\in\mathbb C[\![t ]\!]$ is a formal power series of the form $F(t)=1+F_1t+F_2t^2+\cdots$, then the matrix $($computed w.r.t.\ the chosen Schubert basis $\sigma_0$, $\sigma_1$, $\sigma_{2}$, $\sigma_{1,1}$, $\sigma_{2,1}$, $\sigma_{2,2})$ representing the endomorphism
\[
\lambda_F\cup (-)\colon \ H^\bullet(\mathbb G;\mathbb C)\to H^\bullet(\mathbb G;\mathbb C),
\]
where $\lambda_F\in H^\bullet(\mathbb G;\mathbb C)$ is such that
\[
\widehat F(T\mathbb G)\cup\lambda_F=\widehat F(T^*\mathbb G),
\]
is an element of $\mathcal C_0(\eta,\mu, R)$. Here $\widehat F(V)$ denotes the Hirzebruch multiplicative characteristic class of the vector bundle $V\to\mathbb G$ associated {with} the formal power series $F(t)$ $($see~{\rm \cite{hirz})}.
\end{Proposition}
\begin{proof}
The equations defining the group $\mathcal C_0(\eta,\mu, R)$ are obtained by direct computation from the requirement that $P(z):=z^\mu z^R\cdot C\cdot z^{-R}z^{-\mu}$ is a polynomial of the form $P(z)=\mathbbm 1+A_1 z+A_2z^2+\cdots$, together with the orthogonality condition $P(-z)^{\rm T}\eta P(z)=\eta$. Notice that the polynomial for the generic matrix of the {above} form is equal to
\[P(z)=\left(
\begin{matrix}
 1 & 0 & 0 & 0 & 0 & 0 \\
 z \alpha_1 & 1 & 0 & 0 & 0 & 0 \\
 z^2 \alpha_2 & z \alpha_1 & 1 & 0 & 0 & 0 \\
 z^2 \alpha_3 & z \alpha_1 & 0 & 1 & 0 & 0 \\
 z^3 \alpha_4 & z^2  (\alpha_2+\alpha_3 ) & z \alpha_1 & z \alpha_1 & 1 & 0 \\
 z^4 \alpha_5 & z^3 \alpha_4 & z^2 \alpha_3 & z^2 \alpha_2 & z \alpha_1 & 1 \\
\end{matrix}
\right).
\]We leave as an exercise to show that such a matrix group is abelian. Let $\delta_1,\dots,\delta_6$ be the Chern roots of $T\mathbb G$. Then, for some complex constants $a_{i,j}\in\mathbb C$, we have
\begin{gather*}
\widehat F(T\mathbb G):=\prod_{j=1}^6F(\delta_j)=1+a_1\sigma_1+a_2\sigma_2+a_{1,1}\sigma_{1,1}+a_{2,1}\sigma_{2,1}+a_{2,2}\sigma_{2,2},
\\
\widehat F(T^*\mathbb G):=\prod_{j=1}^6F(-\delta_j)=1-a_1\sigma_1+a_2\sigma_2+a_{1,1} \sigma_{1,1}-a_{2,1}\sigma_{2,1}+a_{2,2}\sigma_{2,2}.
\end{gather*}
Thus, if
\[\lambda_F=1+x_1\sigma_1+x_2\sigma_2+x_3\sigma_{1,1}+x_4\sigma_{2,1}+x_5\sigma_{2,2},
\]from the condition $\widehat F(T\mathbb G)\cup\lambda_F=\widehat F(T^*\mathbb G)$ we obtain the constraints
\[\begin{cases}
x_1=-2a_1,\\
x_2=2a_1^2,\\
x_3=2a_1^2,\\
x_4=2a_1(a_2+a_{1,1})-4a_1^3-2a_{2,1},\\
x_5=4a_1a_{2,1}-4a_1^2(a_2+a_{1,1})+4a_1^4.
\end{cases}
\]From this it is immediately seen that $x_1^2-x_2-x_3=0$ and $x_2^2+x_3^2-2x_1x_4+2x_5=0$.
\end{proof}

\subsubsection[Stokes rays and computation of $\Xi_{\rm left}$, $\Xi_{\rm right}$]{Stokes rays and computation of $\boldsymbol{\Xi_{\rm left}}$, $\boldsymbol{\Xi_{\rm right}}$}\label{30luglio2016-4}

According to Theorem \ref{mainisoth}, monodromy data of $QH^{\bullet}(\mathbb G)$ can be computed starting from a point $\big(0,t^2,0,\dots ,0\big)$ of the small quantum cohomology. Moreover, thanks to the isomonodromy theorems, it suffices to do the computation at $t^2=0$, i.e., $q=1$, where the canonical coordinates~(\ref{28luglio2016-12}) are
\[
u_1=u_2=0,\qquad u_3=-4{\rm i}\sqrt{2},\qquad u_4=4{\rm i}\sqrt{2},\qquad u_5=-4\sqrt2,\qquad u_6=4\sqrt2.
\]
The Stokes rays (\ref{28luglio2016-13}) are {easily} seen to be
\begin{gather*}R_{13}=R_{23} = \{-\rho\colon\rho\geq0 \},\qquad
R_{14}=R_{24}=R_{34} =\{\rho\colon\rho\geq0 \},\\
R_{15}=R_{25} = \{-{\rm i}\rho\colon\rho\geq 0 \},\qquad
R_{16}=R_{26}=R_{56} = \{{\rm i}\rho\colon\rho\geq 0 \},
\\
R_{35}= \big\{\rho {\rm e}^{-{\rm i}\frac{\pi}{4}}\colon \rho\geq 0 \big\},\qquad
R_{36} = \big\{\rho {\rm e}^{{\rm i}\frac{\pi}{4}}\colon \rho\geq 0 \big\},\\
R_{45}= \big\{{-}\rho {\rm e}^{{\rm i}\frac{\pi}{4}}\colon \rho\geq 0 \big\},\qquad
R_{46} = \big\{{-}\rho {\rm e}^{-{\rm i}\frac{\pi}{4}}\colon \rho\geq 0 \big\},
\qquad R_{ji}=-R_{ij}.
\end{gather*}

We fix the admissible line $\ell$
\[\ell:=\big\{\rho {\rm e}^{{\rm i}\frac{\pi}{6}}\colon\rho\in\mathbb R\big\},
\]
so that the sectors for the asymptotic expansion, containing $\Pi_{\rm left/right}$ and extending up to the nearest Stokes rays are \begin{gather*}
 \mathcal{S}_{\rm right}=\{z\colon -\pi<\arg z<\pi/4\},
 \qquad \mathcal{S}_{\rm left}=\{z\colon -0<\arg z<\pi+\pi/4\}.
\end{gather*}
For such a choice of the line, according to Theorem~\ref{proprstok}, the structure of the Stokes matrix is
\begin{gather}\label{formstok}
S=\begin{pmatrix}
1&0&*&0&0&*\\
0&1&*&0&0&*\\
0&0&1&0&0&*\\
*&*&*&1&0&*\\
*&*&*&*&1&*\\
0&0&0&0&0&1
\end{pmatrix}.
\end{gather}
We use the following notation for fundamental matrices
\begin{align*}
\Xi_{\rm right}=\begin{pmatrix}
\xi_{(1),1}^R&\xi_{(2),1}^R&\xi_{(3),1}^R&\xi_{(4),1}^R&\xi_{(5),1}^R&\xi_{(6),1}^R\vspace{1mm}\\
\xi_{(1),2}^R&\xi_{(2),2}^R&\xi_{(3),2}^R&\xi_{(4),2}^R&\xi_{(5),1}^R&\xi_{(6),2}^R \\
\vdots&\vdots&\vdots&\vdots&\vdots&\vdots   \\
\xi_{(1),6}^R&\xi_{(2),6}^R&\xi_{(3),6}^R&\xi_{(4),6}^R&\xi_{(5),6}^R&\xi_{(6),6}^R
\end{pmatrix},
\\
 \Xi_{\rm left}=\begin{pmatrix}
\xi_{(1),1}^L&\xi_{(2),1}^L&\xi_{(3),1}^L&\xi_{(4),1}^L&\xi_{(5),1}^L&\xi_{(6),1}^L\vspace{1mm}\\
\xi_{(1),2}^L&\xi_{(2),2}^L&\xi_{(3),2}^L&\xi_{(4),2}^L&\xi_{(5),1}^L&\xi_{(6),2}^L\\
\vdots&\vdots&\vdots&\vdots&\vdots&\vdots \\
\xi_{(1),6}^L&\xi_{(2),6}^L&\xi_{(3),6}^L&\xi_{(4),6}^L&\xi_{(5),6}^L&\xi_{(6),6}^L
\end{pmatrix}.
\end{align*}
Note, in particular, that~(\ref{formstok}) implies that the fifth columns of $\Xi_{\rm right}$ and $\Xi_{\rm left}$ coincide. Then~$\xi_{(5),1}^{L}$ is the analytical continuation of $\xi_{(5),1}^{R}$ on $\mathcal{S}_{\rm left}$. Moreover, the exponential ${\rm e}^{zu_5}$ dominates all others ${\rm e}^{zu_j}$'s in the sector between the rays~$R_{45}$ and~$R_{46}$, i.e., for $-\pi -\pi/4<\arg z<-\pi+\pi/4$. This implies that the asymptotics
\[
\xi_{(5),1}^L=\xi_{(5),1}^R= \frac{c^{\frac{1}{2}}}{2\sqrt{2}}{\rm e}^{zu_5}\left(1+O\left(\frac{1}{z}\right)\right)
\]
 is valid in the whole sector $-\pi-\pi/4<\arg z<\pi+\pi/4$. By Lemma~\ref{asympt},
\[
\frac{c^{\frac{1}{2}}}{2\pi^{\frac{3}{2}}}z^2\vp_1(z) =
\frac{c^{\frac{1}{2}}}{2\sqrt{2}}{\rm e}^{zu_5}\left(1+O\left(\frac{1}{z}\right)\right),\qquad \text{for }-\frac{\pi}{2}<\arg z<\frac{\pi}{2}.
\]
Since the exponential ${\rm e}^{zu_5}$ is dominated by all others exponentials ${\rm e}^{zu_j}$ in the region between~$R_{35}$ and~$R_{36}$, namely for $-\pi/4<\arg z <\pi/4$, we conclude necessarily that
\[
\frac{c^{\frac{1}{2}}}{2\pi^{\frac{3}{2}}}z^2\vp_1(z)=\xi_{(5),1}^{L/R}(z).
\]
This determines the 5-{th} column of $\Xi_{\rm right}$ and $\Xi_{\rm left}$ in terms of $\vp_1$, using equations~\eqref{reconstruction1},~\eqref{reconstruction4}. We also obtain an improvement of Lemma~\ref{asympt}:

\begin{Lemma}\label{asympt-revised} $\vp_1$ and $\vp_2$ have the following asymptotic behaviour for $w\to\infty$ in the sectors below
\begin{alignat*}{3}
 & \vp_1(w)=(2\pi)^{\frac{3}{2}}\frac{{\rm e}^{-4\sqrt{2}w}}{4w^2}\left(1+O\left(\frac{1}{w}\right)\right),
 \qquad& &\text{for } -\pi-\frac{\pi}{4}<\arg w<\pi+\frac{\pi}{4}, &
\\
& \vp_2(w) = \frac{{\rm i}\pi^{\frac{5}{2}}}{2w^2}\left(1+O\left(\frac{1}{w}\right)\right),
\qquad&&\text{for }-\frac{\pi}{2}<\arg w<\pi.&
\end{alignat*}
\end{Lemma}

We are ready to determine the other columns of $\Xi_{\rm left/right}$. By Lemma \ref{asympt-revised},
\begin{alignat}{3}
\label{28luglio2016-14-1}
& -\frac{c^{\frac{1}{2}}}{2\pi^{\frac{3}{2}}}z^2\vp_1\big(z{\rm e}^{{\rm i}\frac{\pi}{2}}\big )=\frac{c^{\frac{1}{2}}}{2\sqrt{2}}{\rm e}^{zu_3}\left(1+O\left(\frac{1}{z}\right)\right),
\quad && \text{for }-2\pi+\frac{\pi}{4}<\arg z<\frac{3\pi}{4},& \\
\label{28luglio2016-14-2}
& \frac{c^{\frac{1}{2}}}{2\pi^{\frac{3}{2}}}z^2\vp_1\big(z{\rm e}^{{\rm i}\pi}\big)=\frac{c^{\frac{1}{2}}}{2\sqrt{2}}{\rm e}^{zu_6}\left(1+O\left(\frac{1}{z}\right)\right),
\quad && \text{for }-2\pi-\frac{\pi}{4}<\arg z<\frac{\pi}{4}.&
\end{alignat}
We consider first (\ref{28luglio2016-14-1}). Being solutions of a differential equation, the following holds:
\[
-\frac{c^{\frac{1}{2}}}{2\pi^{\frac{3}{2}}}z^2\vp_1\big(z{\rm e}^{{\rm i}\frac{\pi}{2}}\big)=\text{linear combination of the } \xi_{(1),i}^R,\qquad 1\leq i \leq 6.
 \]
 On the other hand, ${\rm e}^{zu_3}$ is dominated by all other ${\rm e}^{zu_i}$'s in the sector $-\pi+\pi/4<\arg z<-\pi/2$ between~$R_{45}$ and~$R_{35}$. This requires that the linear combination necessarily reduces to
\[
-\frac{c^{\frac{1}{2}}}{2\pi^{\frac{3}{2}}}z^2\vp_1\big(z{\rm e}^{{\rm i}\frac{\pi}{2}}\big)=\xi_{(3),1}^{R}.
\]
Now we consider~(\ref{28luglio2016-14-2}). As above, since ${\rm e}^{zu_6}$ is dominated by all the other ${\rm e}^{zu_i}$'s in the sector $-5\pi/4<\arg z  < -3\pi/4$ between $R_{46}$ and $R_{45}$, we conclude that
\[\frac{c^{\frac{1}{2}}}{2\pi^{\frac{3}{2}}}z^2\vp_1\big(z{\rm e}^{{\rm i}\pi}\big)=\xi_{(6),1}^{R}.
\]
Analogously we find that
\begin{alignat*}{3}
& -\frac{c^{\frac{1}{2}}}{2\pi^{\frac{3}{2}}}z^2\vp_1\big(z{\rm e}^{-{\rm i}\frac{\pi}{2}}\big)=\frac{c^{\frac{1}{2}}}{2\sqrt{2}}{\rm e}^{zu_4}\left(1+O\left(\frac{1}{z}\right)\right)\qquad && \text{for }-\frac{3\pi}{4}<\arg z<\pi+\frac{3\pi}{4},&\\
& \frac{c^{\frac{1}{2}}}{2\pi^{\frac{3}{2}}}z^2\vp_1\big(z{\rm e}^{-{\rm i}\pi}\big)=\frac{c^{\frac{1}{2}}}{2\sqrt{2}}{\rm e}^{zu_6}\left(1+O\left(\frac{1}{z}\right)\right)\qquad && \text{on }-\frac{\pi}{4}<\arg z<2\pi+\frac{\pi}{4}.&
\end{alignat*}
By dominance considerations as above, we conclude that
\[\xi_{(4),1}^{L}=-\frac{c^{\frac{1}{2}}}{2\pi^{\frac{3}{2}}}z^2\vp_1\big(z{\rm e}^{-{\rm i}\frac{\pi}{2}}\big),\qquad \xi_{(6),1}^{L}=\frac{c^{\frac{1}{2}}}{2\pi^{\frac{3}{2}}}z^2\vp_1\big(z{\rm e}^{-{\rm i}\pi}\big).
\]
The above results reconstruct (using identities \eqref{reconstruction1},\eqref{reconstruction4}) three columns of matrices $\Xi_{\rm right}$ and $\Xi_{\rm left}$ respectively. As far as the first two columns are concerned, we invoke again Lemma~\ref{asympt-revised} for $\vp_2$, which yields
\begin{alignat*}{3}
& \frac{c^{\frac{1}{2}}}{\pi^{\frac{5}{2}}}z^2\vp_2\big(z{\rm e}^{{\rm i}\frac{\pi}{2}}\big)=-\frac{{\rm i} c^{\frac{1}{2}}}{2}\left(1+O\left(\frac{1}{z}\right)\right)\qquad && \text{on }-\pi<\arg z<\frac{\pi}{2},&\\
& \frac{c^{\frac{1}{2}}}{\pi^{\frac{5}{2}}}z^2\vp_2\big(z{\rm e}^{-{\rm i}\frac{\pi}{2}}\big)=-\frac{{\rm i} c^{\frac{1}{2}}}{2}\left(1+O\left(\frac{1}{z}\right)\right)\qquad && \text{on }0<\arg z<\frac{3\pi}{2}.&
\end{alignat*}
Exactly as before, dominance relations of the exponentials ${\rm e}^{zu_i}$ yield
\[\frac{c^{\frac{1}{2}}}{\pi^{\frac{5}{2}}}z^2\vp_2\big(z{\rm e}^{{\rm i}\frac{\pi}{2}}\big)=\xi_{(1),1}^R=\xi_{(2),1}^R,\qquad
\frac{c^{\frac{1}{2}}}{\pi^{\frac{5}{2}}}z^2\vp_2\big(z{\rm e}^{-{\rm i}\frac{\pi}{2}}\big)=\xi_{(1),1}^L=\xi_{(2),1}^L.
\]
Using \eqref{reconstruction1}, \eqref{reconstruction2}, \eqref{reconstruction3}, the first two columns are constructed. Summarizing, we have determined the following columns in terms of $\vp_1$ and $\vp_2$.
\begin{gather*}
\Xi_{\rm right}=\begin{pmatrix}
\xi_{(1),1}^R&\xi_{(2),1}^R&\xi_{(3),1}^{R}&\text{\rm unknown}&\xi_{(5)}^{R}&\xi_{(6),1}^{R}\\
\vdots&\vdots&\vdots&\vdots&\vdots&\vdots
\end{pmatrix},\\
\Xi_{\rm left}=\begin{pmatrix}
\xi_{(1),1}^L&\xi_{(2),1}^L&\text{\rm unknown}&\xi_{(4),1}^{L}&\xi_{(5),1}^{L}&\xi_{(6),1}^{L}\\
\vdots&\vdots&\vdots&\vdots&\vdots&\vdots
\end{pmatrix}.
\end{gather*}

In Section \ref{28luglio2016-15} we show that the above partial information and the constraint (2) in Theorem~\ref{constraint} are {sufficient} to determine the Stokes and central connection matrices simultaneously. Since constraint~(2) holds only in case~$S$ and~$C$ are related to Frobenius manifolds, we sketch below~-- for the sake of completeness~-- the general method to obtain the missing columns of~$\Xi_{\rm left/right}$ and~$S$, in a pure context of asymptotic analysis of differential equations.

 We observe that
\[
-\frac{c^{\frac{1}{2}}}{2\pi^{\frac{3}{2}}}z^2\vp_1\big(z{\rm e}^{-{\rm i}\frac{\pi}{2}}\big)=\frac{c^{\frac{1}{2}}}{2\sqrt{2}}{\rm e}^{u_4 z}(1+O(1/z)),\qquad\text{for }
-\pi +\frac{\pi}{4}<\arg z<\pi+\frac{3\pi}{4}.
\]
The sub-sector $-\pi<\arg z<-3\pi/4$ of $\mathcal{S}_{\rm right}$ is not covered by the sector where the above asymptotic behaviour holds. On the sub-sector, the dominance relation $\big|{\rm e}^{zu_4}|<|{\rm e}^{zu_5}\big|$ holds. Thus,
\begin{gather}\label{28luglio2016-18-1}
\xi_{(4),1}^R=-\frac{c^{\frac{1}{2}}}{2\pi^{\frac{3}{2}}}z^2\vp_1\big(z{\rm e}^{-{\rm i}\frac{\pi}{2}}\big)+v \xi_{(5),1}^R,
\end{gather}
for some complex number $v\in\mathbb C$, to be determined. Analogously, we observe that
\[
-\frac{c^{\frac{1}{2}}}{2\pi^{\frac{3}{2}}}z^2\vp_1\big(z{\rm e}^{-{\rm i}\frac{3\pi}{2}}\big)=\frac{c^{\frac{1}{2}}}{2\sqrt{2}}{\rm e}^{u_4 z}(1+O(1/z)),\qquad\text{for }
-2\pi -\frac{3\pi}{4}<\arg z<-\frac{\pi}{4}.
\]
The sub-sector $-\pi/4<\arg z<\pi/4$ of $\mathcal{S}_{\rm right}$ is not covered by the sector where the asymptotic behaviour holds. Now, the following dominance relations hold: $\big|{\rm e}^{zu_4}\big|<\big|{\rm e}^{zu_i}\big|$, for $i=1,2,3,6$, in $0<\arg z<\pi/4$; for $i=6$ in $-\pi/4<\arg z<0$. Thus
\begin{gather}\label{28luglio2016-18-2}
\xi_{(4),1}^R=-\frac{c^{\frac{1}{2}}}{2\pi^{\frac{3}{2}}}z^2\vp_1\big(z{\rm e}^{{\rm i}\frac{3\pi}{2}}\big)+\gamma_1\xi_{(1),1}^R+\gamma_3\xi_{(3),1}^R+\gamma_6\xi_{(6),1}^R
\end{gather}
for some complex number $\gamma_1,\gamma_3,\gamma_6\in\mathbb C$, to be determined.\footnote{There is no need to include a term $+\gamma_2\xi_{(2),1}^R$ in the linear combination, since $\xi_{(1),1}^R =\xi_{(2),1}^R$.} The above (\ref{28luglio2016-18-1}) and (\ref{28luglio2016-18-2}) become a 6-terms linear relation between functions $\vp_2\big(z{\rm e}^{{\rm i}\frac{k\pi}{2}}\big)$, as follows
 \begin{gather*}
-\vp_1\big(z{\rm e}^{-{\rm i}\frac{\pi}{2}}\big)+v \vp_1(z)=-\vp_1\big(z{\rm e}^{{\rm i}\frac{3\pi}{2}}\big)+\frac{\gamma_1}{\pi}\vp_2\big(z{\rm e}^{{\rm i}\frac{\pi}{2}}\big)-\gamma_3
\vp_1(z{\rm e}^{{\rm i}\frac{\pi}{2}})+\gamma_6\vp_1\big(z{\rm e}^{{\rm i}\pi}\big),\\
\vp_1(z)=\frac{1}{2\pi}\big[\vp_2(z)+\vp_2\big(z{\rm e}^{{\rm i}\frac{\pi}{2}}\big)\big].
\end{gather*}
Some further information is {needed} in order to determine the unknown constants $v$, $\gamma_1$, $\gamma_3$, $\gamma_6$, as in the following

\begin{Lemma}The solutions of the equation \eqref{eqip} satisfy the {identity}
\begin{gather}\label{28luglio2016-19}
\vp\big(z{\rm e}^{{\rm i}\frac{5\pi}{2}}\big)-5\vp\big(z{\rm e}^{2\pi{\rm i}}\big)+10\vp\big(z{\rm e}^{{\rm i}\frac{3\pi}{2}}\big)-10\vp\big(z{\rm e}^{{\rm i}\pi}\big)+5\vp\big(z{\rm e}^{{\rm i}\frac{\pi}{2}}\big)-\vp(z)=0.
\end{gather}
\end{Lemma}

\begin{proof} The equation $\Theta^5\vp-1024z^4\Theta\vp-2048z^4\vp=0$ admits the symmetry $z\mapsto z{\rm e}^{{\rm i}\frac{\pi}{2}}$. This means that if $\vp$ is a solution of the equation then also $\vp\big(z{\rm e}^{{\rm i}\frac{\pi}{2}}\big)$ is. Such a symmetry defines a~linear map on the vector space of solutions of the equation defined in a neighborhood of $z=0$. Because of this symmetry, the form~\eqref{near0} can be refined as
\begin{gather}\label{near01}\vp(z)=\sum_{n\geq 0}z^{4n}\big(a_n+b_n\log z+c_n\log^2 z+d_n\log^3 z+e_n\log^4 z\big),
\end{gather}where $a_0$, $b_0$, $c_0$, $d_0$, $e_0$ are arbitrary constants, and successive coefficients can be obtained recursively. In the basis of solutions of the form~\eqref{near01} with $(a_0,b_0,c_0,d_0,e_0)=(1,0,\dots ,0)$, $(0,1,0,\dots ,0)$ and so on, the matrix of the operator
\[ (A\vp )(z):=\vp\big(z{\rm e}^{{\rm i}\frac{\pi}{2}}\big)
\]is of triangular form with 1's on the diagonal. Hence, by Cayley--Hamilton theorem we deduce that
$ (A-\mathbbm 1)^5=0$, namely
\begin{gather*} A^5-5A^4+10A^3-10A^2+5A-\mathbbm 1=0.\tag*{\qed}
\end{gather*}\renewcommand{\qed}{}
\end{proof}

The relation (\ref{28luglio2016-19}) applied to $\vp_2$ determines $v$, $\gamma_1$, $\gamma_3$, $\gamma_6$. For example, $v=6$. This determines $\xi_{(4),1}^R$ through formula~(\ref{28luglio2016-18-1}). The fourth column of $\Xi_{\rm right}$ is then constructed with formula~(\ref{reconstruction1}) applied to $\xi_{(4),1}^R$ (with $h=0$). The value $v=6$ will be determined again in Section~\ref{28luglio2016-15} making use of the constraint~(2) of Theorem~\ref{constraint}.

Proceeding in the same way, we also determine $\xi_{(3),1}^L$. One observes that
\begin{alignat*}{3}
& -\frac{c^{\frac{1}{2}}}{2\pi^{\frac{3}{2}}}z^2\vp_1\big(z{\rm e}^{{\rm i}\frac{3\pi}{2}}\big)=\frac{c^{\frac{1}{2}}}{2\sqrt{2}}{\rm e}^{zu_3}(1+O(1/z)),\qquad && \text{for }
\frac{\pi}{4}<\arg z<\frac{3\pi}{2}+2\pi,& \\
& -\frac{c^{\frac{1}{2}}}{2\pi^{\frac{3}{2}}}z^2\vp_1\big(z{\rm e}^{{\rm i}\frac{\pi}{2}}\big)=\frac{c^{\frac{1}{2}}}{2\sqrt{2}}{\rm e}^{zu_3}(1+O(1/z)),\qquad && \text{for }
-2\pi-\frac{\pi}{2}<\arg z<\frac{3\pi}{4}.&
\end{alignat*}
The first asymptotic relation does not hold in the sub-sector $-\pi/4<\arg z<\pi/4$ of $\mathcal{S}_{\rm left}$. The second {one} does not {hold} in $3\pi/4<\arg z<5\pi/4$. Then, the dominance relations in these sub-sectors generate {a} 6-terms linear relation with unknown coefficients. The coefficients are determined by~(\ref{28luglio2016-19}).

Once $\Xi_{\rm left/right}$ has been determined, $S$ can be computed by direct comparison of the two fundamental matrices (formula~(\ref{28luglio2016-19}) need to be used at some point of the comparison). The final result is the Stokes matrix $S$ of formula~(\ref{29luglio2016-1}) below with $v=6$.

\begin{figure}[t]\centering
\def\svgscale{0.5}
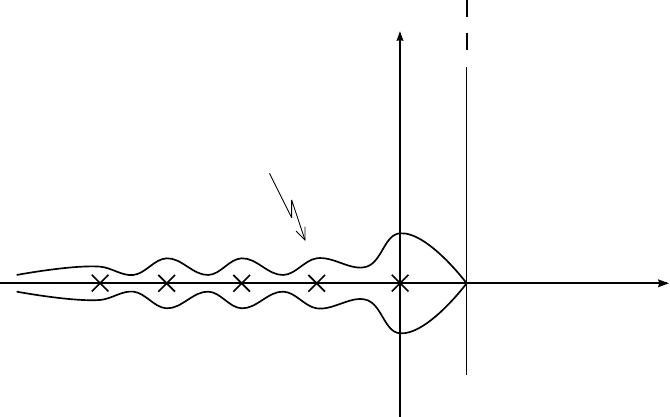
\caption{Deformation of the path $\Lambda_{1/2}$, in order to apply residue theorem. Poles are represented.} \label{10agosto2016-1}
\end{figure}

\subsubsection[Computation of Stokes and central connection matrices, \\ using {constraint} (2) of Theorem \ref{constraint}]{Computation of Stokes and central connection matrices, \\ using {constraint} (2) of Theorem \ref{constraint}}\label{28luglio2016-15}
 We start from formula (\ref{28luglio2016-18-1}):
\[
\xi_{(4),1}^R=-\frac{c^{\frac{1}{2}}}{2\pi^{\frac{3}{2}}}z^2\vp_1\big(z{\rm e}^{-{\rm i}\frac{\pi}{2}}\big)+v \xi_{(5),1}^R\equiv \frac{c^{\frac{1}{2}}}{2\pi^{\frac{3}{2}}}z^2\big({-}\vp_1\big(z{\rm e}^{-{\rm i}\frac{\pi}{2}}\big)+v\vp_1(z)\big).
\]
Though $v$ has already been determined above, we show that constraint (2) of Theorem~\ref{constraint} suffices to determine  the value of~$v$ and reconstruct both the Stokes and the central connection matrices, as follows.

The definition of the central connection matrix $C$ and the transformation \eqref{transformation} imply that
\[\Xi_{\rm right}=\Xi_0C.
\]
The {matrix} $C$ can be obtained by comparing the leading behaviours of $\Xi_{\rm right}$ and $\Xi_0$ near $z=0$.
The leading behaviour of $\Xi_0$ in \eqref{xi0} is $\eta z^\mu z^R$. In order to find the behaviour of $\Xi_{\rm right}$, we need to compute the behaviour of $\vp_1$ and $\vp_2$ near $z=0$. To this end, we consider the integral representations in Lemma~\ref{sol}, and deform both paths~$\Lambda_1$ and $\Lambda_2$ to the left, as shown in Fig.~\ref{10agosto2016-1}. By residue theorem, we obtain a~ representations of $\vp_1$ and $\vp_2$ as a series of residues at the poles $s=0,-1,-2,\dots$. Then, by the  reconstruction dictated by equations \eqref{reconstruction1}, \eqref{reconstruction2}, \eqref{reconstruction3}, \eqref{reconstruction4}, for each entry of the matrix $\Xi_{\rm right}$ we obtain  an expansion in $z$ and $\log z$, converging for small $|z|$.

For example, let us compute the first and second columns of the matrix $C$: by deformation of the path $\Lambda_2$ we obtain that for small $z$ the following series expansions hold:
\begin{align*} \xi_{(1),1}^R=\xi_{(2),1}^R&=\frac{c^{\frac{1}{2}}}{\pi^{\frac{5}{2}}}z^2\vp_2\big(z{\rm e}^{{\rm i}\frac{\pi}{2}}\big) =\frac{c^{\frac{1}{2}}}{\pi^{\frac{5}{2}}}z^2\sum_{n=0}^{\infty}\underset{s=-n}{\text{res}} \left(\Gamma(s)^5\Gamma\left(\frac{1}{2}-s\right){\rm e}^{-{\rm i}\pi s}4^{-s}z^{-4s}\right)\\
&=\alpha_1z^2\log^4z+\alpha_2z^2\log^3z+\alpha_3z^2\log^2 z+\alpha_4z^2\log z+\alpha_5 z^2+ O\big(z^4\big),
\end{align*}
where $\alpha_i$ can be explicitly computed. By comparison with the first row of $\eta z^\mu z^R$ we determine the entries
\begin{alignat*}{3}
& C_{11}=C_{12} =\frac{3}{64c}\alpha_1, \qquad&& C_{21}=C_{22}=\frac{3}{64c}\alpha_2,& \\
& C_{51}=C_{52}=\frac{1}{4c}\alpha_4, \qquad  && C_{61}=C_{62}=\frac{1}{c}\alpha_5.&
\end{alignat*}
For the other entries we have to consider expansions of $\xi^R_{(1),3}$, $\xi^R_{(2),3}$, $\xi^R_{(1),4}$, $\xi^R_{(2),4}$. For example,
\begin{align*}
\xi^R_{(1),3}=\xi^R_{(2),4}&=\frac{c^{\frac{1}{2}}}{\pi^{\frac{5}{2}}}\cdot\frac{1}{32}\big(z\vp'_2\big(z{\rm e}^{{\rm i}\frac{\pi}{2}}\big)+z^2\vp''_2\big(z{\rm e}^{{\rm i}\frac{\pi}{2}}\big)\big)-\frac{c^{\frac{1}{2}}}{2}\\
&=-\frac{c^{\frac{1}{2}}}{2}+\frac{c^{\frac{1}{2}}}{2\pi^{\frac{5}{2}}} \sum_{n=0}^{\infty}\underset{s=-n}{\text{res}}\left(\Gamma(s)^5\Gamma\left(\frac{1}{2}-s\right){\rm e}^{-{\rm i}\pi s}4^{-s}s^2z^{-4s}\right)\\
&=\beta_1\log^2 z+\beta_2\log z+\beta_3+O\big(z^4\big),
\end{align*}where $\beta_i$ can be explicitly computed. So, by comparison of the third row of $gz^\mu z^R$ we obtain
\[C_{31}=C_{42}=\frac{\beta_3}{c}.
\]Analogously one obtains $C_{32}=C_{41}$. Note that the other entries $C_{ij}$, with $j=3,4,5,6$, are uniquely determined only by the expansion of $\xi^R_{(j),i}$ because of \eqref{reconstruction4}.
The compuation for all the other entries of $C$ can be done in the same way, so it will not be repeated here. Due to the length of the result, we write the whole $C$ in Appendix \ref{appc}. As it can be seen in Appendix \ref{appc}, only the fifth column of $C$ is expressed in terms of the constant $v$. This $v$ will now be determined.

Since $S$ and $C$ are associated {with} a Frobenius manifold, the constraint~(2) of Theorem~\ref{constraint} holds:
\begin{gather}\label{29luglio2016-2}
S=C^{-1}{\rm e}^{-\pi{\rm i}R}{\rm e}^{-\pi{\rm i}\mu}\eta^{-1}\big(C^{\rm T}\big)^{-1}.
\end{gather}
Substituting $C$ of Appendix~\ref{appc} with {an indeterminate} $v$ in the {above} constraint, we obtain the Stokes matrix
\begin{gather}
\label{29luglio2016-1}
S=\left(\begin{matrix}
 1 & 0 & 4 & 0 & 0 & 4 \\
 0 & 1 & 4 & 0 & 0 & 4 \\
 0 & 0 & 1 & 0 & 0 & 6 \\
 -4 & -4 & -16 & 1 & 6-v & -6 \\
 4 (v -1) & 4 (v -1) & 16 v -26 & -v & (v -6) v +1 & 6 v -16 \\
 0 & 0 & 0 & 0 & 0 & 1 \end{matrix}\right).
\end{gather}
By a direct comparison with the expected matrix form~\eqref{formstok}, which dictates that~$S_{45}=0$ and $S_{55}=1$, we conclude that necessarily
\begin{gather*}v=6.
\end{gather*}
In this way we have completely determined both the Stokes and central connection matrices as well as the fundamental matrix $\Xi_{\rm right}$. See also~(\ref{29luglio2016-8-1}) below.

\subsection[Monodromy data of the small quantum cohomology and exceptional collections\\ in $\mathcal D^b(\mathbb G)$]{Monodromy data of the small quantum cohomology\\ and exceptional collections in $\boldsymbol{\mathcal D^b(\mathbb G)}$} \label{30luglio2016-7}

The monodromy data $R$ and $C$ computed above can be read as characteristic classes of objects of an exceptional collection in $\mathcal D^b(\mathbb G)$, as it has been conjectured by one of the authors~\cite{dubro0}, though the formulation for the central connection matrix was not well understood {then}. Following~\cite{kkp}
 where the role of the $\widehat \Gamma^\pm$-classes (characteristic classes obtained by the Hirzebruch's procedure starting from the series expansion of the functions $\Gamma(1\pm t)$ near $t=0$) was pointed out, we claim that the central connection matrix (for canonical coordinates in triangular/lexicographical order) can be identified with the matrix of the $\mathbb C$-linear morphism
\begin{gather*}\mathfrak X_{\mathbb G}^\pm\colon \  K_0(\mathbb G) \otimes_{\mathbb Z}\mathbb C
\to H^\bullet(\mathbb G;\mathbb C), \\
\hphantom{\mathfrak X_{\mathbb G}^\pm\colon \  K_0(\mathbb G) \otimes_{\mathbb Z}}{} E \mapsto \frac{1}{(2\pi)^2c^{\frac{1}{2}}}\widehat \Gamma^\pm(\mathbb G)\cup\operatorname{Ch}(E),
\\
\widehat\Gamma^\pm(\mathbb G):=\prod_j \Gamma(1\pm\delta_j), \qquad\text{where }\delta_j\text{'s are the Chern roots of }T\mathbb G,\\
\operatorname{Ch}(V):=\sum_k {\rm e}^{2\pi{\rm i} x_k},\qquad\text{ $x_k$'s are the Chern roots of a vector bundle $V$},
\end{gather*}
 expressed w.r.t.
\begin{itemize}\itemsep=0pt
\item an exceptional basis $(\varepsilon_i)_i$ of $K_0(\mathbb G)\otimes_{\mathbb Z}\mathbb C$, i.e., satisfying $\chi(\varepsilon_i,\varepsilon_i)=1$, and the Grothendieck--Euler--Poincar\'e orthogonality conditions $\chi(\varepsilon_i,\varepsilon_j)=0$ for $i>j$, obtained by projection of a full exceptional collection $(E_i)_i$ in $\mathcal D^b(\mathbb G)$;
\item a basis in $H^\bullet (\mathbb G;\mathbb C)$ related to $(\sigma_0,\sigma_1,\sigma_{2},\sigma_{1,1},\sigma_{2,1},\sigma_{2,2})$ (the Schubert basis we have fixed) by a $(\eta,\mu)$-orthogonal-parabolic $G$ endomorphism (as described in Section \ref{3agosto2016-7}) which commutes with the operator of classical $\cup$-multiplication $c_1(\mathbb G)\cup-\colon H^\bullet (\mathbb G;\mathbb C)\to H^\bullet (\mathbb G;\mathbb C)$.
\end{itemize}
By application of the constraint (\ref{29luglio2016-2}) and the Grothendieck--Hirzebruch--Riemann--Roch theorem, one can prove that the Stokes matrix (in triangular/lexicographical order) is equal to the inverse of the Gram matrix:
\[\big(S^{-1}\big)_{ij}=\chi(\varepsilon_i,\varepsilon_j).
\]
See \cite{CDG1} for a proof.

\begin{Remark}As it was formulated in Theorem \ref{6-07-17-2} in the Introduction and in Section~\ref{freedom}, some natural transformations are allowed, such as
\begin{itemize}\itemsep=0pt
\item the left action of the group $\mathcal C_0(\eta,\mu, R)$:
\begin{gather}\label{29luglio2016-5-1}
\text{no action on }S,\qquad C\longmapsto G C,
\end{gather}
where $G\in \mathcal C_0(\eta,\mu, R)$ and has the form prescribed by Proposition~\ref{29luglio2016-4};
\item the right action of the group $(\mathbb Z/2\mathbb Z)^{\times 6}$:
 \begin{gather}\label{29luglio2016-5-2}
S\longmapsto \mathcal{I} S \mathcal{I},\qquad C\longmapsto C \mathcal I,
\end{gather}
where $\mathcal I$ is a diagonal matrix of $1$'s and $-1$'s;
\item the right action of the braid group $\mathcal B_6$:
 \begin{gather}\label{29luglio2016-5-3}
S\longmapsto A^\beta S \big(A^\beta\big)^{\rm T}, \qquad C\longmapsto C \big(A^\beta\big)^{-1},
\end{gather}
 as in formulae (\ref{connbraid2-29luglio}) and (\ref{connbraid2}).
\end{itemize}
The above actions naturally manifest respectively on the space $H^\bullet (\mathbb G;\mathbb C)$, on the set of full exceptional collections in the category $\mathcal D^b(\mathbb G)$, and/or on the set of exceptional bases of the complexified Grothendieck group $K_0(\mathbb G)\otimes_{\mathbb Z}\mathbb C$. More precisely,
\begin{itemize}\itemsep=0pt
\item $\mathcal C_0(\eta,\mu, R)$ acts on $H^\bullet (\mathbb G;\mathbb C)$ as $(\eta,\mu)$-orthogonal-parabolic endomorphisms commuting with the classical $\cup$-product by the first Chern class $c_1(\mathbb G)$;
\item the action of the shift functor $[1]\colon\mathcal D^b(\mathbb G)\to\mathcal D^b(\mathbb G)$ on the objects of a full exceptional collection projects as an action of $(\mathbb Z/2\mathbb Z)^{\times 6}$ on $K_0(\mathbb G)\otimes_{\mathbb Z}\mathbb C$ by changing of signs of the elements of the corresponding exceptional basis;
\item the braid group $\mathcal B_6$ acts on the set of exceptional collections (and the corresponding exceptional bases) as follows: the generator $\beta_{i,i+1}$ ($1\leq i\leq 5$) transforms the collection $(E_1,\dots, E_{i-1},E_i, E_{i+1},E_{i+2},\dots, E_6)$ into $(E_1,\dots, E_{i-1},L_{E_i}E_{i+1}, E_i,E_{i+2},\dots, E_6)$, where the object $L_{E_i}E_{i+1}$ is defined, up to unique isomorphism, by the distinguished triangle
\[L_{E_i}E_{i+1}[-1]\to \Hom^\bullet(E_i,E_{i+1})\otimes E_i\to E_{i+1}\to L_{E_i}E_{i+1}.
\]
\end{itemize}
Notice that our definition of braid mutations of exceptional objects differs {from} the one given, for example, in~\cite{helix} by a shift: this difference is important in order to obtain the coincidence of the braid group action on the matrix representing the morphism $\mathfrak X_{\mathbb G}^\pm$ with the action on the central connection matrix.
\end{Remark}

\begin{Remark}The conjecture we are discussing was also formulated in~\cite{gamma1} in the same time as~\cite{dubro4} for any Fano manifold~$X$. In~\cite{gamma1} the authors stress the relevance of the class $\widehat{\Gamma}^+(X)$, while in~\cite{dubro4} it was conjectured that $\widehat{\Gamma}^-(X)$ is the relevant characteristic class. As we will show below, $\widehat{\Gamma}^+(X)$ and $\widehat{\Gamma}^-(X)$ can be interchanged by the action~\eqref{29luglio2016-5-1} of the group $\mathcal C_0(\eta,\mu, R)$.
\end{Remark}

We now show that the monodromy data computed in Section~\ref{26gennaio2020-2} are of the {above} form for an exceptional collection in the same orbit of the Kapranov collection, under the action of the braid group. The Kapranov exceptional collection for $\mathbb G$ is formed by vector bundles $\mathbb S^\lambda(\mathcal S^*)$ ($\mathcal S$~is the tautological bundle), where $\mathbb S^\lambda$ denotes the Schur functor corresponding to the Young diagram~$\lambda$.\footnote{The reader can find the definition of Schur functors as endo-functors of the category of vector spaces in~\cite{fulton-harris}. The definition easily extends to the category of vector bundles.} In the general case of~$G_{\mathbb C}(k,n)$, the graded Chern character of these bundles is given by
\[\operatorname{Ch}\big(\mathbb S^{\lambda}(\mathcal S^*)\big)=s_\lambda\big({\rm e}^{2\pi \sqrt{-1} x_1},\dots,{\rm e}^{2\pi \sqrt{-1} x_k}\big):=\frac{\det\big({\rm e}^{2\pi \sqrt{-1} x_i(\lambda_j+r-j)}\big)_{1\leq i,j\leq k}}{\prod\limits_{i<j}\big({\rm e}^{2\pi \sqrt{-1} x_i}-{\rm e}^{2\pi \sqrt{-1} x_j}\big)},
\]i.e., the Schur polynomial calculated {at} the Chern roots $x_1,\dots, x_k$ of~$\mathcal S^*$.
In our case we obtain the following classes: posing $a:={\rm e}^{2\pi{\rm i} x_1}$ and $b:={\rm e}^{2\pi{\rm i} x_2}$ with $x_1+x_2=\sigma_1$ and $x_1x_2=\sigma_{1,1}$ we have that
\begin{alignat*}{3}
& \text{for }\lambda=0\qquad && \operatorname{Ch}\big(\mathbb S^{\lambda}(\mathcal S^*)\big)=1,&
\\
& \text{for }\lambda=\text{\tiny $\yng(1)$}, \qquad && \operatorname{Ch}\big(\mathbb S^{\lambda}(\mathcal S^*)\big)=a+b,&
\\
& \text{for }\lambda=\text{\tiny $\yng(2)$}, \qquad&& \operatorname{Ch}\big(\mathbb S^{\lambda}(\mathcal S^*)\big)=(a+b)^2-ab,&
\\
& \text{for }\lambda=\text{\tiny $\yng(1,1)$}, \qquad&& \operatorname{Ch}\big(\mathbb S^{\lambda}(\mathcal S^*)\big)=ab,&
\\
& \text{for }\lambda=\text{\tiny $\yng(2,1)$}, \qquad&& \operatorname{Ch}\big(\mathbb S^{\lambda}(\mathcal S^*)\big)=(a+b)ab,&
\\
& \text{for }\lambda=\text{\tiny $\yng(2,2)$}, \qquad&& \operatorname{Ch}\big(\mathbb S^{\lambda}(\mathcal S^*)\big)=a^2b^2.&
\end{alignat*}
Observing that
\begin{gather*} ab=1+2\pi{\rm i}\sigma_1-2\pi^2(\sigma_2+\sigma_{1,1})-\frac{8}{3}{\rm i}\pi^3\sigma_{2,1}+\frac{4}{3}\pi^4\sigma_{2,2},\\
a+b=2+2\pi{\rm i}\sigma_{1}-2\pi^2\sigma_2+2\pi^2\sigma_{1,1}+\frac{4}{3}{\rm i}\pi^3\sigma_{2,1},
\end{gather*} after some computations one obtains all graded Chern characters. Recalling the value of the $\widehat\Gamma^\mp$-class
\begin{gather*}\widehat\Gamma^\mp(\mathbb  G)=1 \pm4\gamma\sigma_1+\frac{1}{6}\big(48\gamma^2+\pi^2\big)(\sigma_{1,1}+\sigma_2) \pm\frac{4}{3}\big(16\gamma^3+\gamma\pi^2-\zeta(3)\big)\sigma_{2,1}\\
\hphantom{\widehat\Gamma^\mp(\mathbb  G)=}{} +\frac{1}{36} \big(768\gamma^4+96\gamma^2\pi^2-\pi^4-192\gamma\zeta(3)\big)\sigma_{2,2}
\end{gather*}
we can explicitly compute all the classes
\[\frac{1}{4\pi^2c^{\frac{1}{2}}}\big(\widehat\Gamma^\mp(\mathbb G)\cup\operatorname{Ch}\big(\mathbb S^{\lambda}(\mathcal S^*)\big)\big).
\]
We denote by $C^\mp_{\rm Kap}$ the matrix obtained in this way: in Appendix \ref{appc} the reader can find the entries of the matrix $C^-_{\rm Kap}$.

The Stokes matrix can be put in triangular form by a suitable permutation of $(u_1,\dots ,u_6)$, to which a permutation matrix $P$ is associated, according to the transformations (\ref{29luglio2016-9}). There are two permutations which yield $PSP^{-1}$ in triangular form, namely
\begin{gather}\label{6luglio2017-1}
\tau_1\colon \ (u_1,u_2,u_3,u_4,u_5,u_6)\mapsto (u_1^\prime,u_2^\prime,u_3^\prime,u_4^\prime,u_5^\prime,u_6^\prime):=(u_5,u_4,u_2,u_1,u_3,u_6),
\\
\label{6luglio2017-2}
\tau_2\colon \ (u_1,u_2,u_3,u_4,u_5,u_6)\mapsto (u_1^\prime,u_2^\prime,u_3^\prime,u_4^\prime,u_5^\prime,u_6^\prime):=(u_5,u_4,u_1,u_2,u_3,u_6).
\end{gather}
In both cases, the Stokes matrix $S$ in~(\ref{29luglio2016-1}), with $v=6$, becomes
\begin{gather}\label{29luglio2016-8-1}
S \longmapsto PSP^{-1}=
\left(
\begin{matrix}
 1 & -6 & 20 & 20 & 70 & 20 \\
 0 & 1 & -4 & -4 & -16 & -6 \\
 0 & 0 & 1 & 0 & 4 & 4 \\
 0 & 0 & 0 & 1 & 4 & 4 \\
 0 & 0 & 0 & 0 & 1 & 6 \\
 0 & 0 & 0 & 0 & 0 & 1
\end{matrix}
\right).
\end{gather}
The matrix $C$ in Appendix \ref{appc}, with $v=6$, becomes
\begin{gather}
\label{29luglio2016-8-2}
C\mapsto CP^{-1}.
\end{gather}

\begin{Theorem}\label{resultg24}The Stokes and connection matrices at $0\in QH^\bullet(\mathbb G)$ are related to the exceptional block collections obtained from the Kapranov block collection by mutations under the inverse of the braid $\beta_{12}\beta_{56}\beta_{45}\beta_{23}\beta_{34}$ or the braid $\beta_{34}\beta_{12}\beta_{56}\beta_{45}\beta_{23}\beta_{34}$ $($the action of $\beta_{34}$ acting just as a permutation of the third and fourth elements of the block$)$. \end{Theorem}

It is important to remark that the Kapranov $5$-block exceptional collection itself appears {\it neither} at $t=0$ {\it nor} anywhere else along the locus of the small quantum cohomology, see Corollary~\ref{7febbraio2020-7} below.
\begin{proof}Consider the monodromy data of the quantum cohomology of the Grassmannian $\mathbb G$ at $0\in QH^\bullet(\mathbb G)$, as computed in Section~\ref{28luglio2016-15} with respect to an admissible line\footnote{The computations have been done for $\phi=\pi/6$, but nothing changes if $0<\phi<\frac{\pi}{4}$, since the sectors, where the asymptotic behaviours are studied, are the same $\mathcal{S}_{\rm left/right}$.} $\ell=\ell(\phi)$ of slope $0<\phi<\frac{\pi}{4}$ and w.r.t.\ the basis of solutions (\ref{xi0}). These are the matrix~$S$ in formula~(\ref{29luglio2016-1}) and the matrix~$C$ in Appendix~\ref{appc}, with $v=6$. Arrange~$S$ in triangular form as in~(\ref{29luglio2016-8-1}), with~$P$ associated with one of the above permutations $\tau_1$ or $\tau_2$ above, and transform~$C$ as in~(\ref{29luglio2016-8-2}). The data so obtained are related to the Kapranov exceptional collection by a finite sequence of natural transformations~(\ref{29luglio2016-5-1}), (\ref{29luglio2016-5-2}), (\ref{29luglio2016-5-3}). More precisely, the following sequence transforms~$CP^{-1}$ into~$C^-_{\rm Kap}$:
\begin{itemize}\itemsep=0pt
\item[(1)] the change of signs in the normalised idempotents vector fields, determined by the action~(\ref{29luglio2016-5-2}) of the diagonal matrix $\mathcal I:= \operatorname{diag}(1,-1,-1,1,-1,1)$ (if we start from the cell where~$\tau_1$ is lexicographical), or $\mathcal I:=\operatorname{diag}(1,-1,1,-1,-1,1)$ (if we start from the cell where $\tau_2$ is lexicographical),
\item[(2)] change of solution at the origin through the action~(\ref{29luglio2016-5-1}), with $G$ equal to
\[A=\left(
\begin{matrix}
 1 & 0 & 0 & 0 & 0 & 0 \\
 2 {\rm i}\pi & 1 & 0 & 0 & 0 & 0 \\
 -2 \pi ^2 & 2 {\rm i}\pi & 1 & 0 & 0 & 0 \\
 -2 \pi ^2 & 2 {\rm i}\pi & 0 & 1 & 0 & 0 \\
 -\frac{1}{3} \big(8 {\rm i}\pi ^3\big) & -4 \pi ^2 & 2 {\rm i}\pi & 2 {\rm i}\pi & 1 & 0 \\
 \frac{4 \pi ^4}{3} & -\frac{1}{3} \big(8 {\rm i}\pi ^3\big) & -2 \pi ^2 & -2 \pi ^2 & 2 {\rm i}\pi & 1 \\
\end{matrix}
\right)\in\mathcal C_0(\eta,\mu, R),
\]
\item[(3)] the action~(\ref{29luglio2016-5-3}) with either the braid $\beta_{12}\beta_{56}\beta_{45}\beta_{23}\beta_{34}$ (if we start from the cell where $\tau_1$ is lexicographical), or the braid $\beta_{34}\beta_{12}\beta_{56}\beta_{45}\beta_{23}\beta_{34}$ (if we start from the cell where $\tau_2$ is lexicographical).
\end{itemize}
Moreover, $CP^{-1}$ in \eqref{29luglio2016-8-2} is transformed into $C^+_{\rm Kap}$ if, after the sequence of the above transformations (1), (2), (3) above, the following transformation is further applied:
\begin{itemize}\itemsep=0pt
\item[(4)]
the action (\ref{29luglio2016-5-1}), with matrix $G$ equal to
\[
B=
\left(
\begin{matrix}
 1 & 0 & 0 & 0 & 0 & 0 \\
 -8 \gamma & 1 & 0 & 0 & 0 & 0 \\
 32 \gamma ^2 & -8 \gamma & 1 & 0 & 0 & 0 \\
 32 \gamma ^2 & -8 \gamma & 0 & 1 & 0 & 0 \\
 \frac{8}{3} \big(\zeta (3)-64 \gamma ^3\big) & 64 \gamma ^2 & -8 \gamma & -8 \gamma & 1 & 0 \vspace{1mm}\\
 \frac{64}{3} \big(16 \gamma ^4-\gamma \zeta (3)\big) & \frac{8}{3} \big(\zeta (3)-64 \gamma ^3\big) & 32 \gamma ^2 & 32 \gamma ^2 & -8 \gamma & 1
\end{matrix}
\right)\in\mathcal C_0(\eta,\mu, R).
\]
\end{itemize}

The inverse of the Stokes matrix obtained from $PSP^{-1}$ in \eqref{29luglio2016-8-1} by either the sequence (1), (2), (3) or (1), (2), (3), (4) (recall that steps~(2) and~(4) act trivially on~$S$) coincides with the following Gram matrix of the Kapranov exceptional collection
\begin{gather}\label{6luglio2017-3}
\bold{G}_\text{\rm Kap}=\left(
\begin{matrix}
 1 & 4 & 10 & 6 & 20 & 20 \\
 0 & 1 & 4 & 4 & 16 & 20 \\
 0 & 0 & 1 & 0 & 4 & 10 \\
 0 & 0 & 0 & 1 & 4 & 6 \\
 0 & 0 & 0 & 0 & 1 & 4 \\
 0 & 0 & 0 & 0 & 0 & 1
\end{matrix}\right).
\end{gather}
\end{proof}

\begin{Remark}In both cases $C^{+}_{\rm Kap}$ and $C^{-}_{\rm Kap}$, the relation~(\ref{29luglio2016-2}) holds between $C^{\pm}_{\rm Kap}$ and $\bold{G}_\text{\rm Kap}^{-1}$.
\end{Remark}

\begin{Remark}The algebro-geometric meaning of the matrices $A$ and $B$ of Theorem~\ref{resultg24} will be thoroughly explained in our~\cite{CDG1}.
\end{Remark}

\subsection{Reconstruction of monodromy data along the small quantum locus}\label{recSC}
In this section we reconstruct the monodromy data at all other points of the small quantum cohomology of $\mathbb G$, by applying the procedure described in Section~\ref{30luglio2016-3}, and already illustrated in Section~\ref{A_3case}.

We identify the small quantum {cohomology} with the set of {points} $t=\big(0,t^2,0,\dots ,0\big)$. These points can be represented {on} the real plane $\big(\operatorname{Re} t^2,\operatorname{Im} t^2\big)$. At a point $\big(0,t^2,0,\dots ,0\big)$, the cano\-ni\-cal coordinates are~\eqref{28luglio2016-12}, so that the Stokes rays are
\[
R_{ij}\big(t^2\big)={\rm e}^{\overline{t^2}/4}R_{ij}(0)\equiv {\rm e}^{-{\rm i} \operatorname{Im}{t^2}/4}R_{ij}(0),
\]
where $R_{ij}(0)$ are the rays $R_{ij}$ of Section~\ref{30luglio2016-4}. Let $\ell$ be a line of slope $\varphi\in {}]0,\pi/4[$, admissible for $t^2=0$, i.e., for the Stokes rays $R_{ij}(0)$. Then, whenever $\operatorname{Im} t^2 \in \pi \cdot \mathbb{Z} -4\phi$, at least a pair of rays $R_{ij}\big(t^2\big)$ and $R_{ji}\big(t^2\big)$ lie along the line $\ell$, for some $(i,j)$. This means that the small quantum cohomology of $\mathbb G$ is split into the following horizontal bands of the $\big(\operatorname{Re} t^2,\operatorname{Im} t^2\big)$-plane:
\[
\mathcal{H}_k:=\big\{t^2\colon k\pi-4\phi< \operatorname{Im} t^2< (k+1)\pi-4\phi\big\},\qquad k\in \mathbb Z.
\]
If $t^2$ varies along a curve connecting two neighbouring bands, at least a pair of opposite rays $R_{ij}\big(t^2\big)$ and $R_{ji}\big(t^2\big)$ cross~$\ell$ in correspondence with~$t^2$ crossing the border between the bands.

A point $\big(0,t^2,0,\dots ,0\big)$, such that $t^2$ is interior to a band, is a semisimple coalescence point, where Theorem~\ref{mainisoth} applies. The polydisc $\mathcal{U}_{\epsilon_1}\big(u\big(0,t^2,\dots ,0\big)\big)$ is split into two $\ell$-cells. Each cell corresponds, through the coordinate map $p\mapsto u(p)$, to the closure of an open connected subset of an $\ell$-chamber of $QH^{\bullet}(\mathbb G)$, as explained in Section~\ref{30luglio2016-3}. Therefore, each band $\mathcal{H}_k$ precisely belongs to the boundary of two $\ell$-chambers corresponding to the two cells, while each line $ \operatorname{Im} t^2 =k\pi-4\phi$ between two bands $\mathcal{H}_{k-1}$ and $\mathcal{H}_k$ belongs to the intersection of the boundaries of four neighbouring chambers of $QH^{\bullet}(\mathbb G)$. As explained in Section~\ref{30luglio2016-3}, the monodromy data computed via Theorem~\ref{mainisoth} in $\mathcal{U}_{\epsilon_1}\big(u\big(0,t^2,\dots ,0\big)\big)$ are the data of the two chambers sharing the boundary $\mathcal{H}_k$. In particular, as a necessary consequence of Theorem~\ref{mainisoth}, these data are the data at each point of~$\mathcal{H}_k$. This means that {\it the monodromy data are constant in each band~$\mathcal{H}_k$}.

In order to compute the monodromy data in {every} chamber of $QH^{\bullet}(\mathbb G)$ is sufficient to apply the procedure of Section~\ref{30luglio2016-3} starting from the data~$C$, $S$ computed at $t=0$ in Section~\ref{28luglio2016-15}. Preliminarily, by a permutation~$P$, we have {obtained} upper triangular $PSP^{-1}$ and the corresponding $CP^{-1}$ in~\eqref{29luglio2016-8-1} and~\eqref{29luglio2016-8-2}, which are the monodromy data in the cell of $\mathcal{U}_{\epsilon_1}(u^\prime(0,0,\dots ,0))$ where $u_1^\prime(0,0,\dots ,0),\dots , u_{6}^\prime (0,0,\dots ,0)$ are in lexicographical order as in \eqref{6luglio2017-1} or~\eqref{6luglio2017-2}. Thus, they are the data of the band $\mathcal{H}_0$. Then, the braid group actions~(\ref{connbraid2-29luglio}) and~(\ref{connbraid2}) can be applied. In particular, we have computed the action of those braids which allow to pass from the chamber (with lexicographical order) whose boundary contains $\mathcal{H}_0$, to the chambers whose boundary contains $\mathcal{H}_k$, for $k=1,2,\dots ,8$. The values of $S$ and $C$ so obtained are, as explained above, the constant monodromy data for $\mathcal{H}_0$, $\mathcal{H}_1$, \dots, $\mathcal{H}_8$. They are reported in Table~\ref{30luglio2016-6}. From the table, we can read the monodromy data for the whole small quantum cohomology, since for any $k\in\mathbb Z$, the data for $\mathcal{H}_{k+8}$ are the same as for $\mathcal{H}_k$, as will be clear from the explanation below.

We need to determine the braid connecting neighbouring bands $\mathcal{H}_k$'s, from $\mathcal{H}_0$ to $\mathcal{H}_1$, from~$\mathcal{H}_1$ to~$\mathcal{H}_2$, and so on. The passage from $\mathcal{H}_k$ to $\mathcal{H}_{k+1}$ is achieved by increasing $\operatorname{Im}(t^2)$, to which a~clockwise rotation of the Stokes rays corresponds. In order to identify the corresponding braid, we have to keep track of the rays which cross a fixed~$\ell$ with slope $\phi\in {} ]0,\pi/4[$. Equivalently, we can consider a fixed configuration $u_1\big(0,t^2,\dots ,0\big), \dots , u_{6}\big(0,t^2,\dots ,0\big)$ in lexicographical order, corresponding to a fixed $t=\big(0,t^2,0,\dots,0\big)$ in $\mathcal{H}_0$, so that the corresponding rays $R_{ij}\big(t^2\big)$ are fixed. Then, we let $\ell$ rotate counter-clockwise increasing $ \phi$, with the consequent gliding of the $\ell$-horizontal bands towards $\operatorname{Im}\big(t^2\big)\to -\infty$, and we keep track of the rays which are crossed by~$\ell$.

In order to apply the procedure explained in Section \ref{30luglio2016-3}, we actually need to start with a~fixed configuration $u_1(t), \dots , u_{6}(t)$ of {\it distinct} canonical coordinates. This is achieved by taking~$t$ slightly away from $ \mathcal{H}_0$, in the interior of one of the two chambers whose boundaries contain~$\mathcal{H}_0$, so that the two coalescing canonical coordinates (which equal~$0$ in any~$\mathcal{H}_0$) slightly split. It is after this splitting that we let~$\ell$ rotate and keep track of the rays which are crossed by~$\ell$.

\begin{figure}[t]\centering
\def\svgscale{0.7}
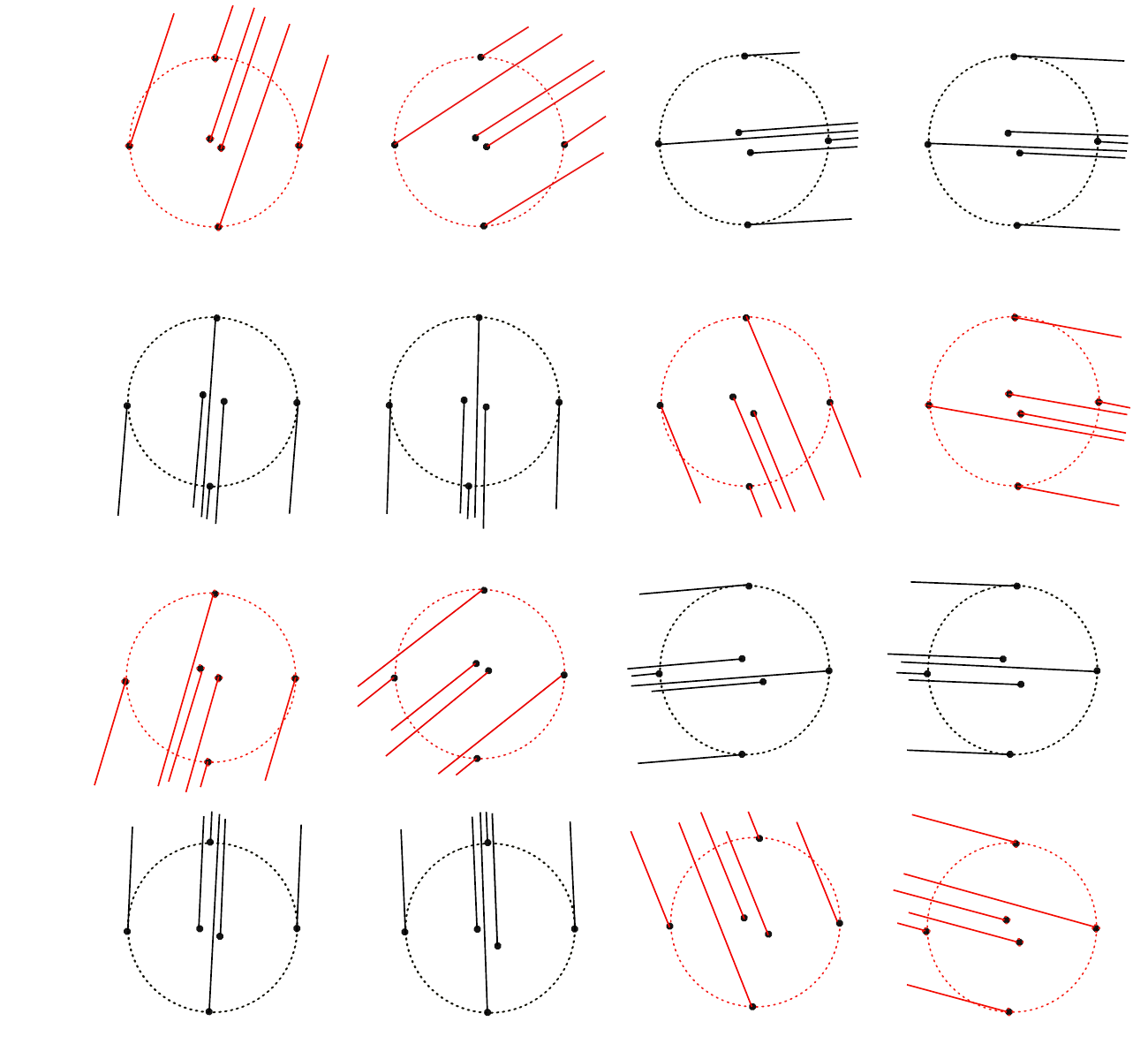
\caption{The picture, to be read in boustrophedon order, shows the braids corresponding to the passage from one band $\mathcal{H}_k$ to $\mathcal{H}_{k+1}$ . Starting from the configuration of the canonical coordinates at $0\in QH^\bullet(\mathbb G)$, we slightly split the coalescence as described in the first red picture in the first line. The numbers represent the lexicographical order of the canonical coordinates w.r.t.\ the admissible line. Letting the admissible line $\ell$ continuously rotate by increasing its slope, we determine all elementary braids acting in the mutation up to the next red configuration (the counter-clockwise rotation of $\ell$ is visualised by the clockwise rotation of the rays $L_j(\phi)$ defined in (\ref{3giugno2019-1})). By coalescence of the points $u_3$, $u_4$ in a red picture we obtain a configuration of canonical coordinates in the locus of small quantum cohomology (i.e., in a strip $\mathcal{H}_k$). Thus we deduce that successive bands of the small quantum cohomology are related by alternate compositions of the braids $\omega_1:=\beta_{12}\beta_{56}$, $\omega_{2}:=\beta_{23}\beta_{45}\beta_{34}\beta_{23}\beta_{45}$ and $\hat\omega_1:=\beta_{12}\beta_{34}\beta_{56}$.} \label{mutazionerettag24}
\end{figure}

This process is shown in Fig.~\ref{mutazionerettag24}. The two canonical coordinates close to~$0$ (the centres of the circles) come from the splitting of the two coalescing eigenvalues~$0$. The rays $L_j(\phi)$ defined in~(\ref{3giugno2019-1}) are represented. Their clockwise rotation corresponds to the counter-clockwise rotation of $\ell$. A Stokes ray $R_{j,j+1}(t)$ is crossed by $\ell$ every time a ray $L_j(\phi)$ aligns with $L_{j+1}(\phi)$. This yields the braid $\beta_{j,j+1}$. The order $u_j$, $u_{j+1}$ is $\ell$-lexicographical with reference to $\ell$ just before the crossing. Each time $\ell$ has just crossed a ray $R_{j,j+1}(t)$ (i.e., $L_j(\phi)$ has just crossed~$L_{j+1}(\phi)$), the coordinates $u_j$'s must be {\it relabelled} in {the} $\ell$-lexicographical order corresponding to $\ell$ just after the crossing.

As it can be read in Fig.~\ref{mutazionerettag24} and Table~\ref{30luglio2016-6}, the passage from $\mathcal{H}_0$ to $\mathcal{H}_{k}$, for $k=1,2,\dots ,8$, is obtained by composition of the braids
\[\omega_{1}:=\beta_{12}\beta_{56},
\qquad \omega_{2}:=\beta_{23}\beta_{45}\beta_{34}\beta_{23}\beta_{45},
\qquad \widehat{\omega}_{1}:=\beta_{12}\beta_{34}\beta_{56},
\]
in the form of products of increasing length
\begin{gather*}\omega_1, \qquad \omega_1\omega_2,\qquad
\omega_1\omega_2\widehat{\omega}_1, \qquad \omega_1\omega_2\widehat{\omega}_1\omega_2,
\qquad \omega_1\omega_2\widehat{\omega}_1\omega_2\omega_1,
\\
\omega_1\omega_2\widehat{\omega}_1\omega_2\omega_1\omega_2,
\qquad \omega_1\omega_2\widehat{\omega}_1\omega_2\omega_1\omega_2\widehat{\omega}_1,
\qquad
\omega_1\omega_2\widehat{\omega}_1\omega_2\omega_1\omega_2\widehat{\omega}_1\omega_2.
\end{gather*}
The result of Table~\ref{30luglio2016-6} is obtained applying formula~\eqref{connbraid2-29luglio}.

The last braid $\omega_1\omega_2\widehat{\omega}_1\omega_2\omega_1\omega_2\widehat{\omega}_1\omega_2$ corresponds to a complete $2\pi$-rotation of the admissible line~$\ell$. It is not difficult to show that
\begin{gather*}\omega_1\omega_2\widehat{\omega}_1\omega_2\omega_1\omega_2 \widehat{\omega}_1\omega_2=(\beta_{12}\beta_{23}\beta_{34}\beta_{45}\beta_{56})^6,\end{gather*}
 the right hand-side being the generator of the center of the braid group $\mathcal B_6$ in Lemma \ref{centerbraidlemma} (the proof can be done by graphically representing the two braids in the l.h.s.\ and r.h.s.\ respectively, and noticing that they can be deformed the one into the other). This corresponds to the cyclical repetition of the same Stokes matrix in $\mathcal{H}_k$ and $\mathcal{H}_{k+8}$ (while the central connection matrix $C$ is transformed to $M_0^{-1}C$).

 Notice that the action~\eqref{stokesbraid1} of $\beta_{34}$ in $\widehat{\omega}_1$ is a permutation of the third and fourth rows and columns of the Stokes matrix, because the entry~$(3,4)$ is zero (see~\eqref{9dicembre2016-2}).

It follows from Table \ref{30luglio2016-6} and the explicit Gram matrix of the Kapranov exceptional collection~\eqref{6luglio2017-3} that the following holds:
\begin{Corollary}\label{7febbraio2020-7}
The Kapranov $5$-block exceptional collection is not associated with the small quantum cohomology locus.
\end{Corollary}

\begin{Remark}
There is a remarkable similarity between the above cyclical repetition and the fact that exceptional collections are organised in algebraic structures called \emph{helices}, introduced in \cite{gore0,goru}, and extensively developed in \cite{gore1,gore2,helix}. This will be thoroughly explained in our paper~\cite{CDG1}.
\end{Remark}

\begin{table}[t]\centering

\caption{List of Stokes matrices {for all} bands decomposing the small quantum cohomology of $\mathbb G$: the computation is done at a point $\big(0,t^2,0,\dots,0\big)$ w.r.t.\ a line $\ell$ of slope $\phi\in{} ]0,\pi/4[$, admissible for $t=0$. The starting matrix $S_\text{lex}$ in $\mathcal{H}_0$ is $ PSP^{-1}$ of formula \eqref{29luglio2016-8-1}, with signs changed by~(\ref{29luglio2016-5-2}) with $\mathcal{I}=\operatorname{diag}(-1,1,1,-1,1,-1)$. The {braids} acting on the monodromy data are $\omega_1:=\beta_{12}\beta_{56}$, $\omega_2:=\beta_{23}\beta_{45}\beta_{34}\beta_{23}\beta_{45}$ and $\widehat{\omega}_1:=\beta_{12}\beta_{34}\beta_{56}$. The table is computed using \eqref{connbraid2-29luglio}, namely by successively applying \eqref{stokesbraid1} for the elementary braids.}\label{30luglio2016-6}

\vspace{1mm}

\begin{tabular}{|c||c|c|}
\hline
Band $\mathcal{H}_k$, $0\leq k \leq 8$ & $S_\text{lex}$ & Braid\\
\hline
\hline
&&\\[-2ex]
$0<\operatorname{Im}\big(t^2\big) +4\phi<\pi$& $\left(
\begin{matrix}
 1 & 6 & -20 & 20 & -70 & 20 \\
 0 & 1 & -4 & 4 & -16 & 6 \\
 0 & 0 & 1 & 0 & 4 & -4 \\
 0 & 0 & 0 & 1 & -4 & 4 \\
 0 & 0 & 0 & 0 & 1 & -6 \\
 0 & 0 & 0 & 0 & 0 & 1
\end{matrix}
\right)$& ${\rm id}$\\
& & \\[-2ex]
\hline
& & \\[-2ex]
$\pi<\operatorname{Im}\big(t^2\big)+4\phi<2\pi$ & $\left(
\begin{matrix}
 1 & -6 & -4 & 4 & 6 & 20 \\
 0 & 1 & 4 & -4 & -16 & -70 \\
 0 & 0 & 1 & 0 & -4 & -20 \\
 0 & 0 & 0 & 1 & 4 & 20 \\
 0 & 0 & 0 & 0 & 1 & 6 \\
 0 & 0 & 0 & 0 & 0 & 1
\end{matrix}
\right)$ & $\omega_1$\\
& & \\[-2ex]
\hline
& & \\[-2ex]
$2\pi<\operatorname{Im}\big(t^2\big)+4\phi<3\pi$ &
$\left(
\begin{matrix}
 1 & 6 & 20 & -20 & -70 & 20 \\
 0 & 1 & 4 & -4 & -16 & 6 \\
 0 & 0 & 1 & 0 & -4 & 4 \\
 0 & 0 & 0 & 1 & 4 & -4 \\
 0 & 0 & 0 & 0 & 1 & -6 \\
 0 & 0 & 0 & 0 & 0 & 1
\end{matrix}
\right)$
& $\omega_1\omega_2$\\
& &\\[-2ex]
\hline
& &\\[-2ex]
$3\pi<\operatorname{Im}\big(t^2\big)+4\phi<4\pi$ & $\left(
\begin{matrix}
 1 & -6 & -4 & 4 & 6 & 20 \\
 0 & 1 & 4 & -4 & -16 & -70 \\
 0 & 0 & 1 & 0 & -4 & -20 \\
 0 & 0 & 0 & 1 & 4 & 20 \\
 0 & 0 & 0 & 0 & 1 & 6 \\
 0 & 0 & 0 & 0 & 0 & 1
\end{matrix}
\right)$ & $\omega_1\omega_2\widehat{\omega}_1$\\
& & \\[-2ex]
\hline
& & \\[-2ex]
$4\pi<\operatorname{Im}\big(t^2\big)+4\phi<5\pi$ & $\left(
\begin{matrix}
 1 & 6 & 20 & -20 & -70 & 20 \\
 0 & 1 & 4 & -4 & -16 & 6 \\
 0 & 0 & 1 & 0 & -4 & 4 \\
 0 & 0 & 0 & 1 & 4 & -4 \\
 0 & 0 & 0 & 0 & 1 & -6 \\
 0 & 0 & 0 & 0 & 0 & 1
\end{matrix}
\right)$ & $\omega_1\omega_2\widehat{\omega}_1\omega_2$\\
& & \\[-2ex]
\hline
& & \\[-2ex]
$5\pi<\operatorname{Im}\big(t^2\big)+4\phi<6\pi$ & $\left(
\begin{matrix}
 1 & -6 & 4 & -4 & 6 & 20 \\
 0 & 1 & -4 & 4 & -16 & -70 \\
 0 & 0 & 1 & 0 & 4 & 20 \\
 0 & 0 & 0 & 1 & -4 & -20 \\
 0 & 0 & 0 & 0 & 1 & 6 \\
 0 & 0 & 0 & 0 & 0 & 1 \\
\end{matrix}
\right)$ & $\omega_1\omega_2\widehat{\omega}_1\omega_2\omega_1$\tsep{35pt}\\
\hline
\end{tabular}

\end{table}

\begin{table}[t]\centering

\rightline{Table 3. {\it Continued from the previous page.}}

\vspace{1mm}

\begin{tabular}{|c||c|c|}
\hline
Band $\mathcal{H}_k$, $0\leq k \leq 8$ & $S_\text{lex}$ & Braid\\
\hline
\hline
& &\\[-2ex]
$6\pi<\operatorname{Im}\big(t^2\big)+4\phi<7\pi$ & $\left(
\begin{matrix}
 1 & 6 & -20 & 20 & -70 & 20 \\
 0 & 1 & -4 & 4 & -16 & 6 \\
 0 & 0 & 1 & 0 & 4 & -4 \\
 0 & 0 & 0 & 1 & -4 & 4 \\
 0 & 0 & 0 & 0 & 1 & -6 \\
 0 & 0 & 0 & 0 & 0 & 1 \\
\end{matrix}
\right)$ & $\omega_1\omega_2\widehat{\omega}_1\omega_2\omega_1\omega_2$
\\
& &\\[-2ex]
\hline
& &\\[-2ex]
$7\pi<\operatorname{Im}\big(t^2\big)+4\phi<8\pi$ & $\left(
\begin{matrix}
 1 & -6 & 4 & -4 & 6 & 20 \\
 0 & 1 & -4 & 4 & -16 & -70 \\
 0 & 0 & 1 & 0 & 4 & 20 \\
 0 & 0 & 0 & 1 & -4 & -20 \\
 0 & 0 & 0 & 0 & 1 & 6 \\
 0 & 0 & 0 & 0 & 0 & 1
\end{matrix}
\right)$ & $\omega_1\omega_2\widehat{\omega}_1\omega_2\omega_1\omega_2\widehat{\omega}_1$\\
& &\\[-2ex]
\hline
& &\\[-2ex]
$8\pi<\operatorname{Im}\big(t^2\big)+4\phi<9\pi$ & $\left(
\begin{matrix}
 1 & 6 & -20 & 20 & -70 & 20 \\
 0 & 1 & -4 & 4 & -16 & 6 \\
 0 & 0 & 1 & 0 & 4 & -4 \\
 0 & 0 & 0 & 1 & -4 & 4 \\
 0 & 0 & 0 & 0 & 1 & -6 \\
 0 & 0 & 0 & 0 & 0 & 1
\end{matrix}
\right)$ & $\omega_1\omega_2\widehat{\omega}_1\omega_2\omega_1\omega_2\widehat{\omega}_1\omega_2$\bsep{35pt}\\
\hline
\end{tabular}
\end{table}

\section{A note on the topological solution for Fano manifolds}\label{topsolfano}

For quantum cohomologies of smooth projective varieties, a fundamental system of solutions of the equation for gradients of deformed flat coordinates
 \begin{gather}\label{16-06-16}
 \begin{cases}
 \partial_\alpha\zeta=z\mathcal C_\alpha\zeta,\\
 \partial_z\zeta=\left(\mathcal U+\dfrac{1}{z}\mu\right)\zeta,
 \end{cases}
 \end{gather}
can be expressed in \emph{enumerative-topological} terms, namely the genus 0 correlations functions.
\begin{Proposition}\label{topsol}For a sufficiently small $R>0$, it is defined an analytic function
\[\Theta\colon \ B_{\mathbb C}(0;R)\times\Omega\to\operatorname{End}(H^\bullet(X;\mathbb C))
\]with series expansion
\begin{gather*}\Theta(z,t): =\operatorname{Id}+\sum_{\alpha=0}^N \llangle[\bigg]\frac{z\cdot(-)}{1-z\psi},T_\alpha\rrangle[\bigg]_0(t)T^\alpha
 =\operatorname{Id}+\sum_{n=0}^\infty z^{n+1}\sum_{\alpha=0}^N\llangle\tau_n(-),T_\alpha\rrangle_0(t)T^\alpha.
\end{gather*}
This function $\Theta$ satisfies the following properties:
\begin{enumerate}\itemsep=0pt
\item[$1)$] for any $\phi\in H^\bullet(X;\mathbb C)$, the vector field
\begin{gather*}\Theta_\phi:=\Theta(z,t)\phi =\phi+\sum_{\alpha=0}^N\llangle[\bigg] \frac{z\phi}{1-z\psi},T_\alpha\rrangle[\bigg]_0(t)T^\alpha
 =\phi+\sum_{n=0}^\infty z^{n+1}\sum_{\alpha=0}^N\llangle\tau_n\phi,T_\alpha\rrangle_0(t)T^\alpha
\end{gather*}
satisfies the equations
\[\partial_\alpha\Theta_\phi=z\partial_\alpha*\Theta_\phi;
\]
\item[$2)$] when restricted to the small quantum locus $\Omega\cap H^2(X;\mathbb C)$, i.e., $t^i=0$ for $i=0,r+1,\dots, N$, then
\[\Theta_\phi={\rm e}^{z\delta}\cup\phi+\sum_{\beta\neq 0}\sum_{\alpha=0}^N{\rm e}^{\int_\beta\delta}\bigg\langle \frac{z{\rm e}^{z\delta}\cup\phi}{1-z\psi},T_\alpha\bigg\rangle^X_{0,2,\beta}T^\alpha,\qquad\delta:=\sum_{i=1}^rt^iT_i\in H^2(X;\mathbb C);
\]
\item[$3)$] for any $\phi_1,\phi_2\in H^\bullet(X;\mathbb C)$ we have
\[\eta\left(\Theta(-z,t)\phi_1,\Theta(z,t)\phi_2\right)=\eta\left(\phi_1,\phi_2\right);
\]
\item[$4)$] for any $\phi\in H^\bullet(X;\mathbb C)$, the vector field
\[\widetilde{\Theta}_\phi:=\big(\Theta(z,t)\circ z^\mu z^{c_1(X)\cup(-)}\big)\phi
\]is a solution of the system \eqref{16-06-16}, i.e.,
\[\partial_\alpha\widetilde{\Theta}_\phi=z\partial_\alpha*\widetilde{\Theta}_\phi,
\qquad
\partial_z\widetilde{\Theta}_\phi=\left(\mathcal U+\frac{1}{z}\mu\right)\widetilde{\Theta}_\phi.
\]
\end{enumerate}Thus, the vector fields $\widetilde{\Theta}_{T_\alpha}$'s are gradients of deformed flat coordinates: if $\big(\Theta^\alpha_\beta\big)_{\alpha,\beta}$, $\big(\widetilde{\Theta}^\alpha_\beta\big)_{\alpha,\beta}$ are the matrices representing the two $\operatorname{End}(H^\bullet(X;\mathbb C))$-valued functions~$\Theta$ and~$\widetilde{\Theta}$ w.r.t.\ the basis~$(T_\alpha)_\alpha$, i.e.,
\[\Theta(z,t)T_\beta=\sum_{\alpha=0}^N\Theta^\alpha_\beta(z,t)T_\alpha,\qquad \widetilde{\Theta}(z,t)T_\beta=\sum_{\alpha=0}^N\widetilde{\Theta}^\alpha_\beta(z,t)T_\alpha,
\]then there exist analytic functions $\big(\widetilde{t}_\alpha(z,t)\big)_\alpha$, $(h_\alpha(z,t))_\alpha$ on $ B_{\mathbb C}(0;R)\times\Omega$ such that
\begin{gather*} \widetilde{\Theta}^\alpha_\beta(z,t)=\big(\operatorname{grad} \tilde{t}_\beta(z,t)\big)^\alpha,\qquad \big(\tilde{t}_0,\tilde{t}_1,\dots,\tilde{t}_N\big)=(h_0,h_1,\dots, h_N)\cdot z^\mu z^{R},\\
\Theta^\alpha_{\beta}(z,t)=(\operatorname{grad} h_\beta(z,t))^\alpha,\qquad  \Theta^{\rm T}(-z,t)\eta\Theta(z,t)=\eta,\\
h_\alpha(z,t):=\sum_{p=0}^\infty h_{\alpha,p}(t)z^p,\qquad h_{\alpha,0}(t)=t_\alpha\equiv t^\lambda\eta_{\lambda\alpha}.
\end{gather*}
\end{Proposition}

\begin{proof} Notice that
\[Y(z,t):=H(z,t)z^\mu z^R,\qquad H(z,t)=\sum_{p=0}^\infty H_p(t)z^p,\qquad H_0(t)\equiv\mathbbm 1
\]is a fundamental solution of \eqref{16-06-16} if and only if $H(z,t)$ {satisfies} the system
\[\begin{cases}
\partial_\alpha H=z\mathcal C_\alpha H,\\
\partial_zH=\mathcal UH+\dfrac{1}{z}[\mu, H]-HR.
\end{cases}
\]Because of the symmetry of $c_{\alpha\beta\gamma}$, the columns of $H$ are the components w.r.t.\ $(\partial_\alpha)_\alpha$ of the gradients of some functions: \begin{gather*}h_\alpha(z,t):=\sum_{p=0}^\infty h_{\alpha,p}(t)z^p,\qquad h_{\alpha,0}(t)= t_\alpha,\\
H_\beta^\alpha(z,t)=(\operatorname{grad}h_\beta)^\alpha,\qquad H_{\beta,p}^\alpha(z,t)=(\operatorname{grad}h_{\beta,p})^\alpha.
\end{gather*}
The above system for $H$ is equivalent to the following recursion relations on $h_{\alpha,p}$'s functions:
\begin{gather}\label{recursion1}\partial_\alpha\partial_\beta h_{\gamma,p}(t)=c_{\alpha\beta}^\nu\partial_\nu h_{\gamma,p-1}(t),\qquad p\geq1,
\\
\mathfrak L_E(\operatorname{grad}h_{\alpha,p})=\left(p+\frac{\dim_{\mathbb C}X-2}{2}+\mu_\alpha\right)\operatorname{grad}h_{\alpha,p}+\sum_{\beta=0}^N (\operatorname{grad}h_{\beta,p-1})R^\beta_\alpha,\qquad p\geq1.\nonumber
\end{gather}
The last equation is equivalent to the recursion relations on the differentials
\begin{gather}\label{recursion3}
\mathfrak L_E({\rm d}h_{\alpha,p})=\left(p-\frac{\dim_{\mathbb C}X-2}{2}+\mu_\alpha\right) {\rm d}h_{\alpha,p}+\sum_{\beta=0}^N {\rm d}h_{\beta,p-1}R^\beta_\alpha,\qquad p\geq1.
\end{gather}
In our case we have
\[H(z,t)=\big(\Theta^\alpha_{\beta}(z,t)\big)_{\alpha,\beta},\qquad \partial_\alpha h_{\beta,p}(t)=\llangle\tau_{p-1}T_\beta,T_\alpha\rrangle_0(t).
\]
The recursion relations \eqref{recursion1} then reads{\samepage
\begin{gather*} \llangle T_\alpha, T_\beta, T_\gamma\rrangle_0=\llangle T_\alpha, T_\beta, T^\nu\rrangle_0 \eta_{\nu\gamma}\qquad\text{for }p=1,\\
\llangle T_\alpha,\tau_{p-1} T_\gamma, T_\beta\rrangle_0=\llangle T_\alpha, T_\beta, T^\nu\rrangle_0\llangle\tau_{p-2}T_\gamma, T_\nu\rrangle_0\qquad\text{for }p\geq2.
\end{gather*}
These are exactly the \emph{topological recursion relations in genus~$0$}.}

Let us now prove that also the recursion relations \eqref{recursion3} hold. K.~Hori~\cite{hori} (see also \cite{eguchi}) proved that, for any $\omega\in H^2(X;\mathbb C)$, we have the following constraint on the genus $g$ free energy
\begin{gather}
\omega^\alpha\frac{\partial\mathcal F^X_g}{\partial t^{\alpha,0}}=
\sum_{\beta\in \operatorname{Eff}(X)}\left(\int_\beta\omega\right)\mathcal F^X_{g,\beta}+\sum_{n,\alpha,\atop \sigma,\nu}\omega^\sigma c_{\sigma\alpha}^\nu t^{\alpha, n}\frac{\partial\mathcal F^X_g}{\partial t^{\nu,n-1}}\nonumber\\
\hphantom{\omega^\alpha\frac{\partial\mathcal F^X_g}{\partial t^{\alpha,0}}=}{} +\frac{\delta^0_g}{2}\sum_{\alpha,\nu,\sigma}\omega^\sigma c_{\sigma\alpha\nu}t^{\alpha,0}t^{\nu,0}-\frac{\delta^1_g}{24}\int_X\omega\cup c_{\dim X-1}(X),\label{hori1}
\end{gather}
where
$\mathcal F^X_{g,\beta}$ is the $(g,\beta)$-free energy
\[\mathcal F^X_{g,\beta}:=\sum_{n=0}^\infty\frac{1}{n!}\langle\underbrace{\gamma,\dots,\gamma}_{n\text{ times}}\rangle_{g,n,\beta}^X.
\]By dimensional consideration, one obtains also the selection rule
\begin{gather}\label{hori2}\sum_{n,\alpha}(n+q_\alpha-1)t^{\alpha,n}\frac{\partial\mathcal F^X_g}{\partial t^{\alpha,n}}=\sum_{\beta\in \operatorname{Eff}(X)}\left(\int_\beta\omega\right)\mathcal F^X_{g,\beta}+(3-\dim X)(g-1)\mathcal F^X_g.
\end{gather}
If we introduce the perturbed first Chern class
\[\mathcal E(\bold t):=c_1(X)+\sum_{m,\sigma}(1-q_\sigma-m)t^{\sigma,m}\tau_m(T_\sigma)-\sum_{m,\sigma}t^{\sigma,m}\tau_{m-1}(c_1(X)\cup T_\sigma),
\]and using the selection rule \eqref{hori2}, the Hori's constraint \eqref{hori1} (specialized to $g=0$ and $\omega=c_1(X)$) can be reformulated as
\[\llangle\mathcal E\rrangle_0=(3-\dim X)\mathcal F^X_0+\frac{1}{2}t^{\sigma,0}t^{\rho,0}\int_X c_1(X)\cup T_\sigma\cup T_\rho.
\]Taking the derivative w.r.t.\ $t^{\alpha,n}$, $t^{\beta,0}$ we obtain
\begin{gather*}\llangle\mathcal E,\tau_nT_\alpha,T_\beta\rrangle_0 -(n+q_\alpha+q_\beta-2)\llangle\tau_nT_\alpha,T_\beta\rrangle_0-\llangle\tau_{n-1}(c_1(X)\cup T_\alpha), T_\beta\rrangle_0\\
\qquad{} =(3-\dim X)\llangle\mathcal \tau_nT_\alpha,T_\beta\rrangle_0+\delta_{n,0}\int_X c_1(X)\cup T_\alpha\cup T_\beta.
\end{gather*}
These recursion relations, restricted to the small phase space, are easily seen to be equivalent to~\eqref{recursion3}. This proves~(1),~(4) and the convergence of~$\Theta(z,t)$ for~$|z|$ small enough, because of the regular feature of the singularity $z=0$. The proof of~(2) can be found in~\cite{cox}. Condition~(3) follows from WDVV and string equation, as shown in~\cite{giv}.
\end{proof}

In the case of Fano manifolds, we have the following analytic characterization of the fundamental solution $\big(\widetilde{\Theta}^\alpha_\beta\big)_{\alpha,\beta}$. Furthermore, because of Proposition~\ref{6marzo2017}, we obtain another proof of~(3) in the previous proposition.

\begin{Proposition}\label{2luglio2017-3}If $X$ is a Fano manifold, among all fundamental matrix solutions of the system~\eqref{16-06-16} for deformed flat coordinates,\footnote{Throughout the paper, $Y(z,t)=H(z,t)z^\mu z^R$ has been denoted $Y(z,t)=\Phi(z,t)z^\mu z^R$.} there exists a unique solution such that, on the small quantum locus $($i.e., $t^i=0$ for $i=0,r+1,\dots, N)$ the function $z^{-\mu}H(z,t)z^\mu$ is holomorphic at $z=0$, with series expansion
\[z^{-\mu}H(z,t)z^\mu= {\rm e}^{t\cup}+z K_1(t)+z^2 K_2(t)+\cdots,\qquad t^i=0 \ \text{for} \ i=0,r+1,\dots, N,
\]
This solution coincides with the solution $\big(\widetilde{\Theta}^\alpha_\beta(z,t)\big)_{\alpha,\beta}$.
\end{Proposition}

\begin{proof} We already know from Proposition~\ref{6marzo2017} that such a solution is unique. Let us now prove the main statement. In what follows, we will denote the degree $\deg T_\alpha$ just by~$|\alpha|$ for brevity. By point~$(2)$ of Proposition~\ref{topsol}, we have that
\[z^{-\mu}\big(\Theta_{(z^\mu\phi)}\big)=z^{-\mu}\left({\rm e}^{z\delta}\cup z^\mu\phi+\sum_{\beta\neq 0}\sum_{\alpha=0}^N{\rm e}^{\int_\beta\delta}\bigg\langle \frac{z{\rm e}^{z\delta}\cup z^\mu\phi}{1-z\psi},T_\alpha\bigg\rangle^X_{0,2,\beta}T^\alpha\right),\]
with $\delta:=\sum\limits_{i=1}^rt^iT_i\in H^2(X;\mathbb C)$. Specialising to $\phi=T_\sigma$, we have
\[z^{-\mu}\left(\Theta_{(z^\mu T_\sigma)}\right)={\rm e}^{\delta}\cup T_\sigma+\sum_{\beta\neq 0}\sum_{\alpha,\lambda=0}^N\sum_{n,k=0}^\infty \frac{{\rm e}^{\int_\beta\delta}}{k!}z^{n+1+k+\mu_\sigma-\mu_{\lambda}}\big\langle\tau_n\big(\delta^{\cup k}\cup T_\sigma\big), T_\alpha\big\rangle^X_{0,2,\beta}\eta^{\alpha\lambda}T_\lambda.
\]In the second addend, we have non-zero terms only if
\begin{itemize}\itemsep=0pt
\item $|\alpha|+|\lambda|=2\dim_\mathbb C X$,
\item $2n+2k+|\sigma|+|\alpha|=\operatorname{vir\ dim}_{\mathbb R} X_{0,2,\beta}$.
\end{itemize}
By putting together these conditions, we obtain
\[n+1+k+\frac{1}{2}(|\sigma|-|\lambda|)=-\int_{\beta}\omega_X.
\]
The assumption of being Fano is equivalent to the requirement that the functional $\beta\mapsto -\int_{\beta}\omega_X$ is positive on the closure of the effective cone. This proves the proposition, the l.h.s.\ being exactly the exponents of $z$ which appear in the {above} series expansion.
\end{proof}

\begin{Example}Notice that the solution \eqref{xi0} that we considered in the previous section for the computation of the monodromy data fo $QH^\bullet(\mathbb G)$ satisfies the condition
\begin{gather*} z^{-\mu}\big(\eta^{-1}S(0,z)\eta \big)z^\mu \qquad \text{is holomorphic near }z=0,
\\
z^{-\mu}\big(\eta^{-1}S(0,z)\eta \big)z^\mu=\left(
\begin{matrix}
 1-2 z^4 & 2 z^4 & - z^4 & - z^4 & z^4 & z^8 \\
 0 & 4 z^4+1 & - z^4 & - z^4 & 0 & z^4 \\
 0 & 0 & 1 & 0 & - z^4 & z^4 \\
 0 & 0 & 0 & 1 & - z^4 & z^4 \\
 0 & 0 & 0 & 0 & 1-4 z^4 & 2 z^4 \\
 0 & 0 & 0 & 0 & 0 & 2 z^4+1
\end{matrix}
\right)+O\big(z^9\big).
\end{gather*} This means that $\big(\eta^{-1}S(0,z)\eta \big)z^\mu z^R$ coincides with the topological solution $\widetilde{\Theta}(0,z)$.
\end{Example}

\appendix

\section{Proofs of Lemmata \ref{sol} and \ref{asympt}}\label{applemma}

In this Appendix we prove Lemmata \ref{sol} and~\ref{asympt}.
Before giving the proof of Lemma~\ref{sol}, we recall the following well-known results (see, e.g., \cite{luke1,nist, whitt-wats}).

\begin{Theorem}[Stirling]\label{stirl}The following estimate holds
\[\log\Gamma(s)=\left(s-\frac{1}{2}\right)\log s-s+\frac{1}{2}\log(2\pi)+O\left(\frac{1}{|s|}\right)
\]for $s\to\infty$ and $|\arg s|<\pi$, and where $\log$ stands for the principal {branch} of the complex logarithm.
\end{Theorem}

\begin{Corollary}\label{corstirl} For $|t|\to +\infty$ we have
\[
|\Gamma(\sigma+{\rm i}t)|=\sqrt{2\pi}|t|^{\sigma-\frac{1}{2}}{\rm e}^{-\frac{\pi}{2}|t|}\left(1+O\left(\frac{1}{|t|}\right)\right),
\]
uniformly on any strip of the complex plane $\sigma_1\leq\sigma\leq\sigma_2$.
\end{Corollary}

\begin{proof}[Proof of Lemma \ref{sol}]First of all let us prove that the functions $\vp_1(w)$, $\vp_2(w)$ are well defined on the sectors
\[
-\frac{\pi}{2}<\arg w<\frac{\pi}{2},\qquad -\frac{\pi}{2}< \arg w<\pi,
\] respectively.
We denote by $\mathcal I_1$ and $\mathcal I_2$ the integrands in $\vp_1$ and $\vp_2$ respectively, and $s=\kappa+{\rm i} t$. By Corollary \ref{corstirl} we have that
\[|\mathcal I_1|\sim (2\pi)^2|t|^{5\left(\kappa-\frac{1}{2}\right)}{\rm e}^{-\frac{5\pi}{2}|t|}|t|^{-\kappa}{\rm e}^{\frac{\pi}{2}|t|}{\rm e}^{-\kappa\log4}{\rm e}^{-4\kappa\log|w|+4t\arg w}.
\]
The dominant part is
\[{\rm e}^{-\frac{5\pi}{2}|t|}{\rm e}^{\frac{\pi}{2}|t|}{\rm e}^{4t\arg w}.
\]In order to have $|\mathcal I_1|\to 0$ for $t\to +\infty$ we must impose
\[-\frac{5\pi}{2}+\frac{\pi}{2}+4\arg w<0,\qquad\text{i.e.,} \qquad \arg w<\frac{\pi}{2};
\]analogously, for $t\to -\infty$ we have to impose
\[\frac{5\pi}{2}-\frac{\pi}{2}+4 \arg w>0,\qquad\text{i.e.,} \qquad \arg w>-\frac{\pi}{2}.
\]
Analogously, from Corollary \ref{corstirl} we deduce that
\[|\mathcal I_2|\sim(2\pi)^3|t|^{5(\kappa-\frac{1}{2})}{\rm e}^{-\frac{5\pi}{2}|t|}|{-}t|^{-\kappa}{\rm e}^{-\frac{\pi}{2}|-t|}{\rm e}^{-\pi t}{\rm e}^{-\kappa\log4}{\rm e}^{-4\kappa\log|w|+4t\arg w},
\]and now the dominant part is
\[{\rm e}^{-\frac{5\pi}{2}|t|}{\rm e}^{-\frac{\pi}{2}|{-}t|}{\rm e}^{-\pi t}{\rm e}^{4t\arg w}.
\]
In order to have $|\mathcal I_2|\to 0$ for $t\to\pm\infty$, we find
\begin{gather*}
-\frac{5\pi}{2}-\frac{\pi}{2}-\pi+4\arg w<0,\qquad\text{i.e.,} \qquad \arg w<\pi,
\\
\frac{5\pi}{2}+\frac{\pi}{2}-\pi+4\arg w>0,\qquad\text{i.e.,} \qquad \arg w<-\frac{\pi}{2}.
\end{gather*}
Let us now prove that $\vp_1$ and $\vp_2$ are solutions of equation~\eqref{eqip}. We have that
\begin{gather*}
\Theta^5\vp_1(w) =\frac{4^5}{2\pi{\rm i}}\int_{\Lambda_1}-s^5\frac{\Gamma(s)^5}{\Gamma\big(s+\frac{1}{2}\big)}4^{-s}w^{-4s}\,{\rm d}s
 =\frac{4^5}{2\pi{\rm i}}\int_{\Lambda_1}-\left(s+\frac{1}{2}\right)\frac{\Gamma(s+1)^5} {\Gamma\big(s+\frac{3}{2}\big)}4^{-s}w^{-4s}\,{\rm d}s
\end{gather*}because $z\Gamma(z)=\Gamma(z+1)$. Changing variable $t:=s+1$, and consequently shifting the line of integration $\Lambda_1$ to $\Lambda_1+1$, we have
\begin{align*}
\Theta^5\vp_1(w)&=\frac{4^5}{2\pi{\rm i}}\int_{\Lambda_1+1}-\left(t-\frac{1}{2}\right)\frac{\Gamma(t)^5}{\Gamma\big(t+\frac{1}{2}\big)}4^{-t}\cdot4 w^{-4(t-1)}\,{\rm d}t\\
&=\frac{4^5}{2\pi{\rm i}}\int_{\Lambda_1+1} (-4t)\frac{\Gamma(t)^5}{\Gamma\big(t+\frac{1}{2}\big)}4^{-t} w^{-4(t-1)}\,{\rm d}t
+\frac{2\cdot4^5}{2\pi{\rm i}}\int_{\Lambda_1+1}\frac{\Gamma(t)^5}{\Gamma\big(t+\frac{1}{2}\big)}4^{-t} w^{-4(t-1)}\,{\rm d}t.
\end{align*}
Note that in the region between $\Lambda_1$ and $\Lambda_1+1$ the two last integrands have no poles; so $\int_{\Lambda_1+1}=\int_{\Lambda_1}$ by Cauchy theorem. This shows that
\[\Theta^5\vp_1=4^5w^4\Theta\vp_1+2\cdot4^5w^4\vp_1.
\]
Analogously we have
\begin{align*}
\Theta^5\vp_2&=\frac{4^5}{2\pi{\rm i}}\int_{\Lambda_2}-s^5\Gamma(s)^5\Gamma\left(\frac{1}{2}-s\right){\rm e}^{{\rm i}\pi s}4^{-s}w^{-4s}\,{\rm d}s\\
&=\frac{4^5}{2\pi{\rm i}}\int_{\Lambda_2}-\left(s+\frac{1}{2}\right)\Gamma(s+1)^5\Gamma\left(-\frac{1}{2}-s\right){\rm e}^{{\rm i}\pi (s+1)}4^{-s}w^{-4s}\,{\rm d}s,
\end{align*}
where the second identity follows from $z\Gamma(z)=\Gamma(z+1)$. Note that the integrand function is holomorphic at $s=-\frac{1}{2}$: indeed we have
\[\lim_{s\to{-\frac{1}{2}}}\left(s+\frac{1}{2}\right)\Gamma\left(-\frac{1}{2}-s\right)=-1.
\]So in the strip of the complex plane $-1<\operatorname{Re}s<\frac{1}{2}$ there are no poles, and by Cauchy theorem, we can change path of integration by shifting $\Lambda_2$ to $\Lambda_2-1$:
\[\Theta^5\vp_2=\frac{4^5}{2\pi{\rm i}}\int_{\Lambda_2-1}-s^5\Gamma(s)^5\Gamma\left(\frac{1}{2}-s\right){\rm e}^{{\rm i}\pi s}4^{-s}w^{-4s}\,{\rm d}s.\] Posing now $t=s+1$, we can rewrite
\begin{gather*}
\Theta^5\vp_2=\frac{4^5}{2\pi{\rm i}}\int_{\Lambda_2}\!\!-\left(t-\frac{1}{2}\right)\Gamma(t)^5\Gamma\left(\frac{1}{2}-t\right){\rm e}^{{\rm i}\pi t}4^{-(t-1)}w^{-4(t-1)}\,{\rm d}t
=
4^5w^4\Theta\vp_2+2\cdot 4^5w^4\vp_2.
\end{gather*}
This shows that effectively $\vp_1$ and $\vp_2$ are solutions.\end{proof}

\begin{proof}[Proof of Lemma \ref{asympt}] By Stirling's formula we have that
\[\vp_1(w)=\frac{(2\pi)^2}{2\pi{\rm i}}\int_{\Lambda_1}{\rm e}^{\phi(s)}\,{\rm d}s,
\]where
\[\phi(s)=-5s+5\left(s-\frac{1}{2}\right)\log s+s+\frac{1}{2}-s\log\left(s+\frac{1}{2}\right)-s\log4-4s\log w+O\left(\frac{1}{|s|}\right)
\]for $s\to\infty$ and where $\log$ stands for the principal {branch} of logarithm.
Let us find stationary points of $\phi(s)$ for large values of $|s|,|w|$. The derivative $\phi'$ is
\[\phi'(s)=-4+5\log s+\frac{10s-5}{2s}-\log\left(s+\frac{1}{2}\right)-\frac{s}{s+\frac{1}{2}}-\log4-4\log w+O\left(\frac{1}{|s|}\right).
\]For $|s|$ large enough, we have
\begin{gather*}
\frac{10s-5}{2s}\sim 5-\frac{5}{2s},\qquad \frac{s}{s+\frac{1}{2}}\sim 1-\frac{1}{2s}, \\
\log\left(s+\frac{1}{2}\right)=\log s+\log\left(1+\frac{1}{2s}\right)\sim \log s+\frac{1}{2s}.
\end{gather*}
Substituting these identities in $\phi'$, we find that the critical point $\bar s(w)$ in functions of $w$ (for $|w|$ large)
\begin{gather*}
\bar s(w)=\sqrt{2}w+\frac{5}{8}+O\left(\frac{1}{|w|}\right).
\end{gather*}
Note that for $-\frac{\pi}{2}<\arg w<\frac{\pi}{2}$, the point $\bar s(w)$ is in the half-plane $\operatorname{Re}s>0$, region in which there are no poles of the integrand functions in $\vp_1$. So we can shift the line $\Lambda_1$ in order that it passes through $\bar s$.
In this way we obtain
\[\vp_1(w)=\frac{(2\pi)^2}{2\pi{\rm i}}{\rm e}^{\phi(\bar s)}\int_{\Lambda_1}{\rm e}^{\phi(s)-\phi(\bar s)}\,{\rm d}s
 \sim\frac{(2\pi)^2}{2\pi{\rm i}}{\rm e}^{\phi(\bar s)}\int_{\Lambda_1}{\rm e}^{\frac{\phi''(\bar s)}{2}(s-\bar s)^2}\,{\rm d}s.
\]
The computation of this Gaussian integral shows that
\[\vp_1(w)\sim\frac{(2\pi)^2}{2\pi }{\rm e}^{\phi(\bar s)}\frac{\sqrt{2\pi}}{\sqrt{\phi''(\bar s)}}=(2\pi)^{\frac{3}{2}}\frac{{\rm e}^{\phi(\bar s)}}{\sqrt{\phi''(\bar s)}},
\]where $\operatorname{Re}\sqrt{\phi''(\bar s)}>0$.
An explicit series expansion shows that
\[\phi(\bar s(w))\sim -4\sqrt{2}w-\frac{5}{2}\log w-\frac{5}{8}\log4+O\left(\frac{1}{|w|}\right),
\]
whereas
\[\phi''(\bar s(w))\sim \frac{2\sqrt{2}}{w}+O\left(\frac{1}{|w|^3}\right)
\]and from this we deduce that
\[\vp_1(w)\sim(2\pi)^{\frac{3}{2}}\frac{{\rm e}^{-4\sqrt{2}w}}{4w^2}\left(1+O\left(\frac{1}{w}\right)\right).
\]

Let us now focus on $\vp_2(w)$. From Theorem \ref{stirl} we deduce that
\[\Gamma(-s)={\rm e}^{-\left(s+\frac{1}{2}\right)\log s}{\rm e}^{-{\rm i}\pi s}{\rm e}^s\big({-}{\rm i}\sqrt{2\pi}\big)\left(1+O\left(\frac{1}{|s|}\right)\right)
\]
for $s\to\infty$ and $s\notin\mathbb R_+$.
So,
\[\vp_2(w)=\frac{(2\pi)^3}{2\pi{\rm i}}\int_{\Lambda_2}{\rm e}^{\phi(s)}\,{\rm d}s,
\]
where
\[\phi(s)=5\left(s-\frac{1}{2}\right)\log s-5s-s\log\left(s-\frac{1}{2}\right)+s-\frac{1}{2}-s\log4-4s\log w+O\left(\frac{1}{w}\right),
\]for $w\to\infty$. By computations analogous to those of the previous case, we find that $\phi$ has a~critical point at
\[\bar s(w)=\sqrt{2}w+\frac{5}{4\sqrt{2}}+O\left(\frac{1}{w}\right)
\]
for large values of $|w|$. Note explicitly that for $-\frac{\pi}{2}<\arg w<\frac{\pi}{2}$ this critical point is in the half-plane $\operatorname{Re}s>0$.

\begin{figure}\centering
\def\svgscale{0.7}
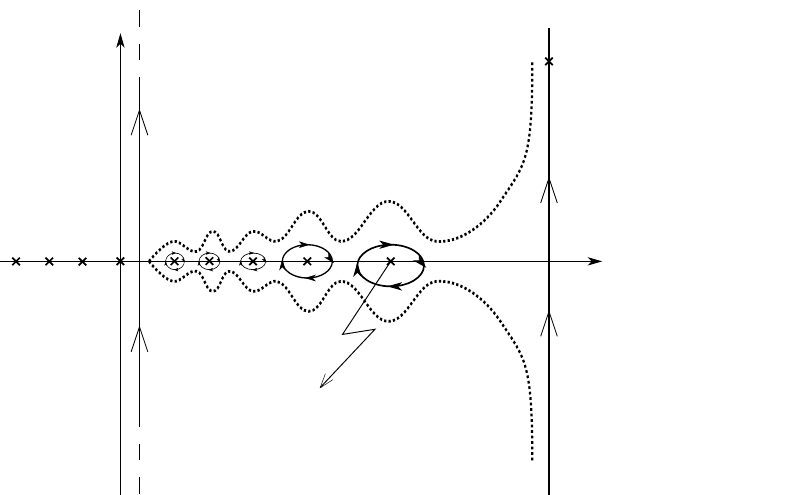
\caption{Deformation of the path $\Lambda_2$.}\label{defor}
\end{figure}

By modifying the path of integration as in Fig.~\ref{defor}, in order that it passes through the critical point, by Cauchy theorem we have
\[\vp_2(w)= \frac{(2\pi)^3}{2\pi{\rm i}}\int_{\Lambda_2^\prime}{\rm e}^{\phi(s)}\,{\rm d}s-
\sum_{p\in P}\underset{s=p}{\text{res}}\ \mathcal{I}_2(s),
\]
where $P$ stands for the set of poles in the region between $\Lambda_2$ and $\Lambda_2'$ and
\begin{gather*}
\mathcal{I}_2(s)= \Gamma(s)^5\Gamma\left(\frac{1}{2}-s\right){\rm e}^{{\rm i}\pi s}4^{-s}w^{-4s}
\end{gather*}
is the integrand in $\vp_2(w)$.
For the first summand we have an asymptotic behavior as before
\[
 \frac{(2\pi)^3}{2\pi{\rm i}}\int_{\Lambda_2^\prime}{\rm e}^{\phi(s)}\,{\rm d}s
 \sim \alpha\frac{{\rm e}^{-4\sqrt{2}w}}{w^2}
\]
with $\alpha$ a certain constant we will not need to specify. For the second summand, keeping into account that $\underset{z=n}{\text{res}}\ \Gamma(-z)=(-1)^{n+1}/n!$ for $n\in\mathbb N$, and that $\Gamma\big(z+\frac{1}{2}\big) =\pi^{1/2}2^{-2z+1}\Gamma(2z)\Gamma(z)^{-1}$, we receive
\begin{align*}\underset{s=n+\frac{1}{2}}{\text{res}}\ \mathcal I_2(s)&=\frac{(-1)^{n+1}}{n!}\Gamma\left(n+\frac{1}{2}\right)^5{\rm e}^{{\rm i}\pi\left(n+\frac{1}{2}\right)}4^{-n-\frac{1}{2}}w^{-4n-2}\\
&=-\frac{i}{n!}\left(\frac{(2n-1)!!}{2^n}\pi^{\frac{1}{2}}\right)^54^{-n-\frac{1}{2}}w^{-4n-2}.
\end{align*}
Therefore,
 \[\sum_{p\in P}\underset{s=p}{\text{res}}\ \mathcal I_2(s)=-\frac{{\rm i}\pi^{\frac{5}{2}}}{2w^2}-\frac{{\rm i}\pi^{\frac{5}{2}}}{256w^6}+O\left(\frac{1}{w^{10}}\right).
\]
In conclusion,
\[\vp_2(w)\sim\frac{{\rm i}\pi^{\frac{5}{2}}}{2w^2}\left(1+O\left(\frac{1}{w}\right)\right),
\qquad \text{for $-\frac{\pi}{2}<\arg w<\frac{\pi}{2}$}.
\]
Let us now use the identity \eqref{rel1} in the following form:
\[
\vp_2(w)=2\pi\vp_1\big(w {\rm e}^{-{\rm i}\frac{\pi}{2}}\big)-\vp_2\big(w{\rm e}^{-{\rm i}\frac{\pi}{2}}\big),\qquad -\frac{\pi}{2}<\arg \big(w {\rm e}^{-{\rm i}\frac{\pi}{2}}\big)<\frac{\pi}{2}.
\]
It implies that
\begin{gather*}
\vp_2(w)\sim\frac{{\rm i}\pi^{\frac{5}{2}}}{2w^2}\left(1+O\left(\frac{1}{w}\right)\right)\qquad\text{on the whole sector }-\frac{\pi}{2}<\arg w<\pi.\tag*{\qed}
\end{gather*}\renewcommand{\qed}{}
\end{proof}

\section{Computation of the central connection matrix}\label{appc}

Here we summarize the explicit values for the columns of the central connection matrix \mbox{$C\!=\!(C_{ij})$}, computed in Section~\ref{28luglio2016-15}, where $v$ is indicated. The correct value is $v=6$ ($v$~was first introduced in~\eqref{28luglio2016-18-1}).
\begin{gather*}
C_{i1}=\left(
\begin{matrix}
 \frac{1}{2 \sqrt{c} \pi ^2} \\
 \frac{4 \gamma +{\rm i}\pi }{2 \sqrt{c} \pi ^2} \\
 \frac{48 \gamma ^2+24 {\rm i} \gamma \pi -5 \pi ^2}{12 \sqrt{c} \pi ^2} \\
 \frac{48 \gamma ^2+24 {\rm i} \gamma \pi +7 \pi ^2}{12 \sqrt{c} \pi ^2} \\
 \frac{64 \gamma ^3+48 {\rm i} \gamma ^2 \pi +4 \gamma \pi ^2+3 {\rm i}\pi ^3-4 \zeta (3)}{6 \sqrt{c} \pi ^2} \\
 \frac{768 \gamma ^4+768 {\rm i} \gamma ^3 \pi +96 \gamma ^2 \pi ^2+144 {\rm i} \gamma \pi ^3-\pi ^4-48 (4 \gamma +{\rm i}\pi ) \zeta (3)}{72 \sqrt{c} \pi ^2}
\end{matrix}
\right),\\
C_{i2}=\left(
\begin{matrix}
 \frac{1}{2 \sqrt{c} \pi ^2} \\
 \frac{4 \gamma +{\rm i}\pi }{2 \sqrt{c} \pi ^2} \\
 \frac{48 \gamma ^2+24 {\rm i} \gamma \pi +7 \pi ^2}{12 \sqrt{c} \pi ^2} \\
 \frac{48 \gamma ^2+24 {\rm i} \gamma \pi -5 \pi ^2}{12 \sqrt{c} \pi ^2} \\
 \frac{64 \gamma ^3+48 {\rm i} \gamma ^2 \pi +4 \gamma \pi ^2+3 {\rm i}\pi ^3-4 \zeta (3)}{6 \sqrt{c} \pi ^2} \\
 \frac{768 \gamma ^4+768 {\rm i} \gamma ^3 \pi +96 \gamma ^2 \pi ^2+144 {\rm i} \gamma \pi ^3-\pi ^4-48 (4 \gamma +{\rm i}\pi ) \zeta (3)}{72 \sqrt{c} \pi ^2}
\end{matrix}
\right),\\
C_{i3}=\left(
\begin{matrix}
 -\frac{1}{4 \sqrt{c} \pi ^2} \\
 \frac{-2 \gamma -{\rm i}\pi }{2 \sqrt{c} \pi ^2} \\
 \frac{-48 \gamma ^2-48 {\rm i} \gamma \pi +11 \pi ^2}{24 \sqrt{c} \pi ^2} \\
 \frac{-48 \gamma ^2-48 {\rm i} \gamma \pi +11 \pi ^2}{24 \sqrt{c} \pi ^2} \\
 \frac{2 \zeta (3)-(2 \gamma +{\rm i}\pi ) (4 \gamma +{\rm i}\pi ) (4 \gamma +3 {\rm i}\pi )}{6 \sqrt{c} \pi ^2} \\
 \frac{-768 \gamma ^4-1536 {\rm i} \gamma ^3 \pi +1056 \gamma ^2 \pi ^2-23 \pi ^4+96 {\rm i}\pi \zeta (3)+96 \gamma \left(3 {\rm i}\pi ^3+2 \zeta (3)\right)}{144 \sqrt{c} \pi ^2}
\end{matrix}
\right),\\
C_{i4}(v)=\left(
\begin{matrix}
 \frac{v -1}{4 \sqrt{c} \pi ^2} \\
 \frac{2 \gamma (v -1)+{\rm i}\pi }{2 \sqrt{c} \pi ^2} \\
 \frac{48 \gamma ^2 (v -1)+48 {\rm i} \gamma \pi +(v +11) \pi ^2}{24 \sqrt{c} \pi ^2} \\
 \frac{48 \gamma ^2 (v -1)+48 {\rm i} \gamma \pi +(v +11) \pi ^2}{24 \sqrt{c} \pi ^2} \\
 \frac{32 \gamma ^3 (v -1)+48 {\rm i} \gamma ^2 \pi +2 \gamma (v +11) \pi ^2-3 {\rm i}\pi ^3-2 (v -1) \zeta (3)}{6 \sqrt{c} \pi ^2} \\
 \frac{768 \gamma ^4 (v -1)+1536 {\rm i} \gamma ^3 \pi +96 \gamma ^2 (v +11) \pi ^2-(v +23) \pi ^4-96 {\rm i}\pi \zeta (3)+96 \gamma \left(-3 {\rm i}\pi ^3-2 (v -1) \zeta (3)\right)}{144 \sqrt{c} \pi ^2}
\end{matrix}
\right),\\
C_{i4}(6)=\left(
\begin{matrix}
 \frac{5}{4 \sqrt{c} \pi ^2} \\
 \frac{10 \gamma +{\rm i}\pi }{2 \sqrt{c} \pi ^2} \\
 \frac{240 \gamma ^2+48 {\rm i} \gamma \pi +17 \pi ^2}{24 \sqrt{c} \pi ^2} \\
 \frac{240 \gamma ^2+48 {\rm i} \gamma \pi +17 \pi ^2}{24 \sqrt{c} \pi ^2} \\
 \frac{160 \gamma ^3+48 {\rm i} \gamma ^2 \pi +34 \gamma \pi ^2-3 {\rm i}\pi ^3-10 \zeta (3)}{6 \sqrt{c} \pi ^2} \\
 \frac{3840 \gamma ^4+1536 {\rm i} \gamma ^3 \pi +1632 \gamma ^2 \pi ^2-288 {\rm i} \gamma \pi ^3-29 \pi ^4-960 \gamma \zeta (3)-96 {\rm i}\pi \zeta (3)}{144 \sqrt{c} \pi ^2}
\end{matrix}
\right),\\
C_{i5}=\left(
\begin{matrix}
 \frac{1}{4 \sqrt{c} \pi ^2} \\
 \frac{\gamma }{\sqrt{c} \pi ^2} \\
 \frac{48 \gamma ^2+\pi ^2}{24 \sqrt{c} \pi ^2} \\
 \frac{48 \gamma ^2+\pi ^2}{24 \sqrt{c} \pi ^2} \\
 \frac{-\zeta (3)+16 \gamma ^3+\gamma \pi ^2}{3 \sqrt{c} \pi ^2} \\
 -\frac{192 \gamma \zeta (3)-768 \gamma ^4+\pi ^4-96 \gamma ^2 \pi ^2}{144 \sqrt{c} \pi ^2}
\end{matrix}
\right),\\
C_{i6}=\left(
\begin{matrix}
 \frac{1}{4 \sqrt{c} \pi ^2} \\
 \frac{\gamma +{\rm i}\pi }{\sqrt{c} \pi ^2} \\
 \frac{48 \gamma ^2+96 {\rm i} \gamma \pi -47 \pi ^2}{24 \sqrt{c} \pi ^2} \\
 \frac{48 \gamma ^2+96 {\rm i} \gamma \pi -47 \pi ^2}{24 \sqrt{c} \pi ^2} \\
 \frac{(\gamma +{\rm i}\pi ) (4 \gamma +3 {\rm i}\pi ) (4 \gamma +5 {\rm i}\pi )-\zeta (3)}{3 \sqrt{c} \pi ^2} \\
 \frac{768 \gamma ^4+3072 {\rm i} \gamma ^3 \pi -4512 \gamma ^2 \pi ^2-2880 {\rm i} \gamma \pi ^3+671 \pi ^4-192 (\gamma +{\rm i}\pi ) \zeta (3)}{144 \sqrt{c} \pi ^2}
\end{matrix}
\right).
\end{gather*}

We write now entries of the matrix $C^-_{\rm Kap}$ whose columns are given by the components of the characteristic classes
\[\frac{1}{4\pi c^\frac{1}{2}}\widehat\Gamma^-(\mathbb G)\cup\operatorname{Ch}\big(\mathbb S^\lambda(\mathcal S^*)\big);
\]the order of the column is given by $\lambda=0$, $\lambda=1$, $\lambda=2$, $\lambda=(1,1)$, $\lambda =(2,1)$ and $\lambda=(2,2)$.
\begin{gather*}
\big(C^-_{\rm Kap}\big)_0=\left(
\begin{matrix}
 \frac{1}{4 \sqrt{c} \pi ^2} \\
 \frac{\gamma }{\sqrt{c} \pi ^2} \\
 \frac{\frac{1}{24}+\frac{2 \gamma ^2}{\pi ^2}}{\sqrt{c}} \\
 \frac{\frac{1}{24}+\frac{2 \gamma ^2}{\pi ^2}}{\sqrt{c}} \\
 \frac{-\zeta (3)+16 \gamma ^3+\gamma \pi ^2}{3 \sqrt{c} \pi ^2} \\
 -\frac{192 \gamma \zeta (3)-768 \gamma ^4+\pi ^4-96 \gamma ^2 \pi ^2}{144 \sqrt{c} \pi ^2}
\end{matrix}
\right),\\
\big(C^-_{\rm Kap}\big)_{\text{\tiny $\yng(1)$}}=\left(
\begin{matrix}
 \frac{1}{2 \sqrt{c} \pi ^2} \\
 \frac{4 \gamma +{\rm i}\pi }{2 \sqrt{c} \pi ^2} \\
 \frac{\frac{2 \gamma (2 \gamma +{\rm i}\pi )}{\pi ^2}-\frac{5}{12}}{\sqrt{c}} \\
 \frac{\frac{2 \gamma (2 \gamma +{\rm i}\pi )}{\pi ^2}+\frac{7}{12}}{\sqrt{c}} \\
 \frac{64 \gamma ^3+48 {\rm i} \gamma ^2 \pi +4 \gamma \pi ^2+3 {\rm i}\pi ^3-4 \zeta (3)}{6 \sqrt{c} \pi ^2} \\
 \frac{768 \gamma ^4+768 {\rm i} \gamma ^3 \pi +96 \gamma ^2 \pi ^2+144 {\rm i} \gamma \pi ^3-\pi ^4-48 (4 \gamma +{\rm i}\pi ) \zeta (3)}{72 \sqrt{c} \pi ^2} \\
\end{matrix}
\right),
\\
\big(C^-_{\rm Kap}\big)_{\text{\tiny $\yng(2)$}}=\left(
\begin{matrix}
 \frac{3}{4 \sqrt{c} \pi ^2} \\
 \frac{3 (2 \gamma +{\rm i}\pi )}{2 \sqrt{c} \pi ^2} \\
 \frac{\frac{6 \gamma (\gamma +{\rm i}\pi )}{\pi ^2}-\frac{19}{8}}{\sqrt{c}} \\
 \frac{\frac{6 \gamma (\gamma +{\rm i}\pi )}{\pi ^2}+\frac{13}{8}}{\sqrt{c}} \\
 \frac{32 \gamma ^3+48 {\rm i} \gamma ^2 \pi -6 \gamma \pi ^2+5 {\rm i}\pi ^3-2 \zeta (3)}{2 \sqrt{c} \pi ^2} \\
 \frac{-6 \gamma ^2+\frac{16 \gamma ^4}{\pi ^2}+\frac{32 {\rm i} \gamma ^3}{\pi }+10 {\rm i} \gamma \pi +\frac{7 \pi ^2}{48}-\frac{2 (2 \gamma +{\rm i}\pi ) \zeta (3)}{\pi ^2}}{\sqrt{c}}
\end{matrix}
\right),
\\
\big(C^-_{\rm Kap}\big)_{\text{\tiny $\yng(1,1)$}}=\left(
\begin{matrix}
 \frac{1}{4 \sqrt{c} \pi ^2} \\
 \frac{2 \gamma +{\rm i}\pi }{2 \sqrt{c} \pi ^2} \\
 \frac{\frac{2 \gamma (\gamma +{\rm i}\pi )}{\pi ^2}-\frac{11}{24}}{\sqrt{c}} \\
 \frac{\frac{2 \gamma (\gamma +{\rm i}\pi )}{\pi ^2}-\frac{11}{24}}{\sqrt{c}} \\
 \frac{(2 \gamma +{\rm i}\pi ) (4 \gamma +{\rm i}\pi ) (4 \gamma +3 {\rm i}\pi )-2 \zeta (3)}{6 \sqrt{c} \pi ^2} \\
 \frac{768 \gamma ^4+1536 {\rm i} \gamma ^3 \pi -1056 \gamma ^2 \pi ^2-288 {\rm i} \gamma \pi ^3+23 \pi ^4-96 (2 \gamma +{\rm i}\pi ) \zeta (3)}{144 \sqrt{c} \pi ^2}
\end{matrix}
\right),
\\
\big(C^-_{\rm Kap}\big)_{\text{\tiny $\yng(2,1)$}}=\left(
\begin{matrix}
 \frac{1}{2 \sqrt{c} \pi ^2} \\
 \frac{4 \gamma +3 {\rm i}\pi }{2 \sqrt{c} \pi ^2} \\
 \frac{\frac{2 \gamma (2 \gamma +3 {\rm i}\pi )}{\pi ^2}-\frac{29}{12}}{\sqrt{c}} \\
 \frac{\frac{2 \gamma (2 \gamma +3 {\rm i}\pi )}{\pi ^2}-\frac{17}{12}}{\sqrt{c}} \\
 \frac{(4 \gamma +{\rm i}\pi ) (4 \gamma +3 {\rm i}\pi ) (4 \gamma +5 {\rm i}\pi )-4 \zeta (3)}{6 \sqrt{c} \pi ^2} \\
 \frac{768 \gamma ^4+2304 v \gamma ^3 \pi -2208 \gamma ^2 \pi ^2-720 {\rm i} \gamma \pi ^3+47 \pi ^4-48 (4 \gamma +3 {\rm i}\pi ) \zeta (3)}{72 \sqrt{c} \pi ^2}
\end{matrix}
\right),
\\
\big(C^-_{\rm Kap}\big)_{\text{\tiny $\yng(2,2)$}}=\left(
\begin{matrix}
 \frac{1}{4 \sqrt{c} \pi ^2} \\
 \frac{\gamma +{\rm i}\pi }{\sqrt{c} \pi ^2} \\
 \frac{\frac{2 \gamma (\gamma +2 {\rm i}\pi )}{\pi ^2}-\frac{47}{24}}{\sqrt{c}} \\
 \frac{\frac{2 \gamma (\gamma +2 {\rm i}\pi )}{\pi ^2}-\frac{47}{24}}{\sqrt{c}} \\
 \frac{(\gamma +{\rm i}\pi ) (4 \gamma +3 {\rm i}\pi ) (4 \gamma +5 {\rm i}\pi )-\zeta (3)}{3 \sqrt{c} \pi ^2} \\
 \frac{768 \gamma ^4+3072 {\rm i} \gamma ^3 \pi -4512 \gamma ^2 \pi ^2-2880 {\rm i} \gamma \pi ^3+671 \pi ^4-192 (\gamma +{\rm i}\pi ) \zeta (3)}{144 \sqrt{c} \pi ^2}
\end{matrix}
\right).
\end{gather*}

By application of the constraint \[S=\big(C^-_{\rm Kap}\big)^{-1}{\rm e}^{-\pi{\rm i}R}{\rm e}^{-\pi{\rm i}\mu}\eta^{-1}\big((C^-_{\rm Kap})^{\rm T}\big)^{-1},
\]we find
\[S_{\rm Kap}=\left(
\begin{matrix}
 1 & -4 & 6 & 10 & -20 & 20 \\
 0 & 1 & -4 & -4 & 16 & -20 \\
 0 & 0 & 1 & 0 & -4 & 6 \\
 0 & 0 & 0 & 1 & -4 & 10 \\
 0 & 0 & 0 & 0 & 1 & -4 \\
 0 & 0 & 0 & 0 & 0 & 1
\end{matrix}
\right),
\qquad
S_{\rm Kap}^{-1}=\left(
\begin{matrix}
 1 & 4 & 10 & 6 & 20 & 20 \\
 0 & 1 & 4 & 4 & 16 & 20 \\
 0 & 0 & 1 & 0 & 4 & 10 \\
 0 & 0 & 0 & 1 & 4 & 6 \\
 0 & 0 & 0 & 0 & 1 & 4 \\
 0 & 0 & 0 & 0 & 0 & 1
\end{matrix}
\right).
\]
Now, $S_{\rm Kap}^{-1}$ coincides with the Gram matrix $\bold{G}_{\rm Kap}=\big(\chi\big(\mathbb S^\lambda\mathcal S^*,\mathbb S^\mu\mathcal S^*\big)\big)_{\lambda,\mu}$ of the Kapranov exceptional collection.

\subsection*{Acknowledgements} We would like to thank Marco Bertola, Ugo Bruzzo, Barbara Fantechi, Claus Hertling, Claude Sabbah, Maxim Smirnov, Jacopo Stoppa, Ian Strachan and Di Yang for several discussions and helpful comments. We also would like to thank the anonymous referees, whose valuable comments have improved the paper. The first author is grateful to the Max-Planck Institut f\"ur Mathematik in Bonn, for hospitality and support. The third author is a member of the European Union's H2020 research and innovation programme under the Marie Sk\l{l}odowska-Curie grant No.~778010 {\it IPaDEGAN}.

\LastPageEnding

\end{document}